%% file: Franc-RR-9445.tex
\documentclass[twoside]{book}
\usepackage[a4paper]{geometry}
\usepackage{RR}

\usepackage[T1]{fontenc}
\usepackage[utf8]{inputenc}
\usepackage[square,sort,comma,numbers]{natbib}
\usepackage{amsmath,amssymb,amscd,amsfonts,amsthm}
\usepackage{mathrsfs}
\usepackage{stmaryrd,bm}
\usepackage{enumitem}
\usepackage{manfnt} 
\usepackage{algorithm}
\usepackage{algorithmic}
\usepackage{listings}
\usepackage{fancybox}
\usepackage{array}
\usepackage{longtable}
\usepackage{datetime}
\usepackage{tikz}
\usetikzlibrary{shapes,arrows,cd}
\usepackage{pgfplots}
\usepackage{tensor}

\usepackage[colorlinks=true, allcolors=blue]{hyperref}
\def\A{\mathcal{A}}
\def\B{\mathcal{B}}
\def\C{\mathbb{C}}
\def\F{\mathcal{F}}
\def\G{\mathcal{G}}
\def\H{\mathcal{H}}
\def\I{\mathbb{I}}
\def\K{\mathbb{K}}
\def\L{\mathcal{L}}
\def\M{\mathcal{M}}
\def\N{\mathbb{N}}
\def\O{\mathcal{O}}
\def\P{\mathbb{P}}
\def\Q{\mathbb{Q}}
\def\R{\mathbb{R}}
\def\S{\mathbb{S}}
\def\T{\mathcal{T}}
\def\U{\mathcal{U}}

\def\Z{\mathbb{Z}}
\def\AS{\mathscr{A}}

\def\EC{\mathcal{E}}
\def\ES{\mathscr{E}}
\def\FS{\mathscr{F}}
\def\GL{\mathbb{GL}}

\def\IC{\mathcal{I}}
\def\JC{\mathcal{J}}

\def\NC{\mathcal{N}}

\def\PC{\mathcal{P}}
\def\PS{\mathscr{P}}
\def\QS{\mathscr{Q}}

\def\SF{\mathfrak{S}}

\def\Abold{\mathbf{A}}
\def\Bbold{\mathbf{B}}
\def\Cbold{\mathbf{C}}
\def\Dbold{\mathbf{D}}
\def\Ebold{\mathbf{E}}

\def\Gbold{\mathbf{G}}

\def\Pbold{\mathbf{P}}

\def\Tbold{\mathbf{T}}

\def\Xbold{\mathbf{X}}
\def\Ybold{\mathbf{Y}}

\def\abold{\mathbf{a}}
\def\bbold{\mathbf{b}}
\def\cbold{\mathbf{c}}
\def\dbold{\mathbf{d}}
\def\ebold{\mathbf{e}}
\def\fbold{\mathbf{f}}
\def\gbold{\mathbf{g}}
\def\ibold{\mathbf{i}}
\def\jbold{\mathbf{j}}

\def\mbold{\mathbf{m}}
\def\nbold{\mathbf{n}}
\def\rbold{\mathbf{r}}
\def\tbold{\mathbf{t}}
\def\ubold{\mathbf{u}}
\def\vbold{\mathbf{v}}
\def\wbold{\mathbf{w}}
\def\xbold{\mathbf{x}}
\def\ybold{\mathbf{y}}
\def\zbold{\mathbf{z}}
\def\phibold{\bm{\phi}}
\def\a{\textsc{a}}
\def\b{\textsc{b}}
\def\c{\textsc{c}}
\def\e{\textsc{e}} 
\def\f{\textsc{f}}
\def\g{\textsc{g}}

\def\i{\textsc{i}}
\def\j{\textsc{j}}
\def\k{\textsc{k}}

\def\m{\textsc{m}}
\def\n{\textsc{n}}
\def\p{\textsc{p}}
\def\q{\textsc{q}}

\def\t{\textsc{t}}

\def\rfrak{\mathfrak{r}}

\def\brank{\mathrm{brank}\:}
\def\det{\mathrm{det}\:}

\def\Diag{\mathrm{Diag}\:}
\def\dim{\mathrm{dim}\:}

\def\Im{\mathrm{Im}\:}

\def\Log{\mathrm{Log}\:}
\def\ones{\mathbf{1}}

\def\rank{\mathrm{rank}\:}
\def\rr{\mathfrak{r}}
\def\span{\mathrm{span}\:}
\def\st{\quad\mathrm{s.t.}\quad}
\def\Tr{\mathrm{Tr}\:}
\def\ttrank{\mbox{TT-rank}\:}
\def\vec{\mathrm{vec}\:}

\def\nB{\vspace*{5mm}\noindent}
\def\nS{\vspace*{5mm}\noindent$\bullet\quad$}
\def\nT#1{\vspace*{5mm}\noindent$\bullet\quad$\textbf{#1:} }
\def\nP#1{\vspace*{5mm}\noindent#1$\quad$}

\numberwithin{equation}{section}
\def\notes{\paragraph{References and notes:}}

\usepackage{makeidx}
\makeindex
\def\kw#1{\emph{#1}\marginpar{\small{\emph{#1}}}\index{#1}}
\def\kwnind#1{\emph{#1}\marginpar{\small{\emph{#1}}}}

\RRNo{9445}
\RRdate{December 2021}

\RRauthor{%
  Alain Franc%
  \thanks{Pleiade team and INRAE, Biogeco, University of Bordeaux,
    69, route d'Arcachon, 33610, Cestas}%
  \thanks{Correspondence: \texttt{alain.franc@inrae.fr}}
}

\authorhead{Franc}
\RRtitle{%
 Rangs des tenseurs : bases pour la réduction de
 dimension et la séparation des variables
}
\RRetitle{%
  Tensor Ranks for the Pedestrian \\
  for Dimension Reduction \\
  and Disentangling Interactions}
\titlehead{Tensor Ranks for the Pedestrian}
\RRresume{%
  Un tenseur est notamment un tableau à plusieurs entrées
  qui peut représenter, outre un jeu de données,
  l'expression d'une loi jointe ou d'une fonction
  multivariée. Il contient alors la description des
  interactions entre les variables correspondant à chacune
  des entrées. Le rang d'un tenseur étend à des tableaux à
  plus de deux entrées la notion de rang d'une matrice,
  sachant qu'il existe plusieurs approches pour construire
  une telle extension. Lorsque le rang vaut un, les
  variables sont séparées, et lorsqu'il est faible, les
  variables sont faiblement couplées. Bien des calculs sont
  plus simples sur des tenseurs de rang faible. Aussi,
  approcher un tenseur donné par un tenseur de rang faible
  permet de les rendre possibles pour calculer certaines
  caractéristiques d'un tableau, comme par exemple la
  fonction de partition quand il s'agit d'une loi
  jointe. Dans cette note, nous présentons en détail une
  approche intégrée et progressive pour approcher un tenseur
  donné par un tenseur de rang plus faible, par une
  utilisation systématique de l'algèbre tensorielle. La
  notion de tenseur est définie rigoureusement, puis des
  opérations élémentaires mais utiles sur les tenseurs sont
  présentées.  Après avoir rappelé plusieurs notions
  différentes pour le rang d'un tenseur, nous montrons
  comment ces opérations élémentaires peuvent être combinées
  pour construire des algorithmes d'approximation de rang
  faible. Le dernier chapitre est consacré à appliquer cette
  approche aux tenseurs construits comme la discrétisation
  d'une fonction multivariée, pour montrer que sur une
  grille cartésienne, le rang de tels tenseurs est en
  général faible. 
}

\RRabstract{%
  A tensor is a multi-way array that can represent, in
  addition to a data set, the expression of a joint law or a
  multivariate function. As such it contains
  the description of the interactions between the variables
  corresponding to each of the entries. The rank of a tensor
  extends to arrays with more than two entries the notion of
  rank of a matrix, bearing in mind that there are several
  approaches to build such an extension. When the rank is
  one, the variables are separated, and when it is low, the
  variables are weakly coupled. Many calculations are
  simpler on tensors of low rank. Furthermore, approximating
  a given tensor by a low-rank tensor makes it possible to
  compute some characteristics of a table, such as the
  partition function when it is a joint law. In this note,
  we present in detail an integrated and progressive
  approach to approximate a given tensor by a tensor of
  lower rank, through a systematic use of tensor
  algebra. The notion of tensor is rigorously defined, then
  elementary but useful operations on tensors are presented.
  After recalling several different notions for extending
  the rank to tensors, we show how these elementary
  operations can be combined to build best low rank
  approximation algorithms. The last chapter is devoted to
  applying this approach to tensors constructed as the
  discretisation of a multivariate function, to show that on
  a Cartesian grid, the rank of such tensors is expected to
  be low.
}

\RRmotcle{%
  Tenseur,
  Produit tensoriel,
  Algèbre tensorielle,
  rang d'un tenseur,
  Candecomp,
  Parafac,
  modèle de Tucker,
  Tensor Train,
  Réduction de dimension,
  Séparation de variables
}

\RRkeyword{%
  Tensors,
  Tensor Product,
  Tensor Algebra,
  Tensor rank,
  Candecomp,
  Parafac,
  Tucker model,
  Tensor Train,
  Dimension reduction,
  Variable separation
}

\RRprojet{Pleiade}

\RCBordeaux

\begin{document}
\makeRR


\tableofcontents



\tikzstyle{tensor}=[draw,circle,fill=blue!25]
\tikzstyle{tensor_ell}=[draw,ellipse,fill=blue!25]
\tikzstyle{tensor_box}=[draw,ellipse,fill=blue!25]

\include{introduction}
\part{Tensor products}\label{part:tensors}
\include{tensor_product}

\part{A primer in tensor algebra}\label{part:primer}
\include{primer}

\part{Tensor ranks}\label{part:ranks}
\include{tensor_rank}

\include{intermezzo}
\part{Best approximations of low ranks}\label{part:blra}

\include{best_low_rank} 
\part{Multivariate functions as tensors}\label{part:funct}
\include{quantization}

\clearpage
\bibliographystyle{alpha}
\bibliography{tensors}

\clearpage
\printindex
\end{document}

%% file: introduction.tex
\begin{center}
    \textbf{\Huge Table of notations}

\vspace*{2cm}

\begin{tabular}{cl}
 $\otimes$ & tensor product \\
 $\bullet$ & contraction \\
 $\llbracket m,n \rrbracket$ & the set of integers $k$ with $m \leq k \leq n$ \\
 $\{m,\ldots,n\}$ & the set of integers $k$ with $m \leq k \leq n$ \\
 $\abold, \bbold, \ldots$ & vectors \\
 $a_\ibold$ & a coordinate of $\Abold$ in basis $\EC_\mu$ for each space $E_\mu$ \\
 $A,B,\ldots, M,N, \ldots$ & matrices in $E \otimes F$ \\
 $\Abold, \Bbold, \ldots, \Tbold, \ldots$ & tensors in $E_1 \otimes \ldots \otimes E_d$ \\
 $\alpha,\beta, \ldots$ & scalars in $\R$ \\
 $d$ & order of a tensor \\
 $B(E,F)$ & bilinear maps on $E \times F$ \\
 $\ebold_{\mu,i_\mu}, \: \ebold_{i_\mu}^{(\mu)}$ & an element of $\EC_\mu$ \\
 $E, F, \ldots$ & a vector space on a field $\K$ \\
 $E^*$ & the dual of $E$ \\
 $E_\mu$ & a vector space in a family $\ES$ of $d$ vector spaces \\
 $\Ebold_\ibold$ & element $e_{1i_1} \otimes \ldots \otimes e_{di_d}$ of basis $\bm{\EC}$ \\
 $\EC$ & $(e_1,\ldots,e_n)$, a basis of $E$ \\
 $\EC^*$ & $(e_1^*,\ldots,e_n^*)$, the basis of $E^*$ dual of $\EC$\\
 $\EC_\mu$ & an orthonormal basis of $E_\mu$ \\
 $\ES$ & a family $(E_1, \ldots,E_d)$ of $d$ vector spaces indexed by $\mu$ \\  
 $\displaystyle \bigotimes_{\mu=1}^d E_\mu^*$ & vector space of $d-$linear forms on $E_1 \times \ldots \times E_d$ \\
 $\bm{\EC}$ & basis of $E_1 \otimes \ldots \otimes E_d$ from $(\EC_1, \ldots,\EC_d)$\\
 $\GL(E)$ & linear group of $E$ \\
 $\ibold$ & $\in \IC$, a multiindex $\ibold = (i_1, \ldots, i_d)$ \\
 $I$ & $\llbracket 1,d \rrbracket \subset \N$ \\
 $I_\mu$ & set of indices for components in $E_\mu$;$I_\mu = \{1,n_\mu\}$\\
 $\IC$ & $I_1 \times \ldots \times I_d$ \\
 $\K$ & a field (usually $\R$, or $\C$)\\
 $\K^{m \times n}$ & set of matrices with $m$ rows and $n$ columns with elements in $\K$ \\
 $\K^{n_1 \times \ldots \times n_d}$ & equivalent to $\K^{n_1} \otimes \ldots \otimes \K^{n_d}$\\ 
 $\L(E,F)$ & linear maps from $E$ to $F$ \\
 $\mu$ & index running over $I$, $\mu \in I$ \\
 $n_\mu$ & dimension of space $E_\mu: \: n_\mu=\dim E_\mu$ \\
 $\T^d(E)$ & vector space of $d-$linear forms on a vector space $E$ \\
 $\T^d(E_1, \ldots, E_d)$ & vector space of $d-$linear forms on $E_1 \times \ldots \times E_d$ \\
 $\xbold,\ybold,\ldots$ & vectors in $E, F, \ldots$
\end{tabular}

\end{center}

\clearpage

\section*{About notations}

\noindent Notations will be introduced gradually when
necessary. There is no consensus for notations among
textbooks and research papers dealing with tensors and
tensor rank, and I will point out some differences to
facilitate connections between different works. For example,
the symbol $\otimes$ for tensor product itself (chosen here)
is not universally adopted, some classical papers proposing
$\circ$ (like \cite{Kolda2009}).

%
\section*{Conventions for fonts}
%

\nS The integers, mainly indices, will be denoted by lower
case Greek or Latin letters, as in $\alpha,\beta, \ldots$ or
$a,b,\ldots$, $i,k,\ldots$ or $m,n,\ldots$.

\nS Scalars in $\R$ will be denoted in the same
manner. Scalars will be denoted by a Greek letter especially
when they act as a factor before a vector or a matrix, as in
$\alpha\, \xbold$.

\nS Vectors will be denoted by boldface small Latin letters,
like $\xbold, \ybold, \ldots$, or $\ubold, \vbold, \ldots$
or $\abold, \bbold, \ldots $

\nS Matrices will be denoted by capital Latin letters, like
$A,B, \ldots, M, \ldots, X, \ldots$

\nS Tensors will be denoted by capital bold Latin letters,
like
$\Abold, \Bbold, \ldots, \Tbold, \ldots, \Xbold, \ldots$

\nS Vector spaces will be denoted by capital letters, like
$E,F,G,\ldots$

\nS Sets will be denoted by script capital letters like $\A$
or $\AS$, like a basis $\EC$ for a vector space $E$ with
$\EC=(\ebold_1,\ldots,\ebold_n)$.

%
\section*{Notations for vector spaces}
%

\nS Vector spaces will be denoted by $E,F, \ldots$. Let us
present some conventions adopted here for a vector space
denoted $E$.

\nS Its dimension will be denoted $m$ or $n$. It is endowed
with a basis denoted $\EC$, defined as
\[
\EC = (\ebold_1, \ldots,\ebold_n)
\]
A vector $\xbold \in E$ is decomposed on this basis as
\[
\xbold = \sum_{i=1}^n\, x_i\, \ebold_i
\]

\nS A family of vector spaces will be denoted $\ES$ with
\[
\ES = (E_1,\ldots,E_\mu, \ldots, E_d)
\]
The number of vector spaces in the family will usually be
denoted by $d$, and the running index over the family by
$\mu$. The dimension of $E_\mu$ is $n_\mu$.
The basis of $E_\mu$ is denoted $\EC_\mu$, with
\[
\EC_\mu = (\ebold_{\mu,1}, \ldots,\ebold_{\mu,d})
\]
or
\[
\EC_\mu = \left(\ebold_1^{(\mu)}, \ldots, \ebold_d^{(\mu)}\right)
\]
A vector $\xbold_\mu \in E_\mu$ is decomposed as
\[
\xbold_\mu = \sum_{i=1}^{n_\mu}\, x_{\mu,i}\, \ebold_{\mu,i}
\]
or
\[
\xbold_\mu = \sum_{i=1}^{n_\mu}\, x_i^{(\mu)}\, \ebold_i^{(\mu)}
\]

\nS The same conventions are adopted for vector spaces
denoted $E,F,G$ with following dimensions, indices, and
basis:\\
\begin{center}
\begin{tabular}{cccc}
     space & dimension & indices & basis \\
     \hline
     $E$ & $m$ & $i$ & $(\ebold_i)_i$ \\ 
     $F$ & $n$ & $j$ & $(\fbold_j)_j$ \\ 
     $G$ & $p$ & $k$ & $(\gbold_k)_k$ \\ 
\end{tabular}
\end{center}

%
\section*{Notations for multi-indices}
%

Multi-indices are everywhere, observed or hidden, in tensor algebra. 

\nS Let us have a family $\ES$ of $d$ vector spaces indexed
by $\mu$ as above. I denote
\[
I = \llbracket 1,d \rrbracket \subset \N
\]
We have
\[
\mu \in I
\]

\nS If $n_\mu = \dim E_\mu$, I denote 
\[
I_\mu = \llbracket 1,n_\mu\rrbracket \subset \N
\]
and
\[
i_\mu \in I_\mu
\]

\nS I denote
\[
\IC = I_1 \times \ldots \times I_d
\]
A multi-index is a tuple
\[
\ibold = (i_1, \ldots,i_d) \in \IC
\]

%
\section*{Notations for Tensors}
%

\nS If
\[
\Tbold \in E_1 \otimes \ldots \otimes E_d
\]
with $\EC_\mu = \left(\ebold_1^{(\mu)}, \ldots,\ebold_d^{(\mu)}\right)$, we have
\[
\Tbold = \sum_{i_1=1}^{n_1} \ldots \sum_{i_d=1}^{n_d} \, t_{i_1 \ldots i_d}\, \ebold_{1,i_1} \otimes  \ldots \otimes \ebold_{d,i_d}
\]
or (both will be used)
\[
\Tbold = \sum_{i_1=1}^{n_1} \ldots \sum_{i_d=1}^{n_d} \, t_{i_1 \ldots i_d}\, \ebold_1^{(i_1)} \otimes  \ldots \otimes \ebold_d^{(i_d)}
\]
I denote
\[
t_\ibold = t_{i_1\ldots i_d}, \qquad \Ebold_\ibold =  \ebold_{1,i_1} \otimes  \ldots \otimes \ebold_{d,i_d}
\]
Then, more compactly
\[
\Tbold = \sum_{\ibold \in \IC}\, t_\ibold \, \Ebold_\ibold
\]

%
\chapter{Introduction}\label{chap:intro}
%

This book is about dimension reduction in tensors and
disentangling interactions in systems modeled by tensors.

\nB Many problems arise in spaces whose large size can
exclude the possibility of finding a solution in a
reasonable time. This can happen in two situations at
least. First, with a few variables but in high dimensional
spaces, where solutions with algorithms with polynomial
complexity become out of scope. One classical solution is to
look for dimension reduction. Second, with many-variable
entangled systems, where building a solution comes up
against the curse of dimensionality. One classical solution
is to approximate the system with a system with local
interactions only, with the ultimate goal of reaching a
separation of variables. A solution to both problems can be
looked for with the following approach: select families of
small dimensional subspaces (linear or not) whose solution
is tractable, project into this family, and solve the
problem in the projected subspaces, hoping the solution will
be close to (or not too far from) the solution of the
initial problem. A classical example of dimension reduction
is Principal Component Analysis (PCA) which is the iconic
linear dimension reduction technique. A classical example of
separation of variables is mean-field approximation.

\nB A remarkable property of tensors is that both guises of
these problems have been rigorously formulated (and
sometimes solved) with the notion of \emph{rank} in
tensors. The candidate family for dimension reduction in
tensors is the set of tensors of a given (small) rank. So,
we study the best approximation of a given tensor by a
tensor of lower rank. The simplification of interactions can
be addressed by selecting families of tensors with a given
structure, like Tensor Trains with low ranks, which lead to
separation of variables. The possibility of addressing both
issues, dimension reduction and disentangling interactions,
is part of the beauty of tensors for modeling complex
problems. In many situations, there is a tension between
selecting a very low dimension for the problem to solve to
be simple, and a high one, for the solution in lower
dimension space to be an accurate approximation of the
initial one.

\nB Tensors are encountered in many domains. They have been
proposed (without the concept of vector, much less of
tensors, using identities in determinants and linear changes
of coordinates: see \cite{Elliott13,Grace1903}) by the
middle of the 19$^{\mathrm{th}}$ century for the study of
invariants in forms under a change of basis. At this time,
the equivalent terms \emph{form} or \emph{quantics} emerged
in algebra for what is now a symmetric tensor representing
an homogeneous polynomial in several variables (see notes
page \pageref{page:order} for further references). The term
\emph{tensor} with the current meaning was coined later by
the christallographer W. Voigt in 1898 (the term was
introduced by Hamilton in 1846 with a different meaning:
what we would call nowadays the norm of a Quaternion). This
notion was developed around 1900 for studying the geometry
of Riemannian manifolds, especially curvature, or transport,
by an Italian school around Levi-Civita, and
Ricci-Curbastro. Einstein used them for setting his theory
of general relativity and popularized them. A tensor can be
seen as a multilinear form, and it is this approach that is
adopted in this book. It is a linear combination of
elementary tensors, an elementary tensor being a tensor
product of vectors, in the same way a matrix is a linear
combination of rank on matrices (elementary matrices).

\nB Here is a partial overview of the different areas where
tensors are routinely found:
\begin{itemize}
    \item in data science, they formalize multiway
      arrays;
    \item in the study of numerical solutions of PDEs,
      especially linear ones, they are used as multiway
      arrays for discretization of multivariate functions on
      a Cartesian grid;
    \item in statistics, they represent either the
      distribution of a joint law on a product of discrete
      spaces, or the array of cumulants of a given order;
    \item in many-body quantum systems, they represent
      entangled states where the wave function lives in a
      tensor product of Hilbert spaces, finite dimensional
      or not.
\end{itemize}

\nB A tensor at this stage can be looked at as an array with
several entries (a matrix is a tensor with two entries), or
a box in a multidimensional space. The number of entries or
the dimension of the space in which the box lives is called
the \emph{order} of the tensor. The dimensions are the
lengths of the edges of the box. A tensor with $d$ orders
and dimension $n$ (to keep it simple) in each order has
$n^d$ elements, which becomes huge as soon as $d$ is
significant. Both order and dimensions cannot be large
simultaneously for the problem to be addressable with
current techniques. Dimension reduction is relevant when the
order is reasonable (the number of entries can be counted on
the fingers of one hand) and disentangling a large number of
variables in a system (one per order) is feasible when
dimensions for each order are small. Here is a cursory
presentation of the complexity in dimension/order of
different areas where tensors play a role:
\\
\begin{center}
    \begin{tabular}{l|ll}
       Area & order & dimension \\
       \hline
       Multiways array & low & medium\\
       Discretization of functions & low & high \\
       Statistics & high & low \\
       Quantum systems & very high & low \\
    \end{tabular}
\end{center}

\nB Dimension reduction in matrices is classically done by
truncation of the Singular Value Decomposition, which yields
its best approximation by another matrix of prescribed
rank. The notion of rank of a matrix has been extended to
different notions of ranks for tensors, like
Candecomp-Parafac, Tucker, and Tensor-Train; and dimension
reduction often comes down to an approximation of a tensor
by a tensor of lower rank. An ultimate goal when possible is
an approximation by a tensor of rank one, which leads to
separation of variables, which is a key issue e.g.~for
solving analytically some PDEs like $\Delta\,u=\lambda u$
when $u(x,y)=f(x)g(y)$. This shows that both issues of
dimension reduction and disentangling interactions are
linked.

\nB Let us come to some algebra, the domain in which the
concept of the rank of a tensor belongs. This will be the
common thread for the implementation of these ideas. Given a
family of vector spaces, say $(E,F,G)$, a \emph{tensor
  product} is an algebraic operator defined between them and
their elements, usually denoted $\otimes$. Vector spaces can
be finite or infinite dimensional. One can define
$E \otimes E$, $F \otimes G$, $E \otimes F \otimes G$,
etc...  and $\xbold \otimes \xbold'$,
$\ybold \otimes \zbold$,
$\xbold \otimes \ybold \otimes \zbold$ if
$\xbold, \xbold' \in E$, $\ybold \in F$ and $\zbold \in
G$. Without entering here into details, it is possible to
define a tensor algebra built from a family of vector
spaces. An algebra is a vector space endowed with a product,
here $\otimes$, that is associative
\[
\xbold \otimes (\ybold \otimes \zbold) = (\xbold \otimes \ybold) \otimes \zbold = \xbold \otimes \ybold \otimes \zbold  
\]
and distributive over $+$
\[
\xbold \otimes (\ybold + \ybold') = (\xbold \otimes \ybold) + (\xbold \otimes \ybold') 
\]
For example, the space of square matrices of dimension $n$
is an algebra with matrix product, because $A(BC)=(AB)C=ABC$
and $A(B+B')=AB+AB'$.

\nB Tensor products express the separation of variables,
which is better seen when vector spaces are spaces of
multivariate functions. If $f,g$ are two univariate
functions, $f \otimes g$ is a bivariate function with
$(f \otimes g)(x,y)=f(x)g(y)$. On says that $f \otimes g$
has rank one. If functions are defined not on $\R$, but on
finite set $\llbracket 1,m \rrbracket$, we have
$h_{ij}= f_ig_j$, or, shifting notations,
$(\xbold \otimes \ybold)_{ij} = x_iy_j$, which is a
definition of the tensor product of vectors
$(\xbold,\ybold)$: a matrix denoted $\xbold \otimes
\ybold$. But not all bivariate functions are in this form,
like $h(x,y)=x^2+2xy+y^2$; $h$ is, however, the sum of 3
functions of rank one. The smallest number $r$ (here $r=3$)
for which such a sum exists is the rank of $h$. This extends
to tensor products of vectors: any matrix is a linear
combination of rank one matrices. In many application areas,
being of low rank offers some advantages.

\nB First, it is expected that operations on tensors of low
rank have a lower complexity than operations on tensors of
large rank. This is a driving force for approximating, as
much as possible, a generic tensor by a tensor of low
rank. For example, in matrix numerical calculus, if $A,B$
are each of dimension $n \times n$, computing their product
$C=AB$ requires $n^3$ multiplications. However, if both have
rank $r$, it can be done with $n^2r$ operations only which
makes a significant improvement when $r \ll n$.

\nB Second, being of low rank spares much storage
capacity. Let us assume that the reader lives on Sirius,
to which communication is very expensive, and that I
want to convey to you the elements of the matrix $A = \xbold
\otimes \ybold$; i.e. $a_{ij}=x_iy_j$ for any $(i,j)$. The
cost for sending $A$ is $n^2$ (the number of elements in
$A$, and sending $(\xbold,\ybold)$ is $2n$. So, I decide to
send the two elements $(\xbold,\ybold)$ and the program to
build $A$ from them:\\
\begin{algorithm}[H]
\begin{algorithmic}[1]
\STATE \textbf{input} $\xbold, \ybold$ 
\FOR{$i \in \{1,n\}$}
\FOR{$j \in \{1,n\}$}
\STATE $a_{ij} = x_iy_j$
\ENDFOR
\ENDFOR
\RETURN $A=(a_{ij})_{ij}$
\end{algorithmic}
\end{algorithm}
\noindent If I send both vectors and the program to build
$A$ from them, I can send $A$ to Sirius with a low
cost. This can be quantified by the Kolmogorov complexity of
rank one matrix $A$, which is low.

\nB Approximation of a tensor by another one in a given
family requires a notion of \emph{distance}, which is
induced here by a norm. The choice of the norm can make
things easy, or difficult. Unless otherwise stated, the
selected norm is Frobenius or $\ell^2$ norm
$(\|\abold\|^2 = \sum_i a_i^2)$. The use of other norms is
evoked in chapter \ref{chap:inter}.

\nB An objective of this work is to present, in a consistent
way, several possibilities for approximating a tensor by a
tensor of lower prescribed rank (or precision). The main
possibilities rely on 
\begin{itemize}[label=$\rightarrow$]
    \item Candecomp-Parafac decomposition (CP)
    \item Tucker model with its two guises: HOSVD and HOOI
    \item Tensor-Train decomposition
\end{itemize}
This book is organized in several parts, of unequal length,
the logic of which is as follows:
\begin{itemize}
    \item The key part is part \ref{part:ranks} which
      presents simultaneously the different ranks, each an
      extension of the rank of a matrix to tensors with more
      than two modes and different norms which are different
      ways to compute a distance between two tensors. This
      part presents furthermore a solution to best low rank
      approximation for matrices: the Singular Value
      Decomposition and Eckart-Young theorem.
    \item The operational part is part \ref{part:blra} which
      presents different algorithms for best low rank
      approximation of a tensor, one per rank.
    \item These algorithms (and their understanding) rely
      all on a combination of elementary operations on
      tensors, which are presented in part \ref{part:primer}.
    \item The reason why "real world  tensors" are of low
      rank, much lower than what is expected for random
      tensors, still is mysterious. However, it is possible
      to show that tensors that are discretizations of
      smooth multivariate functions on a Cartesian grid have
      a low rank, and this is presented in part
      \ref{part:funct} .
    \item And, of course, the definition of a tensor is
      given at the very beginning of this story, in part
      \ref{part:tensors}.
\end{itemize}

\nB The reader is expected to have basic notions in linear
algebra, but basic notions only. Every time a more elaborate
notion is required, it will be recalled, like SVD recalled
in section \ref{sec:inter:svd}. I have avoided as much as
possible the use of too high levels of abstraction, knowing
that a minimum of abstraction is required because it permits
often to link notions which, beforehand, seemed to be
separate.

\nB It is expected that chapters can be read independently,
even if there is a sort of progression throughout the book,
provided the notations at the beginning of the book, which
will be used throughout it, have been understood, as well as
chapter \ref{chap:dmodestens}.

\paragraph{Acknowledgements:} It is a pleasure to thank the
members of the "working group" on tensors at Inria Bordeaux
Sud-Ouest with whom this project has been shared. Several
parts of this work have been presented piece wise and
discussed in the framework of this group, especially the one
on tensor-trains in October 2020. Many thanks to
Mohammed-Anwar Abouabdallah, Olivier Coulaud, Luc Giraud,
Martina Iannacito, Oguz Kaya, Nathalie Peyrard and Jean-René
Poirier who have read some parts of it or with whom I have
had many discussions over some years. Chapter
\ref{chap:quant} on tensors as discretization of
multivariate functions owes a lot to discussions with OC, LG
and JRP.

\nB It is obligatory that this text contains errors, I hope
simply of absent-mindedness. Please send to
\href{mailto:alain.franc@inrae.fr}{alain.franc@inrae.fr} any
criticisms, suggestions, or corrections.

%% file: tensor_product.tex
The (undeserved) reputation of tensors is that they are difficult to understand and master. The reason is probably multiple. First, the understanding of a multi-way array is very intuitive, but tensors are as well objects in tensor algebra, which is richer than linear algebra. Second, in many applications, tensors appear as arrays of coordinates with many indices, which are difficult to grasp together. The aim of this part is to make the concept of tensor accessible, because it is accessible and does not require concepts in very abstract algebra, while remaining rigorous. However, a minimum of incursion in abstract algebra is necessary, and useful. This will be done step by step.

\nB Hence, this part is at the same time the most intuitive and the most abstract of the whole book. Its objective is to give a definition of tensors (there can be several of them, all equivalent) with a sound algebraic basis. For this, I have followed \cite{Charles1984, Deheuvels1993, Lang1984}. But abstract algebra is not very intuitive, although it can illuminate and connect archipelagos of apparently different applied domains. Hence, I have tried to propose a quiet walk based on intuition for a progressive access to the rigorous definition of a tensor as a universal property. Much of the material used in the first chapter (chapter \ref{chap:tens_prod}) will not be used in the rest of the book or, if it is useful, will be recalled. \\
\\
It is organized as follows:
\begin{description}
\item[Chapter \ref{chap:tens_prod}] A definition of a tensor as a multilinear map is given, with smooth progression, from motivation and intuition to setting rigorous definitions. 
\item[Chapter \ref{chap:matelem}] To get accustomed to the notion of tensor if it is "exotic" to the reader, main elementary notions for matrices which will be extended to tensors are rewritten with tensor notations.
\item[Chapter \ref{chap:dmodestens}] Finally, notations which will be used throughout the book for $d-$modes tensors are presented. The chapter is short, because any required notion has been explained beforehand. 
\end{description}

%
\chapter{Tensors, tensor spaces, and tensor product}\label{chap:tens_prod}
%

What are tensors useful for in applied mathematics?

\nB Probably the most concise answer is: for calculations on multidimensional arrays avoiding the burden of tracking multi-indices. What are such calculations useful for? A multidimensional array, be it in numerical calculus for solutions of PDEs, statistics, statistical physics, quantum systems, represents often a set of entangled variables. Because of these dependencies, very few systems are analytically solvable, and when they are, it is often due to a possibility to separate the variables. Variational approach is a way to approximate as accurately as possible a complex system with entangled variables by another system with separate variables, as in mean-field approximation in statistical physics. Such an approximation is made much easier to derive and understand with tensors.

\nB This requires a detour through algebra, not necessarily very abstract, in order to derive a toolbox of concise formulas. An example of such formulas for matrices (which are tensors as well) is given in equation (\ref{eq:elemat}), section \ref{sec:matelem:useful}. This is a classical and basic approach in matrix calculation. Using the matrix $\times$ vector product $\xbold = A.\ybold$ instead of component-wise calculation $x_i = \sum_j\, a_{ij}y_j$ where $\xbold \in \R^m$, $\ybold \in \R^n$ and $A \in \R^{m \times n}$ is a simple and familiar example ($\R^{m \times n}$ is the space of real matrices with $m$ rows and $n$ columns, see section \ref{sec:tens_prod:notations} for notations). Associativity of matrix product
\[
A(BC)=(AB)C=ABC
\]
is better understood with such a compact notation than with coordinates
\[
\sum_j\, a_{ij}\left(\sum_k\, b_{jk}c_{k\ell}\right)= \sum_k\left(\sum_j\, a_{ij}b_{jk}\right)c_{k\ell} = \sum_{j,k}\, a_{ij}b_{jk}c_{k\ell}
\]

\nB As well as there are vectors and vector spaces, they are tensors and tensor spaces. There is as well a product, called tensor product, denoted $\otimes$ which acts between tensors (a tensor product of two tensors is a tensor), and on tensor spaces (a tensor product of two tensor spaces is a tensor space). If $\xbold \in E$ and $\ybold \in F$ are vectors, $\xbold \otimes \ybold \in E \otimes F$.

\nB This chapter presents basic definitions of a tensor, of tensor space, of tensor product of vectors and of vector spaces, as well as it introduces main notations. It is organized as follows:
\begin{description}
\item[section \ref{sec:tens_prod:detour}] As a motivation, a concrete example is given for which such a detour through algebra is useful.
\item[section \ref{sec:tens_prod:nutshell}] A gentle and intuitive introduction to what a tensor is, written in a loose way.
\item[section \ref{sec:tens_prod:def}] A rigorous definition of a tensor, of tensor product, and of tensor spaces.
\item [section \ref{sec:tens_prod:cocontra}] A presentation of covariant and contravariant vectors and tensors, which will not be used further in this book but is standard in other domains.
\item[section \ref{sec:tens_prod:notations}] A presentation of main notations which will be used throughout this book.
\item[section \ref{sec:tens_prod:otherdefs}] Some other equivalent definitions of a tensor product, be it between vectors or vector spaces, are given. It is useful to consider the definition of a tensor product of vector spaces as a universal property, even if more abstract.
\end{description}

%
\section{Why a detour through algebra?}\label{sec:tens_prod:detour}
%

As tensor notions have not yet been presented, some keys will be given along the presentation of the example. 

\nS Let us accept that a tensor is a multidimensional array (this will be stated more precisely in section \ref{sec:tens_prod:nutshell}). Let us denote by $\Abold$ a $4-$ways tensor given by its coefficients
\begin{equation*}
    \Abold = (a_{ijk\ell})_{i,j,k,\ell} \quad \mbox{with} \quad 
    \begin{cases}
     i & \in \llbracket 1,m\rrbracket \\
     j & \in \llbracket 1,n\rrbracket \\
     k & \in \llbracket 1,p\rrbracket \\
     \ell & \in \llbracket 1,q\rrbracket \\
    \end{cases}
\end{equation*}

\nT{Notation} In the same way that for a vector $\xbold=\sum_ix_i\ebold_i$, one denotes the set of coefficients by $(x_i)_i$ and each component by
\begin{equation*}
    x_i \quad \mbox{or} \quad \xbold[i]
\end{equation*}
the set of coefficients for a multi-array is denoted $(a_{ijk\ell})_{i,j,k,\ell}$ and each component by
\begin{equation*}
      a_{ijk\ell} \quad \mbox{or} \quad \Abold[i,j,k,\ell] \quad \mbox{or sometimes} \quad \Abold_{ijk\ell}  
\end{equation*}

\nS It can be "seen" as a four dimensional box, with edges of length respectively $m,n,p$ and $q$. So, $\Abold$ has $mnpq$ elements. We shall write
\begin{equation*}
    \Abold \in \R^{m \times n \times \times q}
\end{equation*}
Let us accept that if
\begin{equation*}
    \xbold \in \R^m, \qquad \ybold \in \R^n, \qquad \zbold \in \R^p, \qquad \tbold \in \R^q
\end{equation*}
the tensor product $\xbold \otimes \ybold \otimes \zbold \otimes \tbold$ is defined as the $4-$ways tensor $\Xbold \in \R^{m \times n \times p \times q}$ with
\begin{equation*}
    \Xbold_{ijk\ell} = x_iy_jz_kt_\ell
\end{equation*}
$\Xbold$ is called an elementary tensor, or a rank one tensor. Such tensors are essential, because the variables $x_i,y_j,z_k,t_\ell$ are separate. A question of interest is solving the following optimisation problem:\\
\begin{center}
    \ovalbox{
    \begin{tabular}{ll}
       given & a tensor $\Abold \in \R^{m \times n \times p \times q}$ \\
       find & an elementary tensor $\Tbold = \xbold \otimes \ybold \otimes \zbold \otimes \tbold$ \\
       with & $\xbold \in \R^m, \: \ybold \in \R^n, \: \zbold \in \R^p, \: \tbold \in \R^q$\\
       such that & $\|\Abold - \Tbold\|$ is minimal
    \end{tabular}
    }
\end{center}
where $\|.\|$ is the Frobenius norm extended naturally to tensors as $\|\Abold\|^2 = \sum_{ijk\ell} \, a_{ijk\ell}^2$. If $\Abold$ describes a joint law between 4 discrete variables, $\Tbold$ is a joint law of 4 independent variables, each taking respectively $m,n,p$ and $q$ values. I now show how to build a solution without and with notions and notations specific to tensor algebra. Let us note first that the problem is overparametrized, because if $(\xbold,\ybold,\zbold,\tbold)$ is a solution, then $(\alpha\xbold, \beta\ybold, \gamma\zbold, \delta\tbold)$  is another solution provided $\alpha\beta\gamma\delta=1$. This can be solved by setting
\begin{equation*}
    \Tbold = \alpha \, \xbold \otimes \ybold \otimes \zbold \otimes \tbold \qquad \mbox{with} \quad \|\xbold\|=\|\ybold\|=\|\zbold\|=\|\tbold\|=1, \quad \alpha \in \R
\end{equation*}

\nT{Setting the problem} Let us denote
\begin{equation*}
\begin{array}{lcl}
    \Delta &=& \|\Abold - \Tbold\|^2 \\
    &=& \displaystyle \sum_{i,j,k,\ell}\left(a_{ijk\ell}-\alpha\,x_iy_jz_kt_\ell\right)^2
\end{array}
\end{equation*}
We then have to solve
\begin{center}
    \ovalbox{
    \begin{tabular}{ll}
       given & a multiarray $(a_{ijk\ell})_{i,j,k,\ell}$ \\
       find &  $\alpha, \xbold, \ybold, \zbold, \tbold$ \\
       with & $\|\xbold\|^2=\|\ybold\|^2=\|\zbold\|^2=\|\tbold\|^2=1$ \\
       such that & $\Delta$ is minimal
    \end{tabular}
    }
\end{center}


\nT{Solution with indices} We have
\begin{equation*}
    \Delta = \Delta(\alpha, x_1, \ldots,x_m,y_1,\ldots,y_n,z_1,\ldots,z_p,t_1,\ldots,t_q)
\end{equation*}
with
\begin{equation*}
    \Delta = \sum_{i,j,k,\ell}\, a_{ijk\ell}^2 + \alpha^2\, \sum_{i,j,k,\ell}\, x_i^2y_j^2z_k^2t_\ell^2 - 2\alpha\,\sum_{i,j,k,\ell}\, a_{ijk\ell}\,x_iy_jz_jt_\ell 
\end{equation*}
So, for example for index $i$
\begin{equation*}
    \forall \: i, \quad \frac{\partial \Delta}{\partial x_i} = 2\alpha^2\,\left(\sum_{j,k,\ell}\, y_j^2z_k^2t_\ell^2\right)\, x_i - 2 \, \alpha\,\sum_{j,k,\ell}\,a_{ijk\ell}\, y_jz_kt_\ell
\end{equation*}
Using Lagrange multpliers for condition $\sum_ix_i^2=1$ yields
\begin{equation*}
    2\alpha^2\,\sum_{j,k,\ell}\, y_j^2z_k^2t_\ell^2\, x_i - 2 \, \alpha\,\sum_{j,k,\ell}\,a_{ijk\ell}\, y_jz_kt_\ell -2 \lambda_x\, x_i=0
\end{equation*}
or
\begin{equation*}
    x_i = \frac
    {\alpha\, \sum_{j,k,\ell}\,a_{ijk\ell} \, y_jz_kt_\ell}
    {\alpha^2\,\sum_{j,k,\ell}y_j^2z_k^2t_\ell^2 - \lambda_x}
\end{equation*}
Let us now recall that $\sum_jy_j^2=\sum_kz_k^2=\sum_\ell t_\ell^2=1$. Then $\sum_{j,k,\ell}y_j^2z_k^2t_\ell^2=1$, and
\begin{equation}\label{eq:tens_prod:detour:xi}
    x_i = \frac
    {\alpha\,\sum_{j,k,\ell}\,a_{ijk\ell} \, y_jz_kt_\ell}
    {\alpha^2- \lambda_x}
\end{equation}

\nB Similar equations can be derived for indices $j,k,\ell$ yielding a fixed point solution
\begin{equation*}
    \left\{
    \begin{array}{lcl}
        \xbold &=& f_x(\ybold,\zbold,\tbold) \\
        \ybold &=& f_y(\xbold,\zbold,\tbold) \\
        \zbold &=& f_z(\xbold,\ybold,\tbold) \\
        \tbold &=& f_t(\xbold,\ybold,\zbold) \\
    \end{array}
    \right.
\end{equation*}

\nB Finally, this set is completed by setting $\partial \Delta/\partial \alpha=0$, or
\begin{equation*}
    \begin{array}{lcl}
         \displaystyle \frac{\partial \Delta}{\partial \alpha} &=& \displaystyle 2\alpha \left(\sum_{i,j,k,\ell}\, x_i^2y_j^2z_k^2t_\ell^2\right) - 2 \left(\sum_{i,j,k,\ell}\,a_{ijk\ell}\, x_iy_jz_kt_\ell\right) \\
         &=& 0
    \end{array}
\end{equation*}
or
\begin{equation*}
    \alpha = \sum_{i,j,k,\ell}\,a_{ijk\ell}\, x_iy_jz_kt_\ell
\end{equation*}
because $\sum_{i,j,k,\ell}\, x_i^2y_j^2z_k^2t_\ell^2=\left(\sum_{i}\, x_i^2\right)\left(\sum_{j}\, y_j^2\right)\left(\sum_{k}\, z_k^2\right)\left(\sum_{\ell}\, t_\ell^2\right) = 1$.


\nT{Solution with tensor notations} We rewrite a path towards the solution by using more compact tensor notation, labelled by a letter and explained just below. We have
\begin{equation}\label{eq:tensor_prod:detour:formal}
\begin{array}{lclll}
\Delta &=& \|\Abold - \Tbold\|^2 &\\
&=& \|\Abold\|^2 + \|\Tbold\|^2 - 2 \langle \Abold \, , \, \Tbold\rangle & \qquad & (a) \\
&=&\|\Abold\|^2 + \alpha^2\, \|\xbold\|^2 \|\ybold\|^2  \|\zbold\|^2 \|\tbold\|^2 - 2 \alpha \, \langle \Abold \, , \, \xbold \otimes \ybold \otimes \zbold \otimes \tbold\rangle & \quad & (b)\\
&=& \|\Abold\|^2 + \alpha^2 \|\xbold\|^2 \|\ybold\|^2  \|\zbold\|^2 \|\tbold\|^2 - 2 \alpha \, \langle \Abold.(\ybold \otimes \zbold \otimes \tbold) \, , \,  \xbold\rangle & \qquad & (c) \\
&=& \|\Abold\|^2 + \alpha^2 \|\xbold\|^2 \|\ybold\|^2  \|\zbold\|^2 \|\tbold\|^2 - 2 \alpha \, \langle \bbold \, , \,  \xbold\rangle & \qquad & (d) 
\end{array}
\end{equation}
with following explanations:\\
\begin{description}
\item[$(a)$] the inner product and norms are extended to tensors, because the space of tensors is a vector space (tensors can be added and multiplied by a scalar)

\item[$(b)$]it can be shown that $\|\xbold \otimes \ybold \otimes \zbold \otimes \tbold\|^2= \|\xbold\|^2 \, \|\ybold\|^2 \, \|\zbold\|^2 \, \|\tbold\|^2$. Indeed
\begin{equation*}
    \sum_{i,j,k,\ell} x_i^2y_j^2z_k^2t_\ell^2 = \left(\sum_ix_i^2\right)\left(\sum_jy_j^2\right)\left(\sum_kz_k^2\right)\left(\sum_\ell t_\ell^2\right)
\end{equation*}

\item[$(c)$] This implies a key tensor operation: the contraction, presented in chapter \ref{chap:contract}. A simple example of contraction is the matrix $\times$ vector product: $\xbold=A.\ybold$, with $\xbold \in \R^m$, $A \in \R^{m \times n}$, $\ybold \in \R^m$, written component-wise:
\begin{equation*}
    x_i = \sum_j\, a_{ij}y_j, \qquad \mbox{with} \quad 
    \begin{cases}
      i & \in \llbracket 1,m\rrbracket \\
      j & \in \llbracket 1,n\rrbracket
    \end{cases}
\end{equation*}
We have $\langle A,\xbold \otimes \ybold\rangle=\langle A\ybold,\xbold\rangle$ which can be shown with coordinates by
\begin{equation*}
    \sum_{i,j}\, a_{ij}\, x_iy_j = \sum_i\left(\sum_j\, a_{ij}y_j\right)x_i
\end{equation*}
So, because $\otimes$ is associative, hence $\xbold \otimes \ybold \otimes \zbold \otimes \tbold=\xbold \otimes (\ybold \otimes \zbold \otimes \tbold)$, this implies 
\begin{equation*}
    \langle \Abold \, , \, \xbold \otimes \ybold \otimes \zbold \otimes \tbold\rangle =  
    \langle \Abold.(\ybold \otimes \zbold \otimes \tbold) \, , \, \xbold\rangle \\
\end{equation*}
which shows $(c)$.
\item[$(d)$] $\bbold=\Abold.(\ybold \otimes \zbold \otimes \tbold) \in \R^m$
\end{description}
If $f\: : \: \R^n \longrightarrow \R$, let us denote 
\begin{equation*}
    \nabla_\xbold f = \frac{\partial \;f}{\partial \; \xbold}
\end{equation*}
the gradient of $f$. Then
\begin{equation*}
    \nabla_\xbold\Delta = 2\alpha\left(\|\ybold\|^2  \|\zbold\|^2 \|\tbold\|^2 \right)\xbold -2\alpha\bbold
\end{equation*}
The gradient for constraint $\|\xbold\|^2=1$ is $2\xbold$. Hence
\begin{equation*}
    2\alpha^2\left(\|\ybold\|^2  \|\zbold\|^2 \|\tbold\|^2 \right)\xbold -2\alpha\bbold = 2\lambda \xbold
\end{equation*}
or
\begin{equation*}
    \xbold = \frac{\alpha\, \bbold}{\alpha^2\left(\|\ybold\|^2  \|\zbold\|^2 \|\tbold\|^2\right) -\lambda}
\end{equation*}
Knowing that $\|\ybold\|=\|\zbold\|=\|\tbold\|=1$ yields
\begin{equation}
    \xbold = \frac{\alpha\, \bbold}{\alpha^2 -\lambda}
\end{equation}
which is equation (\ref{eq:tens_prod:detour:xi}) as
\begin{equation*}
    b_i = \sum_{j,k,\ell}\, a_{ijk\ell}\, y_jz_kt_\ell
\end{equation*} 

\nB The solution for $\ybold, \zbold,\tbold$ is given with a similar calculation after applying a same permutation to $\Abold$ and $\Tbold$ (see section \ref{sec:tenselem:permut}).

\nB Finally, writing $\partial \Delta / \partial \alpha=0$ in (\ref{eq:tensor_prod:detour:formal}(d)) yields $\alpha^2\,\|x\|^2\|y\|^2\|z\|^2\|t\|^2=\langle \bbold,\xbold\rangle$ or
\begin{equation*}
    \alpha = \langle \bbold,\xbold\rangle
\end{equation*}

\nT{Extension to tensors with many modes} The mode of a tensor is, in this example, the number of indices. Here, the mode is 4. Usefulness of tensor notations becomes clearer when the number of modes increases. Let us consider a tensor with $d$ modes, with $d \in \N$. It can be written
\begin{equation*}
    \Abold = \left(a_{i_1\ldots i_d}\right)_{i_1,\ldots,i_d} \qquad \mbox{with} \quad i_\mu \in \llbracket 1,n_\mu\rrbracket
\end{equation*}
Let us define
\begin{equation*}
    \Tbold = \alpha \, \xbold_1 \otimes \ldots \otimes \xbold_d
\end{equation*}
which can be written
\begin{equation*}
    \Tbold = \alpha \, \bigotimes_{\mu=1}^d\, \xbold_\mu
\end{equation*}
We then have, extending (\ref{eq:tensor_prod:detour:formal})
\begin{equation*}
    \begin{array}{lcl}
        \Delta &=& \|\Abold - \Tbold\|^2 \\
        &=& \|\Abold\|^2 - 2 \langle \Abold \, , \, \Tbold\rangle + \|\Tbold\|^2 \\
        &=& \displaystyle \left\| \Abold\right\| + \alpha^2 \prod_{\mu=1}^d\, \|\xbold_\mu\|^2 - 2 \alpha \, \left\langle \Abold \, , \, \bigotimes_{\mu=1}^d\, \xbold_\mu\right\rangle \\
         &=& \displaystyle \left\| \Abold\right\| + \alpha^2 \prod_\mu\, \|\xbold_\mu\|^2 - 2 \alpha \, \left\langle \Abold. \left(\bigotimes_{\mu>1}\, \xbold_\mu\right)\, , \, \xbold_1\right\rangle \\
         &=& \displaystyle \left\| \Abold\right\| + \alpha^2 \prod_\mu\, \|\xbold_\mu\|^2 - 2 \alpha \, \langle \bbold \, , \, \xbold_1\rangle \\
    \end{array}
\end{equation*}
where $\bbold$ does not depend on $\xbold_1$. Then
\begin{equation*}
    \nabla_{\xbold_1}\Delta = 2\alpha\, \left(\prod_{\mu>1}\, \|\xbold_\mu\|^2\right)\xbold_1 - 2\alpha\, \bbold
\end{equation*}
Including the constrain through Lagrange multipliers and acknowledging that $\|\xbold_\mu\|=1$ for $\mu > 1$ yields
\begin{equation*}
    \xbold_1 = \frac{\alpha\, \bbold}{\alpha^2 - \lambda} \qquad \mbox{with} \quad \bbold = \Abold. \left(\bigotimes_{\mu>1}\, \xbold_\mu\right)
\end{equation*}

\nS It is possible to develop such a calculation with indices only, but it would be more tricky. Indices of indices would be required.

%
\section{Tensors in a nutshell}\label{sec:tens_prod:nutshell}
%

Tensors are often presented as multidimensional arrays. But this is a limited approach\marginpar{\dbend}, as if a vector would be an array of coordinates only. For vectors, an array of coordinate is the expression of a vector in a given basis, and a vector is an element of a vector space, which has a structure. The emergence of the structure of a vector space from calculations on arrays has coincided with the rise of abstract algebra. Many operations based on the structure of a vector space are simplified by coordinate free notations, like the definition of an eigenvector and an eigenvalue: $L\ubold=\lambda\ubold$. In linear algebra, one distinguishes a linear map, $L \in \L(E,F)$ from a vector space $E$ on a vector space $F$, from its expression as a matrix $A \in \R^{m \times n}$ once a basis has been selected in $E$ and in $F$. Similarly, a multidimensional array is the expression of a tensor in a given basis, and many operations based on the structure of a tensor space are simplified by coordinate free notations. One distinguishes a tensor $\Abold$, which is a multilinear form on a product of vector spaces, from its expression as a multiway array once a basis has been selected for each vector space on which the form operates. The multiway set of coordinates of a tensor in called an \kw{hypermatrix}. 

\nB The set of multidimensional arrays of given dimensions is a vector space. A multidimensional array is the expression of the tensor in a given basis. A  tensor being given, if the basis for its expression with coordinates changes, the array changes. The array is a representation of a tensor, not the tensor itself. This is illustrated by the following diagram\\
\\
\begin{center}
\begin{tikzcd}
& \llbracket \Tbold \rrbracket_{\EC} \arrow[dd, "(c)"]\\
\Tbold \arrow[ur, "(a)"] \arrow[dr, "(b)"]& \\
& \llbracket \Tbold\rrbracket_{\EC'}
\end{tikzcd}
\end{center}
where $\EC,\EC'$ are two different basis of a tensor space, with
\begin{description}
\item[$(a)$] expression of tensor $\Tbold$ in basis $\EC$
\item[$(b)$] expression of tensor $\Tbold$ in basis $\EC'$
\item[($c$)] basis change
\end{description}

\nT{Connecting numerical calculations on multiway arrays and tensor algebra} This simple observation is a basic motivation for this book: connecting numerical calculations on multidimensional arrays with tensor algebra. This implies to develop tools to work with tensors rather than with arrays, and to dialogue between tensors and arrays. Many notions which do not depend on the basis selected for the expression of a tensor as an array will be made simple. Conversely, some tensors will be better expressed in a specific basis for some optimal properties to be simple to detect. Such a choice is a very natural extension of similar approaches with matrices: rather than opposing numerical calculations and algebra, both are intermingled to develop matrix calculus. Such a calculus has a facet as formal calculation, like $A(B+C)=AB + AC$, and another one as numerical calculation, like algorithms to compute the inverse of a given matrix. For example, the notion of eigenvector or eigenvalue is an intrinsic property of a matrix, independent on the basis in which the matrix is expressed, whereas the expression of a matrix as diagonal in a basis of eigenvectors permits to simplify dramatically many problems. Similarly, the decomposition of a matrix $A$ into products of simpler matrices, like $LU$ decomposition, Cholesky decomposition, or $QR$ decomposition, is a key ingredient in numerical matrix calculation. All those matrix decomposition dramatically simplify numerical calculations for inverting a matrix, computing its determinant, etc ... One objective of this book is to intermingle formal and numerical calculations for tensors.

\nS Let us consider a representation of a tensor $\Tbold$  in a given basis as a $d-$dimensional array, with coordinates
\begin{equation}
 \Tbold_{i_1\ldots i_d} \in \K
\end{equation}
(here and elsewhere, $\K$ is a field, usually $\R$ or $\C$) where each $i_\mu$ with $\mu \in \llbracket 1,d \rrbracket$ is an index along a given axis
\begin{equation}
 i_\mu \in \llbracket 1,n_\mu \rrbracket 
\end{equation}
Let us denote
\begin{equation*}
    I_\mu = \llbracket 1, n_\mu \rrbracket
\end{equation*}
with
\begin{equation}
    i_\mu \in I_\mu
\end{equation}
The expression of $\Tbold$ can be given using a multi-index notation
\begin{equation}
 \ibold = (i_1, \ldots,i_d) 
\end{equation}
If
\begin{equation}
 \IC = I_1 \times \ldots, \times I_d, 
\end{equation}
then
\begin{equation}
    \ibold \in \IC
\end{equation}
The elements of $\Tbold$ are indexed by $\ibold$, and the representation of  $\Tbold$ is a map 
\begin{equation}
 \begin{CD}
   \IC @>>> \K
 \end{CD}
\end{equation}
where usually $\K$ is $\R$ or $\C$. As this concerns a representation of a tensor $\Tbold$, the map depends on the basis. 

\nS The number $d$ of entries of a tensor is called its \kw{order} (see notes on page \pageref{page:order} for the origin of this word). A vector is a tensor of order one. A matrix is a tensor of order two. 

\nS Let $E,F$ be finite dimensional vector spaces of dimension $m,n$ respectively. Let 
\begin{equation}
 \xbold = (x_1, \ldots,x_m) \in E, \qquad \ybold = (y_1,\ldots,y_n) \in F
\end{equation}
Loosely speaking, one can define the tensor product $x \otimes y$ as the matrix in $\R^{m \times n}$ defined by
\begin{equation}
 (\xbold \otimes \ybold)_{i,j} = x_iy_j \qquad \mbox{with} \quad 
 \begin{cases}
  i &\in \llbracket 1,m\rrbracket \\
  j &\in \llbracket 1,n\rrbracket \\
 \end{cases}
\end{equation}
This tensor product of two vectors is called a tensor. In such a case, it is called an elementary tensor. This can be extended to more than two vectors, defining elementary tensors (or rank one tensors) of order 3 or more.  Let
\begin{equation}
    \zbold = (z_1,\ldots, z_p)
\end{equation}
Then
\begin{equation}
(\xbold \otimes \ybold \otimes \zbold)_{i,j,k} = x_iy_jz_k \qquad \mbox{with} \quad 
\begin{cases}
 i &\in \llbracket 1,m\rrbracket \\
  j &\in \llbracket 1,n\rrbracket \\
  p &\in \llbracket 1,p\rrbracket \\
 \end{cases}
\end{equation}
and more generally with a shift in notations
\begin{equation}
 (\xbold_1 \otimes \xbold_2 \otimes \ldots \otimes \xbold_d)_{i_1i_2\ldots i_d} = x_{1,i_1}x_{2,i_2} \ldots x_{d,i_d} \qquad \mbox{with} \quad 
 \left\{
 \begin{array}{lll}
  & \mu &\in \llbracket 1,d \rrbracket \\
  \forall \: \mu, & i_\mu & \in \llbracket 1,n_\mu \rrbracket
 \end{array}
 \right.
\end{equation}
or even more compactly
\begin{equation}
 \left(\bigotimes_{\mu=1}^d \xbold_\mu\right)_{i_1,\ldots,i_d} = \prod_{\mu=1}^d \, x_{i_\mu}^{(\mu)}
\end{equation}

\nS Such an elementary tensor is given as a multidimensional array, i.e. is known by a representation in a given basis. I indicate here loosely how definitions of general tensors can be built from these elementary tensors. This is one way among several others to define a tensor. A rigorous definition of a tensor will be given in section \ref{sec:tens_prod:def}. The approach here leads to: a tensor can be defined as a linear combination of elementary tensors. It can be shown that, for finite dimensional vector spaces, such a construction spans the whole space of multidimensional arrays, i.e., for any tensor $\Tbold$, there is an integer $r$ and $r$ elementary tensors $\Tbold_a$ with
\begin{equation}
 \Tbold = \sum_{a=1}^r\,\Tbold_a 
\end{equation}
with
\begin{equation}
 \Tbold_a = \xbold_1^{(a)} \otimes \ldots \otimes \xbold_d^{(a)}, \qquad \mbox{with} \quad
 \begin{cases}
   \xbold_\mu^{(a)} \in E_\mu \\
   1 \leq \mu \leq d
 \end{cases}
\end{equation}
The integer $d$ is called the \kw{order} of the tensor $\Tbold$. If $\dim E_\mu=n_\mu$, the set $(n_1,\ldots,n_d)$ is the set of \kw{dimensions} of the tensor. A vector space $E_\mu$ is called a \kw{mode} of the tensor. The smallest possible integer $r$ for which such a decomposition holds is called the \kw{rank} of the tensor, or CP-rank. In the same way that a matrix can be approximated at best by another matrix of prescribed low rank, a key issue in numerical tensor algebra is to find a best approximation of a given tensor by a tensor of prescribed low rank.

\nT{What will be shown} It is elementary to see that the set of tensors of a given order with given dimensions is a vector space (they can be added and multiplied by a scalar). Basis spanning this vector space can be built from basis spanning each mode. Some elementary operations like permutation, matricization, slicing etc. can be defined leading to nice and elegant algebraic formula, coordinate free, which considerably ease the calculations in some situations. A tensor product can be defined between tensors of any order and dimensions. The contraction is a key operation which can be defined as the dual of the tensor product provided the modes are endowed with an inner product (the matrix $\times$ vector product is an example of a contraction, through $(\xbold \otimes \ybold).\abold = \langle \abold,\ybold\rangle \, \xbold$). It would be impossible to understand (and prove well-posedness or ill-posedness of) classical best low rank approximations of tensors without this algebraic toolbox.

%
\section{Definition of tensors, tensor spaces and tensor products}\label{sec:tens_prod:def}
%

A tensor is a multilinear form on a family of vector spaces.

\nB This key section is organised as follows to develop these notions:
\begin{enumerate}
\item First, some preliminaries are recalled
\begin{enumerate}
\item the notion of dual space which will be used 
\item the definition of a multilinear form on a product of vector spaces, some examples of its expression on a basis for linear, bilinear and trilinear forms are given
\end{enumerate}
\item Second, the definitions of tensor notions are given: 
\begin{enumerate}
\item the definition of a tensor and of a tensor space 
\item the tensor product between two multilinear forms or tensors
\item which induces a tensor product between tensor spaces.
\end{enumerate}
\item Finally, the tensor product between vectors is defined by two ways:
\begin{enumerate}
\item elementary tensors are defined as tensor product between linear forms; then the canonical isomorphism between a vector space and its bidual is used (a vector is a linear form acting on linear forms)
\item in a less abstract way, when an inner product is defined on vector spaces, it can be used to define an isomorphism between a vector space and its dual, which permits to extend the tensor product between linear forms to the tensor product between vectors. But this is not canonical, and depends on the inner product. So, it is defined up to an isomorphism only.
\end{enumerate}
\end{enumerate}

%
\subsection*{Preliminaries}
%

\nT{Linear form and dual space} In this section, I will use the notion of dual of a vector space. If $E$ is a vector space on a field $\K$, a \kw{linear form}\index{form!linear} $\varphi$ on $E$ is a map 
\[
 \begin{CD}
  E @>\varphi>>\K
 \end{CD}
\]
such that
\[
 \left\{
    \begin{array}{lcl}
      \varphi(\xbold+\ybold) &=& \varphi(\xbold)+\varphi(\ybold) \\
      \varphi(\lambda \xbold) &=& \lambda \varphi(\xbold)
    \end{array}
 \right.
\]
The \kw{dual} of $E$ is the vector space of the linear forms on $E$. It is classically denoted $E^*$.


\nS \textbf{$d-$linear forms:} Let us denote by $E$ a vector space on a field $\K$. A \kw{$d-$linear form} $T$ on $E$ is a map 
\[
 \begin{CD}
  \underbrace{E \times \ldots \times E}_{d \: \mathrm{times}} @>T>> \R
 \end{CD}
\]
which is linear on each of its components, i.e.
\begin{multline}
 \forall \: \mu, \quad T(\xbold_1, \ldots, \xbold_{\mu-1}, \lambda \xbold_\mu + \lambda'\xbold_\mu', \xbold_{\mu+1}, \ldots,\xbold_d) = \lambda T(\xbold_1, \ldots, \xbold_{\mu-1}, \xbold_\mu , \xbold_{\mu+1}, \ldots,\xbold_d) \\
 + \lambda' T(\xbold_1, \ldots, \xbold_{\mu-1}, \xbold_\mu', \xbold_{\mu+1}, \ldots,\xbold_d)
\end{multline}
(here, $\xbold_\mu \in E)$. It is easy to check that $d-$linear forms on $E$ form a vector space, of dimension $n^d$ if $n= \dim E$. It is denoted $\T^d(E)$. By convention if $d=0$, $\T^0(E)=\K$.


\nT{Expression in a basis} Let us select a basis for $E$, and let $\xbold = (x_1,\ldots,x_n)$ in this basis, $\ybold = (y_1,\ldots,y_n)$ and $\zbold = (z_1,\ldots,z_n)$. A \kwnind{linear form}\index{form!linear} on $E$ is a map given by
\begin{equation}
 \begin{CD}
  \xbold @>>> \displaystyle \sum_i\alpha_ix_i
 \end{CD}
\end{equation}
and defined in the given basis by the coefficients $\bm{\alpha} = (\alpha_1,\ldots,\alpha_n)$. Any linear form can be written like this. A \kw{bilinear form}\index{form!bilinear} on $E$ is a map 
\begin{equation}\label{eq:tens_prod:def:1}
 \begin{CD}
  (\xbold,\ybold) @>>> \displaystyle \sum_{ij}\alpha_{ij}x_iy_j
 \end{CD}
\end{equation}
defined by the coefficients $\bm{\alpha} = (\alpha_{ij})_{i,j}$, and a \kw{trilinear form} is a  
\begin{equation}
 \begin{CD}
  (\xbold,\ybold,\zbold) @>>> \displaystyle \sum_{ijk}\alpha_{ijk}x_iy_jz_k
 \end{CD}
\end{equation}
defined by the coefficients $\bm{\alpha} = (\alpha_{ijk})_{i,j,k}$. More generally,  a \kw{$d-$linear form} can be seen as a homogeneous polynomial of degree one for each variable:
\begin{equation}
 \begin{CD}
  (\xbold_1, \ldots,\xbold_d) @>>> \displaystyle \sum_{i_1\ldots,i_d}\alpha_{i_1\ldots,i_d}x_{1i_1}\ldots x_{di_d}
 \end{CD}
\end{equation}
defined by the coefficients $\bm{\alpha} = (\alpha_{i_1\ldots i_d})_{i_1,\ldots, i_d}$. One can see that the number of indices grows with $d$, and tensor algebra enables to simplify notations. 

%
\subsection*{Definitions}
%

This is easily extended to more general situation of multilinear forms on families of vector spaces.

\nT{Tensor} Let us have a family 
\[
 E_1, \ldots, E_\mu, \ldots, E_d
\]
of vector spaces on a same field $\K$. A \kw{tensor} $\Tbold$ is a $d-$linear form on $E_1 \times \ldots \times E_d$:
\begin{equation}
    \begin{CD}
        E_1 \times \ldots \times E_d @>\Tbold>> \K
    \end{CD}
\end{equation}
or
\begin{equation}
    \Tbold \in \T(E_1, \ldots, E_d)
\end{equation}
$\T(E_1, \ldots, E_d)$ is the set of $d-$linear forms on  $E_1 \times \ldots \times E_d$. 


\nT{Tensor space}The space $\T(E_1, \ldots, E_d)$ is called a \kw{tensor space}. It will be denoted
\begin{equation}
    E_1 \otimes \ldots \otimes E_d := \T(E_1, \ldots, E_d)
\end{equation}
If all spaces $E_\mu$ are the same, I will denote
\begin{equation}
    \T^d(E) := \T\underbrace{(E \times \ldots \times E)}_{d \:\: \mathrm{times}}
\end{equation}
or even $\T^d$ when there is no ambiguity on space $E$, and
\begin{equation}
    E^{\otimes d} := \T^d(E).
\end{equation}


\nT{Tensor product of two $d-$linear forms} Let $d,d' \in \N$ be two integers (we can have $d=d'=0)$. Let 
\begin{equation}
 \left\{ 
    \begin{array}{llcll}
      \xbold_\mu & \in E_\mu & \qquad & \mbox{for} & \mu \in \{1,\ldots,d\} \\
      \xbold'_\mu & \in E'_\mu & \qquad & \mbox{for} & \mu \in \{1,\ldots,d'\} \\
    \end{array}
 \right.
\end{equation}
and $\Tbold \in \T^d$, $\Tbold' \in \T^{d'}$. The \kw{tensor product} $\otimes$ between $\Tbold$ and $\Tbold'$ 
\[
 \begin{CD}
  \T^d \times \T^{d'} @>\otimes>> \T^{d+d'}\\
  (\Tbold,\Tbold') @>>> \Tbold \otimes \Tbold'
 \end{CD}
\]
is the $(d+d')-$linear map $\Tbold \otimes \Tbold' \in \T^{d+d'}$ on $E_1 \times \ldots \times E_d \times E'_1 \times \ldots \times E'_{d'}$ defined by
\begin{equation}\label{eq:tp}
 (\Tbold\otimes \Tbold')(\bm{x},\bm{x}')=\Tbold(\bm{x}).\Tbold(\bm{x}')
\end{equation}
for any
\[
\left\{
  \begin{array}{lll}
    \bm{x} & = (\xbold_1,\ldots,\xbold_d) &\in E_1 \times \ldots \times E_d \\
    \bm{x}'&=(\xbold'_1,\ldots,\xbold'_{d'}) &\in E'_1 \times \ldots \times E'_{d'} 
  \end{array}
\right.
\]


\nT{Key observation} The tensor product $\otimes$ defines an operation which separates the variables $\bm{x}$ and $\bm{x}'$. To illustrate this, let $d=d'=1$, $\xbold = (x_1, \ldots,x_n)$, $\xbold'=(x'_1,\ldots,x'_n)$ and define 
\begin{equation}
 \Tbold(\xbold) = \sum_i\, \alpha_ix_i, \qquad \Tbold'(\xbold')=\sum_j\, \alpha'_jx'_j
\end{equation}
Then
\begin{equation}
 \begin{array}{lcl}
    (\Tbold \otimes \Tbold')(\xbold,\xbold') &=& \displaystyle \left(\sum_i\, \alpha_ix_i\right)\left(\sum_j\, \alpha'_jx'_j\right)\\
    &=& \displaystyle \sum_{i,j} \alpha_i\alpha'_j\, x_ix'_j \\
    &=& \displaystyle \sum_{i,j} \bm{\alpha}_{ij}x_ix'_j, \qquad \mbox{see equation (\ref{eq:tens_prod:def:1})}
 \end{array}
\end{equation}
which yields
\begin{equation}
 \bm{\alpha}_{ij} = \alpha_i\alpha'_j
\end{equation}
Let $\dim E=n$, and $\Tbold''$ be a tensor defined component-wise in a given basis by
\begin{equation}
 \Tbold''(\xbold,\xbold')= \sum_{i=1}^n \sum_{j=1}^n \bm{\alpha}_{ij}x_ix'_j
\end{equation}
The calculation of $\Tbold''(\xbold,\xbold')$ requires the summation of $n^2$ terms in $\sum_{i=1}^n\sum_{j=1}^n [\ldots]$. If $\Tbold'' = \Tbold \otimes \Tbold'$, then it requires $2n$ additions and one multiplication only. The gain will be dramatic for tensors of higher order.


\nS It is easy to show that
\begin{itemize}[label=$\rightarrow$]
 \item $\Tbold \otimes \Tbold'$ is a $(d+d')-$ linear form, hence belongs to $\T^{d+d'}$
 \item The map $(\Tbold,\Tbold') \longrightarrow \Tbold \otimes \Tbold'$ is bilinear, i.e.
\begin{equation}
 \Tbold \otimes (\alpha'\Tbold' + \alpha '' \Tbold'') = \alpha ' \Tbold \otimes \Tbold' + \alpha'' \Tbold \otimes \Tbold''
\end{equation}
and other equalities for bilinearity hold
\item The operation $\otimes$ is associative, i.e.
\begin{equation}
 \Tbold \otimes (\Tbold' \otimes \Tbold'') = (\Tbold \otimes \Tbold') \otimes \Tbold''
\end{equation}
or
\begin{equation}
 (\Tbold \otimes \Tbold' \otimes \Tbold'')(\xbold,\xbold',\xbold'') = \Tbold(\xbold)\Tbold'(\xbold')\Tbold''(\xbold'')
\end{equation}

\end{itemize}


\nT{Tensor product of tensor spaces} Let $\T(E_1, \ldots, E_d)$ be the set of $d-$linear forms on $E_1 \times \ldots \times E_{d}$, and $\T(E'_1, \ldots, E'_d)$ the set of $d'-$linear forms on $E'_1 \times \ldots \times E'_{d'}$. If $\Tbold \in \T(E_1, \ldots, E_d)$ and $\Tbold' \in \T(E'_1, \ldots, E'_d)$, the \kw{tensor product} $\Tbold \otimes \Tbold'$ belongs to the set of $d+d'-$linear forms on $E_1 \times \ldots \times E_d \times E'_1 \times \ldots \times E'_{d'}$, which is spanned by the elementary tensors of the form $(\Tbold\otimes \Tbold')(\bm{x},\bm{x}')=\Tbold(\bm{x}).\Tbold(\bm{x}')$ as in (\ref{eq:tp}). Or
\begin{equation}
\left.
\begin{array}{l}
\Tbold \in \T(E_1, \ldots, E_d) \\
\Tbold' \in \T(E'_1, \ldots, E'_d)
\end{array}
\right\}
\quad \Longrightarrow \quad 
\Tbold \otimes \Tbold' \in \T(E_1 \times \ldots \times E_d \times E'_1 \times \ldots \times E'_{d'})
\end{equation}
We define
\begin{equation}
    \T(E_1 \times \ldots \times E_d \times E'_1 \times \ldots \times E'_{d'}) = \T(E_1, \ldots, E_d) \otimes \T(E'_1, \ldots, E'_d)
\end{equation}
with the same symbol $\otimes$ for tensor product between tensor spaces. This justifies the notation
\begin{equation}
    \T(E_1, \ldots, E_d) := E_1 \otimes \ldots \otimes E_d
\end{equation}
by repeating this operation, which can be compacted when useful as
\begin{equation}
    \T(E_1, \ldots, E_d) := \bigotimes_{\mu=1}^d E_\mu
\end{equation}

\subsection*{Tensor product between vectors}

Tensor product has been defined between multilinear forms. Next step is to extend the notion of tensor product to vectors. This can be done in a straightforward way with a few words. As there is a canonical isomorphism between a vector space $E$ and its bidual $E^{**}$ (the dual of its dual), a vector $\abold \in E$ is an element of $(E^*)^*$ with $\abold(\varphi)=\varphi(\abold)$ if $\varphi \in E^*$. Hence, any vector $\abold$ is a linear form on $E^*$, and tensor product between linear forms induces a tensor product between vectors through this isomorphism. It is however useful to give some more details. Here, I use another way by defining an isomorphism between a vector space $E$ and its dual $E^*$ when $E$ is endowed with an inner product. The notion of inner product and how it can define an isomorphism between a space $E$ and its dual $E^*$ is developed in section \ref{sec:matelem:duality}.

\nT{Tensor product of linear forms} Let us consider the case $d=d'=1$, i.e. the tensor product of linear forms. Let $\langle .,.\rangle$ be an inner product on $E$. Let us recall that if $\abold \in E$, it is customary to denote by $\abold^*$ the linear forms on $E$ defined by
\begin{equation}
 \forall \: \xbold \in E, \quad \abold^*(\xbold)= \langle \abold,\xbold\rangle
\end{equation} 
and the same in $F$ with $\bbold \mapsto \bbold^*$. Let
\[
 \abold^* \in E^*, \qquad \bbold^* \in  F^*\]
Then, $\abold^* \otimes \bbold^* \in E^* \otimes F^*$ is the bilinear form on $E \times F$ such that
\begin{equation}\label{eq:tens:tens:def:1}
 \begin{array}{lcl}
 \forall \: (\xbold,\ybold) \in E \times F, \quad (\abold^* \otimes \bbold^*)(\xbold,\ybold) &=& \abold^*(\xbold).\bbold^*(\ybold) \\
 &=& \langle \abold,\xbold\rangle \langle \bbold,\ybold\rangle
 \end{array}
\end{equation}
Let now $\abold_\mu^* \in E_\mu^*$ and $\xbold_\mu \in E_\mu$. Then
\begin{equation}
 \begin{array}{lcl}
   (\abold_1^* \otimes  \ldots \abold_d^*)(\xbold_1,\ldots,\xbold_d) &=& \abold_1^*(\xbold_1) \times \ldots \times \abold_d^*(\xbold_d) \\
   &=& \langle \abold_1,\xbold_1 \rangle \ldots \langle \abold_d,\xbold_d\rangle
 \end{array}
\end{equation}
or
\begin{equation}\label{eq:tens:tens:def:11}
 \left(\bigotimes_{\mu=1}^d \abold_\mu^*\right) (\xbold_1,\ldots,\xbold_d) = \prod_{\mu=1}^d \abold_\mu^*(\xbold_\mu)
\end{equation}


\nS \textbf{Tensor product between vectors:} Let us recall that, if $E$ is a vector space, the vector space of linear forms on $E$ is called its \kw{dual}, and is denoted $E^*$. We have defined the tensor product of linear forms in equation(\ref{eq:tens:tens:def:11}). Each space $E_\mu$ is embedded with an inner product. This induces an isomorphism between each space $E_\mu$ and its dual $E_\mu^*$ through
\begin{equation}
 \begin{CD}
  \abold_\mu \in E @>>> \abold_\mu^* \in E^*
 \end{CD}
\end{equation}
with
\begin{equation}
\begin{CD}
 \xbold_\mu \in E @>\abold_\mu^*>> \langle \abold_\mu,\xbold_\mu\rangle \in \K
\end{CD}
\end{equation}
These isomorphisms enable to defines the tensor product between vectors $\abold_1, \ldots,\abold_d$ as 
\begin{equation}
 \begin{array}{lcl}
   (\abold_1 \otimes \ldots \otimes \abold_d) (\xbold_1, \ldots,\xbold_d) &=& (\abold_1^* \otimes \ldots \otimes \abold_d^*) (\xbold_1, \ldots,\xbold_d) \\
   &=& \displaystyle \prod_\mu \langle \abold_\mu,\xbold_\mu\rangle
 \end{array}
\end{equation}
This induces an inner product in $E_1 \otimes \ldots \otimes E_d$, with 
\begin{equation}
 \langle \abold_1 \otimes \ldots \otimes \abold_d \; , \; \xbold_1 \otimes \ldots \otimes \xbold_d \rangle = \langle \abold_1,\xbold_1\rangle \ldots \langle \abold_d,\xbold_d\rangle
\end{equation}
which can be extended to the whole space by linearity. Let us note that such a tensor product is not canonical\marginpar{\dbend}, as it depends on the inner product defined on $E$. If the inner product changes, the outcome of the tensor product changes.

\notes\label{page:tensor_history} Multilinear algebra has been developed with seminal work of H. Grassman on multivectors, (\emph{Ausdehnungslehre}, theory of Extension, 1844) to address geometry in spaces of more than three dimensions, and of G. Boole who at the same time developed invariant theory which has been further shifted and developed by A. Cayley (see \cite{Wolfson2008,Elliott13}). Grassmann introduced the wedge product $\wedge$ between vectors and his work is intimately linked with the later emergence of the notion of vector and vector space (see \citep{Dieudonne1979}). Cayley and Clifford made a link with Riemann's geometry (see \cite{Struik1987}). Such a link had been addressed as well by Riemann, Christoffel and Lipschitz (see \cite{Farwell1990}). Tensors emerged as mathematical objects in the framework of tensor calculus for differential geometry. The idea of tensor has been "in the air" before the term and definition have been coined. The name \emph{tensor} has been coined by the crystallographer Voigt in 1898\footnote{Voigt, W. (1898) Die fundamentalen physikalischen Eigenschaften der Krystalle in elementarer Darstellung, Verlag von Veit \& Comp., Leipzig.}. Tensors as a mathematical objects in differential geometry have been introduced by Ricci and Levi-Civita in 1900\footnote{Ricci-Curbastro, Gregorio; Levi-Civita, Tullio (1900), Méthodes de calcul différentiel absolu et leurs applications, \emph{Mathematische Annalen}, \textbf{54(1):}125-201, doi:10.1007/BF01454201, ISSN 1432-1807} with the name \emph{systems} in their theory of differential calculus. They have been thoroughly studied and used in the development of differential geometry, under the name of Ricci Calculus (see \cite{Schouten1954,berger03}). Einstein used them extensively with the name of tensors in his theory of relativity, coining the term \emph{tensor analysis}, which made them successful in mathematical physics (see e.g. \cite{Frankel1997,jost09}). The notation $\otimes$ is due to Murray and von Neumann in 1938\footnote{Murray, F. J. \& von Neumann, J. (1936) On Rings of Operators, \emph{Annals Math.} \textbf{37:}116-229}. The extension of tensor product to the tensor product of two Abelian groups is due to Whitney in 1938\footnote{Whitney, H. (1938) Tensor Products of Abelian Groups, \emph{Duke Math. Journal}, \textbf{4:}495-528.}, with the notation $A \circ B$, and not $A \otimes B$. The extension of tensor product of modules over a commutative ring is due to Bourbaki in 1948\footnote{Bourbaki, N. (1948) Livre II Algèbre Chapitre III (état 4) Algèbre Multilinéaire, \\http://archives-bourbaki.ahp-numerique.fr/archive/files/97a9fed708bdde4dc55547ab5a8ff943.pdf.}. These notes are borrowed from a note by Keith Conrad found on the web at http://www.math.uconn.edu/~kconrad/blurbs/linmultialg/tensorprod.pdf, of which we could not determine whether they had been published or not.\\
\\
As soon as these origins, two developments could be recognized: a calculation approach, with manipulations of indices, or a geometric insight of structures behind these calculations. \cite{Grinfeld2013} recalls that famous mathematicians of the beginning of 20th century (Elie Cartan, Herman Weyl) have warned against an extreme use of any approach. Some questions are better addressed with coordinates and indices calculations, some others with coordinate free approaches. The definition of tensor product of linear forms is adapted from \cite[section II.3]{Deheuvels1993}. It is presented in a more general framework of bilinear forms on modules in \cite[section 7D]{Charles1984}. The definition of tensors as multilinear forms can be found in \cite[sect. 1.3.1]{Landsberg2019}.\\
\\
\marginpar{\dbend}Definitions are not always consistent across literature, especially between different application domains. What is called here \emph{order} is sometimes called \emph{dimension}, or \emph{rank}, or \emph{degree}.

%
\section{Contravariant and covariant tensors}\label{sec:tens_prod:cocontra}
%

Covariance and contravariance are key notions for tensors when they are used in differential geometry or in physics. To keep the story short, if $E$ is a vector space, a contravariant tensor lives in $E \otimes \ldots \otimes E$, whereas a covariant tensor lives in $E^* \otimes \ldots \otimes E^*$. I will not make such a distinction in this book because a space $E$ and its dual $E^*$ will be made isomorphic through an inner product in $E$ (it will be presented in section \ref{sec:matelem:duality}). However, it is important to have these notions in mind. This section can be omitted: it will not be used afterwards. It will however facilitate the connection with other texts on tensors which may rely on such notions.

\nB A tensor can be given as a multidimensional array once a basis has been selected. If the basis changes, the components of the tensor change accordingly. This generalizes the change of components of a vector with a change of a basis of the vector space it belongs to. This leads to the notions of contravariant and covariant tensors, which extend the same notions for vectors, presented first.

\nS Let us select a simple case with $E$ being a one-dimensional space with basis $\ubold$. A vector $ \xbold \in E$ is given as
\begin{equation}
 \xbold = x \, \ubold \qquad \mbox{with} \quad 
 \begin{cases}
  x & \in \K \\
  \ubold, \xbold & \in E
 \end{cases}
\end{equation}
$\xbold$ is a \kw{vector}, $\ubold$ a \kw{basis} or a \kw{coordinate system}, and $x$ a \kw{component} in coordinate system given by $
\ubold$. Let us now change the basis of $E$ by
\begin{equation}
 \ubold' = \lambda \, \ubold \qquad \mbox{with} \quad \xbold = x' \, \ubold'
\end{equation}
It is easy to check by equating $\xbold = x\ubold=x'\ubold'=x'\lambda\ubold$  that
\begin{equation}\label{eq:tens_prod:cocontra:1}
 x' = x\, \lambda^{-1}
\end{equation}
Hence, a basis change $\ubold \longrightarrow \ubold' = \lambda \, \ubold$ induces an opposite change $x \longrightarrow x' = \lambda^{-1}\, x$ for the components of $\xbold$ in the new basis. The vector $\xbold$ is said \emph{contravariant}\marginpar{\emph{contravariant vector}}\index{vector!contravariant} (the components vary in the direction opposite to the basis change).

\nS Let us consider now the \emph{dual}\kw{dual space}\index{space!dual} $E^*$ of $E$, i.e. the vector space of linear forms on $E$. Let $\bm{\varphi} \in E^*$:
\begin{equation}
 \begin{CD}
   E @>\bm{\varphi}>> \K
 \end{CD}
\end{equation}
with
\begin{equation}
 \bm{\varphi}(\alpha \xbold)= \alpha \, \varphi(\xbold)
\end{equation}
It is possible to associate to $\ubold$, basis of $E$, its dual basis $\ubold^* \in E^*$ such that
\begin{equation}
 \ubold^*(\ubold) = 1
\end{equation}
In such a base, a linear form $\bm{\varphi} \in E^*$ can be written
\begin{equation}
 \bm{\varphi} = \varphi\, \ubold^*
\end{equation}
Hence
\begin{equation}
    \begin{array}{lcl}
      \bm{\varphi}(\xbold) &=& (\varphi\, \ubold^*)(x \, \ubold) \\
      &=& (\varphi\,x)\; \ubold^*(\ubold) \\
      &=& \varphi\,x
    \end{array}
\end{equation}
Let us now change the basis in $E$ as
\begin{equation*}
    \begin{CD}
      \ubold @>>> \ubold' = \lambda \, \ubold
    \end{CD}
\end{equation*}
This induces a dual $\ubold^{'*}$ such that $\ubold^{'*}(\ubold')=1$, and we have
\begin{equation}
 \bm{\varphi} = \varphi' \, \ubold^{'*}
\end{equation}
And, similarily, for same vector $\xbold = x' \, \ubold'$
\begin{equation}
   \bm{\varphi}(\xbold) = \varphi' \,x'
\end{equation}
So
\begin{equation}
 \varphi\,x = \varphi' \,x'
\end{equation}
But, as $x' = x\, \lambda^{-1}$ (see equation (\ref{eq:tens_prod:cocontra:1})), this leads to
\begin{equation}
 \varphi' = \lambda \, \varphi
\end{equation}
Hence, a basis change $\ubold \longrightarrow \ubold' = \lambda \, \ubold$ induces an similar change $\varphi \longrightarrow \varphi' = \lambda\, \varphi$ for the components of $\bm{\varphi}$ in the new basis. The linear form $\bm{\varphi}$ is said \emph{covariant}\marginpar{\emph{covariant vector}}\index{covariant vector}\index{vector!covariant} (the components vary in the same direction as the basis change). It is called as well a \kw{covector}.


\nT{Einstein convention} The choice of $\dim E=1$ has been made for sake of simplicity of formulas, to stick to the duality between covariance and contravariance. It is true for spaces of any dimension. Let $\xbold \in E$ with $\dim E=n$. If $(\ubold_1, \ldots,\ubold_n)$ is a basis of $E$, we can write
\begin{equation}
 \xbold = \sum_{i=1}^n\, x_i \, \ubold_i
\end{equation}
Here comes a convention in 2 steps to simplify notations:
\begin{itemize}[label=$\rightarrow$]
 \item the indices of contravariant vectors (the $\ubold_i$) are written as subscripts
 \item the indices of contravariant components are written as superscripts $x^i$ (and powers $k$ will be denoted $(x^i)^k$)
\end{itemize}
Hence, we have
\begin{equation}
 \xbold = \sum_{i=1}^n\, x^i \, \ubold_i
\end{equation}
The opposite is selected for covariant vectors: the indices of a covector if any are denoted by superscripts, and the indices of its components as subscripts. This can be summarized in following table:\\
\begin{center}
    \begin{tabular}{l|cc}
      & components & vector \\
      \hline 
      contravariant & $x^i$ & $\ubold_i$ \\
      covariant & $x_i$ & $\ubold^i$ \\
    \end{tabular}
\end{center}
The second step is \kw{Einstein convention}: when in a product $x^i\ubold_i$ a same index $i$ is superscript for one term and subscript for another, it should be understood as a summation, i.e. the formula is simplified as
\begin{equation}
 \begin{CD}
  \displaystyle \sum_{i}\, x^i \, \ubold_i @>>> x^i\, \ubold_i
 \end{CD}
\end{equation}
So, we have
\begin{equation}
 \xbold = x^i\, \ubold_i
\end{equation}


\nT{Components of a matrix}If $A \in\R^{m \times n}$ is a $m \times n$ matrix with coefficients $a_{ij}$, Einstein convention for component-wise formula of $\ybold = A\xbold$ reads
\begin{equation}
 \begin{CD}
  y_i = \sum_j\, a_{ij}x_j @>>> y^i = a^{i}_{\; j}\,x^j
 \end{CD}
\end{equation}
The components of a contravariant vector can be thought of as a vertical set of  components with superscripts, and of a covariant vector as an horizontal set of components with subscript. Vectors are written with subscripts and covectors with superscripts, as in\\
\begin{center}
 \begin{tabular}{ll}
   vector & $\xbold = x^i\, \ubold_i$ \\
   covector & $\xbold = x_i\, \ubold^i$ \\
 \end{tabular}
\end{center}
A matrix $A \in \R^{m \times n}$ has $m$ row covectors $\abold^1, \ldots,\abold^m$, and each row covector $\abold^i$ has components $(a^i_1, \ldots,a^i_n)$. The same matrix has $n$ column vectors $\abold_1, \ldots,\abold_n$, and each column vector $\abold_j$ has components $(a_j^1, \ldots,a_j^m)$.  Hence, the components of a matrix in a given basis are denoted $a^i_j$ with $i$ the row index and $j$ the column index. 

\nS This notation permits to extend to vector spaces with $d > 1$ the change of basis for vectors and covectors.\\
\\
\nP{$\rightarrow$} Let us have a vector
\begin{equation}
 \xbold = x^i\, \ubold_i
\end{equation}
and a new basis (or coordinate system)
\begin{equation}
 \ubold'_j = a^i_{\; j}\, \ubold_i \qquad \mbox{with} \quad \xbold = x^{'j}\, \ubold'_j 
\end{equation}
So
\begin{equation}
  \xbold = x^{'j}\, a^i_{\; j}\, \ubold_i
\end{equation}
(here, two indices, $i$ and $j$ are repeated; so there are two sums in this simple formula!). Then, by equating both expressions of $\xbold$
\begin{equation}\label{sec:tens_prod:cocontra:eq1}
  x^i = a^i_{\; j}\, x^{'j}
\end{equation}
and (if $\ubold_i=b_i^j\ubold'_j$)
\begin{equation}\label{sec:tens_prod:cocontra:eq2}
  x^{'j}  = b^{\; j}_i\,x^i 
\end{equation}
where
\begin{equation}\label{sec:tens_prod:cocontra:eq3}
  a^i_{\; j}b^{\; j}_k = \delta_k^j 
\end{equation}
where $\delta_k^j$ is Kronecker symbol, i.e. $B=A^{-1}$. If the basis vectors are transformed by application of matrix $A$, components are transformed by application of matrix $B=A^{-1}$. Hence, $\xbold$ is a \kw{contravariant vector}.

\nP{$\rightarrow$} The same calculation can be done for covariant vectors. It relies on the notion of dual basis. If $(\ubold_i)_i$ is a basis of $E$, the \kwnind{dual basis}\index{dual!basis} of it is the basis $(\ubold^{*i})_i$ defined by
\begin{equation}
 \ubold^{*i}(\ubold_j) = 
 \begin{cases}
  0 & \mbox{if}\quad i \neq j \\
  1 & \mbox{if}\quad i = j \\
 \end{cases}
\end{equation}
which can be written with Kronecker $\delta$ as well
\begin{equation}
  \ubold^{*i}(\ubold_j) = \delta^i_j
\end{equation}
Let $\bm{\varphi} \in E^*$. It can be written in basis $(\ubold^{*i})_i$
\begin{equation}
 \bm{\varphi} = \varphi_i\ubold^{*i}
\end{equation}
and, if $\xbold = x^k\ubold_k \in E$, we have
\begin{equation}
 \begin{array}{lcl}
   \bm{\varphi}(\xbold) &=& \left(\varphi_i\ubold^{*i}\right)\left(x^k\ubold_k\right) \\
   &=& \varphi_ix^k\, \ubold^{*i}(\ubold_k) \\
   &=& \varphi_ix^i
 \end{array}
\end{equation}
Let us now have the same basis change in $E$ defined by
\begin{equation}
 \ubold'_j = a^j_i\ubold^i 
\end{equation}
This induces a new dual basis $(\ubold^{'j*})_j$, and we have
\begin{equation}
 \bm{\varphi} = \varphi'_j\, \ubold^{'j*}
\end{equation}
We have as well
\begin{equation}
 \bm{\varphi}(\xbold) = \varphi'_jx^{'j}
\end{equation}
Then, we have (from equation \ref{sec:tens_prod:cocontra:eq2})
\begin{equation}
  \varphi_ix^i = \varphi'_jx^{'j} = \varphi'_j \, b^{\; j}_i\,x^i 
\end{equation}
So
\begin{equation}
  \varphi_i = b^{\; j}_i \, \varphi'_j
\end{equation}
or (from equation \ref{sec:tens_prod:cocontra:eq3})
\begin{equation}
  \varphi'_j = a^i_{\; j} \varphi_i
\end{equation}
The components are transformed with the same transformation than basis vector, so $\bm{\varphi}$ is covariant.

\nB Let us note that such a calculation did not require the calculation of the dual basis.

\nT{Components of a tensor} Let us just give an elementary example on $\Abold \in E \otimes E^* \otimes E$. Let $(\ubold_i)_i$ be a basis of contravariant vectors in $E$, and $(\ubold^{*j})_j$ of covariant vectors in $E^*$. Then, in these basis, one can write
\begin{equation}
    \Abold = a\indices{^i_j^k} \; \ubold_i \otimes \ubold^{*j} \otimes \ubold_k
\end{equation}

%
\section{Notations}\label{sec:tens_prod:notations}
%

Basic notations are recalled here:\\
\\
\begin{center}
\ovalbox{
\begin{tabular}{cl}
$\Abold$ & a tensor in $E_1 \otimes \ldots \otimes E_d$ \\
$A$ & a matrix in $E \otimes F$ \\
$\abold$ & a vector in $E$ \\
$d$ & order of a tensor \\
$I$ & $\llbracket 1,d \rrbracket \subset \N$ \\
$\mu$ & index running over $I$, $\mu \in I$ \\
%
\hline
%
$E$ & a finite dimensional vector space \\
$\EC$ & an orthonormal basis of $E$ \\
$\ES$ & a family $(E_1, \ldots,E_d)$ of $d$ vector spaces indexed by $\mu$ \\  
$E_\mu$ & a vector space in a family $\ES$ of $d$ vector spaces \\
$\EC_\mu$ & an orthonormal basis of $E_\mu$ \\
$\ebold_{\mu,i_\mu}$ & an element of $\EC_\mu$ (a basis vector)\\
$\bm{\EC}$ & basis of $E_1 \otimes \ldots \otimes E_d$ from $(\EC_1, \ldots,\EC_d)$\\
$\Ebold_\ibold$ & element $e_{1i_1} \otimes \ldots \otimes e_{di_d}$ of basis $\bm{\EC}$ \\
%
\hline
%
$n_\mu$ & dimension of space $E_\mu: \: n_\mu=\dim E_\mu$ \\
$I_\mu$ & set of indices for components in $E_\mu$;$I_\mu = \{1,n_\mu\}$\\
$\IC$ & $I_1 \times \ldots \times I_d$ \\
$\ibold$ & $\in \IC$, a multi-index $\ibold = (i_1, \ldots, i_d)$ \\
$a_\ibold$ & a component of $\Abold$ in basis $\EC_\mu$ for each space $E_\mu$ \\
\end{tabular}
}
\end{center}

%
\section{Some other definitions of tensor product}\label{sec:tens_prod:otherdefs}
%

Here are some other equivalent definitions of the tensor product.

\nS (see \citep[sect. XI-5]{Grillet2007}). This definition implements the observation that $\xbold \otimes \ybold = \xbold \ybold^\t$ where $\xbold$ is a column vector and $\ybold^\t$ is a row vector. Let $\dim E=m$ and $\EC = (\ebold_1,\ldots,\ebold_m)$ be a basis for $E$, $\dim F=n$ and $\F=(\fbold_1,\ldots,\fbold_n)$ be a basis for $F$. Let
\[
 \left\{
 \begin{array}{lcl}
  \xbold \in E &=& \displaystyle \sum_{i=1}^m x_i\ebold_i \\
  \ybold \in F &=& \displaystyle \sum_{j=1}^n y_j\fbold_j \\
 \end{array}
 \right.
\]
Then, $\xbold \otimes \ybold$ is the matrix $A \in \L(F,E)$ which has as coordinates
\begin{equation}\label{eq:tens_def}
 \xbold \otimes \ybold = \sum_{i,j} \, x_iy_j \; \ebold_i \otimes \ebold_j
\end{equation}
where $\ebold_i \otimes \ebold_j$ is the matrix with 1 at row $i$ and column $j$ and 0 elsewhere. This can be formulated in a more abstract way as follows: $E \otimes F$ is the $mn-$dimensional vector spaces with basis $(\ebold_i \otimes \fbold_j)_{i,j}$ such that $\xbold \otimes \ybold \in E \otimes F$ and equation (\ref{eq:tens_def}) holds. Then, it is easy to show that the map $(\xbold,\ybold) \longmapsto \xbold \otimes \ybold$ is bilinear.

\nB It is necessary to show that $\xbold \otimes \ybold$ does not depend on the selected basis. Moreover, if $\xbold \in \R^{m \times 1}$, $\ybold, \abold \in \R^{n \times 1}$, then (see equation (\ref{eq:matelem:tensors:1}) as well)
\begin{equation}
    \begin{array}{lcl}
         (\xbold \otimes \ybold).\abold &=& (\xbold\ybold^\t).\abold \\
         &=& \xbold(\ybold^\t.\abold) \\
         &=& \langle \ybold,\abold\rangle \, \xbold
    \end{array}
\end{equation}

\nS (see \cite[chap.~1]{Ryan2002}, \cite[sec.~ 3.2]{Garling2011}) Let $E,F$ be two finite dimensional vector spaces, $\B$ the set of a bilinear forms on $E \times F$, and $B \in \B$. Let $\xbold \in E$ and $\ybold \in F$. The \emph{tensor product} $\xbold \otimes \ybold$ is the linear form on $\B$ ("on" $\B$, not "of" $\B$\marginpar{\dbend})
\[
 \begin{CD}
  \B @> \xbold \otimes \ybold >> \K
 \end{CD}
\]
defined by
\begin{equation}
 (\xbold \otimes \ybold)(B) = B(\xbold,\ybold)
\end{equation}
Then, 
\[
 \xbold \otimes \ybold \in \B^*
\]
and is called an \kw{elementary tensor product}. Indeed, it is easy to show that
\begin{equation}
    \begin{array}{lcl}
        (\xbold \otimes \ybold)(\beta B + \beta' B') &=& (\beta B + \beta' B')(\xbold,\ybold) \\ 
        &=& \beta B(\xbold,\ybold) + \beta' B'(\xbold,\ybold)  \\
        &=& \beta (\xbold \otimes \ybold)(B) + \beta' (\xbold \otimes \ybold)(B)
    \end{array}
\end{equation}
The tensor product $E \otimes F$ between spaces $E$ and $F$ is the subspace of $\B^*$ spanned by these elementary tensor products. Hence, an element in $E \otimes F$ can be written as
\[
 A = \sum_k \, \alpha_k \; \xbold_k \otimes \ybold_k
\]
It is possible to show that these elementary tensor products span the entire space $\B^*$ provided $E$ and $F$ are finite dimensional. This is no longer the case if they are infinite dimensional. A \kw{tensor} is a linear combination of elementary tensor products.\\
\\
This can be extended to the tensor product of $d$ vectors and $d$ spaces. Let $(E_\mu)_\mu$ be a set of $d$ finite dimensional vector spaces on a same field $\K$. Let $\xbold_\mu \in E_\mu$. Let 
\[
 \ES = E_1 \times \ldots \times E_d
\]
and $\T^d(\ES)$ be the set of $d-$linear forms on $E_1 \times  \ldots \times E_d$. The \kw{tensor product} $\bigotimes_\mu \xbold_\mu$  between vectors $\xbold_\mu$ is an element of $\left(\T^d(\ES)\right)^*$
\[
 \begin{CD}
  \displaystyle \bigotimes_{\mu=1}^d\xbold_\mu \: : \: \T^d(\ES) @>>> \K
 \end{CD}
\]
such that, for any $\Tbold \in \T^d(\ES)$
\begin{equation}
 \left(\bigotimes_\mu \xbold_\mu\right)(\Tbold) = \Tbold(\xbold_1, \ldots,\xbold_d)
\end{equation}
As $\Tbold$ is multilinear, so is $\bigotimes_\mu \xbold_\mu$ as
\begin{equation}
 \begin{array}{l}
  \left(\xbold_1 \otimes \ldots \otimes \xbold_{\mu-1} \otimes (\alpha' \xbold'_\mu + \alpha'' \xbold''_\mu) \otimes \xbold_{\mu+1} \otimes \ldots \otimes \xbold_d\right)(\Tbold) \\
  \qquad = \Tbold\left(\xbold_1, \ldots , \xbold_{\mu-1}, (\alpha' \xbold'_\mu + \alpha'' \xbold''_\mu), \xbold_{\mu+1}, \ldots , \xbold_d\right) \\
  \qquad = \alpha' \Tbold\left(\xbold_1, \ldots , \xbold_{\mu-1}, \xbold'_\mu, \xbold_{\mu+1}, \ldots , \xbold_d\right) + \alpha'' \Tbold\left(\xbold_1, \ldots , \xbold_{\mu-1}, \xbold''_\mu, \xbold_{\mu+1}, \ldots , \xbold_d\right) \\
  \qquad = \alpha' \left(\xbold_1 \otimes \ldots \otimes \xbold_{\mu-1} \otimes \xbold'_\mu  \otimes \xbold_{\mu+1} \otimes \ldots \otimes \xbold_d\right)(\Tbold) \\
  \qquad\quad + \: \alpha'' \left(\xbold_1 \otimes \ldots \otimes \xbold_{\mu-1} \otimes \xbold''_\mu  \otimes \xbold_{\mu+1} \otimes \ldots \otimes \xbold_d\right)(\Tbold)
 \end{array}
\end{equation}
and as this is true for any $\Tbold$
\begin{equation}
 \begin{array}{l}
  \xbold_1 \otimes \ldots \otimes \xbold_{\mu-1} \otimes (\alpha' \xbold'_\mu + \alpha'' \xbold''_\mu) \otimes \xbold_{\mu+1} \otimes \ldots \otimes \xbold_d \\
  \qquad = \quad \alpha' \left(\xbold_1 \otimes \ldots \otimes \xbold_{\mu-1} \otimes \xbold'_\mu  \otimes \xbold_{\mu+1} \otimes \ldots \otimes \xbold_d\right) \\
  \qquad\quad + \: \alpha'' \left(\xbold_1 \otimes \ldots \otimes \xbold_{\mu-1} \otimes \xbold''_\mu  \otimes \xbold_{\mu+1} \otimes \ldots \otimes \xbold_d\right)
 \end{array}
\end{equation}
It is possible to show that elementary tensors span the entire space $\left(\T^d(\ES)\right)^*$ provided that every space $E_\mu$ is finite dimensional.


\nT{Universal property} The \kw{tensor product} of two vector spaces is often given in abstract algebra as a \kw{universal property}. Without entering in the domain of categories (to which such a definition belongs), let us show that it is quite simple. Let $E$ and $F$ be two vector spaces on a same field $\K$. Let $V$ be a vector space on $\K$ and $B$ a bilinear map
\begin{equation}
 \begin{CD}
   E \times F @>B>> V
 \end{CD}
\end{equation}
Then, there exists a vector space $W$ denoted $W = E \otimes F$, unique up to an isomorphism, such that, for any bilinear map $B$ and vector space $V$, there exists a linear map $L$ such that the following diagram is commutative:\\
\\
\begin{center}
\begin{tikzpicture}
\node (EF) at (0,0) {$E \times F$} ;
\node (V) at (3,1) {$V$} ; 
\node (W) at (3,-1) {$E \otimes F$} ; 
\draw[->, >=stealth] (EF) -- (V) node [pos=0.5, above, sloped] {$B$} ; 
\draw[->, >=stealth] (EF) -- (W) node [pos=0.5, below, sloped] {$\otimes$} ; 
\draw[->, >=stealth] (W) -- (V) node [pos=0.5, right] {$L$} ; 
\end{tikzpicture}
\end{center}
(see section \ref{sec:contract:bilin} for an illustration and a development on this universal property). This can easily be extended to $d-$linear maps with\\
\\
\begin{center}
\begin{tikzpicture}
\node (EF) at (0,0) {$E_1 \times \ldots \times E_d$} ;
\node (V) at (3,1) {$V$} ; 
\node (W) at (3,-1) {$E_1 \otimes \ldots \otimes E_d$} ; 
\draw[->, >=stealth] (EF) -- (V) node [pos=0.5, above, sloped] {$T$} ; 
\draw[->, >=stealth] (EF) -- (W) node [pos=0.5, below, sloped] {$\otimes$} ; 
\draw[->, >=stealth] (W) -- (V) node [pos=0.5, right] {$L$} ; 
\end{tikzpicture}
\end{center}
such that this diagram is commutative.

\notes This presentation of a tensor product of two vector spaces as a universal property is borrowed from \cite[p.~42]{Deheuvels1993} and \cite[p. 555]{Lang1984}. It is the classical way to define tensor product with a slight incursion into more abstract algebra.

%
\chapter{Elementary matrix operations with tensor notations}\label{chap:matelem}
%

Tensors can be seen as extensions of matrices (a matrix is a $2-$modes tensor, and tensors when given as arrays of coordinates are sometimes called hypermatrices), and many elementary operations are better known on matrices. Thus, to familiarize the reader with the concept of tensor, we rewrite in this section some matrix definitions and properties as tensor definitions and properties.

\nB One thing cannot be extended to tensors, and has been a source of much confusion because many works have tried (unsuccessfully) to find such an extension: a matrix is a representation of a linear map from a vector space $E$ on a vector space $F$. This is no longer the case with tensors,\marginpar{\dbend} which are by essence multilinear, and not bilinear, even if some matrices can be built from tensors (this is called matricization, see section \ref{sec:contract:alamat}). 

\nB This has far reaching consequences. For example, the product of matrices is defined from the combination of linear maps
\[
 \begin{CD}
   E @>A>> F @>B>> G
 \end{CD}
\]
and $BA \in \L(E,G)$. Such a product cannot be defined canonically for tensors. Matrix product leads to many classical matrix decomposition as a product of other matrices with elements having simpler structures, like SVD, Schur product, QR decomposition, etc ... At a more elementary level, if $E=F$, one can define $A^2$ or any power $A^k$. One can define eigenvalues $\lambda$ and eigenvectors $\xbold$ with the property $A\xbold = \lambda \xbold$. None of these operations can be extended as such to tensors. 


\nT{Notations}I denote in this chapter
\begin{center}
 \ovalbox{
 \begin{tabular}{cl}
  $E,F,G$ & vector spaces on $\R$ \\
  $\EC, \F, \G$ & orthonormal basis for these spaces \\
  $\abold,\xbold$ & $\in E$ \\
  $\bbold,\ybold$ & $\in F$ \\
  $\cbold,\zbold$ & $\in G$ \\
  $A$ & $\in E \otimes F \simeq \L(F,E)$ \\
  $m$ & $\dim E$ \\
  $n$ & $\dim F$ \\
  $p$ & $\dim G$ \\
  $i$ & $\in \llbracket 1,m\rrbracket := \{1,\ldots,m\}$\\
  $j$ & $\in \llbracket 1,n \rrbracket := \{1,\ldots,n\}$ \\
  $k$ & $\in \llbracket 1,p \rrbracket:= \{1,\ldots,p\}$ \\
 \end{tabular}
 }
\end{center}


\nT{Organization of the chapter}  This chapter presents some elementary properties of matrices as tensor properties. It is organized as follows:

\begin{description}
\item[section \ref{sec:matelem:duality}] An isomorphism is built between the dual $E^*$ of a space $E$ endowed with an inner product and $E$. Such an isomorphism is not canonical. 
\item[section \ref{sec:matelem:tensors}] A matrix is presented as a tensor, which establishes an isomorphism between $E \otimes F$ and $\L(F,E)$.
\item[section \ref{sec:matelem:contract}] The contraction as a key operation in tensor algebra is presented on the example of a matrix $\times$ vector product.
\item[section \ref{sec:matelem:diag}] The way to build diagrams for sketching contractions is presented.
\item[section \ref{sec:matelem:basis}] The decomposition of a matrix on a basis of the tensor space $E \otimes F$ is presented. 
\item[section \ref{sec:matelem:basis}] Some simple formulas for deriving calculation in matrix algebra as calculations in tensor algebra without use of components are given.
\item[section \ref{sec:matelem:iso}] Finally, some canonical and non canonical isomorphisms between various spaces involved in matrix and bilinear algebra are recalled.

\end{description}

%
\section{Duality and inner product}\label{sec:matelem:duality}
%

Here, I will use the notion of \kw{inner product}\index{product!inner} on a vector space $E$ on $\R$. An inner product is a bilinear form $B$ (i.e. a bilinear map from $E \times E$ on $\R$) which is\\
\begin{center}
\begin{tabular}{ll}
symmetric & $B(\xbold,\ybold)=B(\ybold,\xbold)$\\
definite & $B(\xbold,\xbold)=0 \quad \Rightarrow \quad \xbold=0$\\
positive & $B(\xbold,\xbold) \geq 0$\\
\end{tabular}
\end{center}

\nB An inner product will be denoted $\langle . , . \rangle$, i.e. $B(\xbold,\ybold):=\langle \xbold,\ybold\rangle$. It is sometimes denoted $(.,.)$. There is a deep connection between inner product and duality which will be used in this book and is presented here. 

\nS Let $E$ be a finite dimensional vector space on a field $\K$. The dual $E^*$ is the set of linear forms on $E$. As $\dim E=\dim E^*$, they are isomorphic. Let us first consider $\K=\R$ (the case $\K=\C$ will be considered afterwards). We show here that it is equivalent to define an inner product $\langle .,.\rangle$ in $E$ (hence the choice $\K=\R$) and an isomorphism $\varphi$ between $E$ and $E^*$. The reason for using this property is twofold: first, Euclidean distances built from inner products are ubiquitous, and second, this avoids to make a distinction between vectors (in $E$) and covectors (in $E^*$) or between contravariant (in $E$) and covariant (in $E^*$) vectors. As a consequence, Einstein convention will not be used.

\nS Let $(\ebold_i)_i$ be a basis in $E$. Any vector $\xbold \in E$ can be written $\xbold = \sum_ix_i\ebold_i$. The key notion for the equivalence between an inner product and an isomorphism between $E$ and $E^*$ is the notion of \kwnind{dual basis}\index{dual!basis}. The dual basis of $(\ebold_i)_i$ is the basis of $E^*$ defined by
\begin{equation}
 \ebold_i^*(\ebold_j) = 
 \begin{cases}
  0 & \mbox{if} \quad  i \neq j \\
  1 & \mbox{if} \quad  i=j
 \end{cases}
\end{equation}
This defines an isomorphisms between $E$ and $E^*$ once a basis has been selected in $E$ by
\begin{equation}
 \begin{CD}
  \displaystyle \xbold = \sum_ix_i\ebold_i @>>> \displaystyle \xbold^* = \sum_ix_i\ebold_i^*
 \end{CD}
\end{equation}

\nS It induces an inner product in $E$ through 
\begin{equation}
 \langle \abold , \xbold \rangle = \abold^*(\xbold)
\end{equation}
\begin{proof}
If $\abold = \sum_ia_i\ebold_i$, $\abold^*=\sum_ia_i\ebold_i^*$ and
\begin{equation}
 \begin{array}{lcl}
   \langle \abold , \xbold\rangle &=& \abold^*(\xbold) \\
   &=& \displaystyle \left(\sum_ia_i \ebold_i^*\right)\left(\sum_j x_j\ebold_j\right) \\
   &=& \displaystyle \sum_{i,j} a_ix_j \, \ebold_i^*(\ebold_j) \\
   &=& \displaystyle \sum_i a_ix_i
 \end{array}
\end{equation}
Then, the map $(\abold , \xbold) \longrightarrow \langle \abold,\xbold\rangle$ is bilinear and\\
\begin{center}
    \begin{tabular}{ll}
         symmetric & because $\langle \abold , \xbold\rangle = \langle \xbold , \abold\rangle$ \\
         definite &  because if $\langle \abold , \abold\rangle=0$, $\sum_ia_i^2=0$ and $a_i=0$ for any $i$, hence $\abold=0$ \\
         positive &because $\langle \abold , \abold\rangle = \sum_ia_i^2 \geq 0$.
    \end{tabular}
\end{center}

\noindent Then it is an inner product.
\end{proof}

\nS Conversely, if $\langle .,.\rangle$ is an inner product in $E$, the map
\begin{equation}
 \begin{CD}
  E @>>> E^* \\
  \abold @>>> \abold^*
 \end{CD}
\end{equation}
defined by $\abold^*(\xbold)=\langle \abold , \xbold \rangle$ is an isomorphism.
\begin{proof}
It is linear. Its kernel is $\abold=0$ because if $\abold^*=0$, $\langle \abold,\xbold\rangle=0$ for any $\xbold$ and $\abold=0$. It is surjective because the image of a basis in $E$ is its dual basis which spans $E^*$.  
\end{proof}

\nT{Key observation} \marginpar{\dbend} It is important to understand that such an isomorphism is not canonical, but depends on the selected basis $(\ebold_i)_i$. One way to define an inner product in $E$ is to select a basis and set it as an orthonormal basis. So, the choice of a basis permits to define an isomorphism between $E$ and $E^*$ and to define an inner product in $E$.

\nS If $\K=\C$, the same holds but with a sesquilinear product (not developed here). Unless otherwise stated, the field with which I will work is $\R$, but extension to $\C$ will be straightforward \emph{mutatis mutandis}. Most of algebraic properties will be valid for $\K=\C$ as well.

%
\section{A matrix as a tensor}\label{sec:matelem:tensors}
%

Let $E$ be a finite dimensional vector space on $\R$ endowed with an inner product $\langle .,.\rangle$. Let $\xbold \in E$. It can be seen as a vertical vector, i.e. a $n \times 1$ matrix if $\dim E=n$. The transpose of $\xbold$ is denoted $\xbold^\t$ and is a $1 \times n$ matrix. It represents a linear form on $E$, and is a covector.

\nS The translation of classical matrix identities into the notation of tensor algebra can be done by simply recalling that
\begin{equation}
 \left\{
 \begin{array}{lcl}
   \langle \xbold,\ybold\rangle &=& \xbold^\t \ybold \\
   \xbold \otimes \ybold &=& \xbold\ybold^\t
 \end{array}
 \right.
\end{equation}
$\xbold \otimes \ybold$ can be visualized as
\begin{center}
 \begin{tikzpicture}
   \draw (0,0) rectangle (0.4,2) ;
   \draw (0.8,2.4) rectangle (2.8,2.8) ;
   \draw (0.8,0) rectangle (2.8,2) ;
   \node () at (0.2,2.6){$\otimes$};
   \node () at (0.2, 1){$x_i$};
   \node () at (1.8, 2.6){$y_j$};
   \node () at (1.8, 1){$x_iy_j$};
   \draw[dashed] (.4, 1) -- (1.4,1) ; 
   \draw[dashed] (1.8, 2.4) -- (1.8,1.2) ; 
 \end{tikzpicture}
\end{center}
If $\xbold,\ybold$ are vectors, $\xbold \otimes \ybold$ is a matrix. We have 
\begin{equation}
 \left\{
    \begin{array}{cl}
      \xbold & \in E \\
      \ybold & \in F \\
      \xbold \otimes \ybold & \in E \otimes F \\
      & \in \L(F,E) \\
    \end{array}
 \right.
\end{equation}

\nS Let us show the isomorphism between $E \otimes F$ and $\L(F,E)$. Let $\abold \in E$ and $\bbold,\ybold \in F$. As $\abold \otimes \bbold = \abold \bbold^\t$, $(\abold \otimes \bbold).\ybold$ is defined by
\begin{equation}
 \begin{array}{lcl}\label{eq:matelem:tensors:1}
   (\abold \otimes \bbold)\ybold &=& (\abold \bbold^\t)\ybold \\
   &=& \abold(\bbold^\t \ybold) \\
   &=& \langle \bbold,\ybold\rangle \abold
 \end{array}
\end{equation}
Hence
\begin{equation}
 (\abold \otimes \bbold)\ybold = \langle \bbold,\ybold\rangle \abold
\end{equation}
It is easy to show that, $\abold$ and $\bbold$ being given,
\[
 \begin{CD}
   \ybold \in F @>\abold\otimes \bbold>> (\abold\otimes \bbold)\ybold \in E
 \end{CD}
\]
is linear, so $\abold \otimes \bbold \in \L(F,E)$. We then have
\begin{equation}
 E \otimes F \simeq \L(F,E)
\end{equation}
because any linear map in $\L(F,E)$ can be written as a linear combination of elementary (i.e. rank one) maps $\abold \otimes \bbold$.

\nS We have as well
\begin{equation*}
 \begin{array}{lcl}
   (\abold \otimes \bbold)(\bbold' \otimes \cbold) &=& (\abold \bbold^\t)(\bbold'\cbold^\t) \\
   &=& \abold(\bbold^\t \bbold')\cbold^\t \\
   &=& \langle \bbold,\bbold'\rangle \, \abold \cbold^\t  \\
   &=& \langle \bbold,\bbold'\rangle \, \abold \otimes \cbold
 \end{array}
\end{equation*}
Hence
\begin{equation}
  (\abold \otimes \bbold)(\bbold' \otimes \cbold) = \langle \bbold,\bbold'\rangle \, \abold \otimes \cbold
\end{equation}

%
\section{Matrix $\times$ vector product and contraction}\label{sec:matelem:contract}
%

Here, I introduce a key notion in tensor algebra which is the contraction, through a simple example. This notion will be much more developed further (see chapter \ref{chap:contract}).

\nB The matrix $\times$ vector product
\[
 \begin{CD}
  (A,\ybold) @>>> A.\ybold
 \end{CD}
\]
is an example of a more general operation, called \kw{contraction}. This can be seen first on elementary tensors as in
\begin{equation}
 \begin{CD}
  (\abold \otimes \bbold, \ybold) @>>> \langle \bbold,\ybold\rangle \abold
 \end{CD}
\end{equation}
where $\abold \in E$ and $\bbold,\ybold \in F$. The pair $(\bbold,\ybold) \in F \times F$ has been contracted  into the scalar $\langle \bbold,\ybold\rangle \in \R$. This defines a contraction over $F$ between an elementary tensor $\abold \otimes \bbold$ and a vector $\ybold$. As the inner product is a bilinear form, and as any matrix is a linear combination of a set of elementary tensors, the contraction matrix $\times$ vector can be extended to any matrix by linearity. This retrieves the matrix $\times$ vector product.
\begin{proof}
Indeed, let $(\ebold_i)_i$ and $(\fbold_j)_j$ be orthonormal basis of respectively $E$ and $F$. We have
\begin{equation}
 \begin{array}{lcl}
  A.\ybold &=& \displaystyle \left(\sum_{i,j}a_{ij}\, \ebold_i \otimes \fbold_j\right)\left(\sum_ky_k\fbold_k\right) \\
  &=& \displaystyle \sum_{i,j,k} a_{ij}y_k \, (\ebold_i \otimes \fbold_j).\fbold_k \\
  &=& \displaystyle\sum_{i,j,k} a_{ij}y_k \, \langle \fbold_j,\fbold_k\rangle \ebold_i \\
  &=& \displaystyle \sum_i\left(\sum_j a_{ij}y_j\right)\ebold_i
 \end{array}
\end{equation}
where the matrix $\times$ vector product $A.\ybold$ can be recognized.
\end{proof}

\nS The contraction (over $F$) can be denoted by
\begin{equation*}
 \begin{CD}
  \begin{array}{c}
    E \otimes F \\
    \times \\
    F
  \end{array}
  @>\mathrm{contraction}>> E
 \end{CD}
\end{equation*}
with
\begin{equation}
 \begin{CD}
  \begin{array}{c}
    \abold \otimes \bbold \\
    \bullet \\
    \ybold
  \end{array}
  @>\mathrm{contraction}>> \langle \bbold,\ybold\rangle \abold
 \end{CD}
\end{equation}
(here, the contraction is denoted by a bullet). 

\nS One can define a contraction over space $E$ as well through
\begin{equation*}
 \begin{CD}
  \begin{array}{c}
    E \otimes F \\
    \times \\
    E
  \end{array}
  @>\mathrm{contraction}>> F
 \end{CD}
\end{equation*}
with
\begin{equation}
 \begin{CD}
  \begin{array}{c}
    \abold \otimes \bbold \\
    \bullet_\e \\
    \xbold
  \end{array}
  @>\mathrm{contraction}>> \langle \abold,\xbold\rangle \bbold
 \end{CD}
\end{equation}
Here, the operator $\bullet_\e$ simply reminds that the contraction is done over space $E$. So, we have
\begin{equation}
 \left\{
    \begin{array}{lcl}
      (\abold \otimes \bbold) \bullet_\e \xbold &=& \langle \abold,\xbold\rangle \bbold \\
      (\abold \otimes \bbold) \bullet_\f \ybold &=& \langle \bbold,\ybold\rangle \abold
    \end{array}
 \right.
\end{equation}
Such a notation is ambiguous if $E=F$. In such a case, the order of the mode $E$ or $F$, can be specified, as in
\begin{equation}
 \left\{
    \begin{array}{lcl}
      (\abold \otimes \bbold) \bullet_1 \xbold &=& \langle \abold,\xbold\rangle \bbold \\
      (\abold \otimes \bbold) \bullet_2 \ybold &=& \langle \bbold,\ybold\rangle \abold
    \end{array}
 \right.
\end{equation}

\nS All these contractions defined on elementary tensors can be extended to the whole spaces by linearity.

%
\section{Diagrams for contraction}\label{sec:matelem:diag}
%

This is an occasion to introduce a visualization of elementary tensor operations which can be implemented whatever the order of the tensor. 

\nS In such a \kw{tensor diagram}\index{tensor!diagram}, a tensor product $A = \xbold \otimes \ybold \in E \otimes F$ is an array which can be sketched as\footnote{Vectors in diagrams are not written with bold latin letters, but plain ones.}\\
\\
\begin{center}
\begin{tikzpicture}
\node[tensor_ell] (ab) at (0,0) {$x \otimes y$} ;  
\draw (ab) -- (-1.5,0) node [pos=0.5, above] {$\e$}  ; 
\draw (ab) -- (1.5,0)  node [pos=0.5, above] {$\f$} ; 
\end{tikzpicture}
\end{center}
Bars with a free tip sketch the modes of the matrix (here 2 modes). Vectors $\xbold \in E$ and $\ybold \in F$ can be sketched as\\
\\
\begin{center}
\begin{tikzpicture}
\node[tensor_ell] (x) at (0,0) {$x$} ;  
\draw (x) -- (1,0)  node [pos=0.5, above] {$\e$} ; 
\node[tensor_ell] (y) at (4,0) {$y$} ;  
\draw (y) -- (5,0) node [pos=0.5, above] {$\f$} ; 
\end{tikzpicture}
\end{center}
There is one mode for a vector, hence one bar with a free tip. A contraction between a matrix and a vector is sketched by sharing a bar, which has no free tip, as in\\
\\
\begin{center}
\begin{tikzpicture}
\node[tensor_ell] (ab) at (0,0) {$a \otimes b$} ;  
\node[tensor_ell] (y) at (1.5,0) {$y$} ;  
\draw (ab) -- (-1.5,0) node [pos=0.5, above] {$\e$}  ; 
\draw (ab) -- (y)  node [pos=0.5, above] {$\f$} ; 
\end{tikzpicture}
\end{center}
If we collapse over the shared bar labelled by $F$, and accept the rule that a collapse is an inner product, this yields\\
\\
\begin{center}
\begin{tikzpicture}
\node[tensor_ell] (ab) at (0,0) {$a \; \langle b,y\rangle$} ;  
\draw (ab) -- (-1.5,0) node [pos=0.5, above] {$\e$}  ; 
\end{tikzpicture}
\end{center}
which is a visualization of
\[
 (\abold \otimes  \bbold).\ybold = \langle \bbold,\ybold\rangle \; \abold
\]

\nS Let us have
\[
 \abold \otimes \bbold \in E \otimes F, \qquad \bbold' \otimes \cbold \in F \otimes G 
\]
This can be sketched as\\
\\
\begin{center}
\begin{tikzpicture}
\node[tensor_ell] (ab) at (0,0) {$a \otimes b$} ;  
\draw (ab) -- (-1.5,0) node [pos=0.5, above] {$\e$}  ; 
\draw (ab) -- (1.5,0)  node [pos=0.5, above] {$\f$} ; 
\node[tensor_ell] (aabb) at (4,0) {$b' \otimes c$} ;  
\draw (aabb) -- (2.5,0) node [pos=0.5, above] {$\f$}  ; 
\draw (aabb) -- (5.5,0)  node [pos=0.5, above] {$\g$} ; 
\end{tikzpicture}
\end{center}
Contracting over $F$ yields\\
\\
\begin{center}
\begin{tikzpicture}
\node[tensor_ell] (ab) at (0,0) {$a \otimes b$} ;  
\node[tensor_ell] (aabb) at (2,0) {$b' \otimes c$} ;  
\draw (ab) -- (-1.5,0) node [pos=0.5, above] {$\e$}  ; 
\draw (ab) -- (aabb)  node [pos=0.5, above] {$\f$} ; 
\draw (aabb) -- (3.5,0)  node [pos=0.5, above] {$\g$} ; 
\end{tikzpicture}
\end{center}
which yields by a collapse over the shared bar as an inner product in $F$:\\
\\
\begin{center}
\begin{tikzpicture}
\node[tensor_ell] (ab) at (0,0) {$\langle b,b'\rangle \, a \otimes c $} ;  
\draw (ab) -- (-2.5,0) node [pos=0.5, above] {$\e$}  ; 
\draw (ab) -- (2.5,0)  node [pos=0.5, above] {$\g$} ; 
\end{tikzpicture}
\end{center}
which is a visualization of
\[
(\abold \otimes \bbold).(\bbold' \otimes \cbold) = \langle \bbold,\bbold'\rangle \, \abold \otimes \cbold 
\]
This gives at the same time, for $(\abold \otimes \bbold)(\bbold' \otimes \cbold)$
\begin{itemize}[label=$\rightarrow$]
    \item the outcome: $\langle \bbold,\bbold' \rangle \, \abold \otimes \cbold$ by looking at the content
    \item the space it belongs to: $E \otimes G$ by looking at the labels of the free tips
\end{itemize}

%
\section{Decomposition of a matrix on a basis}\label{sec:matelem:basis}
%

Let $ \EC = (\ebold_i)_i$ be an orthonormal basis of $E$ and $ \F = (\fbold_j)_j$ an orthonormal basis of $F$. Then
\[
\EC \otimes \F = (\ebold_i \otimes \fbold_j)_{i,j}
\]
is an orthonormal basis of $E \otimes F$, in which
\begin{equation}
 A = \sum_{i,j} \alpha_{ij} \; \ebold_i \otimes \fbold_j
\end{equation}
(such a decomposition is true as well if $\EC$ and $\F$ are not orthonormal). We have
\begin{equation}
 \begin{array}{lcl}
  A\fbold_\ell &=& \displaystyle \left(\sum_{i,j} \alpha_{ij} \; \ebold_i \otimes \fbold_j\right).\fbold_\ell \\
  &=& \displaystyle \sum_{i,j} \alpha_{ij} (\ebold_i \otimes \fbold_j)\fbold_\ell \\
  &=& \displaystyle \sum_{i,j} \alpha_{ij} \langle \fbold_j,\fbold_\ell\rangle \; \ebold_i\\  
  &=& \displaystyle \sum_i \alpha_{i\ell}\ebold_i
 \end{array}
\end{equation}
and
\begin{equation}
\begin{array}{lcl}
 A &=& \displaystyle \sum_{i,j} \alpha_{ij} \; \ebold_i \otimes \fbold_j \\
 &=& \displaystyle \sum_j \left(\sum_i \alpha_{ij}\ebold_i\right) \otimes \fbold_j \\
 &=& \displaystyle \sum_j A\fbold_j \otimes \fbold_j
\end{array}
\end{equation}

\nS If $(\fbold_j)_j$ is not orthonormal, we have
\begin{equation}
 \begin{array}{lcl}
  A\fbold_\ell^* &=& \displaystyle \left(\sum_{i,j} \alpha_{ij} \; \ebold_i \otimes \fbold_j\right).\fbold_\ell^* \\
  &=& \displaystyle \sum_{i,j} \alpha_{ij} (\ebold_i \otimes \fbold_j)\fbold_\ell^* \\
  &=& \displaystyle \sum_{i,j} \alpha_{ij} \langle \fbold_j,\fbold_\ell^*\rangle \; \ebold_i \\
  &=& \displaystyle \sum_i \alpha_{i\ell}\ebold_i
 \end{array}
\end{equation}
and
\begin{equation}\label{eq:matelem_nonortho}
 A = \sum_j A\fbold_j^*\otimes \fbold_j
\end{equation}
where $(\fbold_j^*)_j$ is the dual basis of $(\fbold_j)_j$.

\nS As $(\fbold_j^*)^*=\fbold_j$, this reads
\begin{equation}
 A = \sum_j A\fbold_j \otimes \fbold_j^*
\end{equation}
If $A$ is square $(n \times n)$ and diagonalisable with $A\ubold_i= \lambda \ubold_i$, this reads
\begin{equation}\label{eq:matdecomdiago}
 A = \sum_i \lambda_i \; \ubold_i \otimes \ubold_i^*
\end{equation}
If $A$ is symmetric, $(\ubold_i)_i$ is an orthonormal family, and 
\begin{equation}
 A = \sum_i \lambda_i \; \ubold_i \otimes \ubold_i
\end{equation}

\nS This can be summarized in the following table:\\
\\
\begin{center}
\ovalbox{
\begin{tabular}{rlcl}
$A$ & $ = \displaystyle \sum_j A\fbold_j^*\otimes \fbold_j$ & & if $(\fbold_j)_j$ is a general basis \\
    & $ = \displaystyle \sum_j A\fbold_j \otimes \fbold_j^*$ & & if $(\fbold_j)_j$ is a general basis \\
    & $ = \displaystyle \sum_j A\fbold_j\otimes \fbold_j$ & & if $(\fbold_j)_j$ is an orthonormal basis \\
    & $ = \displaystyle \sum_j \lambda_i \, \ubold_i \otimes \ubold_i^*$ & & if $A$ diagonalisable with $A\ubold_i=\lambda \ubold_i$\\  
    & $ = \displaystyle \sum_j \lambda_i \, \ubold_i \otimes \ubold_i$ & & if moreover $A$ is symmetric\\       
\end{tabular}
}
\end{center}

\nS We have as well
\begin{equation}
   A^\t =\sum_{i,j}\, \alpha_{ij}\, \fbold_j \otimes \ebold_i
\end{equation}
Hence
\begin{equation}
 A^\t e_k^* = \sum_j\, \alpha_{kj}\, \fbold_j
\end{equation}
and
\begin{equation}
 \begin{array}{lcl}
   A &=& \displaystyle \sum_{i,j}\, \alpha_{ij} \, \ebold_i \otimes \fbold_j \\
   &=& \displaystyle\sum_i \ebold_i \otimes \left(\sum_j\alpha_ {ij}\fbold_j\right) \\
   &=& \displaystyle \sum_i \, \ebold_i \otimes A^\t \ebold_i^*
 \end{array}
\end{equation}
which translates into
\begin{equation}
 A = \sum_i \ebold_i \otimes A^\t \ebold_i
\end{equation}
when $(\ebold_i)_i$ is orthonormal. Hence\\
\\
\begin{center}
\ovalbox{
\begin{tabular}{rlcl}
$A$ & $ = \displaystyle \sum_i \, \ebold_i \otimes A^\t \ebold_i^* $ &  & if $(\ebold_i)_i$ is a general basis \\
    & $ = \displaystyle  \sum_i \ebold_i \otimes A^\t \ebold_i $ & & if $(\ebold_i)_i$ is an orthonormal basis \\
\end{tabular}
}
\end{center}

\nT{Note} There is more than one way to do one thing ... For example: if $(\ebold_i)_i$ is an orthonormal basis of $E$, $\I=\sum_i\ebold_i \otimes \ebold_i$ is the identity matrix. Indeed, if $\xbold=\sum_i\ebold_i$, then
$\I \xbold = \left(\sum_i\ebold_i \otimes \ebold_i\right)\xbold=\sum_i \langle \xbold, \ebold_i\rangle \ebold_i = \xbold$. Then, $A=A\I=A\left(\sum_i\ebold_i \otimes \ebold_i\right)=\sum_iA(\ebold_i \otimes \ebold_i)=\sum_iA\ebold_i \otimes \ebold_i$ (see equations \ref{eq:elemat}). The other formula can be shown in the same way.

%
\section{Useful formulas}\label{sec:matelem:useful}
%

I present here some useful formulas, which are elementary (it is useful to see in which space the vectors and matrices are, and half of the work is done). Some have already been shown and are recalled here. This section can be seen as an exercise.

\nS Let $E,F,G,V,W$ be vector spaces, and
\[
  \abold,\xbold \in E, \quad \bbold,\abold',\ybold \in F, \quad \bbold' \in G, \quad A \in E \otimes F, \qquad M \in V \otimes E, \quad Q \in W \otimes F
\]
Then
\begin{equation}\label{eq:elemat}
 \left\{ 
 \begin{array}{lcl}
  (\abold \otimes \bbold)^\t &=& \bbold \otimes \abold \\
  (\abold \otimes \bbold)\ybold &=& \langle \bbold,\ybold\rangle \, \abold \\ 
  \langle A,\xbold \otimes \ybold \rangle &=& \langle A\ybold, \xbold \rangle \\
  (\abold \otimes \bbold)(\abold' \otimes \bbold') &=& \langle \bbold,\abold'\rangle  \; \abold \otimes \bbold'\\
  M(\abold \otimes \bbold) &=& M\abold \otimes \bbold \\
  (\abold \otimes \bbold)Q &=& \abold \otimes Q^\t \bbold \\
  \|\abold \otimes \bbold\| &=& \|\abold\|\,\|\bbold\|
 \end{array}
 \right.
\end{equation}
We have the following diagram\\
\\
\begin{center}
 \begin{tikzpicture}
   \node (V) at (4,4) {$V$} ;
   \node (G) at (0,2) {$G$} ;
   \node (F) at (2,2) {$F$} ;
   \node (E) at (4,2) {$E$} ;
   \node (W) at (2,0) {$W$} ;
   \draw[->, >=latex] (E) -- (V) node[pos=0.5, right]{$M$} ; 
   \draw[->, >=latex] (F) -- (E) node[pos=0.5, above]{$\abold \otimes \bbold$} ; 
   \draw[->, >=latex] (G) -- (F) node[pos=0.5, above]{$\abold' \otimes \bbold'$} ;    
   \draw[->, >=latex] (F) -- (V) node[pos=0.5, sloped, above]{$M(\abold \otimes \bbold)$} ;    
   \draw[->, >=latex] (W) -- (F) node[pos=0.5, left]{$Q$} ;    
   \draw[->, >=latex] (W) -- (E) node[pos=0.5, sloped, below]{$(\abold \otimes \bbold)Q$} ;    
 \end{tikzpicture}
\end{center}

\begin{proof}Let us recall that $(\abold \otimes \bbold)\ybold = (\abold \bbold^\t)\ybold = \abold (\bbold^\t \ybold) = \langle \bbold,\ybold\rangle \abold$. Some proofs are given for elementary matrices, like $M=\ubold \otimes \vbold$ or $Q=\ubold \otimes \vbold$, and can be extended to all matrices by linearity.
\begin{itemize}[label=$\rightarrow$]
 \item for $\langle A,\xbold \otimes \ybold \rangle $. Let $A = \abold \otimes \bbold$. Then, 
 \[
  \begin{array}{lcl}
  \langle A,\xbold \otimes \ybold \rangle &=& \langle \abold \otimes \bbold ,\xbold \otimes \ybold \rangle \\
  &=& \langle \abold,\xbold\rangle \langle \bbold,\ybold\rangle \\
  &=& \langle \langle \bbold,\ybold\rangle \abold \; , \; \xbold \rangle \\
  &=& \langle (\abold \otimes \bbold)\ybold \, , \, \xbold \rangle \\
  &=& \langle A\ybold,\xbold\rangle
  \end{array}
 \]
 \item For computing $(\abold \otimes \bbold)(\abold' \otimes \bbold')$:
 \[
  \begin{aligned}
   \forall \: \ybold, \quad (\abold \otimes \bbold)(\abold' \otimes \bbold')\ybold &= (\abold \otimes \bbold)(\langle \bbold',\ybold\rangle \abold') \\
   &= \langle \bbold',\ybold\rangle (\abold \otimes \bbold)\abold' \\
   &= \langle \bbold',\ybold\rangle \langle \bbold,\abold'\rangle \abold \\
   &= \langle \bbold,\abold'\rangle (\langle \bbold',\ybold\rangle \abold \\
   &= \langle \bbold,\abold'\rangle (\abold \otimes \bbold')\ybold
  \end{aligned}
 \]
 \item For $M(\abold \otimes \bbold)$, let us select $M = \ubold \otimes v$. Then
 \[
  \begin{aligned}
   M(\abold \otimes \bbold) &= (\ubold \otimes \vbold)(\abold \otimes \bbold) \\
   &= \langle \vbold,\abold\rangle \ubold \otimes \bbold \\
   &= (\ubold \otimes \vbold).\abold \otimes \bbold \\
   &= M\abold \otimes \bbold
  \end{aligned}
 \]
 \item For $(\abold \otimes \bbold)Q$, let us select $Q = \ubold \otimes \vbold$. Then
 \[
  \begin{aligned}
  (\abold \otimes \bbold)Q &= (\abold \otimes \bbold)(\ubold \otimes \vbold) \\
  &= \langle \bbold,\ubold \rangle \abold \otimes \vbold \\
  &= \abold \otimes \langle \bbold,\ubold \rangle \vbold \\ 
  &= \abold \otimes (\vbold \otimes \ubold)\bbold \\
  &= \abold \otimes Q^\t \bbold
  \end{aligned}
 \]
\end{itemize}
\end{proof}

%
\section{Useful isomorphisms}\label{sec:matelem:iso}
%

In this section, I summarize the classical isomorphisms between spaces, underlining their status of canonical or non canonical (i.e. which are independent or dependent on the choice of a basis, here an orthonormal basis). Let us recall that selecting an inner product in a vector space $E$ enables to define an isomorphism between $E$ and $E^*$ which depends on the inner product. This however avoids to specify whether vectors are covariant or contravariant, and to distinguish their status (a covariant vector is a vector of $E^*$ and a contravariant vector is a vector of $E$, see section \ref{sec:tens_prod:cocontra}) . 

\nS A $2-$modes tensor in $E \otimes F$ is defined as a bilinear form on $E \times F$, i.e. a form $B$ such that
\begin{equation}
 \left\{
    \begin{array}{lcl}
        B(\lambda \xbold + \lambda' \xbold' \, , \, \ybold) &=& \lambda \, B(\xbold,\ybold) + \lambda' \, B(\xbold',\ybold) \\
        B(\xbold \, , \, \mu \ybold + \mu' \ybold') &=& \mu \, B(\xbold,\ybold) + \mu' \, B(\xbold,\ybold')
    \end{array}
 \right.
\end{equation}

\nS Let $\B(E,F)$ be the space of bilinear forms on $E \times F$. We have the following isomorphism
\begin{equation}
 \B(E,F) \simeq E^* \otimes F^*
\end{equation}
with
\begin{equation}
 (\abold^* \otimes \bbold^*)(\xbold,\ybold) = \abold^*(\xbold)\bbold^*(\ybold)
\end{equation}
which can be extended to $E^* \otimes F^*$ by linearity. 

\nS Selecting an inner product in $E$ and in $F$ defines an isomorphism between $E$ (resp. $F$) and $E^*$ (resp. $F^*$) by
\begin{equation}
 \begin{CD}
  E @>>> E^* \\
  \abold @>>> \abold^*
 \end{CD}
\end{equation}
with
\begin{equation}
 \forall \: \xbold \in E, \quad \abold^*(\xbold) = \langle \abold,\xbold\rangle
\end{equation}
This isomorphism depends on the inner product.

\nS There is a canonical isomorphism between $\L(F,E)$ and $E \otimes F^*$ defined by
\begin{equation}
 \forall \: \ybold \in F, \quad (\abold \otimes \bbold^*)\ybold = \bbold^*(\ybold)\, \abold
\end{equation}
which can be extended to $E \otimes F^*$ by linearity. 

\nS As there is an isomorphism between $F$ and $F^*$, this can be transported as
\begin{equation}
 \forall \: \ybold \in F, \quad (\abold \otimes \bbold)\ybold = \langle \bbold,\ybold\rangle \abold
\end{equation}
which defines an isomorphism between $E \otimes F$ and $\L(F,E)$.

\nS We then have the following isomorphisms:\\
\\
\begin{center}
 \ovalbox{
 \begin{tabular}{cclc}
   between & and & description & status\\
   \hline
   $\B(E,F)$ & $E^* \otimes F^*$ & $B(\xbold,\ybold)=\abold^*(\xbold)\bbold^*(\ybold)$ & canonical \\
   $\L(F,E)$ & $E \otimes F^*$ & $(\abold \otimes \bbold^*)\ybold=\bbold^*(\ybold)\, \abold$ & canonical \\
   $E$ & $E^*$ & $\abold^*(\xbold)=\langle \abold,\xbold\rangle$ & non canonical \\
   $\B(E,F)$ & $E \otimes F$ & $B(\xbold,\ybold)=\langle \abold,\xbold \rangle\langle \bbold,\ybold\rangle$ & non canonical \\
   $\L(F,E)$ & $E \otimes F$ & $(\abold \otimes \bbold).\ybold = \langle \bbold,\ybold\rangle\,\abold$ & non canonical
 \end{tabular}
 }
\end{center}

\nS These are a key ingredients for a smooth development of tensor algebra for $d-$modes tensors as presented in chapter \ref{chap:dmodestens}.

\notes The material in this chapter is standard. There exists several excellent textbooks from which I have been inspired for matrix algebra and analysis, like \cite{Schott1997,Horn2012}. Contraction diagrams have been presented here in one of their simplest guise. The framework where they have been developed is what is called nowadays tensor networks. They have been proposed for the first time, apparently, in \cite{Penrose1971}. Rewriting elementary matrix operations with the toolbox of tensor algebra is standard, and can be found in \cite{Franc1992}.

%
\chapter{$d-$modes tensors}\label{chap:dmodestens}
%

Let $E$, $F$ and $G$ be three vector spaces on a same field $\K$, and let us assume that an orthonormal basis has been selected in $E$, $F$, $G$. Hence, $\K=\R$. Extension to $\K=\C$ is straightforward through sesquilinear forms, and will not be considered explicitly here. \\
\\
I give some details about \kw{tensor product of $d$ vector spaces}. It is easier to figure out what is happening with $d=3$ because a visualization still is possible. This is however impossible for $d>3$. Therefore, some elementary notations are first presented for $d=3$ and then extended to $d>3$ with some algebra. Anything relies on a clear understanding of what order, dimensions, indices are.\\
\\
The key result in this section is, beyond a specification of $d-$modes tensors, an elementary (and old) result on the expression of a given tensor in a new basis.

\nT{Organization of the chapter} This chapter is organized as follows:
\begin{description}
\item[section \ref{sec:dmodestens:3modes}] An introduction of notations for $3-$modes tensors.
\item[section \ref{sec:dmodestens:dmodes}] An introduction of notations for general $d-$modes tensors.
\item[section \ref{sec:dmodestens:basis}] Tucker decomposition is presented as the expression of a given tensor in a new basis.
\end{description}

All spaces in this chapter are finite dimensional.

%
\section{$3-$modes tensors}\label{sec:dmodestens:3modes}
%

Let us start with tensor product of 3 vector spaces.

\nT{Definitions} Let $E,F,G$ be three vector spaces with $\dim E = m$, $\dim F=n$ and $\dim G = p$, endowed each with an inner product.

\begin{itemize}[label=$\rightarrow$]
 \item The tensor product $E^* \otimes F^* \otimes G^*$ of vector spaces $E^*,F^*,G^*$ is the space of trilinear forms on $E \times F \times G$.
 \item it can be transported to $E \otimes F \otimes G$ by (non canonical) isomorphisms $E \leftrightarrow E^*, F \leftrightarrow F^*, G \leftrightarrow G^*$, which is equivalent to selecting an inner product in $E, F, G$.
 \item a trilinear form on $E \otimes F \otimes G$ is called a $3-$modes tensor
\end{itemize}

\nS Let $\abold \in E$, $\bbold \in F$ and $\cbold \in G$. The tensor product $\abold \otimes \bbold \otimes \cbold$ is the trilinear form on $E \otimes F \otimes G$ 
\begin{equation}
 \begin{CD}
  E \times F \times G @>\abold \otimes \bbold \otimes \cbold >> \R 
 \end{CD}
\end{equation}
defined by 
\begin{equation}
  (\abold \otimes \bbold \otimes \cbold) \, (\xbold,\ybold,\zbold) = \langle \abold, \xbold \rangle \, \langle \bbold,\ybold \rangle \, \langle \cbold,\zbold \rangle
\end{equation}
or, component-wise
\begin{equation}
  \begin{aligned}
    \abold \otimes \bbold \otimes \cbold \, (x,y,z) &= \left(\sum_ia_ix_i\right)\left(\sum_jb_jy_j\right)\left(\sum_kc_kz_k\right) \\
    &= \sum_{i,j,k} \alpha_{ijk} \,x_iy_jz_k, \qquad \mbox{with} \quad \alpha_{ijk} = a_ib_jc_k
  \end{aligned}
\end{equation}
The tensor $\abold \otimes \bbold \otimes \cbold$ is called an \kw{elementary tensor}.

\nS The tensor product $E \otimes F \otimes G$ of spaces $E,F,G$ is the set of linear combinations of elementary tensors. It can be shown that every trilinear form on $E \times F \times G$ can be written as a linear combination of elementary trilinear forms in $E,F,G$ are finite dimensional. 

\nS Visualization is given here for $d=3$ with $E=\R^m, F=\R^n,G=\R^p$ even though it will not be available for $d>3$. We assume that $E$ is equipped with a basis $\EC=(\ebold_i)_i$, $F$ with a basis $\F = (\fbold_j)_j$and $G$ with a basis $\G = (\gbold_k)_k$.\\
\\
\begin{center}
\begin{tikzpicture}
\draw (0,0) rectangle (2,3);
\draw (1,4) -- (3,4) -- (3,1) -- (2,0) ; 
\draw (2,3) -- (3,4) ; 
\draw[->] (0,3) -- (2,5) ;
\draw[->] (0,3) -- (0,-1) ;
\draw[->] (0,3) -- (4,3) ;
\draw[dashed] (1,4) -- (1,1) ;
\draw[dashed] (1,1) -- (3,1) ;
\draw[dashed] (0,0) -- (1,1) ; 
\node () at (0.3,3.8 ) {$G$} ;
\node () at (1.5,3.3 ) {$F$} ;
\node () at (-.5, 1.5) {$E$} ;
\node () at (2.3, 5.3) {$k \in \llbracket 1,p \rrbracket$} ; 
\node () at (4,3.3) {$j \in \llbracket 1,n \rrbracket$} ; 
\node () at (0,-1.3) {$i \in \llbracket 1,m \rrbracket$} ; 
\node () at (-2.5,3) {$E \otimes F \otimes G$} ;
\node () at (-2.5,2.5) {$\R^m \otimes \R^n \otimes \R^p$} ; 
\end{tikzpicture}
\end{center}

\nS If we have the (non necessarily orthonormal) basis
\begin{equation}
\begin{array}{ccccc}
\ebold_i \in E=\R^m   & \qquad & \fbold_j \in F=\R^n & \qquad & \gbold_k \in G=\R^p \\
i \in \llbracket 1,n \rrbracket & \qquad & j \in \llbracket 1,m \rrbracket & \qquad & k \in \llbracket 1,p \rrbracket
\end{array}
\end{equation}
then
\begin{equation}
\begin{aligned}
\Abold & = \sum_{i,j,k} \alpha_{ijk} \; \ebold_i \otimes \fbold_j \otimes \gbold_k \\
& \in \R^m \otimes \R^n \otimes \R^p \\
& := \R^{m \times n \times p}
\end{aligned}
\end{equation}
For example,
\begin{equation}\label{eq:x_tens_y_tens_z}
 \Abold = \xbold \otimes \ybold \otimes \zbold \quad \Longrightarrow \quad \alpha_{ijk} = x_iy_jz_k
\end{equation}

\nT{Spaces, indices, dimensions} It will be useful for the rest of this book to keep in mind the following table of correspondence between spaces, indices and dimensions\\
\begin{center}
 \begin{tabular}{ccc}
   \hline
   Space & Indices & Dimension \\
   \hline
   $E$ & $i$ & $m$ \\
   $F$ & $j$ & $n$ \\
   $G$ & $k$ & $p$ \\   
   \hline
 \end{tabular}
\end{center}

\nB $\Abold$ is a tensor of order 3, of dimensions $(m,n,p)$. The modes are $E,F,G$.

%
\section{$d-$modes tensors}\label{sec:dmodestens:dmodes}
%

This section is a technical complement only. All the ingredient for tensor product have been introduced for $d=3$ modes with a geometric intuition of a $3-$modes tensor being a shoebox. Nothing else is required to evolve from 3 modes to $d$ modes with $d \in \N$. Such an extension is presented in this short section.

\nS Let $\mu$ be an index running over $\llbracket 1,d \rrbracket$. Let $(E_\mu)_\mu$ be a family of vector spaces with $\dim E_\mu=n_\mu$, each equipped with an inner product. Let $i_\mu$ be an index running over $\llbracket 1,n_\mu \rrbracket$ (see section \ref{sec:tens_prod:notations} for more details about notations). 

\nT{Tensor product} The \kw{tensor product} 
\[
 E_1 \otimes \ldots \otimes E_d = \bigotimes_\mu E_\mu
\]
is the vector space of $d-$linear forms on $E_1 \times \ldots \times E_d$. A $d-$linear form in $E_1 \otimes \ldots \otimes E_d$ is called a $d-$modes tensor.

\nS An elementary $d-$modes tensor $\abold_1 \otimes \ldots \otimes \abold_d$ is defined by 
\begin{equation}
\begin{array}{lcl}
  \displaystyle \left(\bigotimes_\mu \abold_\mu\right)(\xbold_1, \ldots,\xbold_d) &=& \displaystyle \prod_\mu \langle \abold_\mu,\xbold_\mu\rangle \\
  &=& \displaystyle \prod_\mu \left(\sum_{i_\mu}a_{i_\mu}^{(\mu )} x_{i_\mu}^{(\mu )}\right) \\
  &=& \displaystyle \sum_{i_1} \ldots \sum_{i_d} \left(\prod_{\mu}\alpha_{i_\mu}^{(\mu )}\right)\left(\prod_{\mu}x_{i_\mu}^{(\mu )}\right)
\end{array}
\end{equation}
So, as it is true whatever $\xbold_1, \ldots,\xbold_d$
\begin{equation}
 \bigotimes_{\mu=1}^d\abold_\mu = \sum_\ibold \alpha_\ibold \Ebold_\ibold \qquad \mbox{with} \quad \alpha_{i_1\ldots i_d} = \prod_\mu \alpha_{i_\mu}^{(\mu )}
\end{equation}

\nS A $d-$modes tensor is a linear combination of $d-$modes elementary tensors
\begin{equation}
 \Abold = \sum_a \lambda_a \; \xbold_{1a} \otimes \ldots \otimes \xbold_{da}
\end{equation}
with
\begin{equation}
 \xbold_{\mu a} \in E_\mu
\end{equation}
If $(\ebold_{\mu,i_\mu})_{i_\mu}$ is a basis of $E_\mu$, $\Abold$ can be decomposed on this basis as
\begin{equation}
 \Abold = \sum_{i_1} \ldots \sum_{i_d} \, \alpha_{i_1\ldots i_d} \, \ebold_{1i_1} \otimes  \ldots \otimes \ebold_{di_d}
\end{equation}

\nS We denote
\begin{center}
 \begin{tabular}{l|c}
 \hline
  Space & $E_\mu$ \\
  Dimension & $n_\mu$\\
  Indices & $i_\mu$ \\
  \hline
 \end{tabular}
\end{center}

%
\section{Decomposition on a basis and Tucker decomposition} \label{sec:dmodestens:basis}
%

In this section (and some others), we distinguish
\begin{itemize}
 \item[$\rhd$] a mathematical object: a tensor $\Abold$ element of a tensor space
 \item[$\rhd$] the expression of such an object in a given basis $\EC$, denoted $\llbracket \Abold \rrbracket_\EC$ when useful.
\end{itemize}
The calculations developed here will assemble elementary operations on tensors like permutations (see section \ref{sec:tenselem:permut}) and associated linear maps (see section \ref{sec:contract:alamat}). However, this section is self-contained (any required notion is explained). It is merely a matter of notations.

\nS Let $\Abold \in E \otimes F \otimes G$, with $\EC=(\ebold_i)_i$ an orthonormal basis of $E$, $\F = (\fbold_j)_j$ an orthonormal basis of $F$ and $\G =(g_k)_k$ an orthonormal basis of $G$. Let us have
\begin{equation}
 \Abold = \sum_{i,j,k} \alpha_{ijk} \; \ebold_i \otimes \fbold_j \otimes \gbold_k
\end{equation}
We denote 
\begin{equation}
 [\alpha_{ijk}]_{i,j,k} = \llbracket \Abold \rrbracket_{\EC \otimes \F \otimes G}
\end{equation}
the $3-$modes array of the coordinates of $\Abold$ in basis $(\EC, \F, \G)$. Let
\[
 \EC'=(\ebold_a')_a, \qquad \F' = (\fbold'_b)_b, \qquad \G' = (\gbold'_c)_c
\]
be new orthonormal\footnote{The same calculation can be done for general basis, just by using the dual basis in inner product in (\ref{eq:tucker:newbasis}).} basis of $E,F,G$ respectively. We compute $(\alpha_{abc}')_{a,b,c}$ such that
\begin{equation}
\begin{aligned}
 \Abold &= \sum_{i,j,k} \alpha_{ijk} \; \ebold_i \otimes \fbold_j \otimes \gbold_k \\
 &= \sum_{a,b,c} \alpha'_{abc} \; \ebold'_a \otimes \fbold'_b \otimes \gbold'_c 
 \end{aligned}
\end{equation}
We have
\begin{equation}
  [\alpha'_{abc}]_{a,b,c} = \llbracket \Abold \rrbracket_{\EC' \otimes \F' \otimes G'}
\end{equation}
Then
\begin{equation}\label{eq:tucker:newbasis}
  \begin{array}{lcl}
  \alpha'_{abc} &=& \langle \Abold \, , \; \ebold'_a \otimes \fbold'_b \otimes \gbold'_c \rangle \\
  &=& \displaystyle \left\langle \sum_{i,j,k} \, \alpha_{ijk} \, \ebold_i \otimes \fbold_j \otimes \gbold_k \; , \; \ebold'_a \otimes \ebold'_b \otimes \ebold'_c \right\rangle \\
  &=& \displaystyle \sum_{i,j,k} \, \alpha_{ijk} \langle \ebold_i,\ebold'_a\rangle \, \langle \fbold_j,\fbold'_b\rangle \, \langle \gbold_k,\gbold'_c\rangle 
 \end{array}
\end{equation}

\nS Let us define transition matrices $M,\: N, \: P$ in $E,F,G$ respectively such that 
\begin{equation}
 \ebold'_a = \sum_im_{ia}\ebold_i, \qquad \fbold'_b = \sum_j n_{jb}\fbold_j, \qquad \ebold'_c = \sum_kp_{kc}\gbold_k
\end{equation}
(the column $a$ of $M$ is the set of coordinates of the vector $\ebold'_a$ in basis $\EC$, etc ..). Then
\begin{equation}
 \langle \ebold_i,\ebold'_a\rangle = m_{ia}, \qquad \langle \fbold_j,\fbold'_b\rangle = n_{jb}, \qquad \langle \gbold_k,\gbold'_v\rangle = p_{kc}
\end{equation}
and
\begin{equation}\label{eq:basis-change}
 \begin{aligned}
   \alpha'_{abc} &= \sum_{i,j,k} \alpha_{ijk}\, m_{ia} \, n_{jb} \, p_{kc}
 \end{aligned}
\end{equation}

\nT{Tucker decomposition} Equation (\ref{eq:basis-change}) is called the \kwnind{Tucker decomposition}\index{Tucker!decomposition} of $\Abold$ on basis $\EC \otimes \F \otimes \G$.

\nT{Computation} Let us observe that these quantities can be computed iteratively, as
\begin{equation}
 \sum_{i,j,k} \alpha_{ijk}\, m_{ia} \, n_{jb} \, p_{kc} = \sum_i \left(\sum_j\left(\sum_k \alpha_{ijk}\, p_{kc} \right) n_{jb} \right) m_{ia} 
\end{equation}
The term $\sum_k \alpha_{ijk}\, p_{kc}$ is a term of a matrix product. To see this, let us consider the matrix $A_\textsc{g} \in EF \otimes G$ where $EF$ stands for $E \otimes F$ (by associativity of $\otimes$, $E \otimes F \otimes G = (E \otimes F) \otimes G = EF \otimes G$) by reshaping modes and coordinates (see section  \ref{sec:contract:alamat}). It has $mn$ rows and $p$ columns. Let us write $\alpha_{ijk}= \alpha_{(ij)k}$ (index $(ij)$ for rows and $k$ for columns). Let
\[
 \beta'_{(ij)c} = \sum_k \alpha_{(ij)k}\, p_{kc}
\]
It is the term $\beta'_{(ij)c}$ of matrix
\[
 \underbrace{B_\textsc{g}}_{mn \times p}= \underbrace{A_\textsc{g}}_{mn \times p}.\underbrace{P}_{p \times p}
\]
of dimension $mn \times p$. This matrix can be reshaped into a tensor $\Bbold \in E \otimes F \otimes G$ of dimensions $(m,n,p)$ with
\begin{equation}
 \beta_{ijc} = \beta'_{(ij)c}
\end{equation}
We then have
\begin{equation}
 \sum_{i,j,k} \alpha_{ijk}\, m_{ia} \, n_{jb} \, p_{kc} = \sum_i \left(\sum_j\beta_{ijc} n_{jb} \right) m_{ia} 
\end{equation}

\nS Similarly, The term $\sum_j\beta_{ijc}\, n_{jb}$ is a term of a matrix product. Let us consider the matrix $B_\textsc{f} \in EG \otimes F$ with $mp$ rows and $n$ columns, defined by $\beta''_{(ic)j}= \beta_{ijc}$ (index $(ic)$ for rows and $j$ for columns). Let
\[
 \gamma'_{(ic)b} = \sum_j \beta''_{(ic)j}\, n_{jb}
\]
Then $\gamma'_{(ic)b}$ is the term of matrix
\[
 C_\textsc{f}= B_\textsc{f}N
\]
of dimension $np \times m$ as well. This matrix can be reshaped into a tensor $\Cbold \in E \otimes F \otimes G$ of dimensions $(m,n,p)$ with
\begin{equation}
\gamma_{ibc} = \gamma'_{(ic)b} = \sum_j\beta''_{(ic)j}n_{jb} = \sum_j\beta_{ijc}n_{jb} 
\end{equation}
We then have
\begin{equation}
 \sum_{i,j,k} \alpha_{ijk}\, n_{kc} \, p_{jb} \, q_{ia} = \sum_i \gamma_{ibc} m_{ia} 
\end{equation}

\nS Let us consider finally the term $\sum_i \gamma_{ibc} m_{ia}$, and therefore the matrix $C_\textsc{e} \in FG \otimes E$ with $np$ rows and $m$ columns defined by $\gamma''_{(bc)i}=\gamma_{ibc}$ (index $(bc)$ for rows and $i$ for columns). Let
\[
 \alpha''_{(bc)a} = \sum_i \gamma''_{(bc)i}\, m_{ia}
\]
Then, $\alpha''_{(bc)a}$ us the term of matrix
\[
 A'_\textsc{e}= C_\textsc{e}M
\]
of dimension $mp \times m$ as well. This matrix can be reshaped into a tensor $\Abold' \in E \otimes F \otimes G$ of dimensions $(m,n,p)$ with
\begin{equation}
 \alpha'_{abc} = \alpha''_{(bc)a}
\end{equation}
We then have
\begin{equation}
 \begin{array}{lcl}
   \alpha'_{abc} &=& \alpha''_{(bc)a} \\
   &=& \displaystyle \sum_i \gamma''_{(bc)i}\, m_{ia} \\
   &=& \displaystyle \sum_i \gamma'_{(ic)b}\, m_{ia} \\
   &=& \displaystyle \sum_{i,j} \beta''_{(ic)j}\, n_{jb}\, m_{ia} \\
    &=& \displaystyle \sum_{i,j} \beta'_{(ij)c}\, n_{jb}\, m_{ia} \\
     &=& \displaystyle \sum_{i,j,k} \alpha_{ijk}\, p_{kc} \, n_{jb}\, m_{ia} \\
 \end{array}
\end{equation}

\nS The whole cycle can be written algorithmically as \\
\\
\begin{algorithm}[H]
\begin{algorithmic}[1]
\STATE \textbf{input} $\Abold = \texttt{A[i,j,k]}$, $M=\texttt{M[i,a]}$, $N=\texttt{N[j,b]}$, $P=\texttt{P[k,c]}$
\STATE $A_\textsc{g}$: $(mn \times p)$ matrix of the linear application in $E \otimes F \rightarrow G$ associated to $\Abold$
\STATE $B_\g=A_\textsc{g}P$
\STATE $\Bbold = B_\g.\texttt{reshape}(m,n,p)$
\STATE $B_{\textsc{f}}$: $(mp \times n)$ matrix of the linear application in $E \otimes G \rightarrow F$ associated to $\Bbold$
\STATE $C_\f=B_\f N$
\STATE $\Cbold = C_\f.\texttt{reshape}(m,n,p)$
\STATE $C_\e$: $(np \times m)$ matrix of the linear application in $F \otimes G \rightarrow E$ associated to $\Cbold$
\STATE $A'=C_\e M$
\STATE $\Abold' = A'.\texttt{reshape}(m,n,p)$
\RETURN $\Abold'$
\end{algorithmic}
\caption{Tucker decomposition $\textsc{Tucker-decomposition}(\Abold,M,N,P)$}
\label{algo:tucker-decompose}
\end{algorithm}

\notes Equation (\ref{eq:basis-change}) is formula $(17)$ in \cite[p.~286]{Tucker1966}. There is a lot of confusion about this simple relation, which is nothing more than a basis change.\marginpar{\dbend} 
\begin{itemize}[label=$\rhd$]
 \item Such a relation has been coined as multilinear matrix multiplication in \cite{Silva2008}, and denoted
\begin{equation*}
 \Abold' = (M,N,P).\Abold
\end{equation*}
which has often been accepted as a notation for Tucker decomposition. As basis change in a space $E$ can be written as the action of the linear group $\GL(E)$ on a vector ($\abold = M\abold'$, where $\abold$ and $\abold'$ are the expression of the same vector in two basis, with basis change given by $M$), such a notation is consistent with the notion of action of the linear group (see section \ref{sec:alstruct:group}).
\item It is denoted $\Abold' = (M,N,P) \circ \Abold$ in \cite{Grasedick2010}, where $\Abold$ is called the core tensor of $\Abold'$, and this decomposition is called the Tucker format or representation. In \cite{Grasedick2010}, the matrices $M,N,P$ are called the mode-frames of the Tucker representation. 
\item It is formula $(4.1)$ in \cite[p.~474]{Kolda2009} where it is denoted
\[
 \Abold' = \llbracket \Abold \, ; \, M,N,P \rrbracket
\]
The matrices $M,N,P$ are called factor matrices, and are labeled as principal components for each mode.
\item The equivalence between a Tucker format or decomposition and a change of basis in the modes is presented in \cite{Franc1992}.
\end{itemize}

%% file: primer.tex
Linear algebra has developed by keeping two paths very active and linked:
\begin{itemize}
    \item developing very efficient algorithms for some issues in numerical matrix calculus  
    \item identifying some elementary operations which can be combined to develop and understand these algorithms.
\end{itemize}
Both ways are currently very active research domains in tensor algebra and tensor calculus, like Siamese twins. The operations presented here will be very useful to understand, and control, the algorithms for computing low rank approximations of a given tensor, like CP-decomposition Tucker model, HOSVD, HOOI (see Part \ref{part:blra}), and many other things. This requires a small detour through (slightly) abstract algebra, which is the aim of this chapter. For example, HOSVD is a well accepted method for extending PCA to tensors, even if it is known that it is not the best approximation of a given tensor by a Tucker approximation. The algorithm to build HOSVD of a given tensor is the following sequence (in words):\\
\\
$\rightarrow$ in a loop over all orders:
\begin{enumerate}
    \item matricize the tensor over this order (see section \ref{sec:contract:alamat})
    \item run the SVD of this matrix; store the outcomes (see section \ref{sec:inter:svd}) 
    \item reshape the tensor as  its initial organisation in modes (see section \ref{sec:tenselem:reshape})
    \item run a permutation for next order to be in the front (see section \ref{sec:tenselem:permut})
\end{enumerate}
This is just an example among others. Permutations, matricization and reshaping are elementary bricks of algorithms for best low CP rank approximation, and best Tensor-Train approximation, too. Contraction and Kronecker product (isomorphic to tensor product up to a reshaping) are operations the usefulness of which cannot be underestimated, and deserve each a dedicated section. Contraction over repeated modes within a family of tensors is the "glue" to build tensor networks, which are an approach for disentangling interactions in many-body quantum systems and in joint laws in statistical modeling to facilitate marginalisation.

\nB All notions still are elementary. This chapter is organized as follows:

\begin{description}
\item[Chapter \ref{chap:tenselem}] Elementary algebraic operations on tensors like permutation, symmetrization, antisymmetrization, slicing and are presented. All of them consist in rewriting a given tensor with the same elements but with some reorganization. At the end of the chapter Kronecker product is presented, with a deep connection with vectorization and reshaping. All are key manipulations on tensors.
\item[Chapter \ref{chap:contract}] A key operation, the contraction, is presented. It is the extension to tensors of the matrix $\times$ vector product (which is an example of contraction). Contraction is the dual operation of tensor product.
\item[Chapter \ref{sec:alstruct}] It is organised around the Kronecker product between matrices and the action of the linear group on vector spaces, with presentation of some classical operations between matrices.
\end{description}

%
\chapter{Elementary algebraic tensor operations}\label{chap:tenselem}
%

Algebra is like a lego set: some elementary operations are defined, and its  flavour is about assembling them to build different shapes. Some elementary tensor operations are presented here:\\
\begin{center}
 \begin{tabular}{lll}
   operation & short description & section \\
   \hline
   permutation & permutations of the modes of a tensor & \ref{sec:tenselem:permut} \\
   symmetrization & invariance by a permutation & \ref{sec:tenselem:sym} \\
   antisymmetrization & invariance by a permutation up to its signature & \ref{sec:tenselem:antisym} \\
   slicing & extracting a slice & \ref{sec:tenselem:slice} \\
   reshaping & reorganising modes & \ref{sec:tenselem:reshape} \\
   \hline
 \end{tabular}
\end{center}
All of these operations consist in rewriting a given tensor with some reorganization of its modes.

\nB One elementary operation is a key one, and deserves a chapter by its own: contraction, presented in chapter \ref{chap:contract}. It is an operation dual of tensor product and, in a word, generalizes the matrix $\times$ matrix or matrix $\times$ vector product.

%
\section{Permutation}\label{sec:tenselem:permut}\index{permutation@\textbf{permutation}}
%

This is a generalization of the transposition of a matrix. Let us present it first on an example with $d=3$.

\nS Let $\sigma \in \SF_3$ be a permutation of $(0,1,2)$. For example
\[
 \sigma = 
 \begin{pmatrix}
  0 & 1 & 2 \\
  2 & 0 & 1
 \end{pmatrix}
\]
which means
\[
 \sigma(0)=2, \qquad \sigma(1)=0, \qquad \sigma(2)= 1
\]
This \kw{permutation} is applied on the tensor product as\footnote{Some authors consider $\sigma(\abold_0 \otimes \abold_1 \otimes \abold_2) = \abold_{\sigma^{-1}(0)} \otimes \abold_{\sigma^{-1}(1)} \otimes \abold_{\sigma^{-1}(2)}$}
\begin{equation}
 \sigma(\abold_0 \otimes \abold_1 \otimes \abold_2) = \abold_{\sigma(0)} \otimes \abold_{\sigma(1)} \otimes \abold_{\sigma(2)}
\end{equation}
Then
\begin{equation}
 \begin{CD}
  E_0 \otimes E_1 \otimes E_2 @>\sigma >> E_{\sigma(0)} \otimes E_{\sigma(1)} \otimes E_{\sigma(2)}
 \end{CD}
\end{equation}
or, for $\sigma$ as above
\begin{equation}
 \sigma(\abold \otimes \bbold \otimes \cbold) = \cbold \otimes \abold \otimes \bbold
\end{equation}
Note that here, by a slight abuse of notations, we confound the permutation $\sigma \in \SF_3$ and the operator which acts on $\abold \otimes \bbold \otimes \cbold$.

\nS Let
\[
\abold \in E, \qquad \bbold \in F, \qquad \cbold \in G
\]
The permutation operator $\sigma$ is extended linearly to any tensor $\Abold \in E \otimes F \otimes G$ as (for $\sigma$ as above)
\begin{equation}\label{eq:permut}
 \begin{aligned}
   \Abold_\sigma & = \sigma(\Abold) \\
   &= \sum_{i,j,k} \alpha_{ijk} \; \sigma(\ebold_i \otimes \fbold_j \otimes \gbold_k) \\
   &= \sum_{i,j,k} \alpha_{ijk} \; \gbold_k \otimes \ebold_i \otimes \fbold_j \\
   & \in G \otimes E \otimes F
 \end{aligned}
\end{equation}
or
\begin{equation}
    \begin{CD}
      E \otimes F \otimes G @>\sigma>> G \otimes E \otimes F \\
      \Abold @>>> \Abold_\sigma
    \end{CD}
\end{equation}

\nS Here we see the impact of a permutation on the indices. It is presented on the example above. Permutation $\sigma$ can be transferred to indices as
\begin{equation}
 \sigma = 
 \begin{pmatrix}
  i & j & k \\
  k & i & j
 \end{pmatrix}
 \qquad \mbox{and} \qquad
 \sigma^{-1} = 
 \begin{pmatrix}
   i & j & k \\
   j & k & i
 \end{pmatrix}
\end{equation}
We then have 
\begin{equation}\label{eq:tenselem:permut:1}
 \begin{array}{lcl}
   \sigma(\Abold) &=& \displaystyle \sum_{i,j,k} \alpha_{ijk} \; \gbold_k \otimes \ebold_i \otimes \fbold_j \\
   &=& \displaystyle \sum_{i,j,k} \beta_{kij} \; \gbold_k \otimes \ebold_i \otimes \fbold_j \qquad \mbox{definition of} \: \beta_{kij}
 \end{array}
\end{equation}
with
\begin{equation}
 \begin{array}{lcl}
   \beta[k,i,j] &=& \alpha[i,j,k] \\
   &=& \alpha[\sigma^{-1}(k),\sigma^{-1}(i),\sigma^{-1}(j)]
 \end{array}
\end{equation}


\nS \textbf{Permutation for $d-$tensors:} \index{permutation}This is easily extended to tensor products of $d$ vector spaces. Let $\sigma \in \SF_d$. Then
\begin{equation}
 \sigma\left(\bigotimes_\mu \abold_\mu\right) = \bigotimes_\mu \abold_{\sigma(\mu)}
\end{equation}
or\footnote{Just to convince the reluctant reader that compact formula as above are easier to read ...}
\begin{equation}
 \sigma(\abold_1 \otimes \ldots \otimes \abold_\mu \otimes \ldots \otimes \abold_d) = \abold_{\sigma(0)} \otimes \ldots \otimes \abold_{\sigma(\mu)} \otimes \ldots \otimes \abold_{\sigma(d)}
\end{equation}
which can be extended by linearity
\begin{equation}
 \sigma\left(\sum_\ibold \alpha_\ibold \Ebold_\ibold\right) = \sum_\ibold \alpha_\ibold \; \sigma(\Ebold_\ibold)
\end{equation}
where
\[
 \begin{aligned}
   \sigma(\Ebold_\ibold) &= \sigma(\ebold_{1i_1} \otimes \ldots \otimes \ebold_{di_d}) \\
   &= \ebold_{\sigma(1)i_{\sigma(1)}} \otimes \ldots \otimes \ebold_{\sigma(d)i_{\sigma(d)}}
 \end{aligned}
\]


\nS Let us assume that each space $E_\mu$ is the same space $E$. Let us denote
\begin{equation*}
 E^{\otimes d} = \underbrace{E \otimes \ldots \otimes E}_{d \: \mathrm{times}}
\end{equation*}
Then
\[
\begin{CD}
  E^{\otimes d} @>\sigma >> E^{\otimes d}
\end{CD}
\]
can be seen as the \kw{action of the group of permutations} $\SF_d$ on tensor space $ E^{\otimes d}$:
\begin{equation}
\begin{CD}
  \SF_d \times  E^{\otimes d} @>\sigma>> E^{\otimes d}
\end{CD}
\end{equation}
Indeed
\begin{equation}
 \sigma(\sigma'(\Abold)) = (\sigma\sigma')(\Abold)
\end{equation}
which can be seen easily on elementary tensors and extended to any tensor by linearity. Note that with the convention $\sigma(\bigotimes_\mu \abold_\mu)= \bigotimes \abold_{\sigma^{-1}(\mu)}$ and not $\bigotimes \abold_{\sigma(\mu)}$, we would have 
\begin{equation}
 \begin{array}{lcl}
  \displaystyle\sigma\left(\sigma'\left(\bigotimes_\mu \abold_\mu\right)\right) &=& \displaystyle\sigma\left(\bigotimes_\mu \abold_{\sigma^{'-1}\mu}\right) \\
  &=& \displaystyle\bigotimes_\mu \abold_{\sigma^{-1}\sigma^{'-1}\mu} \\
  &=& \displaystyle\bigotimes_\mu \abold_{(\sigma'\sigma)^{-1}\mu} \\
  &=& \displaystyle (\sigma'\sigma)\left(\bigotimes_\mu \abold_\mu\right)  \\
  & \neq& \displaystyle (\sigma\sigma')\left(\bigotimes_\mu \abold_\mu\right)
 \end{array}
\end{equation}

\nS If in $\bigotimes_\mu E_\mu$ the spaces $E_\mu$ are different, it still is possible to build the action of $\SF_d$, by considering the space
\[
 \mathfrak{E} = \bigoplus_{\sigma \in \SF_d} \: E_{\sigma(1)} \otimes \ldots \otimes E_{\sigma(d)}
\]
and the action
\[
 \begin{CD}
  \SF_d \times \mathfrak{E} @>\sigma>> \mathfrak{E}
 \end{CD}
\]

%
\section{Symmetric tensors}\label{sec:tenselem:sym}
%

This leads to the definition of symmetric and antisymmetric tensors. These definitions make sense if the tensors belong to a tensor power of a given vector space $E$ only. For example, the definition of a symmetric matrix makes sense for square matrices only. Let
\[
 \Abold \in E^{\otimes d} = \underbrace{E \otimes \ldots \otimes E}_{d \: \mathrm{times}}
\]
It is a \kw{symmetric tensor}\index{tensor!symmetric} if it is invariant by any permutation
\begin{equation}
 \forall \: \sigma \in \SF_d, \quad \sigma(\Abold) = \Abold
\end{equation}
The \kw{symmetric part of a tensor} is defined as
\begin{equation}
 \mathrm{Sym\:}\Abold = \frac{1}{n!}\sum_{\sigma \in \SF_d}\sigma(\Abold)
\end{equation}
Then, a symmetric tensor $\Abold$ is such that
\[
  \Abold = \mathrm{Sym\:}\Abold
\]

\nS Let $(\abold_1, \ldots,\abold_d)$ be a family of vectors with $\abold_\mu \in E$. The \kw{symmetric product} of the vectors is the symmetric part of their tensor product. There is no standard notation for it, and we use the dot here
\begin{equation}
 \abold_1. [\ldots] .\abold_d = \frac{1}{n!}\sum_{\sigma \in \SF_d} \abold_{\sigma(1)} \otimes \ldots \otimes \abold_{\sigma(d)}
\end{equation}
For example\footnote{One can see here why tensors must belong to a $d-$tensor power of a given vector space. Indeed, the sum $\abold \otimes \bbold + \bbold \otimes \abold$ is defined if both terms belong to a same vector space, here $E \otimes E$. The sum is not defined if $\abold \in E$ and $\bbold \in F$ with $E \neq F$.}
\[
 \abold.\bbold =\frac{1}{2} (\abold \otimes \bbold + \bbold \otimes \abold)
\]
If $\Abold, \Bbold$ are two tensors, their \kw{symmetric product} is defined as
\begin{equation}
 \Abold.\Bbold = \mathrm{Sym}\:\Abold \otimes \Bbold
\end{equation}
It is possible to show that the symmetric product is associative, and this leads to a structure of a commutative algebra.

%
\section{Antisymmetric tensors}\label{sec:tenselem:antisym}
%

\nT{Parity and signature of a permutation} Any permutation $\sigma \in \SF_d$ can be decomposed as a product of inversions (an inversion is a permutation where two elements are inverted, as in $abcde \longrightarrow adcbe$ where $(b,d)$ are inverted). A permutation is even (resp. odd) if the number of inversions in its decomposition is even (resp. odd). This does not depend on the decomposition. The signature of a permutation is $+1$ if it is even, and $-1$ if it is odd. For example, as $012 \longrightarrow 021 \longrightarrow 201$, the permutation $\sigma(0,1,2)= (2,0,1)$ is even, and its signature is $+1$.

\nS A tensor $\Abold \in E^{\otimes d}$ is an \kw{antisymmetric tensor}\index{tensor!antisymmetric} if
\begin{equation}
 \forall \: \sigma \in \SF_d, \quad \sigma(\Abold) = (-1)^s \Abold
\end{equation}
where $s$ is the signature of $\sigma$. For example, $a \otimes b - b \otimes a$ is antisymmetric. The \kw{antisymmetric part of a tensor} is defined as
\begin{equation}
 \mathrm{Antisym\:}\Abold = \frac{1}{n!}\sum_{\sigma \in \SF_d}(-1)^s\, \sigma(\Abold)
\end{equation}
An antisymmetric tensor $\Abold$ is such that
\[
 \Abold = \mathrm{Antisym\:}\Abold
\]


\nT{Wedge product}Let $(\abold_1, \ldots,\abold_d)$ be a family of vectors with $\abold_\mu \in E$. The \kw{wedge product}\index{product!wedge} of the vectors is the antisymmetric part of their tensor product. The standard notation is $\wedge$. It is called as well the \kw{exterior product}\index{product!exterior}. For example
\[
 \abold \wedge \bbold = \frac{1}{2} (\abold \otimes \bbold- \bbold \otimes \abold)
\]
If
\[
 \abold = (a_1, a_2), \qquad \bbold = (b_1,b_2)
\]
we have
\[
 \abold \otimes \bbold = 
 \begin{pmatrix}
  a_1b_1 & a_1b_2 \\
  a_2b_1 & a_2b_2
 \end{pmatrix}
 , \qquad 
 \bbold \otimes \abold = 
 \begin{pmatrix}
  a_1b_1 & a_2b_1 \\
  a_1b_2 & a_2b_2
 \end{pmatrix}
\]
and
\[
 \abold \wedge \bbold = \frac{1}{2}
 \begin{pmatrix}
  0 & a_1b_2 - a_2b_1 \\
  a_2b_1-a_1b_2 & 0
 \end{pmatrix}
\]

\nS The wedge product is anticommutative
\begin{equation}
 \bbold \wedge \abold = - \abold \wedge \bbold
\end{equation}
If $\abold$ and $\bbold$ are colinear, $\abold \wedge \bbold=0$, because $\abold \wedge \abold = - \abold \wedge \abold =0$, and $(\lambda \abold) \wedge \abold = \lambda (\abold \wedge \abold)=0$.  So, we have 
\begin{equation}
 \begin{array}{lcl}
  \abold \wedge \bbold &=& (a_1\ebold_1 + a_2\ebold_2) \wedge (b_1\ebold_1 + b_2\ebold_2) \\
  &=& (a_1b_1 \ebold_1 \wedge \ebold_1) + (a_1b_2 \ebold_1 \wedge \ebold_2) + (a_2b_1 \ebold_2 \wedge \ebold_1) + (a_2b_2 \ebold_2 \wedge \ebold_2) \\
  &=& (a_1b_2-a_2b_1) \; \ebold_1 \wedge \ebold_2
 \end{array}
\end{equation}
where
\begin{equation}
  \ebold_1 \wedge \ebold_2 =    \frac{1}{2}
 \begin{pmatrix}
  0 & 1\\
  -1 & 0
 \end{pmatrix}
\end{equation}
If two vectors $\abold$ and $\bbold$ are not colinear, they define a plane, and a parallelogram in this plane with two sides $\abold$ and $\bbold$. Then, $\|\abold \wedge \bbold\|$ is the area of this parallelogram.

\nT{Grassmann algebra} Let $E$ be a vector space on a field $\K$ and consider
\begin{equation}
 \Lambda(E) = \K \oplus E \oplus E^{\otimes 2} \oplus \ldots \oplus E^{\otimes n}
\end{equation}
Then, $(\Lambda(E), +, \wedge)$ is an algebra, called the \kw{Grassmann algebra}\index{algebra!Grassmann} or the \kw{exterior algebra}\index{algebra!exterior} of $E$.

%
\section{Slicing}\label{sec:tenselem:slice}\index{slicing@\textbf{slicing}}
%

Slicing is \emph{par excellence} the elementary operation on tensors which is clearly understood visually, and is often presented likewise, sometimes with colors. Here, I develop an algebraic presentation of slicing, and finish by a visual presentation for $3-$modes tensors.

\nS I first present it on matrices. Let 
\[
 A = \sum_{i,j}\, \alpha_{ij} \, \ebold_i \otimes \fbold_j, \qquad \mbox{with} \quad 
 \begin{cases}
   A & \in  E \otimes F  \\
   \dim E &= m \\
   \dim F &= n
 \end{cases}
\]
Slicing in $A$ is selecting a row or a column. Rows are vectors $\abold_{i*} \in F$ and columns are vectors $\abold'_{*j} \in E$. This is a two steps procedure:
\begin{itemize}[label=$\rightarrow$]
 \item selecting the mode: $E$ for columns or $F$ for rows
 \item once the mode has been selected, select an index $i \in \llbracket 1,m\rrbracket$ if rows have been selected or $j \in \llbracket 1,n \rrbracket$ if columns have been selected.
\end{itemize}
This gives all the information needed to build $\abold_{i*}$ or $\abold'_{*j}$. 


\nS Let us extend this to $3-$modes tensors. Let 
\[
 \Abold = \sum_{i,j,k} \alpha_{ijk} \; \ebold_i \otimes \fbold_j \otimes \gbold_k
\]
Fixing an index and letting the others run yields a \kw{slice}. A slice of $\Abold$ is defined by
\begin{itemize}[label=$\rightarrow$]
 \item selecting the mode for which the index is fixed
 \item selecting the value for this index
\end{itemize}
As there as three modes, there are three sets of slices,\\
\\
\begin{center}
\begin{tabular}{ccclc}
Selected mode & Fixed index & Running indices & Slice & dimensions\\
\hline 
$E$ & $i$ & $j,k$ & $ \displaystyle \sum_{j,k} \alpha_{ijk} \; \fbold_j \otimes \gbold_k$ & $n \times p$\\
$F$ & $j$ & $i,k$ & $ \displaystyle \sum_{i,k} \alpha_{ijk} \; \ebold_i \otimes \gbold_k$ & $m \times p$\\
$G$ & $k$ & $i,j$ & $ \displaystyle \sum_{i,j} \alpha_{ijk} \; \ebold_i \otimes \fbold_j$ & $m \times n$\\
\end{tabular} 
\end{center}


\nS It can be easily extended to $d-$tensors. Let 
\[
 \Abold = \sum_{i_1} \ldots \sum_{i_d}\alpha_{i_1 \ldots i_d} \; \ebold_{1i_1}  \otimes \ldots \otimes \ebold_{di_d}
\]
As for $3-$modes tensors, fixing an index and letting the others run yields a \kw{slice}. A slice of $\Abold$ is defined by
\begin{itemize}[label=$\rightarrow$]
 \item selecting the mode $\mu$ for which the index is fixed
 \item selecting the value $i_\mu$ selected for this index
\end{itemize}
Let us define
\[
 \ibold \backslash \mu := (i_1, \ldots, i_{\mu-1},i_{\mu+1}, \ldots , i_d)
\]
i.e. the multi-index $\ibold$ but index $i_\mu$, easily extended to $\mu=1$ and $\mu=d$. Then, fixing the index $i_\mu$ for mode $\mu$ yields the $\mu-$slice $\Abold_{\mu,i_\mu}$ defined by 
\begin{equation}
 \begin{array}{lcl}
   \Abold_{\mu,i_\mu} &=& \displaystyle \sum_{\ibold \backslash \mu} \alpha_\ibold \, \Ebold_{\ibold \backslash \mu } \\
   &=& \displaystyle \sum_{i_1, \ldots, i_{\mu-1},i_{\mu+1}, \ldots , i_d} \: \alpha_{i_1 \ldots i_\mu \ldots i_d}\; \ebold_{1,i_1} \otimes \ldots \otimes \ebold_{\mu-1,i_{\mu-1}} \otimes \ebold_{\mu+1,i_{\mu+1}} \otimes \ldots \otimes \ebold_{di_d} 
 \end{array}
\end{equation}
In $\Abold_{\mu, i_\mu}$, $\mu$ denotes the mode with fixed index, and $i_\mu$ the value of the fixed index.

\nS It can be extended towards fixing indices in a subset of modes. Let
\[
 J \subset I
\]
with $|J|=d' < d$. If $\overline{J}= I\backslash J$, we have
\[
 I = J \sqcup \overline{J}
\]
Let
\[
 \left|
 \begin{array}{cll}
 J &=(j_1, \ldots, j_{d'}) &\\
 \\
 \ibold_\j &\quad \mbox{denote} \quad (i_{j_1}, \ldots, i_{j_{\nu}}, \ldots,  i_{j_{d'}}) \quad &\mbox{for} \quad j_\nu \in J \\
 \\
 \ibold \backslash \j &\quad \mbox{denote} \quad (i_{k_1}, \ldots, i_{k_{d-d'}}) \quad &\mbox{for} \quad k_\nu \notin J
 \end{array}
 \right.
\]
(the set $J$ is the subset of indices $\mu$ in $J$, $\ibold_\j$ the subset of multi-index $\ibold$ with indices in $J$ and $\ibold \backslash \j$ the subset of multi-index with indices not in $J$). Then
\begin{equation}
 \Abold_{\j, \ibold_\j } = \sum_{\ibold \backslash \j } \; \alpha_\ibold \; \Ebold_{\ibold \backslash \j}
\end{equation}
with
\[
 \Abold_{\j, \ibold_\j } \in \bigotimes_{\mu \notin J}E_\mu
\]
$\Abold_{\j, \ibold_\j }$ is obtained by
\begin{itemize}[label=$\rightarrow$]
    \item fixing indices in $J$, a subset of modes, one per mode
    \item summing up coefficients in $A$ over all indices in other modes $(I\backslash J)$
\end{itemize}

\nT{Example} It is easily understood on an example. Let
\[
 d=5, \qquad J = \{2,4\}, \qquad \mbox{fix} \quad i_2,i_4
\]
with
\[
 \Abold = \sum_{i_1,i_2,i_3,i_4,i_5} \, \alpha_{i_1i_2i_3i_4i_5} \; \ebold_{1,i_1} \otimes \ebold_{2,i_2} \otimes \ebold_{3,i_3} \otimes \ebold_{4,i_4} \otimes \ebold_{5,i_5}
\]
Then
\[
 \left|
 \begin{array}{cl}
 J &=\{2,4\} \\
 \\
 \ibold_\j &\quad \mbox{denotes} \quad (i_2,i_4) \\
 \\
 \ibold \backslash \j &\quad \mbox{denotes} \quad (i_1,i_3,i_5)
 \end{array}
 \right.
\]
and
\[
 \Abold_{(2,4),(i_2,i_4)} = \sum_{i_1,i_3,i_5} \; \alpha_{i_1i_2i_3i_4i_5} \; \ebold_{1,i_1} \otimes \ebold_{3,i_3} \otimes \ebold_{5,i_5}
\]

\nT{An even more concrete example} Let $\Abold \in E \otimes F \otimes G$. Let us select
\[
 J = (F,G), \qquad \mbox{fix } j,k
\]
Then, the slice is called a fiber, and is
\begin{equation}
 \Abold_{(\f,\g),(j,k)} = \sum_i \, \alpha_{ijk} \, \ebold_i \in E
\end{equation}

\nS And here is the promise: visualisation of the extraction of slice with index $i_2$ on mode $\mu=2$ from a $3-$modes tensor.\\
\\
\begin{center}
    \begin{tikzpicture}
       \fill[green!40] (0,0,2)--(2,0,2)--(2,2,2)--(0,2,2)--cycle ; 
       \fill[green!20] (0,2,2)--(0,2,0)--(2,2,0)--(2,2,2)--cycle ;
       \fill[green!60] (2,0,2)--(2,0,0)--(2,2,0)--(2,2,2)--cycle; 
       \fill[red!40] (1,0,2)--(1.3,0,2)--(1.3,2,2)--(1,2,2)--cycle;
       \fill[red!20] (1,2,2)--(1.3,2,2)--(1.3,2,0)--(1,2,0)--cycle;
       \draw[->] (3,1,.5) -- (5,1,.5) ; 
       \node () at (4,1.5,.5){$\mu=2$} ;
       \node () at (-.5,1,2){$\mu=1$};
       \node () at (1,-.5,2){$\mu=2$};
       \node () at (-.5,2.2,1){$\mu=3$};
       \node () at (7.75,-.5,2){$i_2$};
       \fill[red!60] (8,0,2)--(8,0,0)--(8,2,0)--(8,2,2) --cycle ; 
       \fill[red!40] (8,0,2)--(8,2,2)--(7.5,2,2)--(7.5,0,2) -- cycle ; 
       \fill[red!20] (8,2,2)--(7.5,2,2)--(7.5,2,0)--(8,2,0)--cycle ;
    \end{tikzpicture}
\end{center}

\nT{Note} Slicing is intimately link with marginalization in statistics. Let us have a bivariate distribution over set $\Lambda_1 \times \Lambda_2$. Let $p_{ij}$ be defined as $p_{ij}= \P(X_1=i,X_2=j)$. Then, the marginal over $\Lambda_2$ for value $i$ in $\Lambda_1$ is defined as $M_1(i)=\sum_j p_{ij}$. The distribution of the marginal is a row of matrix $P$, i.e. a slice. This can be naturally extended to $d-$variate distributions over $\Lambda_1 \times \ldots \times \Lambda_d$ with tensor $\Pbold$ defined by $\Pbold_{i_1\ldots i_d}=\P(X_1=i_1,\ldots,X_d=i_d)$. Marginal distributions are the sums of elements in a slice of $\Pbold$.

%
\section{Kronecker product, vectorization and reshaping}\label{sec:tenselem:reshape}
%

Reshaping is a very simple operation consisting in aggregating modes. It is better understood and formalized with vectorization in mind. Vectorization is better understood and formalized with Kronecker product between vectors, which I present now.

\nT{Kronecker product between two vectors} \index{Kronecker product between vectors}Let
\[
\abold = 
\begin{pmatrix}
  a_1 \\
  \vdots \\
  a_m
\end{pmatrix}
\in \R^m, \qquad 
 \bbold = 
\begin{pmatrix}
  b_1 \\
  \vdots \\
  b_n
\end{pmatrix}
\in \R^n, 
\]
Then, the \kw{Kronecker product} between $\abold$ and $\bbold$, denoted $\otimes_\k$, is defined by \
\begin{equation}
 \abold \otimes_\k \bbold = 
 \begin{pmatrix}
  a_1\bbold \\
  a_2\bbold \\
  \vdots \\
  a_m\bbold
 \end{pmatrix}
\quad \in \R^{mn}
\end{equation}
For example,
\begin{equation}
    \begin{pmatrix}
      a \\
      a'
    \end{pmatrix}
    \otimes_\k
    \begin{pmatrix}
      b \\
      b' \\
      b''
    \end{pmatrix}
    =
    \begin{pmatrix}
      ab \\
      ab' \\
      ab'' \\
      a'b \\
      a'b' \\
      a'b''
    \end{pmatrix}
\end{equation}
The Kronecker product between matrices is presented in section \ref{sec:alstruct:prodmat}. Let us present it here componentwise. Let
\begin{equation}
    A = 
    \begin{pmatrix}
      a_{11} & a_{12} & \ldots & a_{1n} \\
      a_{21} & a_{22} & \ldots & a_{2n} \\
      \vdots & \vdots & \ddots & \vdots \\
      a_{m1} & a_{m2} & \ldots & a_{mn} \\
    \end{pmatrix}
\end{equation}
and $X$ be a $p \times q$ matrix. Then, $A \otimes_\k X$ is a $mp \times nq$ matrix which is written blockwise as
\begin{equation}
    A = 
    \begin{pmatrix}
      a_{11}X & a_{12}X & \ldots & a_{1n}X \\
      a_{21}X & a_{22}X & \ldots & a_{2n}X \\
      \vdots & \vdots & \ddots & \vdots \\
      a_{m1}X & a_{m2}X & \ldots & a_{mn}X \\
    \end{pmatrix}
\end{equation}
One can check that, if dimensions are consistent
\begin{equation}
    (A \otimes_\k B).(\xbold \otimes_\k \ybold) = A\xbold \otimes_\k B\ybold
\end{equation}
Here, a vector in $\R^m$ is considered as a $m \times 1$ matrix, or $\abold \in \R^{m \times 1}$.

\nT{Vectorization of a matrix} Let us consider two vectors $\abold \in \R^m$, $\bbold \in \R^n$. We define the \emph{vec} application on elementary matrices (rank one) by 
\begin{equation}
 \begin{CD}
   \R^m \otimes \R^n @>\mathrm{vec}>> \R^{mn} \\
   \abold \otimes \bbold @>>> \abold \otimes_\k \bbold
 \end{CD}
\end{equation}
and extend it to all matrices in $\R^m \otimes \R^n$ by linearity as
\begin{equation}
 \begin{CD}
  \displaystyle A = \sum_{i,j}\, a_{ij}\, \ebold_i \otimes \fbold_j @>\mathrm{vec}>> \displaystyle \sum_{i,j} a_{ij} \, \ebold_i \otimes_\k \fbold_j \in \R^{mn}
 \end{CD}
\end{equation}
It is easy to show that \emph{vec} is an isomorphism (the vectors $(\ebold_i \otimes_\k \fbold_j)_{i,j}$ form a basis in $\R^{mn}$). Some times, we will write
\begin{equation}
 \R^m \otimes_\k \R^n  \simeq \R^{mn}
\end{equation}
to say explicitly that $\R^{mn}$ is the target space of \emph{vec}. This is a slight abuse of notations, because $\otimes_\k$ is no defined between spaces. This can be solved by defining $E \otimes_\k F$ as the vector space spanned by Kronecker products $\abold \otimes_\k \bbold$ where $\abold \in E$ and $\bbold \in F$. The isomorphisms are as follows:\\
\begin{center}
 \begin{tikzcd}
   \R^m \times \R^n \ar[r, "\otimes"] \ar[ddr, "\otimes_\k"] & \R^m \otimes \R^n \ar[dd, "\mathrm{vec}"]\\
   &\\
   & \R^{mn}
  \end{tikzcd}
\end{center}
So
\begin{equation}
\abold \otimes_\k \bbold = \vec (\abold \otimes \bbold) 
\end{equation}

\nT{Example} Let for example
\begin{equation}
 A = 
 \begin{pmatrix}
  a & b \\
  c & d
 \end{pmatrix}
\end{equation}
Then
\begin{equation}
 A = a \, \ebold_1 \otimes \fbold_1 + b \, \ebold_1 \otimes \fbold_2 + c \, \ebold_2 \otimes \fbold_1 + d \, \ebold_2 \otimes \fbold_2 
\end{equation}
We have
\begin{equation}
  \ebold_1 \otimes_\k \fbold_1 = 
  \begin{pmatrix}
   1 \\
   0 \\
   0 \\
   0
  \end{pmatrix}
  , \quad
  \ebold_1 \otimes_\k \fbold_2 = 
  \begin{pmatrix}
   0 \\
   1 \\
   0 \\
   0
  \end{pmatrix}
   , \quad
  \ebold_2 \otimes_\k \fbold_1 = 
  \begin{pmatrix}
   0 \\
   0 \\
   1 \\
   0
  \end{pmatrix}
   , \quad
  \ebold_2 \otimes_\k \fbold_2 = 
  \begin{pmatrix}
   0 \\
   0 \\
   0 \\
   1
  \end{pmatrix}
\end{equation}
and
\begin{equation}
 \vec A = a \, \ebold_1 \otimes_\k \fbold_1 + b \, \ebold_1 \otimes_\k \fbold_2 + c \, \ebold_2 \otimes_\k \fbold_1 + d \, \ebold_2 \otimes_\k \fbold_2 
\end{equation}
or
\begin{equation}
 \vec A = 
 \begin{pmatrix}
  a \\
  b \\
  c \\
  d
 \end{pmatrix}
\end{equation}
More generally, the vectorization of a $m \times n$ matrix is a $mn$ vector with, as coefficients, the rows of $A$ stacked from top to bottom. This simple example shows the economy of means obtained by using tensor algebra instead of component-wise calculations. 


\nT{Basic algebra of \emph{vec} operator} Some elementary properties of the \emph{vec} operator are given here without proof. See \cite[section 7.5]{Schott1997} for details. If $A,B$ are two matrices
\begin{equation}
 \langle A,B\rangle = \Tr A^\t B = \langle \vec A \, , \, \vec B \rangle
\end{equation}
Let $A \in \R^{m \times n}, B \in \R^{n \times p}, C \in \R^{p \times q}$. Then
\begin{equation}
 \vec ABC = (C^\t \otimes_\k A)\vec B
\end{equation}
This can be extended to more than three matrices (see \cite[p.264]{Schott1997}).

\nT{Vectorization of a tensor}This can easily be extended to tensors, by observing that Kronecker product between vectors is associative. Let us explicit it for the Kronecker product of three vectors. Let
\begin{equation}
\begin{cases}
  \abold & \in \R^m \\
  \bbold & \in \R^n \\
  \cbold & \in \R^p
\end{cases}
\end{equation}
Then
\begin{equation}
 \abold \otimes_\k \bbold \otimes_\k \cbold \in \R^{mnp}
\end{equation}
with
\begin{equation}
  \abold \otimes_\k \bbold \otimes_\k \cbold = \abold \otimes_\k (\bbold \otimes_\k \cbold) = (\abold \otimes_\k \bbold) \otimes_\k \cbold
\end{equation}
This enables to define an isomorphisms \emph{vec} as
\begin{equation}
 \begin{CD}
   \R^{m \times n \times p} @>\mathrm{vec}>> \R^{mnp} \\
   \abold \otimes \bbold \otimes \cbold @>>> \abold \otimes_\k \bbold \otimes_\k \cbold
 \end{CD}
\end{equation}
extended to $\R^{m \times n \times p}$ by linearity. This can be extended to $d-$modes tensors as
\begin{equation}
 \begin{CD}
   \bigotimes_\mu \R^{n_\mu} @>\mathrm{vec}>> \R^{n_1\ldots n_d} \\
   \bigotimes_\mu \abold_\mu @>>> \left(\bigotimes_\k\right)_\mu \abold_\mu
 \end{CD}
\end{equation}
and extended to $\bigotimes_\mu \R^{n_\mu}$ by linearity.


\nT{Reshaping on an example} Reshaping is a sort of intermediate step in vectorization. Let us start with an example. Let $d=5, I = \{1,5\}$ and 
\begin{equation}
 \Abold \in E_1 \otimes E_2 \otimes E_3 \otimes E_4 \otimes E_5
\end{equation}
with $\dim E_\mu = n_\mu$. Let us define
\begin{equation}
 I_1 = (1,2), \qquad I_2 = 3, \qquad I_3 = (4,5)
\end{equation}
(the order of indices of spaces is respected) with
\begin{equation}
 I = I_1 \sqcup I_2 \sqcup I_3
\end{equation}
and
\begin{equation}
 \Ebold_1 = E_1 \otimes E_2, \qquad \Ebold_2 = E_3, \qquad \Ebold_3 = E_4 \otimes E_5
\end{equation}
Then, reshaping is a two steps process: 
\begin{enumerate}
 \item writing
\begin{equation}
 \Abold \in (E_1 \otimes E_2) \otimes E_3 \otimes (E_4 \otimes E_5) := \Ebold_1 \otimes \Ebold_2 \otimes \Ebold_3
\end{equation}
so
\begin{equation}
   \Abold =  \displaystyle \sum_{i_1\ldots i_5}\, a_{i_1\ldots i_5} \, (\ebold_1 \otimes \ebold_2) \otimes \ebold_3 \otimes (\ebold_4 \otimes \ebold_5 )
\end{equation}
\item applying the \emph{vec} operator on each component $\Ebold_1,\Ebold_2,\Ebold_3$. Hence, if $\Abold'$ is the desired reshaping of $\Abold$, \begin{equation}
    \begin{array}{lcl}
        \Abold' &=& \displaystyle \sum_{i_1\ldots i_5}\, a_{i_1\ldots i_5} \, \vec (\ebold_1 \otimes \ebold_2) \otimes \ebold_3 \otimes \vec (\ebold_4 \otimes \ebold_5 )\\
        &=& \displaystyle \sum_{i_1\ldots i_5}\, a_{i_1\ldots i_5} \, (\ebold_1 \otimes_\k \ebold_2) \otimes \ebold_3 \otimes (\ebold_4 \otimes_\k \ebold_5 )
    \end{array}
\end{equation}
\end{enumerate}
Then
\begin{equation}
 \Abold' \in (E_1 \otimes_\k E_2) \otimes E_3 \otimes (E_4 \otimes_\k E_5)
\end{equation}
or, with dimensions
\begin{equation}
 \Abold' \in \R^{n_1n_2} \otimes \R^{n_3} \otimes \R^{n_4n_5}
\end{equation}

\nT{Definition of reshaping} Let us now give the general definition. Let us recall that $d$ is the number of modes and $I = \llbracket 1,d\rrbracket$. Let us consider a partition of $I$
\[
 I = \bigsqcup_\alpha I_\alpha
\]
respecting the order of indices in $I$, like
\[
 \{1,6\} = \{1,2\} \sqcup \{3\} \sqcup \{4,5,6\}
\]
but not
\[
 \{1,6\} = \{1,3\} \sqcup \{4\} \sqcup \{2,5,6\}
\]
This means that all indices in any $I_\alpha$ must be consecutive. 
Let us define
\begin{equation}
 \Ebold_\alpha = \bigotimes_{\mu \in I_\alpha} E_\mu
\end{equation}
Then, reshaping $\Abold \in \bigotimes_\mu E_\mu$ consists in
\begin{enumerate}
    \item writing
    \begin{equation}
        \Abold \in \bigotimes_{\alpha} \Ebold_\alpha
    \end{equation}
    \item applying the operator \emph{vec} in each component $\Ebold_\alpha$ as
    \begin{equation}
        \begin{CD}
            \bigotimes_{\mu \in I_\alpha} \ebold_\mu @>>> \bigotimes_{\mu \in I_\alpha}^{\k} \ebold_\mu
            \end{CD}
    \end{equation}
\end{enumerate}

\nT{Matricization} Matricization \marginpar{\emph{matricization}}\index{matricization} is a reshaping where $I$ is split into 2 subsets. Let $\Abold \in E \otimes F \otimes G$ with
\[
 \Abold = \sum_{i,j,k} \, \alpha_{ijk} \, \ebold_i \otimes \fbold_j \otimes \gbold_k
\]
Let us define
\begin{equation}
 \{1,2,3\} = \{1\} \sqcup \{2,3\}
\end{equation}
Then
\begin{equation}
 \Ebold_1 = E, \qquad \Ebold_2 = F \otimes G
\end{equation}
and
\begin{equation}
 \Abold \in \Ebold_1 \otimes \Ebold_2
\end{equation}
with
\begin{equation}
  \Abold = \sum_{i,jk} \, \alpha_{i,jk} \, \ebold_i \otimes (\fbold_j \otimes \gbold_k)
\end{equation}
where
\begin{itemize}[label=$\rightarrow$]
 \item $jk$ runs over $\llbracket 1,n\rrbracket \times \llbracket 1,p\rrbracket$
 \item $\fbold_j \otimes \gbold_k$ is a basis of $\Ebold_2$
\end{itemize}
Matricization is the following map 
\begin{equation}
 \begin{CD}
  E \otimes F \otimes G @>>> E \otimes (F \otimes G)
 \end{CD}
\end{equation}
or, component-wise 
\begin{equation}
 \begin{CD}
   \R^{m \times n \times p} @>>> \R^{m \times np}
 \end{CD}
\end{equation}

\nT{Note} As much is known about linear algebra of matrices, matricization is an important operation for deciphering the structure of tensors. Hence, it is developed in section \ref{sec:contract:alamat}. For example, it leads directly to HOSVD as a best low rank approximation of a tensor.

%
\chapter{Contraction}\label{chap:contract}
%

Matrix $\times$ vector and matrix $\times$ matrix products are classical operations in matrix algebra. They are two elementary examples of a much wider operator between two tensors: the contraction. Component-wise, contraction is a summation over repeated indices. For example, in matrix $\times$ vector product $C=AB$
\[
\begin{CD}
\left.
   \begin{array}{ll}
     A &=(a_{ij})_{i,j} \\
     B & =(b_{jk})_{j,k} 
    \end{array}
\right\} 
@>>> \displaystyle C = (c_{ik})_{i,k} \qquad \mbox{with} \quad c_{ik}=\sum_j a_{ij}b_{jk}
\end{CD}
\]
The index $j$ is common to $A$ and $B$, and the sum is over $j$ (see section \ref{sec:matelem:contract}). A family of contractions over a set of tensors which share some modes is called a \kw{tensor network}.

\nB There exists a coordinate-free definition of contraction: it is the dual operation of tensor product. This is far from being trivial: tensor product and contraction are two dual ways to navigate between tensors of different orders: tensor product $\otimes$ expands the order (here, $\mathrm{ord}\: \Xbold$ is the order of $\Xbold$), as $\mathrm{ord}\: \Abold \otimes \Bbold = \mathrm{ord}\: \Abold + \mathrm{ord}\: \Bbold$, and contraction shrinks the orders, as $\mathrm{ord}\: \Abold \bullet \Bbold = \mathrm{ord}\: \Abold - \mathrm{ord}\: \Bbold$. This is a reason why it is so important in tensor algebra. Even if it is slightly more abstract, it is important to know that it exists, because it shows that, in our example, even if the matrix $C$ depends obviously on $A$ and $B$, hence the selected basis in $E,F$ and $G$, the result can be transported as such by a simple change of basis if basis in $E,F,G$ change. The invariant behind simply is the composition of linear operators expressed as matrices in a given basis, as $L_\c = L_\a \circ L_\b$.

\nB Let us develop this on a simple example, which is a nice (and not so complicated) approach in abstract algebra. Let $\abold,\xbold \in E$, $\bbold, \ybold \in F$ and $A \in E \otimes F \simeq \L(F,E)$. The dual (or transposed) operator $A^\t$ of $A$ is defined as
\[
\langle A.\ybold \, , \, \xbold \rangle_\e = \langle \ybold \, , \, A^\t.\xbold \rangle_\f 
\]
So
\[
\begin{CD}
A \in \L(F,E) @>>> A^\t \in \L(E,F)
\end{CD}
\]
Let us consider the linear map
\[
\begin{CD}
 F @>L_a>> E \otimes F \\
 \bbold @>>> \abold \otimes \bbold
\end{CD}
\]
It is a linear map associated to $\otimes$. Let us show that its transpose $L_a^\t$ is associated to a contraction. We have, on elementary tensors
\[
\langle L_a(\bbold) \, , \, \xbold \otimes \ybold \rangle = \langle \bbold \, , \, L_a^\t(\xbold \otimes \ybold) \rangle
\]
As 
\[
\begin{array}{lcl}
\langle L_a(\bbold) \, , \, \xbold \otimes \ybold \rangle &=& \langle \abold \otimes \bbold \, , \, \xbold \otimes \ybold \rangle \\
&=& \langle \abold,\xbold\rangle_\e \, \langle \bbold,\ybold\rangle_\f \\
&=& \langle \bbold \, , \, \langle \abold,\xbold\rangle \, \ybold \rangle 
\end{array}
\]
we have
\[ 
\begin{array}{lcl}
L_a^\t(\xbold \otimes \ybold) &=& \langle \abold,\xbold\rangle \, \ybold \\
&=& (\ybold \otimes \xbold).\abold
\end{array}
\]
which is a contraction of $\ybold \otimes \xbold$ (\marginpar{\dbend}not $\xbold \otimes \ybold$) by $\abold$. Indeed, we have seen (see section \ref{sec:matelem:contract}) that 
\[
 (\ybold \otimes \xbold).\abold = \langle \xbold,\abold\rangle \ybold
\]
This defines a \kw{contraction} over $E$ through extension by linearity 
\[
 \begin{CD}
  \begin{array}{c}
   F \otimes E \\
   \times \\
   E 
  \end{array}
  @>\mathrm{contraction}>> F
 \end{CD}
\]
by
\[
\begin{CD}
 (Y,\abold) @>>> Y.\abold
\end{CD}
\]
(matrix $\times$ vector product). So, $\otimes$ is an operator from $E \times F$ to the space of matrices in $E \otimes F$. Its dual operator is the contraction of a matrix in $F \otimes E$ with a vector in $E$.

\nT{Organization of the chapter} This chapter is organized as follows:
\begin{description}
\item[section \ref{sec:contract:3modes}] To go beyond matrix $\times$ vector and matrix $\times$ matrix product, different contractions are defined for $3-$modes tensors. Each one can be sketched by a diagram which is provided.
\item[section \ref{sec:contract:alamat}] The contraction operators in $3-$modes tensors are deeply connected with matricization of a $3-$modes tensor, i.e. how to associate linear maps to a tensor. This is developed in this section, with notations for all 8 possibilities
\item[section \ref{sec:contract:dmodes}] Once the basic notions have been given for $3-$modes tensors, contractions are extended to the general case of $d-$modes tensors, without much changes.
\item[section \ref{sec:contract:hasse}] Contractions for $3-$modes tensors can be organized and piped on Hasse diagram of some posets . A first link with networks diagrams is presented, and a sketch for extending these approaches to $d-$modes tensors are given.
\item[section \ref{sec:contract:bilin}] In the last section, the importance of contraction as an operator is emphasized by showing that any bilinear map, or more generally multilinear map, from a product of vector spaces on a vector space can be written as a combination of a contraction and a tensor product. This is a guise of the universal property of tensor product.
\end{description}

%
\section{Contractions for $3-$modes tensors}\label{sec:contract:3modes}
%

We have seen that algebraically speaking, a matrix $\times$ vector product is a contraction, which is the dual operation of tensor product. This can be extended to $d-$modes $\times$ $d'-$modes tensors products on elementary tensors, and extended to the whole tensor space products by linearity. Let us show it first for $d=3$ and $d'=1$ or $d'=2$. Let
\[
 \Abold \in E \otimes F \otimes G, \qquad \zbold \in G, \qquad Y \in F \otimes G
\]
The mode common to $\Abold$ and $\zbold$ is $G$, indexed by $k$. The modes common to $\Abold$ and $Y$ are $F,G$, indexed by $(j,k)$.Then, a contraction over $G$ (summation over $k$) can be done considering associativity of tensor product: $E \otimes F \otimes G = (E \otimes F) \otimes G$. To emphasize this approach, we denote $EF \overset{.}{=} E \otimes F$, and define the contraction
\[
 \begin{CD}
  \begin{array}{c}
   EF \otimes G \\
   \times \\
   G 
  \end{array}
  @>\mathrm{contraction}>> EF
 \end{CD}
\]
Let us first consider elementary tensors
\[
 \Abold = \abold \otimes \bbold \otimes \cbold, \qquad Y = \ybold \otimes \zbold
\]
This yields 
\begin{equation}
 \begin{array}{lcl}
  (\abold \otimes \bbold \otimes \cbold).\zbold &=&  (\overbrace{\abold \otimes \bbold}^{A \in E \otimes F} \otimes \cbold).\zbold \\
  &=& (A \otimes \cbold).\zbold \\
  &=& \langle \cbold,\zbold\rangle A \\
  &=& \langle \cbold,\zbold\rangle \; \abold \otimes \bbold
 \end{array}
\end{equation}
Similarly, with $E \otimes F \otimes G = E \otimes (F \otimes G)$
\begin{equation}
 \begin{array}{lcl}
  (\abold \otimes \bbold \otimes \cbold).(\ybold \otimes \zbold) &=&  (\abold \otimes \overbrace{\bbold \otimes \cbold}^{B \in F \otimes G}).\overbrace{(\ybold \otimes \zbold)}^{Y \in F \otimes G} \\
  &=& (\abold \otimes B).Y \\
  &=& \langle B,Y\rangle \abold \\
  &=& \langle \bbold \otimes \cbold \, , \, \ybold \otimes \zbold\rangle \; \abold \\
  &=& \langle \bbold,\ybold\rangle \langle \cbold,\zbold\rangle \; \abold
 \end{array}
\end{equation}
So, we have two contractions:
\[
 \begin{CD}
  \begin{array}{c}
   E \otimes F \otimes G\\
   \times \\
   G 
  \end{array}
  @>\mathrm{contraction}>> E \otimes F
 \end{CD}
 \qquad \mbox{and} \qquad 
 \begin{CD}
  \begin{array}{c}
   E \otimes F \otimes G\\
   \times \\
   F \otimes G 
  \end{array}
  @>\mathrm{contraction}>> E
 \end{CD}
\]
Extended by linearity, the contractions are
\begin{equation}
 \left\{ 
 \begin{array}{lcl}
   \Abold.\zbold &=& X \\
   \Abold.Y &=& \xbold
 \end{array}
 \right.
\end{equation}

\nT{Contraction diagram} Each contraction can be sketched by \kw{contraction diagram}\index{diagram!contraction}. Only two of them are provided here. We adopt the same convention for labelling the bars than for labeling the contraction: let $E=E_1$, $F=E_2$ and $G=E_3$. Then, a bar is labeled by the index $\mu$ of a space in $E_1 \otimes E_2 \otimes E_3$. This yields:\\
\\
\begin{center}
\begin{tikzpicture}
\node[tensor_ell] (abc) at (1,0) {$a \otimes b \otimes c$} ;  
\node[tensor_ell] (z) at (3,0) {$z$} ;  
\draw (-.5,0) -- (abc) node[pos=0.5, above] {$1$} ;
\draw (abc) -- (1,1) node[pos=0.5, right] {$2$} ;
\draw (abc) -- (z) node[pos=0.5, above] {$3$} ;
\node () at (4,0) {and} ;
\node[tensor_ell] (abc2) at (7,0) {$a \otimes b \otimes c$} ;  
\node[tensor_ell] (YY) at (9.5,0) {$y \otimes z$} ;  
\draw (5,0) -- (abc2) node[pos=0.5, above] {$1$} ;
\draw  ([yshift=0.5ex]abc2.east) -- ([yshift=.5ex]YY.west) node[pos=0.5, above] {$2$} ;
\draw  ([yshift=-.5ex]abc2.east) -- ([yshift=-.5ex]YY.west)node[pos=0.5, below] {$3$};
\end{tikzpicture} 
\end{center}
which by collapse ($=$ contraction over repeated modes) yields\\
\\
\begin{center}
\begin{tikzpicture}
\node[tensor_ell] (abc) at (1,0) {$a \otimes b \, \langle c, z\rangle$} ;  
\draw (-1,0) -- (abc) node[pos=0.5, above] {$1$} ;
\draw (abc) -- (1,1) node[pos=0.5, right] {$2$} ;
\node () at (4,0) {and} ;
\node[tensor_ell] (abc2) at (7,0) {$a \, \langle b,y \rangle \langle c,z\rangle$} ;  
\draw (5,0) -- (abc2) node[pos=0.5, above] {$1$} ;
\end{tikzpicture} 
\end{center}

%
\section{Associated linear maps and matricization}\label{sec:contract:alamat}
%

Here we introduce matricization which is a key tool in tensor algebra given the immensity of the domain of linear algebra (or matrix algebra, even if matrices are bilinear objects). It is simply another guise of contraction, as its expression component-wise in a given basis. We will use the tools and notations presented in section \ref{sec:tenselem:reshape} for Kronecker product and vectorization. The key ingredient in this section is a clear awareness of the distinction between a linear map and its expression as a matrix, i.e. between what is called here \kw{associated linear map} and \kw{matricization}.

\nS Let us consider the contraction (see section \ref{sec:contract:3modes} for notations)
\begin{equation}
 \begin{CD}
  \begin{array}{c}
   EF \otimes G \\
   \times \\
   G 
  \end{array}
  @>\mathrm{contraction}>> EF
 \end{CD}
\end{equation}
written
\begin{equation}
 \Abold.\zbold = X
\end{equation}
Knowing $\Abold$, the map 
\begin{equation}
 \begin{CD}
  G @>A_\g>> E \otimes F \\
  \zbold @>>> X
 \end{CD}
\end{equation}
is a linear map, called the linear map associated to $\Abold$ by mode $G$ (the mode over which the contraction operates), and denoted $A_\g$. 
\begin{equation}
 A_\g \in \L(G \, , \, E\otimes F)
\end{equation}
As the contraction is defined coordinate free, so is the linear map associated to $\Abold$ by mode $G$. Matricization is building the matrix associated to this linear map, denoted $A_\g$ as well. It can be defined by
\begin{equation}
A_\g.\zbold = \vec X 
\end{equation}
and has $mn$ rows and $p$ columns. This is made possible by the isomorphism between $E \otimes F$ and $E \otimes_\k F$ (see section \ref{sec:contract:3modes})\\
\\
\begin{center}
\begin{tikzcd}
  & E \otimes F \arrow[dd, "\mathrm{vec}"]\\
  G \arrow[ur, "A_\g"] \arrow[dr, "M_\g" ]& \\
  & E \otimes_\k F
\end{tikzcd}
\end{center}
where the matrix associated to $A_\g$ is denoted $M_\g$ in this diagram to emphasize its construction as
\begin{equation}
 M_\g = \mathrm{vec} \circ A_\g
\end{equation}


\nS We define similarily\begin{equation}
 \begin{CD}
  F\otimes G @>>> E \\
  Y @>>> \xbold
 \end{CD}
\end{equation}
denoted
\begin{equation}
 \begin{CD}
   F \otimes G @>A_{\f\g}>> E
 \end{CD}
\end{equation}
which can be written in a matrix form
\begin{equation}
 \xbold = A_{\f\g}(\vec Y)
\end{equation}
where $A_{\f\g}$ is a matrix with $p$ rows and $np$ columns.

\nS Depending on the way $E \otimes F \otimes G$ is split, such an approach enables to define 8 contractions / associated linear maps, presented in the following table:\\
\\
\begin{center}
\begin{tabular}{ccccll}
\hline
 Mode(s) & ALM & from & to & for an elementary tensor &\\
 \hline 
 $\emptyset$ & $A$ & $\K$ & $E \otimes F \otimes G$ & $\abold \otimes \bbold \otimes \cbold \bullet_\emptyset \lambda$ & = $\lambda \abold \otimes \bbold \otimes \cbold $ \\
 $E$ & $A_\e$ & $E$ & $F \otimes G$ & $(\abold \otimes \bbold \otimes \cbold) \bullet_1 \xbold $ & = $\langle \abold,\xbold\rangle \; \bbold \otimes \cbold$\\
 $F$ & $A_\f$ & $F$ & $E \otimes G$ & $(\abold \otimes \bbold \otimes \cbold) \bullet_2 \ybold $ & = $ \langle \bbold,\ybold\rangle \; \abold \otimes \cbold$\\
 $G$ & $A_\g$ & $G$ & $E \otimes F$ & $(\abold \otimes \bbold \otimes \cbold) \bullet_3 z $ & = $\langle \cbold,\zbold\rangle \; \abold \otimes \bbold$\\ 
 $E,F$ & $A_{\e\f}$ & $E \otimes F$ & $G$ & $(\abold \otimes \bbold \otimes \cbold) \bullet_{12} (\xbold \otimes \ybold) $ & = $\langle \abold,\xbold\rangle \langle \bbold,\ybold \rangle  \; \cbold$\\
 $E,G$ & $A_{\e\g}$ & $E \otimes G$ & $F$ & $(\abold \otimes \bbold \otimes \cbold) \bullet_{13} (\xbold \otimes \zbold) $ & = $\langle \abold,\xbold\rangle \langle \cbold,\zbold \rangle  \; \bbold$\\
 $F,G$ & $A_{\f\g}$ & $F \otimes G$ & $E$ & $(\abold \otimes \bbold \otimes \cbold) \bullet_{23} (\ybold \otimes \zbold) $ & = $\langle \bbold,\ybold\rangle \langle \cbold,\zbold \rangle  \; \abold$\\
 $E,F,G$ & $A_{\e\f\g}$ & $E \otimes F \otimes G$ & $\K$ & $(\abold \otimes \bbold \otimes \cbold) \bullet_{123} (\xbold \otimes \ybold \otimes \zbold) $ & = $\langle \abold,\xbold\rangle \langle \bbold,\ybold \rangle \langle \cbold,\zbold \rangle$\\ 
 \hline
\end{tabular}
\end{center}

\nT{Convention for notation} We adopt the convention that a contraction is denoted by the symbol $\bullet_\j$ where $J \subset I$ is the subset of indices of the spaces ( $=$ modes) over which the contraction operates\footnote{One could have adopted the convention that the symbol is $\bullet_\e$ when contraction is made over space $E$ in $E \otimes F \otimes G$. But this leads to a difficulty if $E=F=G$. Moreover, in more than 3 modes, the spaces are denoted $E_1 \otimes \ldots \otimes E_\mu \otimes \ldots \otimes E_d$, and it is natural to chose $\bullet_\mu$.}. 

\nT{Matricization} Any linear map between finite dimensional spaces can be represented by a matrix, once basis have been selected. This is true for any of the above mentioned linear maps. The dimensions are given in the following table, recalling that if $A\: : \: \K^n \longrightarrow \K^m$, $A$ is $m \times n$:\\
\begin{center}
 \begin{tabular}{cccccc} \\
 matrix & \# rows & & \# columns & row index & column index \\\\
 \hline
  $A$ & $mnp$   & $\times$ & $1$ & $ijk$ & $.$ \\
  $A_\e$ & $np$ & $\times$ & $m$ & $jk$ & $i$ \\
  $A_\f$ & $mp$ & $\times$ & $n$ & $ik$ & $j$ \\
  $A_\g$ & $mn$ & $\times$ & $p$ & $ij$ & $k$ \\
  $A_{\e\f}$ & $p$ & $\times$ & $mn$ & $k$ & $ij$ \\
  $A_{\e\g}$ & $n$ & $\times$ & $mp$ & $j$ & $ik$ \\
  $A_{\f\g}$ & $m$ & $\times$ & $np$ & $i$ & $jk$ \\
  $A_{\e\f\g}$ & $1$ & $\times$ & $mnp$ & $.$ & $ijk$ \\
 \end{tabular}
\end{center}


\nT{Matricization and HOSVD} One of the most useful possibility offered by matricization is to build a low rank approximation of the induced matrix through \kw{SVD}. Let us recall that, if $A \in \K^{m \times n}$ with $m \geq n$ and is full rank, SVD is
\begin{equation}
 A = \sum_{a=1}^n \sigma_a \ubold_a \otimes \vbold_a
\end{equation}
with $(\ubold_a)_a$ and $(\vbold_a)_a$ being orthonormal families (see section \ref{sec:inter:svd} for details). Best rank $r$ approximation of $A$ is
\begin{equation}
 A_r = \sum_{a=1}^r \sigma_a \ubold_a \otimes \vbold_a
\end{equation}
Let us give an example with $A_\f \in \L(F,E\otimes G)$, defined by
\begin{equation}
 A_\f.\fbold_j = \sum_{i,k}\alpha_{ijk} \, \ebold_i \otimes \gbold_k
\end{equation}
It is a matrix with $mp$ rows and $n$ columns, with $\alpha_{ijk}$ at row $i'=(i,k)$ and column $j$. Best approximation at rank $r$ reads
\begin{equation}
 A_{\f,r} = \sum_{a=1}^r \sigma_a \, \ubold_a \otimes \vbold_a \qquad \mbox{with} \quad
 \begin{cases}
  \ubold_a \in \R^{mp} \\
  \vbold_a \in \R^n 
 \end{cases}
\end{equation}
Let 
\begin{equation}
 \ubold_a = \sum_{i,k} u_{aik}\, \ebold_{ik}, \qquad \vbold_a = \sum_j v_{aj}\,\ebold_j
\end{equation}
Then, best approximation of $\Abold$ from SVD of $A_\f$ is
\begin{equation}
 \widetilde{A} = \sum_{a=1}^r \sigma_a\,u_{aik}v_{ja}\, \ebold_i \otimes \fbold_j \otimes \gbold_k
\end{equation}
Alternating such a matricization and best low rank approximation over all modes yields \kw{HOSVD} (see further section \ref{sec:blratuck:hosvd3} for details).

%
\section{Contraction between $d-$modes tensors}\label{sec:contract:dmodes}
%

This can be extended to contractions between $d-$modes tensors in a straightforward way. The only difficulty to overcome is to set clear notations. 

\nS Let 
\[
 \Ebold = E_1 \otimes \ldots \otimes E_d = \bigotimes_{\mu=1}^d E_\mu 
\]
and
\[
 \Abold \in \Ebold
\]
Let
\[
 I = \llbracket 1,d \rrbracket \subset \N
\]
and 
\[
 J \subset I
\]
Let
\[
 \Xbold \in \bigotimes_{\mu \in J} E_\mu
\]
We observe that each mode $\mu \in J$ of $\Xbold$ is a mode $\mu \in I$ of $\Abold$ too. Then, the \kw{contraction} of $\Abold$ with $\Xbold$ over $J$, denoted $\bullet_\j$ is defined on elementary tensors 
\begin{equation}
 \left(\bigotimes_{\mu \in I}\abold_\mu\right) \bullet_\j \left(\bigotimes_{\mu' \in J}\xbold_{\mu'}\right) = \left(\prod_{\mu \in J}\langle \abold_\mu \, \, \xbold_\mu\rangle \right) \; \bigotimes_{\mu \notin J}\abold_\mu
\end{equation}
This is extended as a contraction over the whole space
\[
 \begin{CD}
  \begin{array}{c}
   \bigotimes_{\mu \in I} E_\mu\\
   \times \\
   \bigotimes_{\nu \in J} E_\nu 
  \end{array}
  @>\mathrm{contraction}>> \bigotimes_{\mu \notin J} E_\mu
 \end{CD}
\]
by linearity.

\nS Component-wise, a contraction is a summation over shared indices $i_\mu \in J$. Let 
\begin{equation*}
  \left\{
        \begin{array}{lcl}
          \Abold &=& \displaystyle \sum_{\ibold \in \IC}\, \abold_\ibold \, \Ebold_\ibold \\
          \Xbold &=& \displaystyle \sum_{\jbold \in \JC}\, \xbold_\jbold \, \Ebold_\jbold
        \end{array}
  \right.
\end{equation*}
Then
\begin{equation}
 \left(\sum_{\ibold \in \IC} \abold_\ibold \Ebold_\ibold\right) \bullet_\j \left(\sum_{\jbold \in \JC} \xbold_\jbold \Ebold_\jbold\right) = \sum_{\ibold' \notin \JC} \left(\sum_{\jbold \in \JC} \abold_{\ibold'\jbold}\xbold_{\jbold}\right) \Ebold_{\ibold'}
\end{equation}
(this is a formula similar to the formula of matrix $\times$ vector contraction, but over multi-indices).

\nT{Example} This is better understood on an example. Let $d=6$ and $J=\{2,4,5\}$.  Let
\begin{equation*}
    \Abold = \sum_{i_1 \ldots i_6} \, a_{i_1\ldots i_6}\, \ebold_{1,i_1} \otimes \ldots \otimes \ebold_{6,i_6}
\end{equation*}
and
\begin{equation*}
    \Xbold = \sum_{i_2,i_4,i_5}\, x_{i_2i_4i_5}\, \ebold_{2,i_2} \otimes \ebold_{4,i_4} \otimes \ebold_{5,i_5}
\end{equation*}
Then
\begin{equation}
    \Abold \bullet_{2,4,5} \Xbold = \sum_{i_2i_4i_5} \,a_{i_1\ldots i_6}x_{i_2i_4i_5} \, \ebold_{1,i_1} \otimes \ebold_{3,i_3} \otimes \ebold_{6,i_6} 
\end{equation}

%
\section{Hasse diagrams of contractions}\label{sec:contract:hasse}
%

The set $\PC(I)$ of subsets of $I$ is a poset\footnote{Partially Ordered Set} with order relation $\subset$. Subsets in $I$ may be organized as a diagram, the \kw{Hasse diagram}\index{diagram!Hasse}, where an arrow joins $J$ and $J'$ when $J' \subset J$\marginpar{\dbend} (and not when $J \subset J'$). As $\subset$ is transitive, nearest neighbor relations only are drawn. We associate to each subset $J$ the tensor product of vector spaces
\begin{equation}
 \begin{CD}
  J @>>> \displaystyle \bigotimes_{\mu \in J}E_\mu
 \end{CD}
\end{equation}
and extend the order relationship to tensor products of vector spaces as $\prec$. Hence
\begin{equation}
 J' \subset J \quad \Longleftrightarrow \quad \bigotimes_{\mu \in J'}E_\mu \prec \bigotimes_{\mu \in J}E_\mu
\end{equation}
The Hasse diagram can be transferred to $\prec$. For $d=3$ and $I = \{1,2,3\}$, this yields:\\
\\
\begin{center}
\begin{tikzpicture}
\node (EFG) at (2,6) {$E \otimes F \otimes G$} ;  
\node (EF) at (0,4) {$E \otimes F$} ;  
\node (EG) at (2,4) {$E \otimes G$} ;  
\node (FG) at (4,4) {$F \otimes G$} ;  
\node (E) at (0,2) {$E$} ;  
\node (F) at (2,2) {$F$} ;  
\node (G) at (4,2) {$G$} ;  
\node (0) at (2,0) {$\K$} ;  
\draw[->,>=stealth] (E) -- (0) ; 
\draw[->,>=stealth] (F) -- (0) ; 
\draw[->,>=stealth] (G) -- (0) ;
\draw[->,>=stealth] (EF) -- (E) ; 
\draw[->,>=stealth] (EG) -- (E) ; 
\draw[->,>=stealth] (EF) -- (F) ; 
\draw[->,>=stealth] (FG) -- (F) ; 
\draw[->,>=stealth] (EG) -- (G) ; 
\draw[->,>=stealth] (FG) -- (G) ; 
\draw[->,>=stealth] (EFG) -- (EF) ; 
\draw[->,>=stealth] (EFG) -- (EG) ; 
\draw[->,>=stealth] (EFG) -- (FG) ; 
\end{tikzpicture}
\end{center}
An arrow like
\[
 \begin{CD}
    E \otimes F \otimes G @>>> E \otimes G
 \end{CD}
\]
means that there is a contraction from $E \otimes F \otimes G $ on $E \otimes G$ over its complement, here $F$. This is true for any arrow. 

\nS Let us have a path, like
\[
 \begin{CD}
   E \otimes F \otimes G @>\f>> E \otimes G @>\e>> G
 \end{CD}
\]
with indication of the space ( here a mode) over which the contraction is done. Let
\[
 \Abold \in E \otimes F \otimes G, \qquad \ybold \in F, \qquad \xbold \in E
\]
Then\footnote{Let us note that the indices in $J$ in $\bullet_\j$ are always in the same order than in $I$.}
\begin{equation}
 (\Abold \bullet_2 \ybold) \bullet_1 \xbold = \Abold \bullet_{12} (\xbold \otimes \ybold)
\end{equation}
\begin{proof}
Let us show it for elementary tensors, and this can be extended to whole tensor space by linearity. Let 
\[
 \Abold = \abold \otimes \bbold \otimes \cbold
\]
Then
\[
 \begin{array}{lcl}
  ((\abold \otimes \bbold \otimes \cbold) \bullet_2 \ybold) \bullet_1 \xbold &=& (\langle \bbold,\ybold\rangle \; \abold \otimes \cbold) \bullet_1 \xbold \\
  &=& \langle \bbold,\ybold\rangle \; (\abold \otimes \cbold) \bullet_1 \xbold \\
  &=& \langle \bbold,\ybold\rangle \langle \abold , \xbold \rangle \cbold \\
  &=& (\abold \otimes \bbold \otimes \cbold) \bullet_{12} (\xbold \otimes \ybold)
 \end{array}
\]
\end{proof}
\noindent This shows as well that 
\begin{equation}
  (\abold \otimes \bbold \otimes \cbold) \bullet_{12} (\xbold \otimes \ybold) = (\abold \otimes \bbold \otimes \cbold) \bullet_{21} (\ybold \otimes \xbold)
\end{equation}
and that it is always possible to order the indices as in $I$. It simply means that the Hasse diagram is commutative. Let us show it here with the different paths in red and blue, and indication of the spaces over which the contraction operates:\\
\\
\begin{center}
\begin{tikzpicture}
\node (EFG) at (2,6) {\textcolor{red}{$E \otimes F \otimes G$}} ;  
\node (EF) at (0,4) {$E \otimes F$} ;  
\node (EG) at (2,4) {\textcolor{red}{$E \otimes G$}} ;  
\node (FG) at (4,4) {\textcolor{blue}{$F \otimes G$}} ;  
\node (E) at (0,2) {$E$} ;  
\node (F) at (2,2) {$F$} ;  
\node (G) at (4,2) {\textcolor{red}{$G$}} ;  
\node (0) at (2,0) {$\K$} ;  
\draw[->,>=stealth] (E) -- (0) ; 
\draw[->,>=stealth] (F) -- (0) ; 
\draw[->,>=stealth] (G) -- (0) ;
\draw[->,>=stealth] (EF) -- (E) ; 
\draw[->,>=stealth] (EG) -- (E) ; 
\draw[->,>=stealth] (EF) -- (F) ; 
\draw[->,>=stealth] (FG) -- (F) ; 
\draw[->,>=stealth, red] (EG) -- (G) node [pos=0.3, above, sloped] {$\e$} ; 
\draw[->,>=stealth, blue] (FG) -- (G) node[pos=0.5, right] {$\f$} ; 
\draw[->,>=stealth] (EFG) -- (EF) ; 
\draw[->,>=stealth, red] (EFG) -- (EG) node [pos=0.5, right] {$\f$} ; 
\draw[->,>=stealth, blue] (EFG) -- (FG) node [pos=0.5, above, sloped] {$\e$}; 
\draw[->,>=stealth, violet](EFG)--(G) ; 
\end{tikzpicture}
\end{center}
The diagram is commutative: red path (contract first on $F$ and then on $E$) and blue path (contract first on $E$ and then on $F$) yield the same result, drawn in purple (contract on $E \otimes F$).

\nS This can be seen on a \kw{contraction diagram}\index{diagram!contraction}, as follows: the red path above
\[
 \begin{CD}
  E \otimes F \otimes G @>\f>> E \otimes G @>\e>> G
 \end{CD}
\]
can be sketched by the following diagram\\
\\
\begin{center}
\begin{tikzpicture}
\node[tensor_ell] (A) at (2,0) {$a \otimes b \otimes c$} ; 
\node[tensor_ell] (x) at (-0.5,0) {$x$} ; 
\node[tensor_ell] (y) at (2,1.5) {$y$} ; 
\draw (x) -- (A) node [pos=0.5, above] {$1$} ;
\draw (A) -- (y) node [pos=0.5, right] {$2$} ;
\draw (A) -- (4,0) node [pos=0.5, above] {$3$} ;
\end{tikzpicture}
\end{center}
which yields, by collapsing:\\
\\
\begin{center}
\begin{tikzpicture}
\node[tensor_ell] (A) at (2,0) {$\langle a,x \rangle \langle b,y \rangle \; c$} ; 
\draw (A) -- (4,0) node [pos=0.5, above] {$3$} ;
\end{tikzpicture}
\end{center}
where commutativity of the Hasse diagram is obvious.


\nT{Hasse diagrams for $d-$modes tensors} This can be generalized to $d-$tensors as follows. Let us denote 
\[
 \ES_\j = \bigotimes_{\mu \in J}E_\mu
\]
Let $J \subset I$ and $I' = I \backslash J$. Let $K \subset I'$ and $I'' = I' \backslash K$. Hence
\[I = J \sqcup I', \qquad I' = K \sqcup I''
\]
Let us consider the sequence
\begin{equation}
 \begin{CD}
 \ES_{\i} @>>> \ES_{\i'} @>>> \ES_{\i''}
 \end{CD}
\end{equation}
Let us consider
\[
 \Abold \in \ES_{\i}, \qquad \Xbold \in \ES_{\j}, \qquad \Ybold \in \ES_{k}
\]
Then
\begin{equation}
 (\Abold \bullet_{\j} \Xbold) \bullet_{\k} \Ybold = \Abold \bullet_{\j \cup \k} (\Xbold \otimes \Ybold)
\end{equation}
provided the indices are kept in the order of $I$ in $J \cup K$. 
\begin{proof}
This can be shown by above calculation by selecting $E \otimes F \otimes G = \ES_{\i}$, $F= \ES_{\j}$ and $E = \ES_{\k}$ modulo the order of the parameters.\\
\\
Let us give an example with 
\[
d=6, \qquad J= \{3,5\}, \qquad  K = \{2,4\}
\]
so
\[
I' = \{1,2,4,6\}, \qquad I'' = \{1,6\}
\]
Then, we have
\begin{equation*}
\begin{array}{l}
[(\abold_1 \otimes \abold_2 \otimes \abold_3 \otimes \abold_4 \otimes \abold_5 \otimes \abold_6) \bullet_{35} (\xbold_3 \otimes \xbold_5)] \bullet_{24}(\xbold_2 \otimes \xbold_4)\\
\qquad = [\langle \abold_3,\xbold_3\rangle\langle \abold_5,\xbold_5\rangle (\abold_1 \otimes \abold_2 \otimes \abold_4 \otimes \abold_6)] \bullet_{24}(\xbold_2 \otimes \xbold_4) \\
\qquad = \langle \abold_3,\xbold_3\rangle\langle \abold_5,\xbold_5\rangle (\abold_1 \otimes \abold_2 \otimes \abold_4 \otimes \abold_6) \bullet_{24}(\xbold_2 \otimes \xbold_4) \\
\qquad = \langle \abold_3,\xbold_3\rangle\langle \abold_5,\xbold_5\rangle (\langle \abold_2,\xbold_2\rangle \langle \abold_4,\xbold_4\rangle (\abold_1 \otimes \abold_6)) \\
\qquad = \langle \abold_2,\xbold_2\rangle \langle \abold_3,\xbold_3\rangle\langle \abold_4,\xbold_4\rangle\langle \abold_5,\xbold_5\rangle  (\abold_1 \otimes \abold_6) \qquad (\mbox{reordering of the indices as in } I) \\
\qquad = (\abold_1 \otimes \abold_2 \otimes \abold_3 \otimes \abold_4 \otimes \abold_5 \otimes \abold_6) \bullet_{2345} (\xbold_2 \otimes \xbold_3 \otimes \xbold_4 \otimes \xbold_5) 
\end{array} 
\end{equation*}
\end{proof}

%
\section{Cauda: bilinear maps and structure tensor}\label{sec:contract:bilin}
%

Let us present as a perspective for this chapter on contraction a rich consequence of the duality between tensor product and contraction: any bilinear map $B$ on a pair of vector spaces $E \times F$ can be written as a composition between a contraction and a tensor product. This is an outcome of the universal property of tensor product. It is shown here in a less abstract setting, using coordinates. This introduces the \kw{structure tensor}\index{tensor!structure} of a bilinear map.

\nS Let us consider three vector spaces $E,F,G$ on a field $\K$. Let 
\begin{equation}
 \xbold \in E, \qquad \ybold \in F, \qquad \zbold \in G
\end{equation}
Let $B$ be bilinear map 
\begin{equation}
 \begin{CD}
   E \times F @>B>> G
 \end{CD}
\end{equation}
Let us define 
\begin{equation}
 \left\{
    \begin{array}{lcl}
      \xbold &=& \displaystyle \sum_i \, x_i \,\ebold_i \\
      \ybold &=& \displaystyle \sum_j \, y_j \, \fbold_j
    \end{array}
 \right.
\end{equation}
Then
\begin{equation}
  B(\xbold,\ybold) = \sum_{i,j} \, x_iy_j\,B(\ebold_i,\fbold_j)
\end{equation}
We define 
\begin{equation}
 B(\ebold_i,\fbold_j) = \sum_k \, b_{ijk}\, \gbold_k
\end{equation}
The terms $b_{ijk}$ are the \kw{structure coefficients}  of $B$. Let $\Bbold$ be defined as
\begin{equation}
 \Bbold = \sum_{i,j,k}\, b_{ijk} \, \ebold_i \otimes \fbold_j \otimes \gbold_k \qquad \in E  \otimes F \otimes G
\end{equation}
It is the \kw{structure tensor} of the bilinear map $B$. Then
\begin{equation}
 \begin{array}{lcl}
   B(x,y) &=& \displaystyle \sum_{i,j} \, x_iy_j\,B(\ebold_i,\fbold_j) \\
   &=& \displaystyle \sum_{i,j} \, x_iy_j\left(\sum_k \, b_{ijk}\, \gbold_k\right) \\
   &=& \displaystyle \sum_k \left(\sum_{i,j} \, x_iy_jb_{ijk}\right) \gbold_k
 \end{array}
\end{equation}
where we recognize
\begin{equation}
 B(\xbold,\ybold) = \Bbold \bullet_{1,2} (\xbold \otimes \ybold)
\end{equation}
Then, any bilinear map defined by its structure tensor $\Bbold$ can be written as a combination of a contraction and tensor product.

\nT{Link with universal property of tensor product} Let us note that (and this shows the usefulness of tensor algebra with respect to a component-wise calculation is some circumstances)
\begin{equation}
    \begin{CD}
      E \otimes F @>>> G \\
      \xbold \otimes \ybold @>>> \Bbold \bullet_{1,2} (\xbold \otimes \ybold)
    \end{CD}
\end{equation}
is linear. Indeed,
\begin{equation}
    \begin{array}{lcl}
     \Bbold \bullet_{1,2} \left(\alpha \xbold \otimes \ybold + \alpha' \xbold' \otimes \ybold'\right) &=& \displaystyle \sum_k \left(\sum_{i,j} \, b_{ijk}(\alpha x_iy_j + \alpha' x'_iy'_j)\right) \gbold_k \\
     &=& \displaystyle \alpha \sum_k \left(\sum_{ij} b_{ijk}x_iy_j\right) \gbold_k + \alpha' \sum_k \left(\sum_{ij} b_{ijk}x'_iy'_j\right) \gbold_k \\
     &=& \alpha \, \Bbold \bullet_{1,2} (\xbold \otimes \ybold) + \alpha' \, \Bbold \bullet_{1,2} (\xbold' \otimes \ybold')
    \end{array}
\end{equation}
We recognize here the universal property of tensor product which is accepted in algebra as its definition (see section \ref{sec:tens_prod:otherdefs} where $\Bbold \bullet_{1,2}$ is denoted $L$).

\nS Most of standard operations, like contraction itself, are associated with a sparse structure tensor $\Bbold$. Let us give the example of matrix $\times$ vector product with $E \rightarrow E \otimes F, F=F, G \rightarrow E$. We have
\begin{equation}
 B(A,\xbold)=A.\xbold, \qquad \mbox{with }
 \begin{cases}
  A & \in E \otimes F \\
  \xbold & \in F
 \end{cases}
\end{equation}
It is a bilinear map. Then, the structure coefficients are given by
\begin{equation}
  B(\ebold_i \otimes \fbold_j \, , \, \fbold_k) = \sum_{ijk\ell}\, b_{ijk\ell} \, \ebold_\ell
\end{equation}
Then
\begin{equation}
 \begin{array}{lcl}
   b_{ijk\ell} &=& \langle B(\ebold_i \otimes \fbold_j,\fbold_k) \; , \; \ebold_\ell \rangle\\
   &=& \langle (\ebold_i \otimes \fbold_j).\fbold_k \; , \; \ebold_\ell \rangle\\
   &=&  \langle \delta_k^j\, \ebold_i \, , \, \ebold_\ell \rangle \\ 
   &=& \delta_k^j\delta_i^\ell
 \end{array}
\end{equation}
Then
\begin{equation}
    b_{ijk\ell} = 
    \begin{cases}
      1 & \mbox{if} \quad i=\ell \: \mbox{ and } \: j=k \\
      0 & \mbox{otherwise}
    \end{cases}
\end{equation}
It is a sparse tensor because, if $\dim E=\dim F=n$, $\Bbold$ has $n^4$ terms $(\Bbold \in \R^{n^2 \times n \times n})$, among which $n^2$ only are non zero.

\nS This can easily be extended to $d-$linear maps. Let $(E_1,\ldots,E_d)$ be a family of vectors spaces on a field $\K$, and $F$ a vector space on same field, $\xbold_\mu \in E_\mu$. Let us consider a $d-$linear map $T$
\begin{equation}
\begin{CD}
 E_1 \times \ldots \times E_d @>T>> F 
\end{CD}
\end{equation}
Let us define 
\begin{equation}
 \xbold_\mu = \sum_{i_\mu} x_{i_\mu}^{(\mu)}\, \ebold_{i_\mu}^{(\mu)}
\end{equation}
Then
\begin{equation}
  T(\xbold_1, \ldots,\xbold_d) = \sum_{i_1} \ldots \sum_{i_d} \, x_{i_1}^{(1)} \ldots x_{i_d}^{(d)} \,T\left(\ebold_{i_1}^{(1)},\ldots,\ebold_{i_d}^{(d)}\right)
\end{equation}
We define 
\begin{equation}
  T\left(\ebold_{i_1}^{(1)},\ldots,\ebold_{i_d}^{(d)}\right) = \sum_j \, t_{i_1,\ldots,i_d,j}\, \fbold_j
\end{equation}
The terms $t_{i_1,\ldots,i_d,j}$ are the \kw{structure coefficients}  of $T$. Let $\Tbold$ be defined as
\begin{equation}
 \Tbold = \sum_{i_1}\ldots \sum_{i_d}\sum_j\, t_{i_1,\ldots,i_d,j} \, \ebold_{i_1}^{(1)} \otimes \ldots \otimes \ebold_{i_d}^{(d)}\otimes \fbold_j \qquad \in E_1 \otimes \ldots \otimes E_d \otimes F
\end{equation}
It is the \kw{structure tensor} of the multilinear map $T$. Then
\begin{equation}
 \begin{array}{lcl}
   T(\xbold_1, \ldots,\xbold_d) &=& \displaystyle \sum_{i_1} \ldots \sum_{i_d}\, 
   x_{i_1}^{(1)} \ldots x_{i_d}^{(d)} \,T\left(\ebold_{i_1}^{(1)},\ldots,\ebold_{i_d}^{(d)}\right) \\
   &=& \displaystyle \sum_{i_1} \ldots \sum_{i_d}\, 
   x_{i_1}^{(1)} \ldots x_{i_d}^{(d)} \left(\sum_j \, t_{i_1,\ldots,i_d,j}\, \fbold_j\right) \\
   &=& \displaystyle \sum_j \left(\sum_{i_1} \ldots \sum_{i_d}\, 
   x_{i_1}^{(1)} \ldots x_{i_d}^{(d)} \, t_{i_1,\ldots,i_d,j} \right) \fbold_j
 \end{array}
\end{equation}
where we recognize
\begin{equation}
 T(\xbold_1, \ldots,\xbold_d)= \Tbold \bullet_{1,\ldots,d} (\xbold_1 \otimes \ldots, \xbold_d)
\end{equation}
Then, any multilinear map defined by its structure tensor $\Tbold$ can be written as a combination of a tensor product and a contraction.

%
\chapter{Algebraic structures in tensor spaces}\label{sec:alstruct}
%

This chapter, although with an algebraic flavor, is a key chapter with one target. It establishes a precise link between Kronecker product, which is defined originally between matrices, on one hand, and tensor product on the other. In the literature, those two notions are often confounded, which can be understood and accepted because they are equivalent up to an isomorphism, and even denoted by the same symbol $\otimes$. Here, I will use $\otimes$ for the tensor product and $\otimes_\k$ for the Kronecker product.

\nB One of the motivation for Kronecker product and reason for its success is that it extends to operators (linear maps) the tensor product between vectors. One wishes to write something like $(M \otimes_\k Q)(a \otimes b)=Ma \otimes Qb$, the Kronecker product $\otimes_\k$ being a sort of tensor product $\otimes$ between operators $M$ and $Q$. The aim of this chapter is to set this on a rigorous algebraic basis with full details. 

\nT{Organization of the chapter} This chapter is organized as follows:
\begin{description}
\item[section \ref{sec:alstruct:group}] A key preliminary notion is the action of the linear group $\GL(E)$ on $E$, which is given in this section.
\item[section \ref{sec:alstruct:kron}] This is the key section for defining Kronecker product between matrices (seen here as operators, as it is base free) in relation with tensor product. First, the box product $\boxtimes$ between two operators $M,N$ is defined by $(M \boxtimes Q)(\abold \otimes \bbold)= M\abold \otimes Q\bbold$. We show that it extends the action of the linear group on a vector space, and extend it to product of $d$ vector spaces. The link with Kronecker product is presented at the end of this section. 
\item[section \ref{sec:alstruct:prodmat}] To complete this chapter on products, some classical products between matrices are recalled, where the versatility of notations between publications is shown.

\end{description}

%
\section{Action of the linear group}\label{sec:alstruct:group}
%

It is not obvious \emph{a  priori} that the action of a group on a set is a useful concept in tensor algebra. However, it appears in various guises in some elementary operations on tensors in the literature, especially for Tucker model which is a change of basis (see section \ref{sec:dmodestens:basis}). Let us show it for an elementary matrix. Let $A=\abold \otimes \bbold$ with $\abold = M\abold'$ and $\bbold = N\bbold'$ from a change of basis. Then $A=M\abold' \otimes N\bbold'$ which can be denoted $A=(M,N)(\abold' \otimes \bbold')$ where $(.,.)$ denotes the action of the group $\GL(E) \times \GL(F)$ on $E \otimes F$.  

\nT{Action of a group on a set} Let $X$ be a set and $G$ a group. The group $G$ acts on the set $X$ if there is a map
\[
 \begin{CD}
  G \times X @>>> X \\
  (g,x) @>>> gx
 \end{CD}
\]
such that
\begin{equation}
 \left. 
 \begin{array}{cc}
  \forall \: g,g' &\in G \\
   \forall \: x &\in X 
 \end{array}
 \right\},
 \quad g'(gx)=(g'g)x
\end{equation}
The set
\[ 
 Gx = \{ gx \mid g \in G \} 
\]
is called the \kw{orbit} of $x$. It is customary to denote by $gx$ the action of $G$ on $X$.

\nS The group $\GL(E)$ acts on $E$ by
\[
 \begin{CD}
   \GL(E) \times E @>>> E \\
   (M,\xbold) @>>> M.\xbold
 \end{CD}
\]
Let us observe that $\GL(E) \times \GL(F)$ is a group by composition $*$: if $M, M' \in \GL(E)$ and $Q,Q' \in \GL(F)$,
\[
 (M',Q')*(M,Q) = (M'M, Q'Q)
\]
The identity is $(\I,\I)$, and the inverse of $(M,Q)$ is $(M^{-1},Q^{-1})$. 

\nT{Action of $\GL(E) \times \GL(F)$} On can define an action of $\GL(E) \times \GL(F)$ on $E \otimes F$ by defining its action on elementary tensors 
\begin{equation}
 (M,Q)(\abold\otimes \bbold) = M\abold \otimes Q\bbold
\end{equation}
and extend it to $E \otimes F$ by linearity.
\begin{proof}
It defines the action of a group because 
\begin{equation*}
\begin{array}{lcl}
 (M',Q')((M,Q)(\abold \otimes \bbold)) &=& (M',Q')(M\abold \otimes Q\bbold) \\
 &=& M'(M\abold) \otimes Q'(Q\bbold) \\
 &=& (M'M)\abold \otimes (Q'Q)\bbold \\
 &=& (M'M, Q'Q) (\abold \otimes \bbold) \\
 &=& ((M',Q')(M,Q))(\abold \otimes \bbold)
\end{array}
\end{equation*} 
\end{proof}
\noindent Let us note that
\begin{equation}
 (M,Q)(\abold\otimes \bbold) = M(\abold \otimes \bbold)Q^\t
\end{equation}
as matrix product. This as well can be extended by linearity to $E \otimes F$.

\nS If $A = \sum_{i,j}a_{ij}\,\ebold_i \otimes \fbold_j$, we have
\begin{equation}
    \begin{array}{lcl}
      (M,Q)A &=& \displaystyle (M,Q) \left(\sum_{i,j}\, a_{ij}\, \ebold_i \otimes \fbold_j\right) \\
      &=& \displaystyle \sum_{i,j}\, a_{ij}\, (M,Q)(\ebold_i \otimes \fbold_j) \\
      &=& \displaystyle  \sum_{i,j}\, a_{ij}\, M\ebold_i \otimes Q\fbold_j
    \end{array}
\end{equation}

\nS Such an action can be extended in a straightforward way for $d-$modes tensors, with
\begin{equation}
 (M_1, \ldots,M_d)(\abold_1 \otimes \ldots \otimes \abold_d) = (M_1\abold_1) \otimes \ldots \otimes (M_d\abold_d)
\end{equation}
and extended to $\bigotimes_\mu E_\mu$ by linearity.


\nT{Invariant of the action of a group} An \kw{invariant} of the action of a group $G$ on a set $X$ is a function $f$ 
\[
\begin{CD}
 X @>f>> Y
\end{CD}
\]
from $X$ on $Y$ such that 
\begin{equation}
    \forall \: x \in X, \quad \forall \: g \in G, \qquad f(gx) = f(x)
\end{equation}
i.e. a function constant on any orbit. Here is a classical example. Let $E$ be a vector space with $\dim E=n$ and $X$ the space $\R^{n \times n}$ of square matrices of endomorphisms in $\L(E)$ for a given basis in $E$. The action of $\GL(E)$ on $X$ is defined by
\begin{equation}
    \begin{CD}
      P \: : \: A @>>> PAP^{-1}
    \end{CD}
\end{equation}
with $P \in \GL(E)$ and $A \in X$. It is written 
\[
\begin{CD}
 A @>>> P[A]
\end{CD}
\]
to avoid confusion with matrix product. This is an action of $\GL(E)$ as
\begin{equation}
    \begin{array}{lcl}
        (QP)[A] &=& (QP)A(QP)^{-1} \\
        &=& Q(PAP^{-1})Q \\
        &=& QP[A]Q^{-1}\\
        &=& Q[P[A]]
    \end{array}
\end{equation}
Then, we have
\begin{equation}
    \det P[A] = \det A
\end{equation}
i.e. the function \emph{det} is an invariant of this action (here, $Y=\R$). An invariant of the action of $\GL(E) \times \GL(F)$ on $E \otimes F$ is a mapping 
\[
 \begin{CD}
  E \otimes F @>f>> \R
 \end{CD}
\]
such that
\begin{equation}
 \left. 
  \begin{array}{cl}
   \forall \ (M,N) &\in \GL(E) \times \GL(F) \\
   \forall \: A &\in E \otimes F
  \end{array}
 \right\}
 , \quad f((M,N)A)= f(A)
\end{equation}
An invariant of the action of $\GL(\R^2) \times \GL(\R^2) \times \GL(\R^2)$ on $\R^{2 \times 2 \times 2}$ is known: the Cayley hyperdeterminant or Kruskal's polynomial of a tensor $\Abold \in \R^2 \times \R^2 \times \R^2$ (see section \ref{sec:rank3:rank222}, page \pageref{sec:rank:rank3:rank222:page1}). But exhibiting them is quite an arduous task, and very few results are known.

\notes\label{page:order} The action of a group on a set is a classical development of group theory. Its presentation can be found in any textbook on algebra, see e.g. \cite{Gouvea2012}, chapter 4. Groups act on vector spaces through their linear representation. If $E$ is a vector space, and $G$ a group, a representation of $G$ is a group homomorphism $\rho: \:\: G \longrightarrow \GL(E)$. It is standard to write $g\xbold = \rho(g).\xbold$ if $\xbold \in E$ (see \cite{Vinberg89} for further details). A nice and clear presentation of the rich theory of invariants of the action of a group on a set can be found in \cite{Neusel2007}. It is worth noticing that a motivation for the emergence of the concept of tensor in late 19$^{\mathrm{th}}$ century is the quest for invariants in so called \emph{forms} (see \cite{Elliott13,Grace1903,Gurevich1964,Dieudonne1970b,Luque2008}). \cite{Elliott13} and \cite{Grace1903} have an historical interest, presenting the state of the art by the turn of 19th/20th century. \cite{Gurevich1964} is at the same time more modern and very accessible. \cite{Dieudonne1970b} is a very nice, deep and modern synthesis, but rather difficult without acquaintance with abstract algebra. The link with tensors is that in these early works, a form is a $n-$ary homogeneous polynomial of degree $r$, i.e. an expression $\sum_{|\bm{\alpha}|=r}a_{\alpha_1\ldots\alpha_n}x_1^{\alpha_1}\ldots x_n^{\alpha_n}$, like $ax^3+3bx^2y+3cxy^2+dy^3$ for $r=3$ and $n=2$. This example can be written $\langle \Abold, \xbold \otimes \xbold \otimes \xbold \rangle$ with $\xbold = (x,y) \in \R^2$ and $\Abold \in \R^{2 \times 2 \times 2}$ is a symmetric tensor of order $r=3$ and dimension $n=2$ on each mode. Let us note that, at these times, $r$ was called the order of the form (the degree of the homogeneous polynomial), which is the origin of the notion of order of a tensor. In brief, the objective, which has been reached after decades of intensive research (see \cite{Dieudonne1970b}), was to classify all the real or complex functions of the coefficients of a tensor which are invariant under a change of basis, like the trace or the determinant of a matrix (a tensor of order 2). An example of such an invariant for a tensor in $\R^{2 \times 2 \times 2}$ is the Cayley hyperdeterminant given in section \ref{sec:rank3:rank222}.

%
\section{Kronecker product}\label{sec:alstruct:kron}
%

Let us consider $M\abold \otimes Q\bbold$. By analogy with
\[
 (\abold^* \otimes \bbold^*)(\xbold,\ybold) = \abold^*(\xbold)\bbold^*(\ybold)
\]
where $\abold^*$ acts on $\xbold$ and $\bbold^*$ on $\ybold$ one is tempted to write something like
\[
 (M \otimes Q)(\abold,\bbold) = M\abold \otimes Q\bbold
\]
We show that it is not consistent with the tensor product $E \otimes F$:  $\L(E)$ and $\L(F)$ are each vector spaces, and $M \otimes Q$ is a tensor product of matrices in $\L(E) \otimes \L(F)$ which does not act on $(\abold , \bbold)$ as $M\abold \otimes Q\bbold$. 


\nS To see this, let us recall that, if $A \in \L(F,E)$ and $M,Q$ are invertible
\[
 (M,Q)A = MAQ^\t
\]
Let us define
\begin{equation}
 \begin{CD}
   \L(F,E) @>L_{\m,\q}>> \L(F,E)
 \end{CD}
\end{equation}
by
\begin{equation}
 L_{\m,\q}(A) = MAQ^\t 
\end{equation}
which can be illustrated by the following diagram:
\begin{equation}
\begin{CD}
F @>A>> E \\
@A Q^\t AA @VV M V \\
F @>MAQ^\t>> E
\end{CD}
\end{equation}
For $A = \abold \otimes \bbold$, we have (see equations (\ref{eq:elemat}))
\[
 L_{\m,\q} (\abold \otimes \bbold) = M\abold \otimes Q\bbold
\]
$L_{\m,\q}$ is an endomorphism of $\L(F,E) \simeq E \otimes F$, i.e.
\begin{equation}
 L_{\m,\q} \in \L(E \otimes F)
\end{equation}
or
\begin{equation}
 L_{\m,\q} \in E \otimes F \otimes E \otimes F
\end{equation}
However, as $M \in E \otimes E$ and $Q \in F \otimes F$, we have
\begin{equation}
 \begin{array}{lcl}
   M \otimes Q &\in& E \otimes E \otimes F \otimes F \\
   &\neq& E \otimes F \otimes E \otimes F
 \end{array}
\end{equation}
and
\begin{equation}
 L_{\m,\q} \neq M \otimes Q
\end{equation}
as both sides do not belong to the same space.


\nS Let us define\footnote{The notation $\boxtimes$ is not standard, and introduced here for sake of clarity} the linear operator
\begin{equation}
\begin{CD}
 E \otimes F @>M\boxtimes Q >> E \otimes F
\end{CD}
\end{equation}
by
\begin{equation}
 (M \boxtimes Q)(a \otimes b) = Ma \otimes Qb
\end{equation}
and extended to $E \otimes F$ by linearity. If $M,Q$ are invertible, one has
\begin{equation}
 M \boxtimes Q (A) = (M,N)A
\end{equation}
The reason why this notation $\boxtimes$ is introduced in place of $(M,N)$ is twofold:
\begin{itemize}[label =$\rightarrow$]
 \item $M \boxtimes Q$ is defined even if $M,N$ are not invertible 
 \item it will be extended in next paragraph to situation where $M \in \L(E,F)$ with $F \neq E$. 
\end{itemize}
The operation $\boxtimes$ is a bilinear map\footnote{As such, it can be associated with the tensor product $\L(E) \otimes \L(F)$ up to an isomorphism but, as mentioned above, direct association without distinction will not be consistent with the other tensor product $E \otimes F$.}
\begin{equation}
 \begin{CD}
   \L(E) \times \L(F) @>\boxtimes >> \L(E \otimes F)
 \end{CD}
\end{equation}

\nS This operator can be extended to non elementary tensors by linearity:
\begin{equation}
 \begin{aligned}
   (M \boxtimes Q)(A) &=  (M \boxtimes Q)\left(\sum_{i,j} a_{ij} \; \ebold_i \otimes \fbold_j\right) \\
   &= \sum_{i,j} a_{ij} \;  (M \boxtimes Q)(\ebold_i \otimes \fbold_j) \\
   &= \sum_{i,j} a_{ij} \; M\ebold_i \otimes Q\fbold_j
 \end{aligned}
\end{equation}


\nS This can easily be extended to a boxproduct of $d$ linear maps. Let 
\begin{itemize}[label=$\rightarrow$]
 \item $E_1, \ldots, E_d$ be $d$ vector spaces on a same field $\K$,
 \item $\abold_1, \ldots,\abold_d$ be $d$ vectors  with $\abold_\mu \in E_\mu$ 
 \item $M_1, \ldots, M_d$ be $d$ endomorphisms with $ M_\mu \in \L(E_\mu)$.
\end{itemize}
 Then
\begin{equation}
 (M_1 \boxtimes \ldots \boxtimes M_d)(\abold_1 \otimes \ldots \otimes \abold_d) = M_1\abold_1 \otimes \ldots \otimes M_d\abold_d
\end{equation}
This can be extended to $\bigotimes_\mu E_\mu$ by linearity.


\nS Let us extend this scheme by letting the vector spaces to be different. Let $E,\; F, \; V$ and $W$ be four vector spaces on a same field $\K$, with $\dim E=m, \; \dim F=n, \; \dim V=p$ and $\dim W=q$. This is detailed by specifying corresponding indices and basis in table below:
\[
 \begin{array}{l|cccccccc}
  \mbox{space} && E && F && V && W \\
  \hline
  \dim         && m && n && p && q \\
  \mbox{index} && i && j && k && \ell \\
  \mbox{basis} && (\ebold_i)_i && (\fbold_j)_j && (\vbold_k)_k && (\wbold_\ell)_\ell
 \end{array}
\]
Let us consider
\[
 A \in \L(F,E) \simeq E \otimes F, \qquad X \in \L(W,V) \simeq V \otimes W
\]
and let us define
\begin{equation}
 (A \boxtimes X)(\fbold \otimes \wbold) = A\fbold \otimes X\wbold, \qquad \fbold \in F, \quad \wbold \in W
\end{equation}
We have
\begin{equation}
 \begin{CD}
  F \otimes W @> A\boxtimes X >> E \otimes V
 \end{CD}
\end{equation}
and $\boxtimes$ is a bilinear map
\begin{equation}
 \begin{CD}
  \L(F,E) \times \L(W,V) @> \boxtimes >> \L(F \otimes W\, , \, E \otimes V) 
 \end{CD}
\end{equation}
This can be seen on following diagram:\\
\\
\begin{center}
\begin{tikzpicture}
\node (F) at (0,3) {$F$} ;  
\node (E) at (4.5,3) {$E$} ;  
\node (W) at (0,0) {$W$} ;  
\node (V) at (4.5,0) {$V$} ;  
\node (FW) at (1.5,1.5) {$F \otimes W$} ;  
\node (EV) at (6,1.5) {$E \otimes V$} ;  
\draw[->,>=stealth] (F)--(E) node [pos=.5, above] {$A$} ;
\draw[->,>=stealth] (W)--(V) node [pos=.5, above] {$X$} ;
\draw[->,>=stealth] (FW)--(EV) node [pos=.5, above] {$A \boxtimes X$} ;
\draw[->,>=stealth] (F) -- (FW) ;
\draw[->,>=stealth] (W) -- (FW) ;
\draw[->,>=stealth] (E) -- (EV) ;
\draw[->,>=stealth] (V) -- (EV) ;
\end{tikzpicture}
\end{center}

\nT{Link with Kronecker product}As $\dim E \otimes V = mp$ and $\dim F \otimes W =nq$, $A \boxtimes X$ can be represented after reshaping by a $(mp \times nq)$ matrix in a given basis, acting on $\vec (\fbold \otimes \wbold)$. To see this, let 
\[
 A = \sum_{i,j} a_{ij} \; \ebold_i \otimes \fbold_j, \qquad X = \sum_{k\ell} x_{k\ell}\; \vbold_k \otimes \wbold_\ell
\]
Then, $(\ebold_i \otimes \vbold_k)_{i,k}$ is a basis of $E \otimes V$ and $(\fbold_j \otimes \wbold_\ell)_{j,\ell}$ is a basis of $F \otimes W$. We have in these basis
\begin{equation}
 \begin{aligned}
   (A \boxtimes X)(\fbold_j \otimes \wbold_\ell) &= A\fbold_j \otimes X\wbold_\ell \\
   &= \left(\sum_i a_{ij}\, \ebold_i\right) \otimes \left(\sum_k x_{k\ell}\, \vbold_k\right) \\
   &= \sum_{i,k} a_{ij}x_{k\ell}\; \ebold_i \otimes \vbold_k
 \end{aligned}
\end{equation}
As $A \boxtimes X \in \L(F \otimes W \, , \, E \otimes V)$, it can be represented after reshaping (see section \ref{sec:tenselem:reshape} for reshaping) by a matrix the rows of which are indexed par pairs $(i,k)$ and columns by pairs $(j,\ell)$. Then, the matrix of $A \boxtimes X$ has element $a_{ij}x_{k\ell}$ in row $(i,k)$ and column $(j,\ell)$. If double indices are ordered with second index running first and first index running second, i.e.
$(i,k) \rightarrow 11, 12, \ldots, 1p, 21, \ldots, mp$ and $(j,\ell) \rightarrow 11, 12, \ldots, 1q, 21, \ldots, nq$, one recognizes in matrix $A \boxtimes X$ the Kronecker product of matrices $A$ and $X$, which is not their tensor product $\otimes$. Let us denote by $\otimes_\textsc{k}$ the Kronecker product. Then
\begin{equation}
 A \boxtimes X \simeq A \otimes_\textsc{k} X
\end{equation}
Let us note that
\begin{itemize}[label=$\rhd$]
 \item $A \boxtimes X$ is a tensor of order $4$: 
 \[
  A \boxtimes X \in E \otimes V \otimes F \otimes W
 \]
 acting on a matrix in $M \in F \otimes W$ with a matrix in $R \in E \otimes V$ as image
 \item $A \otimes_\textsc{k} X$ is a matrix with $mq$ rows and $mp$ columns, acting on a vector $u = \mathrm{vec} \: M$ with a vector $r = \mathrm{vec}\:R$ as image.
 \item  Thus, $A \otimes_\k X$ is a matricization of $A \boxtimes X$. It is a reshaping of $A \boxtimes X$ built by 
 \begin{equation}
  A \boxtimes X \in \underbrace{E \otimes V}_{\Ebold_1} \otimes \underbrace{F \otimes W}_{\Ebold_2} 
 \end{equation}
 and
 \begin{equation}
  A \otimes_\k X \in (E \otimes_\k V) \otimes (F \otimes_\k W) 
 \end{equation}
\end{itemize}

%
\section{Different products between matrices}\label{sec:alstruct:prodmat}
%

Matrix algebra has been used widely in statistics and data analysis. This has lead to the definition of several products which confer to the set of matrices a structure of algebra. Some of these products are briefly presented here.

\nS Let $\K^{m \times n}$ be the set of $m \times n$ matrices on a field $\K$, which will be omitted when it is not necessary, and $\K^{p \times q}$ the set of $p \times q$ matrices. Let
\[
 A = [a_{ij}]_{i,j} \in \K^{m \times n}, \qquad B = [b_{ij}]_{i,j} \in \K^{p \times q}
\]
Tracy-Singh and Khatri-Rao products assume that the matrices have been partitioned blockwise.

\nT{Hadamard product} Matrices $A$ and $B$ must have the same dimensions, i.e. $m=p$ and $n=q$. \kw{Hadamard product}\index{product!Hadamard} is elementwise product of two matrices of same size. It is called \kw{Schur product}\index{product!Schur} as well.
\begin{equation}
 A \odot B = 
 \begin{pmatrix}
 a_{11}b_{11} & \ldots & a_{1n}b_{1n} \\
 \vdots & & \vdots \\
 a_{m1}b_{m1} & \hdots & a_{mn}b_{mn}
 \end{pmatrix}
\end{equation}
Then
\begin{equation}
 \left. 
 \begin{array}{ll}
  A & \in \K^{m \times n} \\
  B & \in \K^{m \times n}
 \end{array}
 \right\}
 \quad \Longrightarrow \quad A \odot B \in \K^{m \times n}
\end{equation}
Hadamard product is 
\begin{itemize}[label=$\rightarrow$]
 \item commutative: $A \odot B=B \odot A$, 
 \item associative: $A \odot (B \odot C)= (A \odot B) \odot C$
 \item distributive over addition: $A \odot(B+C)= A \odot B + A \odot C$
\end{itemize}
An important theorem is that
\begin{equation}\label{eq:alstruct:prodmat:1}
 \rank A \odot B \leq (\rank A)(\rank B) 
\end{equation}
(see a short proof in \cite[p. 267]{Schott1997}). The Hadamard product $A \odot B$ is a submatrix of Kronecker product $A \otimes_\k B$.

\nB Hadamard product is sometimes denoted $\star$ or $\circ$.


\nT{Kronecker product} There are no constraints on the dimensions. \kw{Kronecker product}\index{product!Kronecker} is defined blockwise. It is denoted here $\otimes_\k$ between matrices to differentiate it from the tensor product, denoted $\otimes$.
\begin{equation}
 A \otimes_\k B = 
 \begin{pmatrix}
 a_{11}B & \ldots & a_{1n}B \\
 \vdots & & \vdots \\
 a_{m1}B & \hdots & a_{mn}B
 \end{pmatrix}
\end{equation}
There are no constraints on the dimension of the matrices. Then
\begin{equation}
 \left. 
 \begin{array}{ll}
  A & \in \K^{m \times n} \\
  B & \in \K^{p \times q}
 \end{array}
 \right\}
 \quad \Longrightarrow \quad A \otimes B \in \K^{mp \times nq}
\end{equation}
Kronecker product is often denoted $\otimes$, and not always distinguished from tensor product (see notes). One of the most fundamental property is the \kw{mixed-product property}: if
\begin{equation}
 A \in \K^{m \times n}, \qquad B \in \K^{p \times q}, \qquad C \in \K^{n \times k}, \qquad D \in \K^{q \times r}
\end{equation}
then
\begin{equation}
 (A \otimes_\k B)(C \otimes_\k D) = AC \otimes_\k BD
\end{equation}
As a corollary, by selecting $C=A^{-1}$ and $D=B^{-1}$
\begin{equation}
 (A \otimes_\k B)^{-1} = A^{-1} \otimes B^{-1}
\end{equation}




\nT{Tracy-Singh product} Let $A$ and $B$ be defined blockwise. Then, \kw{Tracy-Singh product}\index{product!Tracy-Singh} is defined  as a blockwise Kronecker products. It is denoted $\bowtie$. For example, let
\[
 A = 
 \begin{pmatrix}
  A_{11} & A_{12} \\
  A_{21} & A_{22}
 \end{pmatrix}
 , \qquad
 B = 
 \begin{pmatrix}
  B_{11} & B_{12} \\
  B_{21} & B_{22}
 \end{pmatrix} 
\]
Then
\begin{equation}
 A \bowtie B = 
 \begin{pmatrix}
  A_{11} \otimes_\k B_{11} & A_{11} \otimes_\k B_{12} & A_{12} \otimes_\k B_{11} & A_{12} \otimes_\k B_{12} \\ 
  A_{11} \otimes_\k B_{21} & A_{11} \otimes_\k B_{22} & A_{12} \otimes_\k B_{21} & A_{12} \otimes_\k B_{22} \\
  A_{21} \otimes_\k B_{11} & A_{21} \otimes_\k B_{12} & A_{22} \otimes_\k B_{11} & A_{22} \otimes_\k B_{12} \\ 
  A_{21} \otimes_\k B_{21} & A_{21} \otimes_\k B_{22} & A_{22} \otimes_\k B_{21} & A_{22} \otimes_\k B_{22} \\  
 \end{pmatrix}
\end{equation}
If we have one block only in $A$, this yields
\begin{equation}
 A \bowtie B = 
 \begin{pmatrix}
  A \otimes_\k B_{11} & A \otimes_\k B_{12}\\ 
  A \otimes_\k B_{21} & A \otimes_\k B_{22}\\
 \end{pmatrix}
\end{equation}
Hence, when $A$ is given blockwise too
\begin{equation}
 A \bowtie B = 
 \begin{pmatrix}
  A_{11} \bowtie B & A_{12} \bowtie B\\ 
  A_{21} \bowtie B & A_{22} \bowtie B\\
 \end{pmatrix}
\end{equation}
The dimension of $A \bowtie B$ does not depend on the selected block partitions, and is the same than the dimension of $A \otimes_\k B$: $A \bowtie B$ is a rearrangement of rows and columns of $A \otimes_\k B$. Then
\begin{equation}
 \left. 
 \begin{array}{ll}
  A & \in \K^{m \times n} \\
  B & \in \K^{p \times q}
 \end{array}
 \right\}
 \quad \Longrightarrow \quad A \bowtie B \in \K^{mp \times nq}
\end{equation}


\paragraph{Khatri-Rao product:} Both matrices must have the same number of columns, i.e. $n=q$. \kw{Khatri-Rao product}\index{product!Khatri-Rao} is denoted $\ast$. Let 
\[
 A = [\abold_1, \ldots, \abold_n], \qquad B = [\bbold_1, \ldots, \bbold_n]
\]
be given columnwise with
\[
 \abold_j \in \R^m, \qquad \bbold_j \in \R^p, \qquad 1 \leq j \leq n
\]
Then, columnwise
\begin{equation}
 A \ast B = [\abold_1 \otimes_\k \bbold_1, \ldots, \abold_n \otimes_\k \bbold_n]
\end{equation} 
where $\otimes_\k$ is the Kronecker product\footnote{The Kronecker product of two column vectors of dimensions $m$ and $p$ is a column vector of dimension $mp$.}. In more details, we have
\begin{equation}
 A \ast B = 
 \begin{pmatrix}
  a_{11}\bbold_1 & a_{12}\bbold_2 & \ldots & a_{1n}\bbold_n \\
  a_{21}\bbold_1 & a_{22}\bbold_2 & \ldots & a_{2n}\bbold_n \\
  \vdots & \vdots & \ddots & \vdots \\
  a_{n1}\bbold_1 & a_{n2}\bbold_2 & \ldots & a_{mn}\bbold_n
 \end{pmatrix}
\end{equation}
Then
\begin{equation}
 \left. 
 \begin{array}{ll}
  A & \in \K^{m \times n} \\
  B & \in \K^{p \times n}
 \end{array}
 \right\}
 \quad \Longrightarrow \quad A \ast B \in \K^{mp \times n}
\end{equation}

\nT{CP-product} The \kw{CP-product}\index{product!CP}  is useful for CP decomposition of a tensor CP. Let us give first a definition for $3-$modes tensors. Let $X,Y,Z$ be three matrices defined column-wise:
\begin{equation}
    X = [\xbold_a]_{a \leq r}, \qquad Y = [\ybold_a]_{a \leq r}, \qquad Z = [\zbold_a]_{a \leq r}, \qquad 
\end{equation}
Then
\begin{equation}
    \langle X | Y | Z \rangle = \sum_{a=1}^r \xbold_a \otimes \ybold_a \otimes \zbold_a
\end{equation}
If $X \in \K^{m \times r}, Y \in \K^{n \times r}, Z \in \K^{p \times r}$, then $\langle X | Y | Z \rangle \in \K^{m \times n \times p}$. The reason for this notation is that
\begin{equation}
    \langle X | Y \rangle = XY^\t
\end{equation}
(if $r=1$, it is the inner product). I have adapted here the \emph{bra-ket} notation $\langle x|y \rangle$ for inner product in quantum mechanics. It can be naturally extended to more than three matrices as
\begin{equation}\label{eq:sec:alstruct:prodmat:1}
    \langle X_1 | \ldots | X_d \rangle = \sum_{a=1}^r \xbold_1^{(a)}  \otimes \ldots \otimes \xbold_d^{(a)} 
\end{equation}

\nT{Notations} Here is a summary of the notations adopted here and a comparison with those in \cite{Kolda2009,Seber2008,Cichocki2016}.
\begin{center}
 \begin{tabular}{lcccccc}
 \hline
  product & symbol & other name & in \cite{Kolda2009} & in \cite{Seber2008} & in \cite{Cichocki2016}\\
  \hline
  tensor product & $\otimes$ & outer product & $\circ$ &  & $\circ$\\
  Kronecker product & $\boxtimes$ &  & $\otimes$ & $\otimes$ & $\otimes$ \\
  Hadamard product & $\odot$ & Schur product & $\ast$ & $\circ$ & $\circledast$ \\
  Khatri-Rao product & $\ast$ & & $\odot$ & $\odot$ & $\odot$ \\
  CP product & $\langle . | \ldots | . \rangle$ & & $\llbracket ... \rrbracket$ & \\
  Tracy-Singh product & $\bowtie$ & & & \\
  \hline
 \end{tabular}
\end{center}

\notes See \cite{Seber2008,Horn2012} for use of matrix algebra in statistics. See \cite{Liu2008} for definitions of various matrix products, and for notations which we have adopted here. There is no consensus between various authors on notations of even classical products between matrices, and some may clash with other notation systems while extending them to tensor algebra. So, we have been quite cautious in the selection of a standard and consistent notation system which can be extended to tensors. \textbf{Hadamard product:} The Hadamard product is presented in detail in \cite[chapter 5]{Horn2012} where it is denoted $\circ$. The reader interested in this product is encouraged to read this rich chapter, and historical aspects pp. 302 \& \emph{sq.}. According to \cite[p. 298]{Horn2012}, the most basic property of Hadamard product is that the cone of positive semidefinite matrices is closed under this product. It is presented as well in detail in \cite[section 7.6]{Schott1997}, with notation $\odot$, where the assertion of equation (\ref{eq:alstruct:prodmat:1}) is proved p. 267. \textbf{Kronecker product:} This product is studied in detail in \cite[section 4.2]{Horn2012} under the name of Kronecker product, \kw{outer product}\index{product!outer}, or tensor product. Hence, we set a distinction between Kronecker and tensor product, which has not been retained in \cite{Horn2012} (they are identical up to an isomorphism). The mixed product property is shown p. 244. See \cite[p. 243 \emph{\& sq.}]{Horn2012} for historical aspects. According to \cite[section 4.3]{Horn2012}, many linear equations can be conveniently represented with use of Kronecker product because, with adequate dimensions, $
 AXB=C \quad \Longleftrightarrow \quad (B^\t \otimes_\k A) \, \vec X = \vec C$. Kronecker product is developed in \cite[section 7.3]{Schott1997}, with many equalities, among which $\Tr A \otimes_\k B=(\Tr A)(\Tr B)$ (p. 255). See as well \cite{vanLoan2000} for a thorough survey of usefulness of Kronecker product in matrix numerical analysis. \textbf{Tracy-Singh product:} It has been defined by Tracy and Singh in Tracy, D. S. \& Singh, R. P. (1972). A New Matrix Product and Its Applications in Matrix Differentiation. \emph{Statistica Neerlandica},  \textbf{26(4):}143-157. doi:10.1111/j.1467-9574.1972.tb00199.x. \textbf{Khatri-Rao product:} It has been defined in Khatri C. G. \& Rao, C. R.  (1968). Solutions to some functional equations and their applications to characterization of probability distributions. \emph{Sankhya}, \textbf{30:}167-180.

%% file: tensor_rank.tex
Let $A \in \R^{m \times n}$ be a matrix. It is an expression of a linear map $A \in \L(F,E)$ denoted $A$ as well in a given basis for $E$ and $F$ with $\dim E=m$ and $\dim F=n$. The rank of $A$ is the dimension of its image space $A(F) \subset E$. If it is $r$, it is possible to chose a basis in $F$ denoted $(\ubold_1, \ldots, \ubold_r, \ubold_{r+1}, \ldots,\ubold_n)$ such that $A\ubold_k=0$ for $k>r$. In such a basis, the $n-r$ last columns of $A$ are zero. Hence, it is possible to reduce the dimension of a matrix without losing information if its rank is not full $(r < n)$. If, again, as reader you live on Sirius, I can send you exactly matrix $A$ by sending a $m \times r$ matrix in a suitable basis which is less costly than a $m \times n$ matrix. The rank has several properties. First, if $\rank A=r$, there exists a decomposition
\[
A = \sum_{a=1}^{\rank A} \xbold_a \otimes \ybold_a
\]
as a sum of $r$ rank one matrices. The rank can be defined as the smallest integer $r \in \N$ for which such a decomposition exists. Then $(\ybold_a)_{a \leq r}$ is a basis of $F$ and 
\[
\xbold_a = A\ybold_a^*
\]
As this is true for any basis $(\ybold_a)_a$, it is possible to chose an orthonormal basis for $(\ybold_a)_a$ for which $\xbold_a = A\ybold_a$. 

\nB Let us now ask the following question: a full rank matrix $A \in \R^{m \times n}$ being given with $n < m$, and an integer $r<n$, find a rank $r$ matrix $A_r$ such that $\|A-A_r\|$ is minimal, where $\|.\|$ is the Frobenius or $\ell^2$ norm $(\|A\|^2=\sum_{i,j}a_{ij}^2)$. The solution is well known by so called Eckart-Young theorem: let
\[
A = \sum_{a=1}^n \, \sigma_a \, \ubold_a \otimes \vbold_a
\]
be the Singular Value Decomposition of $A$ withy $\sigma_1 \geq \sigma_2 \geq \ldots \geq \sigma_n \geq 0$. Then
\[
A_r = \sum_{a \leq r} \, \sigma_a \, \ubold_a \otimes \vbold_a
\]
This has several consequences. First, the image of $A_r$ is a space $E_r$ with basis $(\ubold_a)_{a \leq r}$. It is part of the definition of SVD that the basis $(\ubold_a)_a$ of $\Im A \subset E$ is orthonormal. Then, the spaces $(E_a)_a$ are nested
\[
E_1 \subset E_2 \subset \ldots  \subset E_r
\]
and orthogonal
\[
\forall \: a \leq r, \quad E_a \perp (E_1 \oplus \ldots \oplus E_{a-1})
\]
Hence, there are at least two possible views on Eckart-Young theorem:
\begin{itemize}[label=$\rightarrow$]
    \item $A_r$ is a solution of a minimization problem: $\|A-A_r\|$ is minimal
    \item a geometric approach: the spaces spanned by this solution are nested and built on a orthonormal basis of $\Im A \subset E$.
\end{itemize}
Such a result is deeply rooted in the property that a matrix $A$ is the expression of a linear map in a given basis, and the properties of the SVD of a matrix, presented in detail in section \ref{sec:inter:svd}. 

\nB It seems natural to look for a generalization of SVD and its use for a best low rank approximation of a tensor, with both an optimization and a geometric approach. However, there is no canonical way to associate a linear map to a tensor\marginpar{\dbend}: several matricizations can be built by choosing different modes (see section \ref{sec:contract:alamat}). The extension of the concept of rank is not straightforward, and several and not equivalent notions of rank can be seen as an extension of the rank of a matrix. For a given tensor $\Abold$, one has to chose between 
\begin{itemize}[label=$\rightarrow$]
    \item the optimization approach: the smallest number of rank one tensors the sum of which is $\Abold$: this leads to CP-rank\\
    \item the nestedness and orthogonality of subspaces which are associated with various (not one \marginpar{$\dbend$}) matricizations of $\Abold$: this leads to Tucker model and multilinear rank, with various guises.
\end{itemize}
There is no closed form like the SVD for solving the best low rank approximation of a tensor, be it for CP or Tucker. Thus is done in an algorithmic way. Such developments are presented in part \ref{part:blra}.

\nB It will be shown that best low rank approximation of a tensor for CP rank may be ill-posed, because the manifold of tensors of a given CP-rank is not closed (so the projection on it may be ill defined). This leads to the notion of border rank. On the contrary, it will be shown that best low rank approximation of Tucker rank is well-posed. Beyond CP and Tucker, a variety of ranks have been proposed in the literature, each built on a given decomposition of a tensor, in the same way that a graphical model in statistics provides a simplification of a joint law by factorization of smaller elements. This has been fully developed in many-body quantum systems with the notion of tensor network. I present later (see chapter \ref{chap:tt}) one of the most accessible and best understood tensor network: the Tensor-Train  decomposition. 

\nB  The best approximation of a tensor by another of prescribed rank, be it for CP or Tucker or TT rank, depends on the norm selected which induces a distance. Different contexts may lead to different choices for a norm. Some elements are given on the diversity of norms in a finite dimensional vector space, which can easily be extended to matrices and tensors. It is highlighted that the reason why "real world matrices" are better represented by a prescribed rank approximation than random matrices is not fully understood, even if some hints can be proposed. 

\begin{description}
\item[Chapter \ref{chap:rank3}] presents the various extensions of the concept of rank to tensors
\item[Chapter \ref{chap:inter}] presents some classical norms in vector spaces, which induce each a distance.
\end{description}

%
\chapter{Tensor ranks for $d-$modes tensors}\label{chap:rank3}
%

Many numerical calculations on large matrices are made possible when scaling the dimension by working on low rank best approximations of them, for a given norm (most of the time Frobenius norm, i.e. $\ell^2$ norm). It has been known for over a century that best low rank approximation of a given matrix is provided by truncating its Singular Value Decomposition. It has then seemed natural to extend the notion of SVD and rank to tensors, and develop a similar theory for best low rank approximation. This has proved to be fallacious. 

\nB Indeed, the singular values and vectors of a matrix $A$ are intimately linked with the singular values and vector of matrix $A^\t A$ which in turn are attached to the linear map associated to $A^\t A$. There is no linear map canonically associated to a tensor. A tensor $\Abold \in E \otimes F \otimes G$ can be associated with a linear map in $\L(G \, , \, E\otimes F)$ or $\L(F\otimes G \, , \, E)$ (see section \ref{sec:contract:alamat}). Each linear map brings its own set singular vectors and singular values.

\nB The tensor rank is a projective property of a tensor, as it is invariant by any transformation $\lambda \Abold$ provided $\lambda \neq 0$. There is an $s$ at rank\emph{s} because there exists several definitions of the rank of a tensor for $d>2$. This is a major difference with $d=2$, where there is one single definition of the rank of a matrix. Let us note that all definitions of the rank of a tensor for $d \geq 2$ do not yield the rank of a matrix when $d=2$: there exists definitions specific to $d-$modes tensors for $d \geq 3$. We start with several ranks of a tensor when $d=3$. Extension to $d>3$ will be straightforward.

\nB This chapter is organized as follows:
\begin{description}
\item [section \ref{sec:rank3:def}] The classical definition of rank as the minimum integer $r$ such that a tensor can be decomposed as a sum of $r$ tensors of rank one is called CP rank, and presented in this chapter. However, topological arguments lead to complement such a basic definition with more elaborate notions, like border rank, or typical rank. These definitions of different ranks are given first for $3-$modes tensors for sake of clarity.
\item[section \ref{sec:rank3:bound}] Finding the maximal rank of a tensor in $\R^{m \times n \times p}$ is far from being simple (whereas the maximum rank of a matrix in $\R^{m \times n}$ with $n \leq m$ is known to be $n$). Such an issue has been thoroughly studied in the 70's, following Strassen seminal paper on the complexity of matrix multiplication the result of which is built on the study of the rank of a given tensor. Elementary results are given in this section.
\item[section \ref{sec:rank3:rank222}] Analytical results are given for typical ranks in $\R^{2 \times 2 \times 2}$ which are known to be $\{2,3\}$.
\item[section \ref{sec:rank3:tucker}] The definition of Tucker rank or multilinear rank is given. As in PCA, the unknowns are vector spaces (modes) of prescribed dimensions. See also section \ref{sec:dmodestens:basis}.
\item[section \ref{sec:rank3:rankd}] Finally, all these definitions can be extended in a straightforward way to $d-$modes tensors.
\end{description} 

\nT{Notations} Here are some notations in this chapter:\\ \\
\begin{center}
    \ovalbox{
        \begin{tabular}{cl}
             $\Abold$ & a tensor \\
             $\brank \Abold$ & border rank of $\Abold$ \\
             $\mathbf{br}(\Abold)$ & border rank of $\Abold$ \\
              $\mathbb{BP}_r$ & set of tensors of border rank $r$ \\
             $\mathbb{BP}_{\leq r}$ & set of tensors of border rank equal to or less than $r$ \\
             $\mathbb{CP}_r$ & set of tensors of CP-rank $r$ \\
             $\mathbb{CP}_{\leq r}$ & set of tensors of CP-rank equal to or less than $r$ \\
             $\mathbf{mr}(\Abold)$ & multilinear rank of tensor $\Abold$ \\
             $r$ & rank of a tensor \\
             $\rank \Abold$ & CP rank of $\Abold$ \\
             $\rbold(\Abold)$ & CP rank of $\Abold$ \\
             $\rbold(m,n,p)$ & rank of a general tensor in $\K^{m \times n \times p}$ \\
              $\rr(m,n,p)$ & maximal typical rank in $\R^{m \times n \times p}$\\
               $\rr(\C;m,n,p)$ & maximal typical rank in $\C^{m \times n \times p}$\\
             $\rr(\R;m,n,p)$ & maximal typical rank in $\R^{m \times n \times p}$\\
             $\mathbf{tr}(\C;m,n,p) $ & typical rank of a tensor in $\C^{m \times n \times p}$ \\
              $\mathbf{tr}(\R;m,n,p) $ & typical rank of a tensor in $\R^{m \times n \times p}$ \\
             $\mathbf{tr}(m,n,p)$ & typical rank of a tensor in $\K^{m \times n \times p}$ \\
        \end{tabular}
    }
\end{center}

%
\section{Definitions for $3-$modes tensors}\label{sec:rank3:def}
%

In this section, the dimension of the spaces we work with will be crucial. Then, without loss of generality, we will use
\[
 E \simeq  \K^m \qquad \mbox{or} \qquad E \quad \mbox{with} \quad \dim E= m
\]
and similarly for $F, G, \ldots$. Here, $\K$ can be $\R$ or $\C$.

\nS \textbf{Rank:} Let $\Abold \in \K^m \otimes \K^n \otimes \K^p$. The \kw{rank} of $\Abold$, denoted $\rank \Abold$, or $\rbold(\Abold)$, is the smallest integer $r$ such that there exists $r$ triplets $(\xbold_a,\ybold_a,\zbold_a)_{a \leq r}$ with $\xbold_a \in E, \ybold_a \in F, \zbold_a \in G$ such that $\Abold = \sum_{a=1}^r \xbold_a \otimes \ybold_a \otimes \zbold_a$
\begin{equation}
  \rbold(\Abold) = \inf \: \left\{r \in \N  \quad \bigg| \quad  \exists \: (\xbold_a,\ybold_a,\zbold_a)_a \st \Abold = \sum_{a=1}^r \xbold_a \otimes \ybold_a \otimes \zbold_a\right\}
\end{equation}
This rank is often called the \kw{CP-rank}\index{rank!CP}. It is a natural extension of the rank of a matrix. A rank 1 tensor is $\Abold = \xbold \otimes \ybold \otimes \zbold$. In such a case, we have $a_{ijk}=x_iy_jz_k$. We will denote by $\mathbb{CP}_r(\K^m \otimes \K^n \otimes \K^p)$ or $\mathbb{CP}_r(n,m,p)$ or $\mathbb{CP}_r$ when non ambiguous the set of tensors of rank $r$
\begin{equation}
 \mathbb{CP}_r = \{\Abold \in \K^m \otimes \K^n \otimes \K^p \mid \rank \Abold = r\}
\end{equation}
and by $\mathbb{CP}_{\leq r}$ the set of tensors of rank $r' \leq r$
\begin{equation}
 \mathbb{CP}_{\leq r} = \{\Abold \in \K^m \otimes \K^n \otimes \K^p \mid \rank \Abold \leq r\}
\end{equation}
We denote by $\rr(n,m,p)$ the smallest integer $r$ such that
\begin{equation}
 \forall \: \Abold \in \K^m \otimes \K^n \otimes \K^p, \quad \rank \Abold \leq \rr
\end{equation}
This bound depends on the field. Hence $\rr(\R;n,m,p)$ is the bound over $\R$ and $\rr(\C;n,m,p)$ over $\C$. As this bound depends on the existence of solutions of a set of polynomial equations, we have
\begin{equation}
    \rr(\C;n,m,p) \leq \rr(\R;n,m,p)
\end{equation}
So, $\rr(n,m,p)$ will denote $\rr(\R;n,m,p)$.


\nS \textbf{Border rank:}  One difficulty with above definition is that the set $\C\P_r$ of tensors of a given rank is not closed, contrary to what is happening for matrices: there exists series of rank $r$ tensors which converge to a tensor of rank $r' \neq r$. Hence, the projection on $\mathbb{CP}_r$ may be ill-defined, which prevents processing best low rank approximations. This induces the notion of border rank. Let $\Abold \in \K^m \otimes \K^n \otimes \K^p$. The \kw{border rank}\index{rank!border} of $\Abold$ ($\brank \Abold$, or $\mathbf{br}(\Abold)$) is the smallest integer $r$ such that there exists a series
\[
 (\Abold_k)_k, \qquad k \in \N, \quad \Abold \in \K^m \otimes \K^n \otimes \K^p
\]
such that 
\begin{equation}
\left\{
 \begin{array}{l}
   \Abold = \displaystyle \lim_{k \rightarrow \infty} \: \Abold_k \\
   \\
   \mbox{and} \\
   \\
   \forall \: k \in \N, \quad \rank \Abold_k = r
 \end{array}
\right.
\end{equation}
The set of tensors of border rank $r$ is denoted $\mathbb{BR}_r$. By definition, $\mathbb{BR}_r$ is closed. But it is not is the closure of $\mathbb{CP}_r$\marginpar{\dbend}: $\mathbb{BR}_r \neq  \overline{\mathbb{CP}_r}$. Such a definition is meaningful because the set $\mathbb{CP}_r$ is not closed. If it were closed, the border rank would be the rank as the limit of $(\Abold_k)_k$ would be in $\mathbb{CP}_r$, hence of rank $r$. There exists tensors $\Abold$ for which $\rank \Abold \neq \brank \Abold$. \cite{Silva2008} have exhibited a simple example of a tensor $\Abold \in \R^{m \times n \times p} $ of rank 3 which can be approximated arbitrarily closely by tensors $\Abold_k$ of rank 2. This example is
\begin{equation}
 \left\{
  \begin{array}{lcl}
  \Abold &=& \xbold_1 \otimes \xbold_2 \otimes \ybold_3 + \xbold_1 \otimes \ybold_2 \otimes \xbold_3 + \ybold_1  \otimes \xbold_2 \otimes \xbold_3 \\
  &&\\
  \Abold_k &=& \displaystyle k\left(\xbold_1+\frac{1}{k}\ybold_1\right) \otimes \left(\xbold_2+\frac{1}{k}\ybold_2\right) \otimes \left(\xbold_3+\frac{1}{k}\ybold_3\right) -k\, \xbold_1 \otimes \xbold_2 \otimes \xbold_3
  \end{array}
 \right.
\end{equation}
Then
\begin{equation}
\Abold - \Abold_k = \frac{1}{k} (\ybold_1 \otimes \ybold_2 \otimes \xbold_3 + \ybold_1 \otimes \xbold_2 \otimes \ybold_3 + \xbold_1 \otimes \ybold_2  \otimes \ybold_3) + \frac{1}{k^2} \ybold_1 \otimes \ybold_2 \otimes \ybold_3
\end{equation}
and
\begin{equation}
 \lim_{k\rightarrow \infty} \: \|\Abold - \Abold_k\| = 0
\end{equation}
It is clear that $\rank \Abold_k=2$, and \cite{Silva2008} have shown that $\rank \Abold=3$ if and only if $(\xbold_i,\ybold_i)$ are linearly independent for $i \in \{1,2,3\}$.


\nS \textbf{Typical rank:} Let $\K = \R$ or $\C$. A rank $r \in \N$ for a $(m,n,p)-$tensor or a tensor in $\K^{m \times n \times p}$ is a \kw{typical rank}\index{rank!typical} if there exists an open set $\Omega \subset \K^{m \times n \times p}$ such that all tensors $\Abold \in \Omega$ have rank $r$, or
\begin{equation}
    \exists \:\: \Omega, \:\: \mbox{open set} \quad \st \quad \forall \: \Abold \in \Omega, \quad \rank \Abold=r
\end{equation}
Such a rank is not unique (contrary to typical rank of matrices in $\mathbb{K}^{m \times n}$ which is unique and equal to $\min \: (m,n)$).The set of all possible typical ranks is denoted $\mathbf{tr}(m,n,p)$. As $\rank \lambda \Abold = \rank \Abold$ for $\lambda \in \K^*$, this can be rephrased as follows. Let $\S^{m \times n \times p}$ be the unit sphere in $\K^{m \times n \times p}$
\[
 \S^{m \times n \times p} = \{\Abold \in \K^{m \times n \times p} \: \mid \: \|\Abold\|=1\}
\]
Then, the \kw{typical rank}\index{rank!typical} in $\K^{m \times n \times p}$, is the the set of values $r$ such that the measure of the set of tensors $\Abold \in \S^{m \times n \times p}$ of rank $r$ is strictly positive
\begin{equation}
 \mathbf{tr}(m,n,p) = \left\{ r \in \N \: : \:  \mu \left( \mathbb{CP}_r \cap \S\right) > 0\right\}
\end{equation}
Computing $\mathbf{tr}(m,n,p)$ is still an open (and difficult) problem. Values are known for a limited set of dimensions $(m,n,p)$ (see some elements in section \ref{sec:rank3:bound}).

\nS We have:\\
\\
\begin{center}
 \begin{tabular}{l|cl|c}
 \hline 
 Name & & Definition & Set\\
 \hline
  rank & $r=\rbold(\Abold)$ & minimum integer $r$ such that & $\mathbb{CP}_r$\\
  & & $\Abold$ is written as a sum of $r$ elementary tensors  \\
  & & &\\
  border rank &  $r=\mathbf{br}(\Abold)$ & minimum integer $r$ such that & $\mathbb{BR}_r$ \\
  & & $\Abold$ is the limit of a sequence of tensors of rank $r$\\
  & & & \\
  typical rank & $r \in \mathbf{tr}(m,n,p)$ & the set of integers $r$ such that & \\
  & & there is an open set of tensors with rank $r$\\
  \hline
 \end{tabular}
\end{center}

\nT{Further notations} I will denote $\rbold(m,n,p)$ the rank of a general tensor in $\K^{m \times n \times n}$. 

\nS The rank in $E \otimes F \otimes G$ behaves very differently from the rank of a matrix in $E \otimes F$. Let $A \in E \otimes F$. A linear map $L_\textsc{a} \in \L(F,E)$ can be attached to $A$, and 
\begin{equation}
 \rank A = \dim (\Im L_\textsc{a})
\end{equation}
Moreover
\begin{equation}
 \dim (\ker L_\textsc{a}) + \dim (\Im L_\textsc{a})= \dim F
\end{equation}
(This is the fundamental theorem of linear algebra). Let $\dim E=m, \dim F = n$ and $n \leq m$. The set $\Gamma_n(E \otimes F)$ of matrices of rank $n$ is an open set, dense in $E \otimes F$. Its complement is a closed set. Indeed, the set of matrices of rank $< n$ is
\begin{equation}
 \{ A \in E \otimes F \: : \: \rank A < n \} = \{A \in E \otimes F \: : \: \det A^\t A = 0\}
\end{equation}
and is the preimage of a closed set $\{0\}$ by a continuous function; hence closed.

\nS The set of matrices $A \in E \otimes F$ of rank $r$ is denoted $\Gamma_r(E \otimes F)$ or, omitting $E$ and $F$
\begin{equation}
 \Gamma_r = \{A \: : \: \rank A=r\}
\end{equation}
$\Gamma_r$ is closed as for each $r$, $\Gamma_r$ is the set of zeros of a family of polynomials: some minors in $A$. It is a finite intersection of preimages of $\{0\}$ (one preimage per minor equal to zero), and as such, a closed set. Nothing of this sort happens for rank of $d-$tensors, for $d \geq 3$. Indeed :
\begin{itemize}[label=$\rhd$]
 \item the rank of a given tensor may depend on the field over which it is computed
 \item The set of tensors of a given rank may be not closed (neither open)
 \item The border rank may be different from the rank
 \item The typical rank for a given field and given dimensions may be not unique.
\end{itemize}
A full understanding of these questions still is a matter of research (see e.g. \cite{Landsberg2012,Landsberg2017}).

\notes The notion of tensor rank has been introduced by Hitchcock in 1927\footnote{Hitchcock, F.L., 1927 : The expression of a tensor or a polyadic as a sum of products. \emph{Journal of Mathematics and Physics}, \textbf{6:}164-189} under the name Gibb's polyadic, or sum of polyads. A polyad is a generalization of the notion of dyads, introduced by Gibbs while defining new tools in vector analysis for physics\footnote{see J. W. Gibbs \& E. B. Wilson - 1926 - Vector Analysis, Yale University Press, p. 264} and for which he developed a new algebra. However, most of the developments in this paper are devoted to define the rank of a $d-$tensor as the rank of any linear applications associated to this tensor (called by Hitchcock \emph{multiplex rank}), although finding the minimal integer $r$ such that a tensor can be written as a sum of polyadics (of elementary tensors) has been announced as the main topic of the paper. A rich presentation of the geometric nature of the rank of tensors is \cite{Landsberg2012}. It makes some connections with elaborate notions in projective geometry like Segre varieties. See as well \cite[sect. 2.1]{Landsberg2017} for rank and border rank.

%
\section{Some bounds for the rank}\label{sec:rank3:bound}
%

The typical ranks of a tensor $\Abold \in \K^{m \times n \times p}$ are the integers $r \in \tbold\rbold(\Abold)$ such that there exists an open set $\Omega \subset \K^{m \times n \times p}$ with $\rbold(\Abold)=r$ for any $\Abold \in \Omega$. Finding the typical ranks for any set of dimension $(m,n,p)$ is an open (and difficult) problem. Here, some elementary results on bounds of typical ranks are given.

\nS Let $\Abold \in \R^m \otimes \R^n \otimes \R^p$ with $m \geq n \geq p$. Then
\begin{equation}
 \rr(m,n,p) \leq np
\end{equation}
Indeed, let $(\ebold_i)_i$ be an orthonormal basis of $\R^m$, $(\fbold_j)_j$ of $\R^n$ and $(\gbold_k)_k$ of $G$. Then
\begin{equation}
 \begin{array}{lcl}
  \Abold &=& \displaystyle \sum_{i,j,k} \, a_{ijk} \, \ebold_i \otimes \fbold_j \otimes \gbold_k \\
  &=& \displaystyle \sum_{j,k} \,  \underbrace{\left(\sum_i \, a_{ijk}\, \ebold_i\right)}_{\xbold_{jk} \in E} \otimes \fbold_j \otimes \gbold_k\\
  &=& \displaystyle \sum_{j,k} \xbold_{jk} \otimes \fbold_j \otimes \gbold_k
 \end{array}
\end{equation}
which is a decomposition with $np$ elementary tensors at most. 

\nT{Dimension analysis}By a simple dimension analysis, it is possible to show that
\begin{equation}\label{eq:rank3:bound:1}
 \rbold(m,n,p) \geq \frac{mnp}{n+m+p-2}
\end{equation}
\begin{proof}
Indeed, let 
\[
\begin{array}{lcl}
  \Abold &=& \displaystyle \sum_{i=1}^m\sum_{j=1}^n\sum_{k=1}^p \, a_{ijk}\ \ebold_i \otimes \fbold_j \otimes \gbold_k \\
  \\
  &=& \displaystyle \sum_{a=1}^r \xbold_a \otimes \ybold_a \otimes \zbold_a
\end{array}
\]
First decomposition $\sum_i\sum_j\sum_k\ldots$ requires $mnp$ terms. Second decomposition $\sum_a\ldots$ requires $r$ components $\xbold_a \otimes \ybold_a \otimes \zbold_a$. One can set $\|\xbold_a\|=\|\ybold_a\|=\|\zbold_a\|=1$ and add a scalar $\lambda$ to write it $\lambda_a \, \xbold_a \otimes \ybold_a \otimes \zbold_a$. Such a component requires $1+(m-1)+(n-1)+(p-1)=m+n+p-2$ terms, hence $r(m+n+p-2)$ terms for $\sum_a \ldots$. For a locally surjective continuous mapping to exist between both forms, one must have
\[
r(m+n+p-2) \geq mnp 
\]
or
\[
 r \geq \frac{mnp}{m+n+p-2}
\]
\end{proof}
\noindent This shows that, if $r$ is a typical, rank for a tensor in $\K^{m \times n \times p}$, equation (\ref{eq:rank3:bound:1}) must be fulfilled. This can be written as a system 
\begin{equation}
 \forall \: i,j,k, \qquad \sum_a\,x_{ia}y_{ja}z_{ka} = a_{ijk}
\end{equation}
of $mnp$ cubic polynomials equations with $r(m+n+p)$ unknowns (dropping $\|\xbold_a\|=\|\ybold_a\|=\|z_a\|=1$ for sake of simplicity). These equations are polynomials, and their solutions may depend on whether the field is $\C$ or $\R$. If there is no solution with $r=mnp/(m+n+p)$ on $\R$, one may have a typical rank larger than this lower bound. Hence, a rank can be typical for $\R$ but not for $\C$ (examples will be given later in $\K^{2 \times 2 \times 2}$, where typical rank in $\C$ is $\{2\}$ and in $\R$ is $\{2,3\}$).

\nS If $m=n=p$ denoted $n$, this leads to the following table, with $r = \large\lceil \frac{n^3}{3n-2} \large\rceil$:\\
\begin{center}
 \begin{tabular}{r|r||r|r}
  \hline
  $n$ & $r$ & $n$ & $r$\\
  \hline
  2 & 2 & 7 & 19 \\
  3 & 4 & 8 & 24 \\
  4 & 7 & 9 & 30 \\
  5 & 10 & 10 & 36 \\  
  6 & 14 & 11 & 43 \\
  \hline
 \end{tabular}
\end{center}
We haven more generally,
\begin{equation}\label{sec:rank3:bound:eq1}
 r(n,n,n) \geq \frac{n^3}{3n+2} \approx \frac{n^2}{3} \quad \mbox{when } n \rightarrow \infty
\end{equation}
So, there is a wide range of possibilities, as this shows that, if $r \in \mathbf{tr}(n,n,n)$
\[
 \frac{n^2}{3} \leq r \leq n^2
\]
We show next that such a gap has been narrowed for some specific cases.

\nT{Typical rank in $\R^{n \times n \times 2}$} Let us have $\Abold \in \R^n \otimes \R^n \otimes \R^2$. Such a tensor can be visualized as\\
\begin{center}
 \begin{tikzpicture}
        \node () at (-.7,1){$\Abold = $};
        \node () at (1.7,0.3){$A_1$};
        \node () at (2.2, 0.8){$A_2$};
        \draw (0,0) rectangle (2,2) ;
        \draw (0.5, 2) -- (0.5, 2.5) -- (2.5,2.5) -- (2.5,0.5) -- (2,0.5) ; 
 \end{tikzpicture}
\end{center}

Application of equation \ref{eq:rank3:bound:1} yields
\begin{equation}
 r(\Abold) \geq \frac{2n^2}{2n+2-2}=n
\end{equation}
So
\begin{equation}
 n \leq r(\Abold) \leq 2n
\end{equation}
Following \cite{Lafon1975}, we explicit here sets of tensors with typical rank $n$. Let $(\ebold_1,\ebold_2)$ be an orthonormal basis of $\R^2$, typically
\begin{equation}
 \ebold_1 = 
 \begin{pmatrix}
  1 \\
  0
 \end{pmatrix}
 , \qquad \ebold_2 = 
 \begin{pmatrix}
  0 \\
  1
 \end{pmatrix}
\end{equation}
Then,
\begin{equation}
 \Abold = A_1 \otimes \ebold_1 + A_2 \otimes \ebold_2
\end{equation}
with
\[
 A_1 = \Abold.\ebold_1, \qquad A_2 = \Abold.\ebold_2
\]
As in general $\rank A_1= \rank A_2 = n$, one expects $n$ terms for $A_1 \otimes \ebold_1$ and for $A_2 \otimes \ebold_2$, hence $\rank \Abold = 2n$. \cite{Lafon1975} has shown that, if $A_1^{-1}A_2$ is diagonalisable, $\rank \Abold = n$.\label{sec:rank:rank3:bound:ref1}
\begin{proof}
Let $A_1^{-1}A_2 \ubold_i = \lambda_i \ubold_i$. Then (see equation \ref{eq:matdecomdiago})
\[
 A_1^{-1}A_2 = \sum_{i=1}^n \lambda_i \, \ubold_i \otimes \ubold_i^*
\]
and
\[
 A_2 = \sum_{i=1}^n \lambda_i \, A_1\ubold_i \otimes \ubold_i^*
\]
(see equation \ref{eq:elemat}). We have
\[
 A_1 = \sum_i A_1\ubold_i \otimes \ubold_i^*
\]
(see equation \ref{eq:matelem_nonortho} with $\fbold_j \equiv \ubold_i^*$ knowing that $\ubold_i^{**}=\ubold_i$). Thus
\[
 \begin{aligned}
  \Abold &= A_1 \otimes \ebold_1 + A_2 \otimes \ebold_2 \\
  &= \left(\sum_{i=1}^n A_1\ubold_i \otimes \ubold_i^*\right) \otimes \ebold_1 + \left(\sum_{i=1}^n \lambda_i \, A_1\ubold_i \otimes \ubold_i^*\right) \otimes \ebold_2 \\
  &= \sum_{i=1}^n (A_1\ubold_i \otimes \ubold_i^* \otimes \ebold_1) +  (\lambda_i \, A_1\ubold_i \otimes \ubold_i^* \otimes \ebold_2)  \\
  &= \sum_{i=1}^n A_1\ubold_i \otimes \ubold_i^* \otimes (\ebold_1 + \lambda_i\ebold_2) 
 \end{aligned}
\]
and 
\[
\rank \Abold \leq n
\]
As by $(\ref{eq:rank3:bound:1})$ one has $r \geq n$, this shows that $r=n$
\end{proof}
\noindent In particular, 2 is a typical rank in $\K^{2 \times 2 \times 2}$ for $K=\R$ or $\C$. One can observe that
\begin{itemize}[label=$\rightarrow$]
    \item This is true for $\K=\C$ because $\lambda_i \in \C$
    \item This is not systematically true for $\K=\R$ or $\Q$, because we can have $\lambda_i \notin \R$ or $\lambda_i \notin \Q$
\end{itemize}
So, the typical rank depends on the field $\K$.

\nT{Eigenspaces of commuting matrices} Let us have two commuting matrices $A,B \in \K^{n \times n}$, i.e.
\[
 AB=BA
\]
Let us assume that each matrix is diagonalisable, with no multiplicity of eigenvalues to keep things technically plain (i.e. the dimension of each eigenspace is one). Let $(\lambda,u)$ be such that
\begin{equation}
 A\ubold=\lambda \ubold
\end{equation}
Then
\begin{equation}
 \begin{array}{lcl}
   A(B\ubold) &=& B(A\ubold) \\
   &=& \lambda B\ubold
 \end{array}
\end{equation}
As the dimension of eigenspace of $A$ associated to $(\lambda,\ubold)$ is one, one has
\begin{equation}
 B\ubold = \mu \ubold
\end{equation}
and $u$ is an eigenvector of $B$ as well: $A$ and $B$ have the same eigenspaces. Let us note that if $AB=BA$ and $A$ is invertible, then $B=A^{-1}BA$ which is the expression of $B$ in another basis. Then, if $B$ is diagonalisable, $A$ is diagonalisable, and both can be expressed in the basis of eigenvectors.

\nT{Simultaneously diagonalisable matrices} A set of matrices $(A_k)_k$ is simultaneously diagonalisable if
\[
 \exists \: P \: : \: \forall \: k, \quad P^{-1}A_kP \quad \mbox{diagonal}
\]
This means that there exists a basis on which all matrices are diagonal (hence the name). A set of matrix is simultaneously diagonalisable if, and only if, $(i)$ each matrix is diagonalisable and $(ii)$ all matrices in the set commute pairwise. Then, they have the same eigenspaces: 
\begin{equation}
 \exists \: (\ubold_i)_i \st A_k\ubold_i= \lambda_{k,i}\ubold_i \quad \forall \: k
\end{equation}


\nS Let $\Abold \in \R^n \otimes \R^n \otimes \R^m$ with $m \leq n$. Let $(\ebold_1,\ldots, \ebold_m)$ be a basis of $\R^m$, and 
\begin{equation}
 \Abold = \sum_{k=1}^m A_k \otimes \ebold_k, \qquad A_k \in \K^{n \times n}
\end{equation}
with
\begin{equation}
    A_k = \Abold.\ebold_k
\end{equation}
The expected rank is $\rank \Abold = nm$. \cite{Lafon1975} has shown that, if the set $(A_1, \ldots,A_m)$ of square $n \times n$ matrices is simultaneously diagonalisable, then $\rank \Abold = n$. 
\begin{proof}
Let 
\[
 \Abold = \sum_{k=1}^m A_k \otimes \ebold_k 
\]
and 
\[
 \exists \: P \in \K^{n \times n} \quad : \quad \forall \: k, \quad P^{-1}A_kP = \Lambda_k, \quad \Lambda_k \: \mbox{diagonal}
\]
Then, there is a set $(\ubold_i)_i$ of eigenvectors such that 
\[
 \forall \: k, \quad \Lambda_k = \sum_{i=1}^n \lambda_{ik} \, \ubold_i \otimes \ubold_i^*
\]
where $(\ubold_i)_i$ is the basis of $\K^n$ on which each $A_k$ is diagonal. We then have
\[
\begin{array}{lcl}
  A_k &=& P\Lambda_kP^{-1} \\
  &=& \displaystyle P \left(\sum_{i=1}^n \lambda_{ik} \, \ubold_i \otimes \ubold_i^*\right)P^{-1} \\
  &=& \displaystyle \sum_{i=1}^n \lambda_{ik} \, P (\ubold_i \otimes \ubold_i^*)P^{-1} \\
  &=& \displaystyle \sum_{i=1}^n \lambda_{ik} \, P\ubold_i \otimes (P^{-1})^\t \ubold_i^*\\
\end{array}
\]
and
\[
 \begin{array}{lcl}
   \Abold &=& \displaystyle \sum_{k=1}^m A_k \otimes \ebold_k \\
   &=& \displaystyle \sum_{k=1}^m \left( \sum_{i=1}^n \lambda_{ik} \, P\ubold_i \otimes (P^{-1})^\t \ubold_i^*\right) \otimes \ebold_k \\
   &=& \displaystyle \sum_{i=1}^n \left(\sum_{k=1}^m \lambda_{ik} \ebold_k\right) \otimes  P\ubold_i \otimes (P^{-1})^\t \ubold_i^*
 \end{array}
\]
which is a decomposition as a sum of $n$ elementary tensors.
\end{proof}

\nS By (\ref{eq:rank3:bound:1}), typical rank in $\K^{m \times n \times n}$ must be equal to or larger than 
\[
  \mathbf{t}_{m,n} = \frac{mn^2}{m+2n-2}
\]
Do we have $n \in \mathbf{tr}(n,n,m)$? in other words: is $\rbold(n,n,m)=n$ typical? This requires that
\begin{equation}
    \frac{mn^2}{m+2n-2} \leq n
\end{equation}
or
\begin{equation}
    mn \leq m+2n-2
\end{equation}
or
\begin{equation}
    mn-m \leq 2n-2 \quad \Leftrightarrow \quad m(n-1) \leq 2(n-1) \quad \Leftrightarrow \quad m \leq 2
\end{equation}
So, $\rbold(n,n,m)=n$ is typical for $m=2$ only. This shows as well that a set of $m$ simultaneously diagonalisable (or pairwise commuting) matrices is not typical for $m >2$.

%
\section{Typical rank in $\R^2 \otimes \R^2 \otimes \R^2$}\label{sec:rank3:rank222}
%

We show in this section that 
\begin{equation}
 \rfrak(2,2,2)=3
\end{equation}

\nS We first show a preliminary, but useful, result. Let $A \in \R^{n \times n}$ be a matrix of rank $n$, and $\xbold,\ybold \in \R^n$ two vectors. Let $\alpha \in \R$ be a scalar, and define $A' \in \R^{n \times n}$ by 
\begin{equation}
 A = \alpha \xbold\otimes \ybold + A'
\end{equation}
Then, if $\langle A^{-1}\xbold,\ybold\rangle \neq 0$, it is possible to select $\alpha$ such that $\rank A'=n-1$.
\begin{proof}
$A,\xbold,\ybold$ are fixed. As $A$ is full rank, it is invertible, and $A^{-1}$ exists. Let us define $\ubold = A^{-1}\xbold$. Then $\xbold=A\ubold$ and
\[
\begin{array}{ll}
  \xbold &= \alpha \langle \ybold,\ubold\rangle \xbold + A'\ubold\\
  &= \alpha \langle \ybold,A^{-1}\xbold\rangle \xbold + A'\ubold
\end{array}
\]
If $\langle A^{-1}\xbold,\ybold\rangle \neq 0$, let us select $\alpha = 1/\langle A^{-1}\xbold,\ybold\rangle$. Then, $\xbold=\xbold+A'\ubold$, and $A'\ubold=0$. As $\ubold \neq 0$, $\rank A' \leq n-1$. As
\[
 A = \frac{1}{\langle A^{-1}\xbold,\ybold\rangle}\, \xbold \otimes \ybold + A'
\]
and $\rank A=n$, we have $\rank A'=n-1$.
\end{proof}


\nT{Tensor in $\R^{2 \times 2 \times 2}$}Let $\Abold \in \R^{2 \times 2 \times 2}$. If $(\ebold_1,\ebold_2)$ is an orthonormal basis of $\R^2$, one can write
\begin{equation}
 \Abold = A_1 \otimes \ebold_1 + A_2 \otimes \ebold_2, \qquad \mbox{with} \quad
 \begin{cases}
  A_1 &= \Abold.\ebold_1 \\
  A_2 &= \Abold.\ebold_2 
 \end{cases}
\end{equation}
We have 
\[
 \rank \Abold \leq \rank A_1 + \rank A_2
\]
Let us assume we are in the "worst" case $\rank A_1 = \rank A_2=2$. Then, $\rank \Abold \leq 4$. We show that $\rank \Abold \leq 3$.
\begin{proof}
Let us assume that both $A_1$ and $A_2$ are full rank, i.e. each invertible. Let us select $\xbold,\ybold \in \R^2$ such that $\langle A_1^{-1}\xbold,\ybold\rangle \neq 0$ and $\langle A_2^{-1}\xbold,\ybold\rangle \neq 0$. Such a pair exists (select first $\xbold$ and second $\ybold$ for the condition to be fulfilled). We then have, with preliminary result,
\[
 \left\{
   \begin{array}{lcl}
    A_1 &=& \alpha_1 \xbold\otimes \ybold + \xbold_1 \otimes \ybold_1 \\
    A_2 &=& \alpha_2 \xbold\otimes \ybold + \xbold_2 \otimes \ybold_2 \\
   \end{array}
 \right.
\]
Then
\[
 \left\{
   \begin{array}{lcl}
     \Abold &=& A_1 \otimes \ebold_1 + A_2 \otimes \ebold_2 \\
     &=& (\alpha_1 \, \xbold\otimes \ybold + \xbold_1 \otimes \ybold_1 )  \otimes \ebold_1 + (\alpha_2 \, \xbold\otimes \ybold + \xbold_2 \otimes \ybold_2)  \otimes \ebold_2 \\
     &=& \xbold \otimes \ybold \otimes (\alpha_1 \, \ebold_1 + \alpha_2 \, \ebold_2) + \xbold_1 \otimes \ybold_1 \otimes \ebold_1 + \xbold_2 \otimes \ybold_2 \otimes \ebold_2
   \end{array}
 \right.
\]
and $\rank \Abold \leq 3$. If one of the matrices $A_1$ or $A_2$ is not full rank, it has rank one, and the result follows. 
\end{proof}

\nT{Some notations} Let us have
\begin{equation}
 \Abold = \sum_{i,j,k=1}^2 \, \alpha_{ijk} \, \ebold_i \otimes \ebold_j \otimes \ebold_k
\end{equation}
with $(\ebold_1,\ebold_2)$ being the standard basis in $\R^2$
\begin{equation}
 \ebold_1 = 
 \begin{pmatrix}
  1 \\
  0
 \end{pmatrix}
 , \qquad \ebold_2 = 
 \begin{pmatrix}
  0 \\
  1
 \end{pmatrix}
\end{equation}
We have $\Abold = A_1.\ebold_1 + A_2.\ebold_2$ with
\begin{equation}
 \left\{
   \begin{array}{lcl}
     A_1 &=& \alpha_{111}\, \ebold_1 \otimes \ebold_1 + \alpha_{121}\, \ebold_1 \otimes \ebold_2 + \alpha_{211}\, \ebold_2 \otimes \ebold_1 + \alpha_{221}\, \ebold_2 \otimes \ebold_2 \\
     A_2 &=& \alpha_{112}\, \ebold_1 \otimes \ebold_1 + \alpha_{122}\, \ebold_1 \otimes \ebold_2 + \alpha_{212}\, \ebold_2 \otimes \ebold_1 + \alpha_{222}\, \ebold_2 \otimes \ebold_2 
   \end{array}
 \right.
\end{equation}
or
\begin{equation}
 A_1 = 
 \begin{pmatrix}
   \alpha_{111} & \alpha_{121} \\ 
   \alpha_{211} & \alpha_{221} \\ 
 \end{pmatrix}
 , \qquad 
  A_2 = 
 \begin{pmatrix}
   \alpha_{112} & \alpha_{122} \\ 
   \alpha_{212} & \alpha_{222} \\ 
 \end{pmatrix}
\end{equation}
which can be written
\begin{equation}
 \Abold = 
 \left(
 \begin{array}{rr}
    \alpha_{111} & \alpha_{121} \\ 
   \alpha_{211} & \alpha_{221} \\ 
\end{array}
\right.
\left|
\begin{array}{rr}
   \alpha_{112} & \alpha_{122} \\ 
   \alpha_{212} & \alpha_{222} \\ 
 \end{array}
 \right)
\end{equation}
In order to facilitate reading of calculations below, we change notations into
\begin{equation}
 \Abold = 
 \left(
 \begin{array}{rr}
    a & b \\ 
    c & d \\ 
\end{array}
\right.
\left|
\begin{array}{ll}
   a' & b' \\ 
   c' & d' \\ 
 \end{array}
 \right)
\end{equation}

\nT{Cayley hyperdeterminant} We know that $\rank \Abold=2$ if $A_1^{-1}A_2$ is diagonalizable (see section \ref{sec:rank3:bound}, page \pageref{sec:rank:rank3:bound:ref1}). We have
\begin{equation}
 A_1^{-1} = 
 \frac{1}{ad-bc}
 \begin{pmatrix}
  d & -b \\
  -c & a
 \end{pmatrix}
 \qquad \Longrightarrow \qquad (ad-bc)A_1^{-1}A_2 = 
 \begin{pmatrix}
   a'd - bc' & b'd-bd' \\
   ac'-a'c & ad'-b'c
 \end{pmatrix}
\end{equation}
Let
\[
 M = \begin{pmatrix}
   a'd - bc' & b'd-bd' \\
   ac'-a'c & ad'-b'c
 \end{pmatrix}
\]
Then, $A_1^{-1}A_2$ is diagonalizable iff $M$ has two real eigenvalues. The eigenvalues are solution of  
\begin{equation}
  \lambda^2 - \Tr M \; \lambda + \det M = 0
\end{equation}
which has two real roots if 
\begin{equation}
 \Delta = (\Tr M)^2 - 4 \, \det M > 0
\end{equation}
$\Delta$ is called the \kw{Cayley hyperdeterminant} of $\Abold$. Let us compute it\footnote{computing $\Delta$ is the reason for adapting notations. We order letters within monomial terms as $a<a'<b<b'<c<c'<d<d'$ to facilitate cancellations.}. We have
\begin{equation}
 \Tr M = ad' + a'd -bc' -b'c
\end{equation}
so
\begin{equation}
 \begin{array}{lcl}
  (\Tr M)^2 &=& (ad' + a'd -bc' -b'c)^2 \\
  &=& a^2d^{'2} + a^{'2}d^2 + b^2c^{'2}  + b^{'2}c^2 \\
  && \qquad + 2aa'dd' - 2abc'd' - 2ab'cd' - 2 a'bc'd - 2 a'b'cd + 2 bb'cc'  
 \end{array}
\end{equation}
and
\begin{equation}
 \begin{array}{lcl}
  \det M &=& (a'd-bc')(ad'-cb') - (ac'-a'c)(b'd-bd') \\
  &=& aa'dd' - a'b'cd - abc'd' + bb'cc' - ab'c'd + abc'd' + a'b'cd - a'bcd' \\
  &=& aa'dd' + bb'cc' - ab'c'd  - a'bcd' \\
 \end{array}
\end{equation}
So, by reordering the monomials
\begin{equation}
 \begin{array}{lcl}
  \Delta &=& (\Tr M)^2 - 4 \; \det M \\
  &=& a^2d^{'2} + a^{'2}d^2 + b^2c^{'2}  + b^{'2}c^2 \\
  && \quad + 2aa'dd' - 2abc'd' - 2ab'cd' - 2 a'bc'd - 2 a'b'cd + 2 bb'cc' \\
  && \qquad -4\, aa'dd' - 4\,bb'cc' + 4\,ab'c'd  + 4\,a'bcd' \\
  &=&  a^2d^{'2} + a^{'2}d^2 + b^2c^{'2}  + b^{'2}c^2 \\
  && \quad - 2abc'd' - 2ab'cd' - 2 a'bc'd - 2 a'b'cd \\
  && \qquad 4\,ab'c'd  + 4\,a'bcd' \\
  &=& a^2d^{'2} + a^{'2}d^2 + b^2c^{'2}  + b^{'2}c^2 \\
  && \quad - 2( abc'd' + ab'cd' + a'bc'd + a'b'cd ) \\
  && \qquad +4(ab'c'd  + a'bcd')
 \end{array}
\end{equation}

\nS We then have the following result: 
\begin{equation}
\rank \Abold = 
 \begin{cases}
  2 &\mbox{if} \quad \Delta > 0 \\
  3 &\mbox{if} \quad \Delta < 0 \\
 \end{cases}
\end{equation}

\nS\label{sec:rank:rank3:rank222:page1} $\Delta$ is known as Cayley hyperdeterminant or Kruskal's polynomial. It is an invariant of $\Abold$, i.e. it is invariant as a polynomial of the coefficients by  change of basis.

\notes The study of typical rank of a tensor has probably begun with seminal paper of Strassen \cite{Strassen1969} who showed that seven multiplications are sufficient, and not eight, for computing the product of two $2 \times 2$ matrices. The map $(A,B) \longmapsto AB$ being bilinear, it can be described by a tensor, the rank of which is seven, and not eight. This has triggered a series of studies on the typical rank of small tensors. Gastinel has shown that the product of two $n \times n$ matrices can be done with $n^3-n+1$ multiplications only \cite{Gastinel1971}. Lafon has shown that this figure is optimal for $2 \times 2$ matrices \cite{Lafon1975}. We denote here by $r(m,n,p)$ the maximum of the typical rank of a $m\times n \times p$ tensor, generally over $\C$, which has been much studied \cite{Lafon1975,Howell1978,Atkinson1979,JaJa1979,Atkinson1980}. Howell has studied the maximum typical rank of a tensor according to the field $\K$ over which it has been built, and even modules. He has shown that $r(n,n,n) \leq 3n^2/4$ on any ring, which yields $3$ for $2 \times 2 \times 2$ tensors \cite{Howell1978} (the rank is 2 over $\C$ and 3 over $\R$ or $\Q$). \cite{Atkinson1979} have improved Howell's bound over $\C$ and shown that $r(n,n,n) \leq n^2/2 + O(n)$. \cite{Strassen1969,Howell1978} have used Gauss elimination algorithm to bound the rank of a tensor. Such an approach is developed in \cite{Franc1992}, section F. See as well \cite{Comon2009b} for a method to compute bonds of typical rank. The typical rank of a tensor $\Abold \in \R^2 \otimes \R^2 \otimes \R^2$ has been studied, e.g. in \cite{Kruskal1989a} and in \cite{Berge1991}. $\R^2 \otimes \R^2 \otimes \R^2$ has dimension $2^3=8$. \cite{Kruskal1989a} give the following results, but without demonstrations: $(i)$Typical rank is $\{2,3\}$ ; $(ii)$ the boundary between $\Gamma_2$ and $\Gamma_3$ is a $7-$dimensional manifold in $\R^8$ ; $(iii)$ the boundary $\overline{\Gamma_2} \cap \overline{\Gamma_3}$ includes tensors of rank 0, 1, 2 and 3. This means that neither $\Gamma_2$ nor $\Gamma_3$ are closed sets in $\R^8$. This has been developed further in \cite{Berge1991} and \cite{Silva2008}.

%
\section{Tucker rank}\label{sec:rank3:tucker}
%

Let $\Abold \in E \otimes F \otimes G$, with $\dim E=m, \: \dim F=n, \: \dim G=p$. Let $E' \subset E, \: F' \subset F, \: G' \subset G$ with $\dim E'=m', \: \dim F'=n', \: \dim G'=p'$. We have
\begin{equation}
 E' \otimes F' \otimes G' \subset E \otimes F \otimes G
\end{equation}

\nS Let $\rbold=(m',n',p')$. We denote by $\T(\rbold)$ the set of tensors $\Abold \in E \otimes F \otimes G$ such that there exists minimum subspaces $E' \subset E, \: F' \subset F, \: G' \subset G$ with $\dim E'=m', \: \dim F'=n', \: \dim G'=p'$ and $\Abold \in E' \otimes F' \otimes G'$:
\begin{equation}
 \T(\rbold) = \left\{ \Abold \in E \otimes F \otimes G \quad \mbox{s.t.} \quad 
 \left|
  \begin{array}{l}
  \exists \quad E' \subset E, \: F' \subset F, \: G' \subset G \\
  \mbox{with (minimum dimensions)} \\
  \dim E'=m', \: \dim F'=n', \: \dim G'=p' \\
  \mbox{and} \\
 \Abold \in E' \otimes F' \otimes G' 
 \end{array}
 \right.
 \right\}
\end{equation}
Then,
\[
\rbold = (m',n',p') 
\]
is the \kw{Tucker rank}\index{rank!Tucker} of $\Abold$. One is tempted to define the Tucker rank as the smallest rank $\rbold=(m',n',p')$ such that such a condition holds. However, the set of possible ranks is a poset (partially ordered set), and not a totally ordered set. So, one must define what is meant by "smallest". One possibility is to set that no dimension can be lowered while the two others do not change or, for $E$, given $(m', n', p')$, 
\begin{equation}
 \nexists \: \: E'' \:\: \mbox{with} \:\: \dim E'' = m'' < m' \st \Abold \in E'' \otimes F' \otimes G' 
\end{equation}
The Tucker rank of tensor $\Abold$ is denoted $\mathbf{mr}(\Abold)$. It is a triplet of integers:
\begin{equation}
 \mathbf{mr}(\Abold) = (m',n',p')
\end{equation}

\notes The notion of Tucker rank is presented in \cite[chapter~8]{Hackbusch2012}. It has been introduced by Tucker and, according to \cite[p.~219]{Hackbusch2012}, by Hitchcock as well\footnote{Hitchcock, F.L., 1927 : The expression of a tensor or a polyadic as a sum of products. \emph{Journal of Mathematics and Physics}, \textbf{6:}164-189}. It is called Tucker rank, \kw{multilinear rank}\index{rank!multilinear}, or \kw{tensor subspace rank}\index{rank!tensor subspace} (see \cite{Hackbusch2019}, remark 8.4). In their classical survey paper, \cite{Kolda2009} define the Tucker rank as the dimensions of the decomposition $\Xbold = \sum_{ijk} g_{ijk}\abold_i \otimes \bbold_j \otimes \cbold_k$ and consider it is a higher order PCA. The more geometric approach is presented in \cite{Hackbusch2012,Hackbusch2019}, and has been fully developed in \cite{Franc1992}.

%
\section{Tensor ranks of $d-$modes tensors}\label{sec:rank3:rankd}
%

Here, the definitions of rank given for $3-$modes tensors given in previous sections are extended in a straightforward way to $d-$modes tensors. It is merely a matter of clear and simple notations.

\nS \textbf{Rank:} Let
\[
 \Abold \in \bigotimes_{\mu=1}^d \, E_\mu, \qquad E_\mu \simeq \K^{n_\mu}
\]
The \kw{rank} of $\Abold$, denoted $\rank \Abold$, is the smallest integer $r$ such that $\Abold$ is a sum of $r$ elementary tensors (rank one tensors). Then, there exists $r$ $d-$uplets $\left(\xbold_1^{(a)},\ldots,\xbold_d^{(a)}\right)$ with $\xbold_\mu^{(a)} \in E_\mu$ such that
\begin{equation}
 \begin{array}{lcl}
   \Abold &=& \displaystyle \sum_{a=1}^r \left(\bigotimes_\mu \xbold_\mu^{(a)}\right) \\
   &=& \displaystyle \sum_{a=1}^r \, \xbold_\mu^{(1)} \otimes \ldots \otimes x_a^{(d)}
 \end{array}
\end{equation}
Let us denote by $X_\mu$ the $n_\mu \times r$ matrix the column $a$ of which is $\xbold_a$
\begin{equation}
    X_\mu = \sum_{a=1}^r \xbold_\mu^{(a)} \otimes \ebold_\mu^{(a)}
\end{equation}
Then
\begin{equation}
    \Abold = \langle X_1 \mid \ldots \mid X_d\rangle
\end{equation}
(see equation \ref{eq:sec:alstruct:prodmat:1}). Let
\[
 \nbold = (n_1, \ldots,n_d), \qquad \mbox{with} \quad n_\mu = \dim E_\mu
\]
We denote by
\[
\mathbb{CP}_r\left(\bigotimes_\mu E_\mu\right)
\]
or $\mathbb{CP}_r(\nbold)$ or $\mathbb{CP}_r$ when non ambiguous the set of tensors of rank $r$
\begin{equation}
 \Gamma_r = \left\{\Abold \in \bigotimes_\mu E_\mu\: : \: \rank \Abold = r\right\}
\end{equation}
The rank is often called the \kw{CP-rank}\index{rank!CP} as well.

\nT{Example} Let $\Abold$ be a tensor such that each element $a[i_1,\ldots,a_d]$ is a sum of functions each of one index, like, for $d=4$
\[
a_{ijk} = i^2 + j^2 + k^2 + \ell^2
\]
This can be written in the general case
\[
 \forall \: \ibold, \quad a_{i_1 \ldots i_d} = \sum_{\mu=1}^d x_\mu(i_\mu)
\]
with
\[
 \begin{CD}
  \llbracket 1, n_\mu \rrbracket @>x_\mu>> \K
 \end{CD}
\]
For example, if $d=4$, 
\begin{equation*}
 \Abold(i_1,i_2,i_3,i_4) = x_1(i_1) + x_2(i_2)+x_3(i_3)+x_4(i_4)
\end{equation*}
Let us denote by $\ones_{n_\mu}$ the vector in $E_\mu$ with ones only (we have $\dim E_\mu = n_\mu$), and
\[
\xbold_\mu =
\begin{pmatrix}
  x_\mu(i_\mu=1) \\
  \vdots \\
  x_\mu(i_\mu = n_\mu)
\end{pmatrix}
\]
(this is an awkward notation, but I have not found a better one ...). Let us denote
\[
 \Xbold_\mu = \underbrace{\ones_{n_1} \otimes \ldots \otimes \ones_{n_{\mu-1}}}_{a < \mu} \otimes   \underbrace{\xbold_\mu}_{a = \mu}  \otimes \underbrace{\ones_{n_{\mu+1}} \otimes \ldots \otimes \ones_{n_d}}_{a > \mu} 
\]
For a simple example, like for $d=2$ and $\dim E=\dim F=n$, let
\begin{equation}
    x_1(i)=i^2, \qquad x_2(j)=j^2
\end{equation}
Then
\begin{equation}
    \xbold_1 = 
    \begin{pmatrix}
      1 \\
      \vdots \\
      i^2 \\
      \vdots \\
      n^2
    \end{pmatrix}
    , \qquad 
    \xbold_2 = 
    \begin{pmatrix}
      1 \\
      \vdots \\
      j^2 \\
      \vdots \\
      n^2
    \end{pmatrix}
\end{equation}
Then 
\begin{equation}
X_1 = \xbold_1 \otimes \ones =
\begin{pmatrix}
 1 & \hdots & 1 \\
 \hdots & \hdots & \hdots \\
 i^2 & \hdots & i^2 \\
 \hdots & \hdots & \hdots \\
 n^2 & \hdots & n^2 \\
\end{pmatrix}
, \qquad
X_2 = \ones \otimes \xbold_2 =
\begin{pmatrix}
 1 & \hdots & j^2 & \hdots & n^2 \\
 \vdots & \vdots & \vdots & \vdots & \vdots \\
 1 & \hdots & j^2 & \hdots & n^2
\end{pmatrix}
\end{equation}
and $X_1+X_2$ is the $n \times n$ matrix of general term $i^2+j^2$. 
We can check that in general
\begin{equation}
 \Abold = \sum_{\mu=1}^d\Xbold_\mu
\end{equation}
and, because $\Xbold_\mu$ is an elementary tensor for each $\mu$ , $\rank \Abold \leq d$.

\nT{Note} Let us note that in such a case, the rank of $\Abold$ does not depend on the dimensions of the tensor $\Abold$, but its order only. This will be developed in chapter \ref{chap:quant} where it will be shown that tensors which are the \kw{discretization} of a smooth multivariate treal function are of low rank, even if the discretization is done with a very fine grain, i.e. with high dimensions on each mode. 

\nS Let 
\[
 \nbold = (n_1, \ldots, n_d)
\]
We denote by $\rr(\nbold)$ the smallest integer $r$ such that
\begin{equation}
 \forall \: \Abold \in \bigotimes_\mu \, \K^{n_\mu}, \quad \rank \Abold \leq r
\end{equation}

\nT{Border rank}  Let $\Abold \in \bigotimes_\mu E_\mu$. The \kw{border rank}\index{rank!border} of $\Abold$ ($\brank \Abold$) is the smallest integer $r$ such that there exists a series
\[
 (\Abold_k)_k, \qquad k \in \N, \qquad \Abold_k \in \bigotimes_\mu E_\mu
\]
such that 
\begin{equation}
\left\{
 \begin{array}{l}
   \Abold = \displaystyle \lim_{k \rightarrow \infty} \: \Abold_k \\
   \\
   \mbox{and} \\
   \\
   \forall \: k \in \N, \quad \rank \Abold_k = r
 \end{array}
\right.
\end{equation}

\nT{Typical rank} Let $\K = \R$ or $\C$ (it has to be a normed space). Let $\S^{n_1 \otimes \ldots \otimes n_d}$ be the unit sphere in $\bigotimes_\mu\K^{n_\mu}$
\[
 \S^{n_1 \times \ldots \times n_d} = \left\{\Abold \in \bigotimes_{\mu}\K^{n_\mu} \: : \: \|\Abold\|=1 \right\}
\]
which can be written
\[
 \begin{array}{ccl}
  \displaystyle \sum_\ibold \, \alpha_{\ibold}^2=1 & & \mbox{if } \K=\R \\
  \displaystyle \sum_\ibold \, |\alpha_{\ibold}|^2=1 & & \mbox{if } \K=\C
 \end{array}
\]
Then, the \kw{typical rank}\index{rank!typical} in $\bigotimes_\mu\K^{n_\mu}$ is the the set of values $r$ such that the measure of the set of tensors $\Abold \in \S^{n_1 \times \ldots \times n_d}$ of rank $r$ is strictly positive
\begin{equation}
 \mu \left\{\Abold \in \S^{n_1 \times \ldots \times n_d} \; : \; \rank \Abold = r \right\}  > 0
\end{equation}
or
\begin{equation}
 \mu \left( \mathbb{CP}_r \cap \S\right) > 0
\end{equation}
It is denoted $\mathbf{tr}(\nbold)$:
\[
 \mathbf{tr}(m,n,p) = \left\{ r \in \N \: : \:  \mu \left( \mathbb{CP}_r \cap \S\right) > 0\right\}
\]
Computing $\mathbf{tr}(\nbold)$ is still an open (and difficult) problem. Nearly nothing is known for $d>3$.

\nS We have, as for $d=3$:\\
\\
\begin{center}
 \begin{tabular}{l|cl|c}
 \hline 
 Name & & Definition \\
 \hline
  rank & $\rbold(\Abold)$ & minimum integer $r$ such that & $\mathbb{CP}_r $\\
  & & $\Abold$ is written as a sum of $r$ elementary tensors  \\
  &&&\\
  border rank &  $\mathbf{br}(\Abold)$ & minimum integer $r$ such that & $\mathbb{BR}_r$ \\
  & & $\Abold$ is the limit of a sequence of tensors of rank $r$\\
  &&&\\
  typical rank & $\mathbf{tr}(n_1,\ldots,n_d)$ & the set of integers $r$ such that & \\
  & & the measure of tensors of rank $r$ and norm 1 is strictly positive\\
  \hline
 \end{tabular}
\end{center}


\nT{Tucker rank} Let $d \in \N$, $\mu \in \llbracket 1,d \rrbracket$, $(E_\mu)_\mu$ a series of vector spaces over $\K$ with $\dim E_\mu=n_\mu$, and $\rbold=(n'_1, \ldots,n'_d)$ with $n'_\mu \leq n_\mu$. We denote by $\T(\rbold)$ the set of tensors $\Abold \in E \otimes \ldots \otimes E_d$ such that there exists minimum subspaces $E'_1 \subset E_1, \: \ldots , \: E'_d \subset E_d$ with $\dim E_\mu = n_\mu$ and $\Abold \in E'_1 \otimes \ldots \otimes E'_d$:
\begin{equation}
 \T(\rbold) = \left\{ \Abold \in E_1 \otimes \ldots  \otimes E_d \quad \mbox{s.t.} \quad 
 \left|
  \begin{array}{l}
  \exists \quad E'_1 \subset E_1, \: \ldots, \: E'_d \subset E_d \\
  \mbox{with (minimum dimensions)} \\
  \dim E'_1=n'_1, \: \ldots, \: \dim E'_d=n'_d \\
  \mbox{and} \\
 \Abold \in E'_1 \otimes \ldots \otimes E'_d 
 \end{array}
 \right.
 \right\}
\end{equation}
Then,
\[
\rbold = (n'_1, \ldots,n'_d) 
\]
is the \kw{Tucker rank}\index{rank!Tucker} of $\Abold$. The Tucker rank of tensor $\Abold$ is denoted $\mathbf{Tr}(\Abold)$. It is a $d-$uplet of integers:
\begin{equation}
 \mathbf{mr}(\Abold) = (n'_1,\ldots,n'_d)
\end{equation}
It is called as well the multilinear rank, or tensor subspace rank.

%% file: intermezzo.tex
%
\chapter{Intermezzo}\label{chap:inter}
%

Here, we deal with matrices, not tensors. So, it is an \emph{intermezzo} in the story about tensors. 

\nB Let $A \in E \otimes F \simeq \L(F,E)$. Its rank is the dimension of the target space, here included in $E$. So, it is $\dim E$ at most. In general, it is the smallest dimension among the dimensions of $E$ and $F$. This is translated in matrix calculus as the existence of a decomposition
\[
A = UW, \qquad \mbox{with} \quad 
\begin{cases}
  A & \in \R^{m \times n} \\
  U & \in \R^{m \times r} \\
  W & \in \R^{r \times n} 
\end{cases}
\]
where $r = \rank A = \dim \Im A$. Such a decomposition requires $(m+n)r$ terms instead of $mn$. For a given $r$, the storage complexity of the decomposition is in $\O(m+n)$ instead of in $\O(mn)$ for the full matrix. As we will see later, this has far reaching consequences for key operations like computing the partition function in statistical physics or the mode of a distribution in statistics. 

\nB A matrix $A$ being given, it is useful and efficient in many situations to find a matrix $A_r$ of same dimension and of rank $r$ as close as possible to $A$, what is called a best rank$-r$ approximation. This requires to have a way to evaluate a distance between matrices, which leads naturally to norms, as the space of matrices is a vector space. One will see that several norms are possible, and that the difficulty of these tasks varies dramatically according to the selected norm. It is well understood for Euclidean or Frobenius norm, and still a challenging problem for $\ell^1$ or $\ell^\infty$ norms which are met naturally in applications ($\ell^1$ for computing the partition function, $\ell^\infty$ for computing the mode of a distribution).

\nB Best low rank approximation of a matrix, or of a linear map in a Hilbert space, can be done with so called SVD (Singular Value Decomposition), for Frobenius norm. Extension of SVD to tensors has been a motivating challenge for decades, and has led to splitting it into different notions of ranks, as seen in chapter \ref{chap:rank3} and as a consequence to splitting of algorithms for best low rank approximation, which will be presented in part \ref{part:blra}.

\nB This chapter is organized as follows:
\begin{description}
\item[section \ref{sec:inter:norm}] The basic notions of normed vector spaces are presented, with an emphasis on links with distances and metric properties
\item[section \ref{sec:inter:equiv}] The equivalence between norms is presented, which will be useful when the dimension of vector spaces increases
\item[section \ref{sec:inter:some}] Some important norms are presented ($\ell^1, \ell^2, \ell^\infty$)
\item[section \ref{sec:inter:key4tens}] Some key consequences of the choice of a norm for tensors are presented
\item[section \ref{sec:inter:notes}] In this key section, some challenging questions on the reason why some matrices are low rank and some other not are presented, and this is still an uncharted continent for tensors.
\item[section \ref{sec:inter:svd}] Finally, SVD for matrices is presented in last section.
\end{description}

%
\section{Norm on a vector space}\label{sec:inter:norm}
%

Here, I recall how a norm leads to a distance. 

\nS A norm on a vector space $E$ is a map 
\begin{equation*}
 \begin{CD}
  E @>\|.\|>> \R^+
 \end{CD}
\end{equation*}
such that
\begin{itemize}[label=$\rightarrow$]
 \item $\|\xbold\|= 0 \quad \Longleftrightarrow \quad \xbold=0$
 \item $\|\alpha \xbold\| = |\alpha|\|\xbold\|$
 \item $\|\xbold+\ybold\| \leq \|\xbold\|+\|\ybold\|$
\end{itemize}
The pair $(E,\|.\|)$ is a normed vector space.

\nS If $\|.\|$ is a norm on $E$, the map 
\begin{equation*}
 \begin{CD}
  E \times E @>d>> \R^+ \\
  (\xbold,\ybold) @>>> \|\xbold-\ybold\|
 \end{CD}
\end{equation*}
is a distance.

\nS The open ball $\B(\xbold,r)$ for $\xbold \in E$ and $r>0$ is defined as
\begin{equation}
 \B(\xbold,r) = \{\ybold \in E \mid \|\ybold-\xbold\| < r\}
\end{equation}
and the closed ball $\overline{\B}(\xbold,r)$ as
\begin{equation}
 \overline{\B}(\xbold,r) = \{\ybold \in E \mid \|\ybold-\xbold\| \leq r\}
\end{equation}
Let us note that
\begin{equation}
 \B(\xbold,r)=\xbold + \B(0,r)
\end{equation}
which is easily seen from $\ybold = \xbold + \ybold - \xbold$.

\nS If $E$ is real and is endowed with an inner product $\langle .,.\rangle$, or complex and endowed with an Hermitian form $\langle . , . \rangle$, the map 
\begin{equation}
 \begin{CD}
  E @>>> \R^+ \\
  \xbold @>>> \sqrt{\langle \xbold,\xbold\rangle}
 \end{CD}
\end{equation}
is a norm.

\nS A norm on a vector space defines a topology on this space. The base of open sets for this topology is the set of open balls of radius $r$, i.e. the sets
\begin{equation}
 \B(\xbold,r) = \{\ybold \in E \mid d(\xbold,\ybold) < r\} = \{\ybold \in E \mid \|\ybold-\xbold\| < r\}
\end{equation}
for $\xbold \in E$ and $r >0$.

\nS We then have the following constructions:
\[
 \begin{CD}
  \mbox{inner product} @>>> \mbox{norm} @>>> \mbox{distance} @>>> \mbox{topology}
 \end{CD}
\]
The topology is useful for studying convergence of series, especially in infinite dimensional vector spaces. This will not be developed here.

%
\section{Equivalence between norms}\label{sec:inter:equiv}
%

Two norms $\|.\|_a$ and $\|.\|_b$ on a same vector space $E$ are said equivalent if
\begin{equation}
  \exists \: c,c' \in \R^{*+} \quad \st \quad \forall \: \xbold \in E, \quad c\|\xbold\|_a \leq \|\xbold\|_b \leq c'\|\xbold\|_a
\end{equation}

\nS Let us show first that it is an equivalence relation.
\begin{proof}Indeed\\
\begin{itemize}[label=$\rightarrow$]
\item $\|.\|_a \equiv \|.\|_a$ simply by selecting $c=c'=1$. 
\item if $\|.\|_a \equiv \|.\|_b$, then $\|.\|_b \equiv \|.\|_a$. Indeed, $\forall \: \xbold \in E$
\[
 \frac{1}{c'}\|\xbold\|_b \leq \|\xbold\|_a \leq \frac{1}{c}\|\xbold_b\|
 \]
 \item if $\|.\|_a \equiv \|.\|_b$ and $\|.\|_b \equiv \|.\|_c$, then $\|.\|_a \equiv \|.\|_c$. Indeed, we have  $c\|\xbold\|_a \leq \|\xbold\|_b \leq c'\|\xbold\|_a$ and $k\|\xbold\|_b \leq \|\xbold\|_c \leq k'\|\xbold\|_b$. Then, $ck\|\xbold\|_a \leq k\|\xbold\|_b \leq \|\xbold\|_c \leq k'\|\xbold\|_b \leq k'c'\|\xbold\|_a$. 
\end{itemize}
\end{proof}

\nT{Main property} The main reason for defining equivalent norms is that two equivalent norms on a same vector space $E$ define the same topology on $E$. We have the following key theorem
\begin{itemize}[label=$\rightarrow$]
 \item two norms on a finite dimensional vector space are equivalent
 \item there exists non equivalent norms on an infinite dimensional vector space
\end{itemize}
which is equivalent to, in vocabulary of topology
\begin{itemize}[label=$\rightarrow$]
 \item The closed unit ball $\overline{\B}(0,1)$ in a vector space is compact if, and only if, $E$ is finite dimensional.
\end{itemize}
Let us recall that
\[
 \overline{\B}(0,1) = \{\xbold \in E \mid \|\xbold\| \leq 1\}
\]

%
\section{Some important norms}\label{sec:inter:some}
%

\nB Let $p$ be an real number in  $[1,+\infty]$. Let 
\[
 \xbold \in E \qquad \mbox{with} \qquad \xbold = (x_1, \ldots,x_n)
\]
The map 
\begin{equation}
 \begin{CD}
  E @>>> \R^+
 \end{CD}
\end{equation}
denoted $\|.\|_p$ defined by
\begin{equation}
 \|\xbold\|_p = \left(\sum_i|x_i|^p\right)^{\frac{1}{p}}
\end{equation}
is a norm, called norm $\ell^p$. It is not a simple exercise to show that it is a norm, especially to show that $\|\xbold+\ybold\|_p \leq \|\xbold\|_p+\|\ybold\|_p$. This is known as \kw{Minkowski inequality}. A nice presentation with a complete demonstration relying on \kw{Hölder inequality} can be found in \cite[chap. 9]{Steele2008}

\nS In particular, we have
\begin{equation}
 \left\{
    \begin{array}{lcl}
     \|\xbold\|_1 &=& \displaystyle \sum_i |x_i| \\
      \|\xbold\|_2 &=& \displaystyle \sqrt{\sum_i x_i^2} \\
      \|\xbold\|_\infty &=& \underset{i}{\max} \: |x_i|
    \end{array}
 \right.
\end{equation}
$\ell^1$ is called \kw{taxicab norm}\index{norm!taxicab} or \kw{Manhattan norm}\index{norm!Manhattan} and $\ell^\infty$ is called \kw{maximum norm}\index{norm!maximum}, whereas $\ell^2$ is the classical \kw{Euclidean norm}\index{norm!Euclidean} or \kw{Frobenius norm}\index{norm!Frobenius}. In $E=\R^2$, these norms can be seen as\\
\\
\begin{center}
\begin{tikzpicture} 
\draw[green] (0,0) circle(2) ; 
\draw[blue] (-2,-2) rectangle(2,2) ;
\draw[red] (0,-2) -- (2,0) -- (0,2) -- (-2,0) -- (0,-2); 
\draw[->] (-3,0) -- (3,0) ;
\draw[->] (0, -3) -- (0,3) ;
\node () at (6,2) {\textcolor{red}{red}: $\{\xbold \mid \|\xbold\|_1=1\}$} ; 
\node () at (6,0) {\textcolor{green}{green}: $\{\xbold \mid \|\xbold\|_2=1\}$} ; 
\node () at (6,-2) {\textcolor{blue}{blue}: $\{\xbold \mid \|\xbold\|_\infty=1\}$} ; 
\end{tikzpicture}
\end{center}

\nS One sees that there is more room in unit ball of $\|.\|_\infty$ than of $\|.\|_2$ and of $\|.\|_1$ which has the smallest area. Indeed, if $\mu(\B)$ denotes the measure (area, volume, ...) of a ball, it is easy to show that
\begin{equation}
 \left\{
  \begin{array}{lcl}
   \mu(\B_1) &=& 2 \\
   \mu(\B_2) &=& \pi \\
   \mu(\B_\infty) &=& 4 \\
  \end{array}
 \right.
\end{equation}
where $\B_p$ is the unit ball for $\|.\|_p$. More generally, for $E=\R^n$, the volume $\mu_n$ of ball $\B_p$ for norm $\ell^p$ is\footnote{This is a nice formula for $p$ with $1 \leq p \leq +\infty$ which can be found at \\ \texttt{https://www.johndcook.com/blog/2010/07/02/volumes-of-generalized-unit-balls/}}
\begin{equation}
 \mu_n(\B_p) = 2^n\frac{\Gamma\left(1+\frac{1}{p}\right)^n}{\Gamma\left(1+\frac{n}{p}\right)}
\end{equation}
recalling that $\Gamma(z)=(z-1)!$ for $z \in \N^*$. We have\\
\begin{center}
\begin{tabular}{l|ccc}
  & $p=1$& $p=2$ & $p=\infty$ \\
  \hline 
  $\Gamma\left(1+\frac{1}{p}\right)$ & 1 & $\Gamma(\frac{3}{2})$ & 1 \\
  $\Gamma\left(1+\frac{n}{p}\right)$ & $n!$ & $\Gamma(\frac{n+2}{2})$ & 1 \\
\end{tabular}
\end{center}
so
\begin{equation}
 \left\{
  \begin{array}{lcl}
   \mu_n(\B_1) &=& 2^n/n! \\
   \mu_n(\B_2) &=& 2^n\frac{\Gamma\left(3/2\right)}{\Gamma\left(\frac{n+2}{2}\right)} \\
   \mu_n(\B_\infty) &=& 2^n \\
  \end{array}
 \right.
\end{equation}

\nS Let  
\begin{equation}
    \rho_n = \frac{\mu_n(\B_\infty)}{\mu_n(\B_1)} = n!
\end{equation}
Then
\begin{equation}
    \lim_{n \rightarrow \infty}\:\rho_n = + \infty
\end{equation}
There is infinitely more room in unit ball with maximum norm than with $\ell^1$ norm in large dimensional spaces. This has far reaching consequences. A tensor of order $r$ with dimension $n$ on each mode lives in $\R^\n$ with $N=n^d$. Very quickly for significant $d$, $N \rightarrow \infty$. It is then much easier to approximate a tensor with maximum norm than with $\ell^1$ norm.

\nS Let us now select 
\[
 \ones = (1,\ldots,1) \in \R^n 
\]
We have
\begin{equation}
 \left\{
   \begin{array}{lcl}
    \|\ones\|_1 &=& n \\
    \|\ones\|_2 &=& \sqrt{n} \\
    \|\ones\|_\infty &=& 1
   \end{array}
 \right.
\end{equation}
This shows that
\begin{equation}
 \exists \: \xbold \in \R^n \quad \st \quad 
 \left\{
   \begin{array}{lcl}
   \|\xbold\|_1 &\geq& \sqrt{n} \; \|\xbold\|_2 \\
   \|\xbold\|_2 &\geq& \sqrt{n} \; \|\xbold\|_\infty \\
      \end{array}
 \right.
\end{equation}
When $n \rightarrow \infty$, these norms are no longer equivalent, which means that a series may converge for $\|.\|_\infty$ and not for $\|.\|_2$, or for $\|.\|_2$ and not for $\|.\|_1$. 

%
\section{Key consequence for tensors and variational approach}\label{sec:inter:key4tens}
%

Let $\Abold \in (\R^n)^{\otimes d}$. It has $n^d$ coefficients, and can be considered as a vector $\abold \in \R^{n^d}$. The dimension of the vector space where the tensor lives is $N=n^d$, which is more than huge when $n$ is significant and $d$ large. The norm $\|.\|_\infty$
is the norm for the component-wise accuracy of an approximation, which is often met in applications. For example, if both $\Abold, \Bbold \in (\R^n)^{\otimes d}$, then, $\|\Abold-\Bbold\|_\infty < \epsilon$ means
\begin{equation}
 \forall \: \ibold = (i_1, \ldots,i_d) \in I^d, \qquad |a_{i_1\ldots i_d} - b_{i_1\ldots i_d} | < \epsilon
\end{equation}
where $\ibold = (i_1,\ldots,i_d)$ runs over all the indices of the different modes. This may occur even if the $\ell^2$ approximation is large. For example, let $\Abold$ be filled with 0's and $\Bbold$ be filled with $\epsilon$. We the have
\begin{equation}
 \left\{
    \begin{array}{lcl}
      \|\Abold - \Bbold\|_\infty &=& \epsilon \\
      \|\Abold - \Bbold\|_2 &=& \epsilon\,\sqrt{N} \\
      \|\Abold - \Bbold\|_1 &=& \epsilon\,N 
    \end{array}
 \right.
 \qquad \mbox{with} \quad N=n^d
\end{equation}

\nB Depending on the problem to solve, one or the other norm is more relevant. We focus here on two standard operations in applications which are cornerstones in statistics and statistical physics:
\begin{itemize}[label=$\rightarrow$]
 \item the map 
 \[
  \begin{CD}
   \Abold @>Z>> Z(\Abold) = \sum_\ibold a_\ibold
  \end{CD}
 \]
 (computation of the partition function $Z$). If all terms in the tensor are non negative $(a_\ibold \geq 0)$, then $Z(\Abold)=\|\Abold\|_1$. If we have a tensor $\Abold$ for which $Z(\Abold)$ is intractable, and a family $\AS$ for which $Z(\Abold')$ is tractable if $\Abold'\in \AS$, then, the best approximation of $Z(\Abold)$ will be given by $Z(\Abold')$ for $\Abold' \in \AS$ such that $\|\Abold-\Abold'\|_1$ is minimum. This is called the variational approach.
 \item the map 
 \[
  \begin{CD}
   \Abold @>\max>> \underset{\ibold}{\max} \:\: a_\ibold
  \end{CD}
 \]
 (the computation of the mode of a distribution). Then, $\max_\ibold \: a_\ibold = \|\Abold\|_\infty$. Similarly, the best approximation when it is intractable for $\Abold$ is to compute $\max\: \Abold'$ for $\Abold'$ being the best approximation of $\Abold$ in $\AS$ for the norm $\|.\|_\infty$.
\end{itemize}

\nS A challenge arises from the following observation: whereas the problem is set in driving applications on best approximations with norms $\ell^1$ (partition function) or $\ell^\infty$ (computational statistics), most of rigorous results and efficient algorithms for the best approximation of a tensor $\Abold$ by a tensor $\Abold'$ living in a set $\AS$ are known for norm $\ell^2$, and are often very difficult for norm $\ell^1$ or $\ell^\infty$. For example, the best rank $r$ approximation of a matrix $A$ (prescribed rank) is well known for norm $\ell^2$ and given by first terms in the SVD of $A$. Such a problem still is open for norm $\ell^1$ (see \cite{Chierichetti2017}).

%
\section{Notes on the norms of matrices}\label{sec:inter:notes}
%

Let us consider the set $\M_n := \R^{n \times n}$ of $n \times n$ real matrices. It is a vector space, with an inner product defined by $\langle A,B\rangle = \Tr A^\t B$ and, as such, endowed with a norm $\ell^2$ with $\|A\|= \sqrt{\langle A,A\rangle}$. Such a norm is known as \kw{Frobenius norm}\index{norm!Frobenius}. 

\nS There exists many norms on $\M_n$, and most common ones are the \kw{spectral norm}\index{norm!spectral} or \kw{operator norm}\index{norm!operator}, and \kw{$\ell^1$ norm}\index{norm!$\ell^1$} and \kw{$\ell^\infty$ norm}\index{norm!$\ell^\infty$}, defined as
\begin{equation}
 \left\{
   \begin{array}{lclcl}
     \|A\|_1 &=& \displaystyle \sum_{i,j} |\alpha_{ij}| &\qquad& \ell^1 \mbox{ norm}\\
     \|A\|_2 &=& \displaystyle \left(\sum_{i,j} |\alpha_{ij}|^2\right)^{1/2} &\qquad& \ell^2 \mbox{ norm, or Frobenius norm, or } \|.\|\\
     \|A\|_\infty &=& \displaystyle \underset{1 \leq i,j \leq n}{\sup} \: |\alpha_{ij}|  &\qquad& \ell^\infty \mbox{ norm, or uniform convergence norm}\\
     \|A\|_{\mathrm{sp}} &=& \displaystyle \max_{\|\xbold\|=1}\;\|A\xbold\| & \qquad & \mbox{spectral norm or operator norm}
   \end{array}
 \right.
\end{equation}
More generally, for $1 \leq p \leq \infty$, we have the \kw{$\ell^p$ norm}\index{norm!$\ell^p$}
\begin{equation}
 \|A\|_p = \left(\sum_{i,j}|\alpha_{ij}|^p\right)^{1/p}
\end{equation}
The Frobenius norm will be denoted $\|.\|$ when there is no ambiguity. It is the norm which is induced by the inner products in $E$ and $F$ (be they finite or infinite dimensional). Let $A=U\Sigma V^\t$ denote the SVD of matrix $A \in \R^{m \times n}$ (see section \ref{sec:inter:svd}). Then
\begin{equation}
    \|A\|^2 = \sum_i\, \sigma_i^2
\end{equation}
i.e. The Frobenius norm of $A$ is the Euclidean norm of the vector of its singular values. The $\ell^1$ norm of the vector of the singular values of $A$ $(\sum_i\sigma_i)$ is called the \kw{nuclear norm}\index{norm!nuclear} of $A$. It is denoted $\|A\|_*$.

\nS If $n$ is finite, all these norms are equivalent (they define the same topology, i.e. series which converge for one norm converge for the others too). This is no longer true if $n$ is infinite, and this has consequences on the choice of the norm for best low rank approximation in finite but large dimensional spaces. The decrease of the singular values of a matrix $A \in \R^{n \times n}$ gives the accuracy with which it can be approximated by a lower rank matrix for norm $\ell^2$. A matrix with worst low rank approximation is a matrix without decrease in its singular values, i.e. all singular values are equal. It is the identity (or proportional to it)
\begin{equation}
 \I_n = \sum_{i=1}^n \ebold_i \otimes \ebold_i
\end{equation}
The best rank $r$ approximation is given by 
\begin{equation}
 \I_{n,r} = \sum_{i=1}^r \ebold_i \otimes \ebold_i
\end{equation}
and 
\begin{equation}
 \|\I_n-\I_{n,r}\|^2 = \sum_{\alpha>r} \sigma_\alpha^2=n-r \qquad (\mbox{as } \sigma_\alpha=1)
\end{equation}
Let us now consider the diagonal matrix
\[
 D_n = \frac{1}{\sqrt{n}}\I_n = \sum_{i=1}^n \, \frac{1}{\sqrt{n}} \; \ebold_i \otimes \ebold_i
\]
We have
\[
 \|D_n\|=1
\]
It is the matrix with unit $\ell^2$ norm which has the worst low rank approximation.  If $D_{n,r}$ is the best rank $r$ approximation of $D_n$
\begin{equation}
 \|D_n-D_{n,r}\| = \sqrt{1-\frac{r}{n}}
\end{equation}

\nS Let us now consider the matrix $D'_{n,r}$ with $1/\sqrt{n}$ as the first $r$ terms in the diagonal and 0 elsewhere in the diagonal. For example
\[
 D'_{5,2} = 
 \begin{pmatrix}
   1/\sqrt{n} & 0 & \hdots & \hdots & 0 \\
   0 & 1/\sqrt{n} & 0 & \ddots & \vdots \\
   \vdots &  0 & 0 & \ddots  &  \vdots\\
   \vdots & \ddots & \ddots & \ddots & \vdots \\
   0 & \hdots & \hdots & 0 & 0
 \end{pmatrix}
\]
Then, $\rank D'_{n,r}=r$, and 
\begin{equation}
 \forall \: r, \quad \|D_n-D'_{n,r}\|_\infty = \underset{i,j}{\sup} \: |D_n[i,j]-D'_{n,r}[i,j]| = \frac{1}{\sqrt{n}}
\end{equation}
This means that, for the norm $\ell^\infty$,  $D$ can be approximated by a matrix of rank $r$ (even $r=1)$ with an accuracy $\epsilon = \O(n^{-1/2})$, which is an excellent accuracy when $n$ is sufficiently large. This is simply due to the fact that $\|D_n\|_\infty= n^{-1/2}$, although $\|D_n\|_2=1$.

\notes A nice review of unusual properties of the rank of big data matrices is \cite{Udell2019} which we have followed for these notes. Mathematically speaking, low rank matrices are sparse among matrices: any randomly selected matrix is expected to be full rank: singular matrices are nowhere dense. Moreover, the singular values of a random Gaussian matrix are high with a high probability. But it has been observed many times and in various situations that a matrix from a real big data set is very well approximated by a low rank matrix with high probability. This simply means that real world is not random. In other words, its singular values decay exponentially fast. It is then interesting to understand which properties of a matrix ensure a fast decay of its singular values. It has been shown recently that matrices with a displacement structure have such fast decaying singular value property \cite{Beckermann2017}. A matrix $X \in \C^{m \times n}$ with $(A,B)-$displacement structure with $A \in \C^{m \times m}$ and $B \in \C^{n \times n}$ with $n \geq m$ is defined as follows: there are a matrix $M \in \C^{m \times r}$ and $N \in \C^{n \times r}$ such that $AX-XB=MN^*$. Such families include Toeplitz, Hankel, Cauchy, Krylov, Vandermonde, Pick, Sylvester and Löwner matrices (with $r=1$ or $r=2$ accordingly). This has been recently expanded to matrices built from a nice latent variables model \cite{Udell2019}. The rank approximation presented here is the approximation with prescribed rank: a rank $r$ is prescribed, and the objective is to find the best rank $r$ approximation. Some problems are related with best low rank approximation with prescribed accuracy $\epsilon$: the objective is to find the minimal rank $r$ such that there exists a rank $r$ approximation with prescribed accuracy. This leads to the $\epsilon-$rank of a matrix $A \in \R^{n \times n}$: it is the smallest integer $r \leq n$ such that there is a matrix $X_r$ of rank $r$ such that $\|A-X_r\|_\infty < \epsilon$, or $\forall \: i,j, \quad |a_{ij}-x_{ij}| < \epsilon$. Still more puzzling \cite[thrm 1.0]{Udell2019} have shown that, for $n$ very large, any $n \times n$ matrix has $\epsilon-$rank $r$ for $r = \lceil 72 \Log (2n+1)/\epsilon^2\rceil$. Knowing that, understanding which properties ensure an accurate low rank approximation of a tensor still is an uncharted continent.  

%
\section{The Singular Value Decomposition (SVD)}\label{sec:inter:svd}
%

I recall here basic notions on Singular Value Decomposition (SVD) of a matrix which will be a cornerstone for computing the best low rank approximation of a tensor especially in Tucker and Tensor Train best approximations. It deserves to be a chapter, but as it is singled-minded and short, it is included as a section in this key part on tensor rank: SVD is a rank revealing tool for matrices.

\nS Let $A \in E \otimes F \simeq \L(F,E) \simeq \R^{m \times n}$ with $m = \dim E$ and $n=\dim F$. Let us assume that $\rank A = r$. Then, SVD reads
\begin{equation}\label{eq:inter:svd:dyadic}
 A = \sum_{\alpha=1}^r\, \sigma_\alpha \, \ubold_\alpha \otimes \vbold_\alpha \qquad \mbox{with} \quad 
 \begin{cases}
  \ubold_\alpha &\in \R^m \\
  \vbold_\alpha &\in \R^n \\
  \sigma_\alpha &\in \R^+
 \end{cases}
\end{equation}
The vectors $(\ubold_\alpha)_\alpha$ and $(\vbold_\alpha)_\alpha$ form an orthonormal system, i.e.
\begin{equation}
 \forall \: \alpha,\alpha', \qquad 
 \begin{cases}
  \langle \ubold_\alpha,\ubold_{\alpha'}\rangle &= \delta^\alpha_{\alpha'} \\
  \langle \vbold_\alpha,\vbold_{\alpha'}\rangle &= \delta^\alpha_{\alpha'} 
 \end{cases}
\end{equation}
If $A$ is full rank, we have $r = \min(m,n)$. Decomposition (\ref{eq:inter:svd:dyadic}) is called the \kw{dyadic decomposition} of matrix $A$. It is easy to see that
\begin{equation}
    \left\{
        \begin{array}{lcl}
             A\vbold_\alpha &=& \sigma_\alpha \, \ubold_\alpha \\
              A^\t \ubold_\alpha &=& \sigma_\alpha \, \vbold_\alpha \\
        \end{array}
    \right.
\end{equation}
The vectors $(\ubold_\alpha)_{\alpha}$ are the \emph{left} \kw{singular vectors} of $A$, and $(\vbold_\alpha)_\alpha$ the \emph{right singular vectors}. The real numbers $(\sigma_\alpha)_\alpha$ are the \kw{singular values}, and are such that
\begin{equation}
    \sigma_1 \geq \sigma_2 \geq \ldots \geq \sigma_r > 0
\end{equation}
We have
\begin{equation}
    \left\{ 
       \begin{array}{lclcl}
            A^\t A \vbold_\alpha &=& \sigma_\alpha A^\t \ubold_\alpha &=& \sigma_\alpha^2 \vbold_\alpha \\
            AA^\t \ubold_\alpha &=& \sigma_\alpha A \vbold_\alpha &=& \sigma_\alpha^2 \ubold_\alpha 
       \end{array}
    \right.
\end{equation}
So, $(\lambda_\alpha)_\alpha$ with $\lambda_\alpha=\sigma_\alpha^2$ are the eigenvalues of symmetric matrices $A^\t A$ and $AA^\t$, the eigenvectors of which are respectively the families $(\vbold_\alpha)_\alpha$ and $(\ubold_\alpha)_\alpha$. This establishes a natural connection between SVD and Principal Component Analysis.

\nS The vectors $(\ubold_\alpha)_\alpha$ (resp. $(\vbold_\alpha)_\alpha$) are concatenated column-wise to build a matrix $U \in \R^{m \times r}$ (resp. $V \in \R^{n \times r}$), with
\begin{equation}
 UU^\t=VV^\t=\I_m
\end{equation}
The SVD of $A$ can be written
\begin{equation}
 A=U\Sigma V^\t
\end{equation}
where $\Sigma = \Diag(\sigma_1, \ldots,\sigma_m)$. It  can be visualized by

\nP{$\rightarrow$} if $n \leq m$
\\
\begin{center}
\begin{tikzpicture}
 \draw (0,0) rectangle (1,3);
\draw (1.25,3.25) rectangle (2.25,4.25) ; 
\draw (2.5,4.5) rectangle (3.5,5.5) ; 
\node () at (0.5,1.5) {$U$} ; 
\node () at (1.75,3.75) {$\Sigma$} ;
\node () at (3,5) {$V^\t$} ;  
\draw[dashed] (1.25,4.25) -- (2.25,3.25) ; 
\node () at (-.3,1.5) {$m$} ; 
\node () at (0.5,3.2) {$n$} ; 
\node () at (1,3.75) {$n$} ; 
\node () at (1.75,4.5) {$n$} ;
\node () at (2.2,5) {$n$} ; 
\node () at (3,4.3) {$n$} ; 
\end{tikzpicture}
\end{center}

\nP{$\rightarrow$} if $n \geq m$
\\
\begin{center}
\begin{tikzpicture}
 \draw (0,0) rectangle (1,1);
\draw (1.25,1.25) rectangle (2.25,2.25) ; 
\draw (2.5,2.5) rectangle (5.5,3.5) ; 
\node () at (0.5,0.5) {$U$} ; 
\node () at (1.75,1.75) {$\Sigma$} ;
\node () at (4,3) {$V^\t$} ;  
\draw[dashed] (1.25,2.25) -- (2.25,1.25) ; 
\node () at (-.3,.5) {$m$} ; 
\node () at (0.5,1.2) {$m$} ; 
\node () at (1,1.75) {$m$} ; 
\node () at (1.75,2.45) {$m$} ;
\node () at (2.2,3) {$m$} ; 
\node () at (4,3.7) {$n$} ; 
\end{tikzpicture}
\end{center}

\noindent where the dashed line represents the diagonal of $m$ non-zeros element in $\Sigma$. 

\nS Here, I will denote
\begin{equation}
 W = \Sigma V^\t, \qquad \mbox{with} \quad W \in \R^{m \times n} 
\end{equation}
Let us denote the components of these matrices as
\begin{equation}
 A(i,j), \qquad U(i,\alpha), \qquad W(\alpha,j), \qquad \mbox{with} \quad
 \begin{cases}
   i & \in \llbracket 1,m \rrbracket \\
   j & \in \llbracket 1,n \rrbracket \\
   \alpha & \in \llbracket 1,r \rrbracket \\
 \end{cases}
\end{equation}
Then
\begin{equation}
 A(i,j) = \sum_{\alpha=1}^r U(i,\alpha)\,W(\alpha,j)
\end{equation}

\nS Let us denote $u(i) \in \R^m$ the row $i$ of $U$, and $w(j) \in \R^m$ the column $j$ of $W$. Then, as $A=UW$,
\begin{equation}
 A(i,j) = \langle u(i) \, , \, w(j) \rangle
\end{equation}
and
\begin{equation}
 \begin{array}{lcl}
   Z(\Abold) &=& \displaystyle \sum_{i,j}\, A(i,j) \\
   &=& \displaystyle \sum_{i,j} \, \langle u(i) \, , \, w(j) \rangle \\
   &=& \displaystyle \left\langle \sum_iu(i) \, , \, \sum_jw(j)\right\rangle
 \end{array}
\end{equation}
which requires $r(m+n)$ additions and $r$ multiplications instead of $mn$ additions. Tensor-train format will extend such a calculation of $d-$modes tensors with $d \gg 2$. Let us note that there is no free lunch: the gain from $mn$ (a product) to $r(m+n)$ (a sum) for $(m,n)$ is compensated by $r$ additional multiplications. This will prove to be stronger for tensors of large order.

\nT{Dyadic decomposition and pseudo-inverse} Let us recall that the \kw{Moore-Penrose pseudo-inverse} of a matrix $A \in \R^{m \times n}$ is the matrix $A^\dagger \in \R^{n \times m}$ satisfying to all four conditions:
\begin{itemize}[label=$\rightarrow$]
    \item $AA^\dagger A=A$
    \item $A^\dagger A A^\dagger = A^\dagger$
    \item $(AA^\dagger)^\t = AA^\dagger$
    \item $(A^\dagger A)^\t = A^\dagger A$
\end{itemize}
It can be shown that, $A$ being given, there exists one, and one only, matrix $A^\dagger$ satisfying these conditions. Let $A$ be defined as in (\ref{eq:inter:svd:dyadic}). Then
\begin{equation}
    A^\dagger = \sum_{\alpha=1}^r\, \frac{1}{\sigma_\alpha}\, \vbold_\alpha \otimes \ubold_\alpha
\end{equation}
or
\begin{equation}
    A^\dagger = V\Sigma^{-1}U^\t
\end{equation}
Checking that this matrix satisfies to all four conditions is left to the reader. It is sufficient to remember that $(\abold \otimes \bbold)(\cbold \otimes \dbold)= \langle \bbold,\cbold\rangle \abold \otimes \dbold$. Or, for the first one
\[
\begin{array}{lcl}
    AA^\dagger A &=& (U\Sigma V^\t)(V\Sigma^{-1} U^\t)(U\Sigma V^\t) \\
    &=& U\Sigma (V^\t V)\Sigma^{-1} (U^\t U)\Sigma V^\t \\
    &=& U\Sigma\Sigma^{-1}\Sigma V^\t \\
    &=& U\Sigma V^\t \\
    &=& A
\end{array}
\]

\nT{CUR decomposition} Let us have $A \in \R^{m \times n}$ with $m,n$ large. SVD is one way to derive an accurate approximation of $A$ in a compact form. However, the interpretation of the values in the singular values or singular vectors may be not obvious. Let us select a rank $r$, and two sets
\[
I \subset \llbracket 1,m \rrbracket \quad \mbox{and} \quad J \subset \llbracket 1,n \rrbracket
\]
with $|I|=|J|=r$. Let us denote by
\begin{itemize}[label=$\rightarrow$]
    \item $\widehat{A}=A[I,J] \in \R^{r \times r}$
    \item $C = A[:,J] \in \R^{m \times r}$
    \item $R = A[I,:] \in R^{r \times n}$
\end{itemize}
Once the subsets $I$ and $J$ have been selected, the best approximation of $A$ obtained by sampling these rows and columns is given by
\begin{equation}
    A \approx C\widehat{A}^{-1}R
\end{equation}
Next question is to find $I$ and $J$ such that this approximation itself as function of $(I,J)$ is optimal. A classical answer is to use a submatrix $\widehat{A}$ with maximal "volume", i.e. determinant in modulus.

\notes Approximating a given matrix $A \in \R^{m \times n}$ by a matrix $\widetilde{A}$ of lower rank $k$ is an overwhelming question in numerical linear algebra which has relentlessly been studied since its early developments by the end of 19th century by Beltrami (1873), Jordan (1874), Sylvester (1890) and Autonne (1915) who discovered as well the polar decomposition in 1902 (see \cite{Stewart1993,Horn2012}). According to \cite{Stewart1999}, the truncated SVD is one of the most elegant solutions, and provides it for both the Frobenius and the spectral norm (the spectral norm of a matrix $A$ is $\max \: \{ \|Ax\| \mid \|x\|=1\}$). There are many textbooks presenting SVD. See e.g. \cite[sect. I.8]{Strang2019} for an accessible presentation and demonstration of its existence (and uniqueness up to some symmetries). We recommend \cite{Stewart1993} and chapter 3 of \cite{Horn2012} for detailed historical remarks. SVD plays an important role in data compression, as well as in unitarily invariant norms. The problem of computing singular values and vectors of a matrix is well conditioned, hence, algorithms are robust.\\
The CUR or skeleton or pseudo-skeleton decomposition of a matrix has been studied in two seminal papers by Goreinov, Tyrtyshnikov and Zamarashkin in \cite{Goreinov1997,Goreinov1997a}. The CUR decomposition has been studied in \cite{Mahoney2009}.

%% file: best_low_rank.tex
Let us start with the example of CP-rank for $d=3$. If a tensor $\Abold \in \R^{m \times n \times p}$ has rank $r$, it is possible to decompose it as a sum of $r$ elementary tensors at most:
\[
 \Abold = \sum_{a=1}^{\rank \Abold} \xbold_a \otimes \ybold_a \otimes \zbold_a
\]
Let us now fix an integer $r' < r$ or an error $\epsilon >0$ and raise two questions:
\begin{itemize}[label=$\rhd$]
 \item Which is the tensor $\Abold'$ of rank $r'$ which is the closest to $\Abold$ ? what is the error?
 \item Which is the smallest integer $r' \leq r$ such that there exists a tensor $\Abold'$ of rank $r'$ with
\[
 \|\Abold - \Abold'\| < \epsilon \, \|\Abold\|
\]
\end{itemize}
"Closest to" is understood in this chapter and the followings as measured with the distance associated to $\ell^2$ norm. The problems can be studied with other norms, but such a domain is much less advanced and understood (see chapter \ref{chap:inter}).

\nB The first problem is known as \emph{best prescribed rank approximation} and the second as \emph{best prescribed accuracy approximation}. The answer to both questions for matrices, i.e. $2-$tensors, is given by Eckart-Young theorem, and is at the root of PCA and SVD. It is an association between an algebraic view (rank minimization) and a geometric view (a point cloud in a Euclidean space is associated to a matrix, and PCA is finding a subspace of a given dimension with optimal projection). Nothing like this happens for $d-$modes tensors with $d \geq 3$:
\begin{itemize}[label=$\rhd$]
 \item there is no equivalent to Eckart-Young theorem, although such a generalization has been looked for for decades,
 \item geometric and algebraic views dissociate: one has to chose between orthogonality and minimum rank
\end{itemize}

\nB As mentioned above, the norm used in this chapter is Frobenius norm, or $\ell^2$ norm:
\begin{equation*}
 \|\Abold\|^2 = \sum_{i_1}\ldots \sum_{i_d}\, a_{i_1\ldots i_d}^2
\end{equation*}
The solution of prescribed accuracy best approximation for norm $\ell^\infty$ is known as \kw{$\epsilon-$rank}\index{rank!$\epsilon$}.

\nB Three types of ranks have been selected in this book:
\begin{itemize}[label=$\rightarrow$]
    \item CP rank
    \item Tucker rank
    \item TT rank (Tensor Train)
\end{itemize}
(Tensor Train decomposition has not been presented before in chapter \ref{chap:rank3} and will be presented at the beginning of chapter \ref{chap:tt}). Best low rank approximation enables to address high-dimensional issues, like dimension reduction of disentangling interactions in a model. Typically, dimension reduction for matrices is PCA (Principal Components Analysis) and its various guises, and this is naturally extended to tensors by best low rank approximation by a Tucker model. This keeps the flavor of Euclidean structures: orthogonal projections, Frobenius norm and, with a selected tuning, nestedness of vector subspaces (in what is called HOSVD, for Higher-Order Singular Value Decomposition). Hence, HOSVD is a well understood and efficient reduction dimension method for tensors with reasonable order and high dimensions along the orders. The other benefit of low rank approximation is to disentangle interactions in complex systems, by selecting structures with local interactions between few orders only. This has been thoroughly studied and developed in many-body quantum systems, as \emph{tensor networks}. Tensor networks model interactions between tensors of a given family by contracting over repeated modes. Interesting models are those where these repetitions are local, between neighbors if the tensors are set as nodes of a graph (a network). A paradigm of such an approach is given by \emph{tensor trains}, where the network is a line, and contractions are between neighbors only. This has been studied and developed in quantum systems as MPS (Matrix Product States), and will be presented here. Knowing that, CP decomposition and approximaton are in-between. They are related to dimension reduction by extending to tensors the most parsimonious way to decompose a tensor as a sum of rank one tensor, and to disentangling interactions between modes as, for example, $\Abold = \sum_{a \leq r} \xbold_a \otimes \ybold_a \otimes \zbold_a$ can be considered as a weak coupling between modes if $r$ is small. 

\nB This part is simply organized as presenting Best Low Rank Approximation methods for different notions of rank, one per chapter:
\begin{description}
 \item[chapter \ref{chap:blracp}] BLRA for CP-rank
 \item[chapter \ref{chap:blratuck}] BLRA for Tucker rank
 \item[chapter \ref{chap:tt}] BLRA for TT rank
\end{description}

%
\chapter{Best low rank approximation of a tensor: CP rank}\label{chap:blracp}
%

In this chapter, elementary results for Best Low rank Approximation of CP-rank of a tensor are given. CP stands for Candecomp-Partafac, or Canonical Polyadic decomposition. It is with Tucker model one of the historical model for rank decomposition or approximation. For a prescribed rank $r$, it is the model with lowest storage complexity, better than Tucker model or Tensor Train model of same rank. It suffers however from a possibility of numerical indeterminacy. It is widely used for Blind Source Separation.

\nB This chapter is short and organized as follows:
\begin{description}
\item[in section \ref{sec:blracp:3modes}] BLRA of CP-rank for 3-modes tensors is presented with a naive ALS scheme
\item[in section \ref{sec:blracp:opt}] An other optimization scheme is presented
\item[in section \ref{sec:blracp:ill}] It is recalled that CP BLRA may be an ill-posed problem
\item[in section \ref{sec:blracp:dmodes}] These notions are extended in a straightforward way to $d-$modes tensors.

\end{description}

%
\section{Prescribed rank Best CP-approximation for $3-$modes tensors}\label{sec:blracp:3modes}
%

We first present prescribed rank best approximation of a tensor on a example of a $3-$modes tensor $\Abold \in E \otimes F \otimes G$. Generalization to $d-$modes tensors is straightforward and presented in section \ref{sec:blracp:dmodes}.

\nT{Setting the problem}
Let $\Abold \in E \otimes F \otimes G$ with
\[
 \dim E = m, \qquad \dim F=n, \qquad \dim G=p
\]
One wishes to solve the following optimization problem:\\
\\
\begin{center}
\ovalbox{
\begin{tabular}{ll}
given & $\Abold \in E \otimes F \otimes G$ \\
& $r \in \N$ \\
&\\
find & $(\xbold_a)_a \in E$ \\
& $(\ybold_a)_a \in F$ \\
& $(\zbold_a)_a \in G$ \\ 
with & $a \in \llbracket 1,r\rrbracket$ \\
&\\
such that & $\displaystyle \left\|\Abold - \sum_a \xbold_a \otimes \ybold_a \otimes \zbold_a\right\|$ is minimal
\end{tabular}
}
\end{center}


\nS \textbf{Cost function:}
Let us denote $X,Y,Z$ the matrices with columns $(\xbold_a)_a$, $(\ybold_a)_a$, $(\zbold_a)_a$
\begin{equation}\label{sec:blracp:3modes:eq1}
\begin{array}{ccccccccccc}
X &=& [\xbold_1|\ldots|\xbold_r] & \qquad & Y &=& [\ybold_1|\ldots|\ybold_r] &\qquad& Z &=& [\zbold_1|\ldots|\zbold_r] \\
& & \in \R^m \otimes \R^r & \qquad & && \in \R^n \otimes \R^r &\qquad& && \in \R^p \otimes \R^r \\
& & (m \times r) & & & & (n \times r) & & & & (p \times r)
\end{array}
\end{equation}
and define
\begin{equation}
 \phi(\Abold,X,Y,Z) = \left\|\Abold - \sum_{a=1}^r \xbold_a \otimes \ybold_a \otimes \zbold_a\right\|^2
\end{equation}
Solving best CP-rank approximation of $\Abold$ is minimizing the cost function $\phi$. Let us note that
\begin{equation}
    \sum_{a=1}^r \xbold_a \otimes \ybold_a \otimes \zbold_a = \langle X \mid Y \mid Z \rangle
\end{equation}

\nS \textbf{Alternate Least Squares:} We present here a naive algorithm by Alternate Least Square on each of the vectors $(\xbold_a)_a$, $(\ybold_a)_a$, $(\zbold_a)_a$. Let us fix $a \in \{1,r\}$ and minimize $\phi$ on $\zbold_a$ only. Let us denote
\[
 \Bbold_a = \Abold - \sum_{i \neq a}  \xbold_i \otimes \ybold_i \otimes \zbold_i
\]
Then
\begin{equation}
 \phi(\Abold,X,Y,Z) = \|\Bbold_a - \xbold_a \otimes \ybold_a \otimes \zbold_a \|^2
\end{equation}
Observing that $\xbold \otimes \ybold \otimes \zbold = \alpha \xbold \otimes \beta \ybold \otimes \gamma \zbold$ if $\alpha\beta\gamma=1$, one can set $\|\xbold_a\|=\|\ybold_a\|=1$. Then\footnote{Let us recall that $(\abold \otimes \bbold \otimes \cbold) \underset{1,2}{\bullet} (\xbold \otimes \ybold) = \langle \abold,\xbold\rangle\langle \bbold,\ybold\rangle \cbold$, which is extended to $E \otimes F \otimes G$ by linearity. See section \ref{sec:contract:3modes}}
\begin{equation}
 \begin{aligned}
   \phi(\Abold,X,Y,Z) &= \left\|\Bbold_a- \xbold_a \otimes \ybold_a \otimes \zbold_a\right\|^2 \\
   &= \|\Bbold_a\|^2 + \|\zbold_a\|^2 - 2\langle \Bbold_a\, , \, \xbold_a\otimes \ybold_a \otimes \zbold_a \rangle \\
   &= \|\Bbold_a\|^2 + \|\zbold_a\|^2 - 2\langle \Bbold_a \underset{1,2}{\bullet} (\xbold_a\otimes \ybold_a)\, , \, \zbold_a \rangle \\
 \end{aligned}
\end{equation}
So
\begin{equation}
 \nabla_{\zbold_a} \, \phi = 2\zbold_a - 2 \, \Bbold_a \underset{1,2}{\bullet} (\xbold_a \otimes \ybold_a)
\end{equation}
Setting $\nabla_{\zbold_a}\, \phi=0$ leads to 
\begin{equation}
 \zbold_a = \Bbold_a \underset{1,2}{\bullet} (\xbold_a \otimes \ybold_a)
\end{equation}

\nP{$\rightarrow$} By setting $\|\xbold_a\|=\|\zbold_a\|=1$, the same can be done on $Y$ as
\begin{equation}
 \begin{aligned}
   \phi(\Abold,X,Y,Z) &= \left\|\Bbold_a- \xbold_a \otimes \ybold_a \otimes \zbold_a\right\|^2 \\
   &= \|\Bbold_a\|^2 + \|\ybold_a\|^2 - 2\langle \Bbold_a\, , \, \xbold_a\otimes \ybold_a \otimes \zbold_a \rangle \\
   &= \|\Bbold_a\|^2 + \|\ybold_a\|^2 - 2\langle \Bbold_a \underset{1,3}{\bullet} (\xbold_a\otimes \zbold_a)\, , \, \ybold_a \rangle \\
 \end{aligned}
\end{equation}
Setting $\nabla_{\ybold_a}\, \phi=0$ leads to 
\begin{equation}
 \ybold_a = \Bbold_a \underset{1,3}{\bullet} (\xbold_a \otimes \zbold_a)
\end{equation}

\nP{$\rightarrow$} and on $X$ by setting $\|\ybold_a\|= \|\zbold_a\|=1$:
\begin{equation}
 \begin{aligned}
   \phi(\Abold,X,Y,Z) &= \left\|\Bbold_a- \xbold_a \otimes \ybold_a \otimes \zbold_a\right\|^2 \\
   &= \|\Bbold_a\|^2 + \|\xbold_a\|^2 - 2\langle \Bbold_a\, , \, \xbold_a\otimes \ybold_a \otimes \zbold_a \rangle \\
   &= \|\Bbold_a\|^2 + \|\xbold_a\|^2 - 2\langle \Bbold_a \underset{2,3}{\bullet} (\ybold_a\otimes \zbold_a)\, , \, \xbold_a \rangle \\
 \end{aligned}
\end{equation}
Setting $\nabla_{\xbold_a}\, \phi=0$ leads to 
\begin{equation}
 \xbold_a = \Bbold_a \underset{2,3}{\bullet} (\ybold_a \otimes \zbold_a)
\end{equation}

\nP{$\rightarrow$} Last step is to combine those three optimizations in a loop with, at step $t$,
\begin{equation}
    \left\{ 
        \begin{array}{lcl}
           \zbold_a^{t+1} &=& \Bbold_a^t \underset{1,2}{\bullet} (\xbold_a^t \otimes \ybold_a^t) \\
           \ybold_a^{t+1} &=& \Bbold_a^{t+1/3} \underset{1,3}{\bullet} (\xbold_a^t \otimes \zbold_a^{t+1}) \\
            \xbold_a^{t+1} &=& \Bbold_a^{t+2/3} \underset{2,3}{\bullet} (\ybold_a^{t+1} \otimes \zbold_a^{t+1})
        \end{array}
    \right.
\end{equation}
with
\begin{equation}
    \left\{ 
        \begin{array}{lcl}
           \Bbold_a^t &=& \Abold - \sum_{i \neq a}  \xbold_i^t \otimes \ybold_i^t \otimes \zbold_i^t \\
           \Bbold_a^{t+1/3} &=& \Abold - \sum_{i \neq a}  \xbold_i^t \otimes \ybold_i^t \otimes \zbold_i^{t+1} \\
           \Bbold_a^{t+2/3} &=& \Abold - \sum_{i \neq a}  \xbold_i^t \otimes \ybold_i^{t+1} \otimes \zbold_i^{t+1} \\
        \end{array}
    \right.
\end{equation}
The $CP-$approximation by ALS simply is iterating this step, with a stopping condition, until a fixed point is reached. Such an algorithm suffers from classical drawbacks of ALS:
\begin{itemize}[label=$\rightarrow$]
    \item the fixed point can be a local minimum
    \item the procedure may diverge.
\end{itemize}

%
\section{Other optimization schemes}\label{sec:blracp:opt}
%

One can easily imagine other ALS procedures, like minimizing on $(\xbold_a, \ybold_a, \zbold_a)$ with loops on index $a$. Several of them are presented in \cite[section~3.4]{Kolda2009}. Here, we present an approach for an optimization of $Z$ globally.

\nS Let us denote
\begin{equation}
 \phi(Z) = \left\| \Abold - \sum_{a=1}^r \xbold_a \otimes \ybold_a \otimes \zbold_a \right\|^2
\end{equation}
with $X,Y,Z$ defined as in equation (\ref{sec:blracp:3modes:eq1}). The objective is to compute $Z$ such that $\phi(Z)$ is minimal. Let us denote
\[
 U_a = \xbold_a\otimes \ybold_a
\]
Then
\begin{equation}
 \phi(Z) = \left\| \Abold - \sum_a U_a \otimes \zbold_a \right\|^2
\end{equation}

\nS Let us recall that, if $A,B \in \R^{m \times n}$
\begin{equation}
 \langle A,B\rangle = \Tr AB^\t = \Tr A^\t B
\end{equation}
Then (provided the dimensions are consistent)
\begin{equation}
 \begin{array}{lcl}
  \langle A,BC \rangle &=& \Tr A(BC)^\t \\
  &=& \Tr AC^\t B^\t \\
  &=& \langle AC^\t \, , B \rangle
 \end{array}
\end{equation}
or
\begin{equation}
 \begin{array}{lcl}
  \langle A,BC \rangle &=& \Tr A^\t BC \\
  &=& \Tr (B^\t A)^\t C \\
  &=& \langle B^\t A, C\rangle
 \end{array}
\end{equation}

\nS Let us now consider $A \in \R^{mn \times p}$, a family $(\ubold_a)_a$ with $\ubold_a \in \R^{mn}$, with $a$ running over $\llbracket 1, r \rrbracket$ and $\|\ubold_a\|=1$, and find the family $(\zbold_a)_a$ such that
\begin{equation}
 \left\| A - \sum_a \ubold_a \otimes \zbold_a\right\| \qquad \mbox{minimal}
\end{equation}
Let us define the matrix $U \in \R^{mn \times r}$ (resp. $Z \in \R^{p \times r}$) with $\ubold_a$ (resp. $\zbold_a$) as column $a$. Then
\begin{equation}
  \sum_a \ubold_a \otimes \zbold_a = UZ^\t 
\end{equation}
and
\begin{equation}
 \begin{array}{lcl}
  \phi(Z) &=& \|A - UZ^\t\|^2 \\
  &=& \|A\|^2 + \|UZ^\t\|^2 -2 \langle A,UZ^\t\rangle 
 \end{array}
\end{equation}

\nP{$\rightarrow$} We have
\begin{equation}
 \begin{array}{lcl}
  \|UZ^\t\|^2 &=& \langle UZ^\t \, , \, UZ^\t \rangle \\
  &=& \langle U^\t UZ^\t \, , \, Z^\t \rangle \\
  &=& \langle ZU^\t U \, , \, Z \rangle \\
 \end{array}
\end{equation}
and
\begin{equation}
  \begin{array}{lcl}
     \langle A,UZ^\t\rangle &=& \langle U^\t A \, , \,Z^\t\rangle \\
     &=& \langle A^\t U\, , \,Z\rangle
  \end{array}
\end{equation}

\nP{$\rightarrow$} Then
\begin{equation}
 \phi(Z) = \|A\|^2 + \langle ZU^\t U \, , \, Z \rangle -2 \langle A^\t U \, , \, Z \rangle
\end{equation}
and
\begin{equation}
 \nabla_Z\phi = 2ZU^\t U  - 2A^\t U
\end{equation}
If $\phi$ is minimal, $\nabla_Z\phi=0$, and
\begin{equation}
 ZU^\t U = A^\t U
\end{equation}
$U \in \R^{m \times r}$. Let us assume that $mn > r$ with $U^\t U$ invertible. Then
\begin{equation}
 Z = A^\t U(U^\t U)^{-1}
\end{equation}
(note that as $U \in \R^{mn \times r}$ and $A \in \R^{mn \times p}$, we have $(U^\t U)^{-1} \in \R^{r \times r}$, so $U(U^\t U)^{-1} \in \R^{mn \times r}$ and $A^\t U(U^\t U)^{-1} \in \R^{p \times r}$).

\nS An ALS scheme on computing $(X,Y,Z)$ with this approach can be implemented, hoping that it converges to a fixed point which is the expected minimum minimorum.

%
\section{Ill-posedness}\label{sec:blracp:ill}
%

The optimization problem of minimizing $\phi$ may be ill-posed. Indeed, let us recall that $\mathbb{CP}_r$ denotes the set of tensors of rank $r$ (see section \ref{sec:rank3:def} for this notation). Finding the best approximation of $\Abold$ by a tensor of rank $r$ is finding a projection of $\Abold$ on $\mathbb{CP}_r$, or 
\begin{equation}
 \mbox{finding} \:\: \widetilde{\Abold} \in \mathbb{CP}_r \st \|\Abold-\widetilde{\Abold}\| \quad \mbox{is minimal}
\end{equation}
For this problem to have a (possibly non unique) solution, the set $\mathbb{CP}_r$ must be closed. However, it may be not closed (see definition of border rank in section \ref{sec:rank3:def}). A classical example is when $E \otimes F \otimes G$ has two typical ranks. For example, the typical rank of a tensor $\Abold \in \R^2 \otimes \R^2 \otimes \R^2$ is 2 or 3. Kruskal \& \emph{al.} have shown that the boundary between $\mathbb{CP}_2$ and $\mathbb{CP}_3$ contains tensors of any rank between 0 and 3 \cite{Kruskal1989a}. As a corollary, neither $\mathbb{CP}_2$ nor $\mathbb{CP}_3$ are closed. This means that there exists infinite Cauchy sequences
\[
 (\Abold_n)_{n \in \N}, \quad \rank A_n = 2 \quad \forall \: n \in \N
\]
with 
\[
 \rank \left(\lim_{n \to \infty} \Abold_n\right) \neq 2
\]
This may cause degeneracy of the computation of the solution. This is a well identified phenomenon, and some remedies are known, like introducing orthogonality constraints (see \cite{Yang2020} for a recent review).

%
\section{CP best approximation for $d-$modes tensors}\label{sec:blracp:dmodes}
%

CP best approximation at prescribed rank can be extended straightfully to $d-$modes tensors, as follows.

\nT{Setting the problem} Let
\[
 \Abold \in E_1 \otimes \ldots \otimes E_d
\]
and $r \in \N$. CP prescribed rank best approximation of $\Abold$ is finding $\left(\xbold_{a}^{(1)}, \ldots, \xbold_{a}^{(d)}\right)$ with $\xbold_{a}^{(\mu)} \in E_\mu$ such that
\begin{equation}
 \left\| \Abold - \sum_{a=1}^r \, \xbold_{a}^{(1)} \otimes \ldots \otimes \xbold_{a}^{(d)}\right\| \qquad \mbox{is minimal}
\end{equation}
It is solving the following optimization problem:\\
\begin{center}
 \ovalbox{
  \begin{tabular}{lll}
    given && $\Abold \in E_1 \otimes \ldots \otimes E_d$ \\
    && $r \in \N$ \\
    &&\\
    find && $r$ $d-$uplets of vectors $\left(\xbold_{a}^{(1)}, \ldots, \xbold_{a}^{(\mu)}, \ldots \xbold_{a}^{(d)}\right)$ \\
    with && $a \in \llbracket 1,r \rrbracket$ \\
    && $\xbold_{a}^{(\mu)} \in E_\mu$\\
    &&\\
    such that && $\displaystyle \left\| \Abold - \sum_{a=1}^r \, \xbold_{a}^{(1)} \otimes \ldots \otimes \xbold_{a}^{(d)}\right\|$ is minimal 
  \end{tabular}
 }
\end{center}

\nT{Solving the problem} There is no recipe for finding the best approximation, but some to find local minima. Let us define matrices
\[
 X_1 \in \R^{n_1 \times r}, \ldots, X_\mu \in \R^{n_\mu \times r}, \ldots, X_d \in \R^{n_d \times r} 
\]
with column $a$ of $X_\mu$ being $\xbold_{a}^{(\mu)}$. Then, we denote 
\begin{equation}
 \sum_{a=1}^r \, \xbold_{a}^{(1)} \otimes \ldots \otimes \xbold_{a}^{(d)} = \langle X_1| \ldots | X_\mu| \ldots |X_d\rangle
\end{equation}
The idea is to iterate on matrices $X_\mu$. Let us start with $\mu=d$: let us fix matrices $X_1, \ldots, X_{d-1}$ and define the cost function
\begin{equation}
 \phi(X_d) = \| \Abold -  \langle X_1| \ldots  | X_{d-1}|X_d\rangle\|^2
\end{equation}
Let us denote 
\begin{equation}
 U_a =  \xbold_{a}^{(1)} \otimes \ldots \otimes \xbold_{a}^{(d-1)} \in E_1 \otimes \ldots \otimes E_{d-1}
\end{equation}
Then
\begin{equation}
 \phi(X_d) = \left\| \Abold - \sum_a U_a \otimes \xbold_{a}^{(d)} \right\|^2
\end{equation}
The calculation developed in section \ref{sec:blracp:opt} for $d=3$ can be developed here \emph{nec varietur} for $d \geq 3$, yielding

\nP{$\rightarrow$} matricize $\Abold$ as
\[
 \begin{array}{ll}
   A & \in (E_1 \otimes \ldots \otimes E_{d-1}) \otimes E_d \\
   & \in \R^{n_1...n_{d-1} \times n_d}
 \end{array}
\]

\nP{$\rightarrow$} vectorize $U_a$ as
\[
 \ubold_a \in \R^{n_1...n_{d-1} \times n_d}
\]

\nP{$\rightarrow$} then
\begin{equation}
 \phi(X_d) = \left\| A - \sum_a \ubold_a \otimes \xbold_{a}^{(d)}\right\|^2
\end{equation}
and the optimization problem is rewritten according to\\
\begin{center}
 \begin{tabular}{|ll}
  given & $N = n_1...n_{d-1} \in \N$ \\
  & $r \in \N $ \\
  & $A \in \R^{N\times n_d}$ \\
  & $\ubold_a \in \R^N$ for $a \in \llbracket 1,r \rrbracket$\\
  &\\
  find & $\xbold_{1}^{(d)}, \ldots, \xbold_{r}^{(d)}$ with $\xbold_{a}^{(d)} \in \R^{n_d}$ \\
  &\\
  such that & $\displaystyle \left\| A - \sum_{a=1}^r \ubold_a \otimes \xbold_{a}^{(d)} \right\|^2$ is minimal
 \end{tabular}
\end{center}

\nP{$\rightarrow$} Let us denote by $U \in \R^{N \times r}$ the matrix with $\ubold_a$ in column $a$. Then
\begin{equation}
 \begin{array}{lcl}
   \phi(X_d) &=& \|A - UX_d^\t\|^2 \\
   &=& \|A\|^2 + \|UX_d^\t\|^2 - 2 \langle A,UX_d^\t \rangle
 \end{array}
\end{equation}
with, as in section \ref{sec:blracp:opt}
\begin{equation}
 \left\{
  \begin{array}{lcl}
    \|UX_d^\t\|^2 &=& \langle X_dU^\t U \, , \, X_d \rangle \\
    \langle A,UX_d^\t \rangle &=& \langle A^\t U \, , \, X_d \rangle 
  \end{array}
 \right.
\end{equation}
and
\begin{equation}
 \nabla \phi = 2 X_dU^\t U - 2A^\t U
\end{equation}
which vanishes when $X_dU^\t U - A^\t U=0$ or
\begin{equation}\label{eq:blra:cp:dmodes:1}
 X_d = A^\t U(U^\t U)^{-1}
\end{equation}
which minimizes (exactly here) $\phi(X_d)$ when all $X_\mu$ but $X_d$ are fixed.

\nP{$\rightarrow$} The solution of the optimization problem can be obtained by iterating successively on all matrices $X_\mu$ until a convergence criteria is met.

\notes PARAFAC (for Parallel Factors) and CANDECOMP (Canonical Decomposition) are two independent procedures for analysis of 3-ways tensors, developed the same year (1970), by Harshman for PARAFAC \cite{Harshman1970} and Carroll \& Chang for CANDECOMP \cite{Carroll1970}. Both procedures proved to be equivalent, and are referred to globally as Parafac/CanDecomp, or even more briefly CP. It is so standard that a new acronym has been derived: CP stands for Canonical Polyadic Decomposition, referring Hitchcok's seminal paper on the rank of a tensor. ALS was proposed in these seminal papers, and still is the classical tool for implementing CP best approximation. This approach is well developed in \cite[section~3]{Kolda2009}, or \cite{Franc1992}. \cite{Morup2011} presents in details CP decomposition and best approximation, with a variety of related models. See as well \cite{Kroonenberg2008}. It is well known that the naive ALS optimisation scheme as presented here is not optimal. It is a tip of an iceberg. For example, several sources of degeneracy have been identified, as when two (bottleneck) or more (swamp) factors are nearly colinear, and  thoroughly studied. See \cite{Comon2009} and literature within for a complete and detailed presentation of these degeneracies, historical remarks and ways to circumvent them.


%
\chapter{Best low rank approximation of a tensor: Tucker rank}\label{chap:blratuck}
%

Let $A\in \K^{m \times n} \simeq E \otimes F$ with here $\K=\R$. Let it be full rank ($\rank A = n$ if $m \geq n$), but this is not compulsory. If $r = \rank A$, it is possible to write
\[
 A = \sum_{a=1}^r \, \sigma_a \;  \ubold_a \otimes \vbold_a, \qquad \mbox{with} \quad 
 \begin{cases}
    A\vbold_a &=\sigma_a\, \ubold_a\\
    A^\t \ubold_a &=\sigma_a\, \vbold_a\\
    A^\t A \vbold_a &= \sigma^2_a \vbold_a \\
    AA^\t \ubold_a &= \sigma^2_a \ubold_a \\
 \end{cases}
\]
with $(\ubold_a)_a$ and $(\vbold_a)_a$ orthonormal families. This is the SVD, which is written as well
\[
 A = U\Sigma V^\t
\]
where $U \in \K^{m \times r}$ (resp. $V \in \K^{n \times r}$) with columns $(\ubold_a)_a$ (resp. $(\vbold_a)_a$). It is tempting to extend such a decomposition beyond $d=2$, with something like
\[
 \Abold = \sum_{a=1}^r \, \sigma_a \; \ubold_a \otimes \vbold_a \otimes \wbold_a \qquad \in E \otimes F \otimes G
\]
with orthonormality of families $(\ubold_a)_a$ etc. But this is not always possible because, if $m \geq n \geq p$ (with $m=\dim E, n=\dim F, p=\dim G$) and $r > p$, no family of $r$ vectors $(\wbold_a)_a$ can be orthonormal.

\nB Trying to build such an extension has been subject of many (unsuccessful) efforts, although it has been clearly acknowledged that there is no decomposition of a $d-$mode tensor with $d \geq 3$ which possesses all the nice properties of the SVD of a matrix ($d=2$) \cite{Franc1989,Zang2001}. SVD is a decomposition of a matrix on a given basis for $E$ and $F$. The extension of this is a decomposition of a tensor on a given basis for each mode. This is Tucker decomposition (see section \ref{sec:rank3:tucker}), which  is the expression of $\Abold$ in a new basis for each of the mode, selected in a way such that some nice properties of SVD can be recovered on each mode, like
\begin{itemize}
 \item[$\rhd$] the slices along the mode are pairwise orthogonal
 \item[$\rhd$] the norm is concentrated as much as possible in the first slices
\end{itemize} 

\nS This chapter is organized in a logical way of nested algorithms as a succession of short sections, each devoted to a problem, as follows:
\begin{description}
 \item[section \ref{sec:blratuck:set}:] Setting the problem in its full generality for $d=3$
 \item[section \ref{sec:blratuck:hosvd3}:] Presenting a non optimal but easy to implement solution, HOSVD, for $3-$modes tensors
 \item[section \ref{sec:blratuck:hooi}:] Presenting a better solution, HOOI, for $3-$modes tensors
 \item[section \ref{sec:blratuck:dmodes}:] Presenting Tucker best low rank approximation for $d-$modes tensors
 \item[section \ref{sec:blratuck:hosvddmodes}:] Presenting HOSVD for $d-$modes tensors
 \item[section \ref{sec:blratuck:hooidmodes}:] Presenting HOOI for $d-$modes tensors
\end{description}

%
\section{Setting the problem for prescribed rank}\label{sec:blratuck:set}
%

Let $E,F,G$ be three vector spaces with $\dim E=m$, $\dim F=n$ and $\dim G=p$. Let us recall that if $\Abold \in E \otimes F \otimes G$, $\mathbf{mr(\Abold)}$ denotes the Tucker or multilinear rank of $\Abold$ (see section \ref{sec:rank3:tucker}). There is one rank per mode. Tucker rank of $\Abold$ is $\rbold = (m',n',p')$ with $m'\leq m, n' \leq n, p' \leq p$ if there exists three subspaces $E' \subset E$, $F' \subset F$ and $G' \subset G$ with $\dim E'=m'$, $\dim F'=n'$, $\dim G'=p'$ and $\Abold \in E' \otimes F' \otimes G'$.

\nB Tucker best approximation at prescribed rank can be stated as follows:\\
\begin{center}
 \ovalbox{
   \begin{tabular}{ll}
     given & a tensor $\Abold \in E \otimes F \otimes G$ \\
     & with $\dim E=m, \dim F=n, \dim G=p$\\
     & a multilinear rank $\rbold = (m',n',p')$ \\
     &\\
     find & a tensor $\widetilde{\Abold} \in E \otimes F \otimes G$ \\
     & with $\mathbf{mr}(\widetilde{\Abold})= \rbold$ \\
     &\\
     such that & $\|\Abold - \widetilde{\Abold}\|$ is minimal
   \end{tabular}
 }
\end{center}

\nB Then,
\begin{equation}
    \exists \: E',F',G' \quad \mbox{with} \quad 
    \begin{cases}
      \dim E' &= m' \\
      \dim F' &= n' \\
      \dim G' &= p' \\
    \end{cases}
    \quad \st \quad
    \widetilde{\Abold} \in E' \otimes F' \otimes G'
\end{equation}
We now rephrase this definition using orthogonal projections on subspaces of given dimension.

\nT{Orthogonal projections on subspaces} Let $\Abold \in E \otimes F \otimes G$. If we consider subspaces 
\[
E' \subset E, \qquad F' \subset F, \qquad G' \subset G 
\]
then
\[
E' \otimes F' \otimes G' \subset E \otimes F \otimes G 
\]
We have the following diagram of orthogonal projections:\\
\begin{center}
  \begin{tikzpicture}
   \node (100) at (0,0) {$E' \otimes F \otimes G$} ;
   \node (110) at (4,0) {$E' \otimes F' \otimes G$} ;
   \node (000) at (0,3) {\textcolor{red}{$E \otimes F \otimes G$}} ;
   \node (010) at (4,3) {$E \otimes F' \otimes G$} ;
   \node (101) at (2,1) {$E' \otimes F \otimes G'$} ;
   \node (111) at (6,1) {\textcolor{red}{$E' \otimes F' \otimes G'$}} ;
   \node (001) at (2,4) {$E \otimes F \otimes G'$} ;
   \node (011) at (6,4) {$E \otimes F' \otimes G'$} ;   
   \draw[->] (000) -- (001) ; 
   \draw[->] (000) -- (010) ; 
   \draw[->] (000) -- (100) ;    
   \draw[->] (001) -- (011) ; 
   \draw[->, dashed] (001) -- (101) ; 
   \draw[->] (010) -- (011) ;  
   \draw[->] (010) -- (110) ;     
   \draw[->, dashed] (100) -- (101) ;  
   \draw[->] (100) -- (110) ;  
   \draw[->] (011) -- (111) ;  
   \draw[->, dashed] (101) -- (111) ;  
   \draw[->] (110) -- (111) ;    
  \end{tikzpicture} 
\end{center}
where $\A \longrightarrow \B$ denotes a projection from $\A$ on $\B$. This diagram is commutative if projections are orthogonal. We can denote as follows the different projections of a given tensor $\Abold$ depending on the subspaces $E',F',G'$:\\
\begin{center}
  \begin{tikzpicture}
   \node (100) at (0,0) {$\Abold_{\e'} \in E' \otimes F \otimes G$} ;
   \node (110) at (6,0) {$\Abold_{\e'\f'} \in E' \otimes F' \otimes G$} ;
   \node (000) at (0,3) {\textcolor{red}{$\Abold \in E \otimes F \otimes G$}} ;
   \node (010) at (6,3) {$\Abold_{\f'} \in E \otimes F' \otimes G$} ;
   \node (101) at (2,1) {$\Abold_{\e'\g'} \in E' \otimes F \otimes G'$} ;
   \node (111) at (8,1) {\textcolor{red}{$\Abold_{\e'\f'\g'} \in E' \otimes F' \otimes G'$}} ;
   \node (001) at (2,4) {$\Abold_{\g'} \in E \otimes F \otimes G'$} ;
   \node (011) at (8,4) {$\Abold_{\f'\g'} \in E \otimes F' \otimes G'$} ;   
   \draw[->] (000) -- (001) ; 
   \draw[->] (000) -- (010) ; 
   \draw[->] (000) -- (100) ;    
   \draw[->] (001) -- (011) ; 
   \draw[->, dashed] (001) -- (101) ; 
   \draw[->] (010) -- (011) ;  
   \draw[->] (010) -- (110) ;     
   \draw[->, dashed] (100) -- (101) ;  
   \draw[->] (100) -- (110) ;  
   \draw[->] (011) -- (111) ;  
   \draw[->, dashed] (101) -- (111) ;  
   \draw[->] (110) -- (111) ;    
  \end{tikzpicture} 
\end{center}

\nT{Setting the problem with projections} Let $\widehat{\Abold}$ be the orthogonal projection of $\Abold$ on $E' \otimes F' \otimes G'$. Then, Tucker best approximation at rank $(m',n',p')$ can be stated as the following problem:\\
\begin{center}
\ovalbox{
\begin{tabular}{lcl}
 Given & & $\Abold \in E \otimes F \otimes G$ \\
 & & with $\dim E = m, \: \dim F = n, \: \dim G = p$ \\ 
 & & ranks $m' \leq m, \: n' \leq n, \: p' \leq p$ \\
 & & \\
 find && a subspace $E' \subset E, \: F' \subset F, \, G' \subset G$ \\
 & & with $\dim E'=m', \: \dim F'=n', \: \dim G' = p'$ \\
 & & \\
 such that & & $\|\Abold - \widehat{\Abold}\|$ minimal\\
 && with $\widehat{\Abold}$ being the projection of $\Abold$ on $E' \otimes F' \otimes G'$ 
\end{tabular}
}
\end{center}

\nT{Well posedness of Tucker best low rank approximation} We show here that, contrary to CP best low rank approximation, Tucker best low rank approximation is well posed. To see this, we use the notion of $m-$frame of a vector space. A $m-$frame of a vector space $E$ is a set of $m$ linearly independent vectors. Let us set $m',n',p'$ and consider the sets $\FS_{m'}(E)$ of orthonormal $m'-$frames of $m'$ vectors in $E$, $\FS_{n'}(F)$ of orthonormal $n'-$frames of $n'$ vectors in $F$ and $\FS_{p'}(G)$ of orthonormal $p'-$frames of $p'$ vectors in $P$. We have $\FS_{m'}(E) \subset E^{m'}$. $E$ is an Euclidean space with norm induced by the inner product. This norm induces a topology such that a base of open sets is the set of open balls with the metric induced by the norm:\\
\[
 \begin{CD}
  \mbox{inner product} @>>> \mbox{norm} @>>> \mbox{distance} @>>> \mbox{open balls} @>>> \mbox{topology}
 \end{CD}
\]
Let us equip $E^{m'}$ with the product topology. The key observation is that $\FS_{m'}(E)$ is compact for this topology.
\begin{proof}
Indeed, it is closed and bounded in a finite dimensional vector space, hence compact. To see this, let $\F = (\xbold_1, \ldots,\xbold_{m'}) \in \FS$ be an orthonormal frame of $m'$ vectors in $E^{m'}$. We have
\[
 \|\F\|^2 = \sum_{i=1}^{m'}\|\xbold_i\|^2
\]
hence 
\[
 \|\F\|^2 = m'
\]
and $\FS$ is bounded. Let us define\footnote{Here, for this demonstration only, $G$ is the Gram matrix, without confusion with vector space $G$}
\[
 \begin{CD}
  E^{m'} @>G>> \K^{m' \times m'} \\
 \end{CD}
\]
with $G(\F)$ being the Gram matrix of the frame $\F =(\xbold_1,\ldots,\xbold_{m'}) \in E^{m'}$. It is the matrix with coefficients $G_{ij}=\langle \xbold_i,\xbold_j\rangle$. The map $G$ is continuous. $G(\F)$ is the identity matrix because $\F$ is orthonormal, and
\[
 \FS_m = G^{-1}(\{\I\})
\]
where $\I$ is the identify matrix. As  $\{\I\}$ is a single point, it is closed, and as the preimage of a closed set by a continuous function is closed, $\FS_m$ is closed.
\end{proof}
Let us now consider\footnote{We come back to $G$ being a vector space ...}
\begin{itemize}[label=$\rightarrow$]
 \item a triplet of frames 
\[
 \FS_{m'}(E) = (\ebold_1, \ldots,\ebold_{m'}), \quad \FS_{n'}(F) = (\fbold_1, \ldots,\fbold_{n'}), \quad \FS_{p'}(G) = (\gbold_1, \ldots,\gbold_{p'}), 
\]
\item the spaces $E'=\span (\FS_{m'}(E))$, $F'=\span (\FS_{n'}(F))$, $G'=\span (\FS_{p'}(G))$,
\item the projection $\widehat{\Abold}$ of $\Abold$ on $E' \otimes F' \otimes G'$
\end{itemize}
Then, $\Abold$ being fixed, the map
\[
 \begin{CD}
   \FS_{m'}(E) \times \FS_{n'}(F) \times \FS_{p'}(G) @>>> \|\Abold - \widehat{\Abold}\|
 \end{CD}
\]
is continuous. As the set of frames is compact, it is a continuous function from a compact set on $\R$. It then has a minimum, and this minimum is reached for at least one frame. It then has at least one solution. This does not guarantee that this solution is unique.

\nS There is no known algorithm yielding the exact solution. There are however two standard heuristics for developing Tucker best approximation :
\begin{itemize}[label=$\rightarrow$]
 \item a simple one, consisting in $d$ SVD of $d$ matricized tensors, called HOSVD
 \item HOOI, which goes beyond HOSVD, but with a more complex procedure
\end{itemize}
Both are presented here. None of them is an exact solution of best Tucker rank approximation of a tensor $\Abold$ with prescribed rank. They are both organized as a same three steps process:
\begin{enumerate}
 \item finding an orthonormal basis for each of the space $E,F,G$ with optimality properties close to SVD,
 \item decomposing the tensor on this basis,
 \item truncating the decomposition to find the best approximation.
\end{enumerate}
The first basis vectors of $E$ (resp. $F$, $G$) span the space $E'$ (resp. $F'$, $G'$).

%
\section{HOSVD for $d=3$}\label{sec:blratuck:hosvd3}
%

Let $\Abold \in E \otimes F \otimes G$, and let us select a mode, e.g. $G$. Let $(\gbold_k)_k$ be an orthonormal basis of $G$. Then (see sections \ref{sec:contract:alamat} for matricization and \ref{sec:matelem:basis} for the decomposition of a matrix on a basis)
\begin{equation}
 \Abold = \sum_k A_k \otimes \gbold_k \qquad \mbox{with} \quad 
 \left\{
  \begin{array}{ll}
   A_k & \in \K^{m \times n} \\
   \gbold_k & \in \K^p
  \end{array}
 \right.
\end{equation}
As $(\gbold_k)_k$ is orthonormal, we have $A_k=\Abold\bullet_3 \gbold_k$ and
\[
 \|\Abold\|^2 = \sum_k \|A_k\|^2
\]
HOSVD is finding an orthonormal basis $(\gbold_k)_k$ of $G$ (and the same for other modes $E$ and $F$) such that the norm of $\Abold$ is concentrated as much as possible in the first slices $(A_k)_k$, with slices $A_k$ being pairwise orthogonal. The concentration of the inertia of the tensor in the first components on each mode permits to build a good low Tucker-rank approximation by truncating the decomposition on its first components, in the same way than a SVD is a decomposition, and truncating it at a given threshold yields a best rank $r$ approximation of a matrix.

\nS Let us consider the matricization $A_\g$ of $\Abold$ defined by
\begin{equation}
A_\g = \sum_kA_k \otimes \gbold_k \qquad \mbox{with} \qquad 
\left\{
\begin{array}{lll}
  A_\g & \in \L(G, E \otimes F) & \simeq \K^{mn \times p}\\
  A_k & \in E \otimes_\k F & \simeq \K^{mn}
\end{array}
\right.
\end{equation}
where the slices $A_k$ have been vectorized (see section \ref{sec:tenselem:reshape}). Let us now run the SVD of $A_\g$, which reads (see section \ref{sec:inter:svd})
\begin{equation}
 A_\g = U_\g \Sigma_\g V_\g^\t \qquad \mbox{with} \quad 
 \left\{
 \begin{array}{ll}
  U_\g & \in \K^{mn \times p} \\
  \Sigma_\g & \in \K^{p \times p} \\
  V_\g & \in \K^{p \times p} \\
 \end{array}
 \right.
\end{equation}
Column of $U_\g$ are orthonormal, as well as columns of $V_\g$. We select
\begin{equation}
 \left\{
   \begin{array}{lcl}
    A_k & \leftarrow & \sigma_kU_{\g,k} \\
    \gbold_k &=& V_{\g,k}
   \end{array}
 \right.
\end{equation}
where\footnote{There is an $\leftarrow$ and not an $=$ to $A_k$ because $A_k$ is not equal to $\sigma_kU_{\g,k}$ but built from it by a reshaping.}
\begin{itemize}[label=-]
 \item $U_{\g,k}$ is the $k-$th column of $U_\g$, which as been reshaped as a $m \times n$ matrix to become a slice of $\Abold$
 \item $V_{\g,k}$ is the $k-$th column of $V_\g$
\end{itemize}
So, we have both properties
\begin{itemize}[label=$\rightarrow$]
 \item slices $(A_k)_k$ are orthogonal as columns of $U_\g$ are orthogonal
 \item the inertia of $\Abold$ is concentrated in the first slices as $\|A_k\|=\sigma_k$ and $(\sigma_k)_k$ are the singular values of $A_\g$.
\end{itemize}

\nS Knowing that, HOSVD decomposition and best approximation is very simple: it is alternating this step over all modes of $\Abold$, i.e. $E$, $F$ and $G$. Let us consider the three linear applications $A_\e, A_\f, A_\g$ associated to each mode of $\Abold$
\begin{equation}
  \left\{
  \begin{aligned}
    \begin{CD}
     E @>A_\e>> F \otimes G & \qquad A_\e \xbold = X & \qquad X_{jk}= \sum_{i}\alpha_{ijk}x_i \\
     F @>A_\f>> E \otimes G & \qquad A_\f \ybold = Y & \qquad Y_{ik}= \sum_{j}\alpha_{ijk}y_j \\
     G @>A_\g>> E \otimes F & \qquad A_\g \zbold = Z & \qquad Z_{ij}= \sum_{k}\alpha_{ijk}z_k \\    
    \end{CD}
  \end{aligned}
  \right. 
\end{equation}
with 
\[
 A_\e \in \K^{np \times m}, \qquad A_\f \in \K^{mp \times n}, \qquad A_\g \in \K^{mn \times p}
\]
HOSVD is building three systems of basis and slices, one per associated linear application or matricization, with following SVDs\\
\begin{center}
 \begin{tabular}{clll}
  Matricized tensor & SVD & Slice & basis \\
  \hline
  $A_\e$ & $A_\e=U_\e\Sigma_\e V_\e^\t$ & $A_i = U_{\e,i}$ & $\ebold_i=V_{\e,i}$ \\
  $A_\f$ & $A_\f=U_\f\Sigma_\f V_\f^\t$ & $A_j = U_{\f,j}$ & $\fbold_j=V_{\f,j}$ \\
  $A_\g$ & $A_\g=U_\g\Sigma_\g V_\g^\t$ & $A_k = U_{\g,k}$ & $\gbold_k=V_{\g,k}$ \\
 \end{tabular}
\end{center}
Once these three basis have been built, the final step is to decompose $\Abold$ on these three basis, along calculations presented in section \ref{sec:dmodestens:basis}. 

\nT{Best approximation at prescribed rank} Once the decomposition has been made, it is possible to truncate the expression of $\Abold$ to a given rank for each mode. This yields a truncated decomposition, expressed as
\begin{equation}
 \widetilde{\Abold} = \sum_{i=1}^{m'}\sum_{j=1}^{n'}\sum_{k=1}^{p'} \; \alpha_{ijk} \; \ebold_i \otimes \fbold_j \otimes \gbold_k
\end{equation}
This yields HOSVD.

%
\section{HOOI for $d=3$}\label{sec:blratuck:hooi}
%

Tucker approximation is finding the best approximation of a given tensor by a tensor of a prescribed Tucker-rank. It is known that truncation of HOSVD decomposition does not yield the best solution. HOOI improves the approximation, but there is no guarantee as well that HOOI yields the best approximation. As usual, a first presentation is made for $3-$modes tensors, and extension to $d-$modes tensors is straightforward.

\nS It is important to understand at this stage that 
\begin{itemize}
 \item[$\rhd$] the unknowns of the problem are the spaces $E',F',G'$, or the frames spanning them
 \item[$\rhd$] HOSVD is not the best solution to such a problem.
\end{itemize}
To see this and understand to which question HOSVD and HOOI respectively answer, it is useful to have in mind the following diagrams of projection between tensor products of spaces and subspaces. HOSVD is computing spaces (basis of) $E',F',G'$, each being a solution of finding the best projection on respectively $E'\otimes F \otimes G$; $E \otimes F' \otimes G$ and $E \otimes F \otimes G'$. It can be sketched as\\
\begin{center}
\begin{tikzcd}
& \Abold_{\e'} \\
\Abold \arrow[ur] \arrow [r] \arrow[dr] & \Abold_{\f'} \\
& \Abold_{\g'} \\
\end{tikzcd}
\end{center}
with three independent matricizations. In HOOI, the projection $\Abold \longrightarrow \Abold_{\e'\f'\g'}$ is given as in the following composition diagram:
\[
 \begin{CD}
  \Abold @>>> \Abold_{\e'} @>>> \Abold_{\e'\f'} @>>> \Abold_{\e'\f'\g'}
 \end{CD}
\]
where for example the projection on $F'$ takes into account the outcome of the projection on $E'$. It is not sure\footnote{It is rather sure that not} that the composition of optimal projections along this chain at each step yields the optimal simultaneous choice of $E',F',G'$ for projection
\[
 \begin{CD}
  \Abold @>>> \Abold_{\e'\f'\g'}
 \end{CD}
\]
Hence, another strategy must be sought for, and Alternate Least Square is presented here. One can say that, knowing that, HOSVD is the greedy algorithm for finding Tucker best approximation at prescribed rank. Let us mention as well that there is no guarantee that HOOI provides the best Tucker approximation. ALS for HOOI until a fixed point is reached can be sketched by the following diagram:\\
\begin{center}
    \begin{tikzcd}
       \left(E^{'t}, F^{'t}, G^{'t}\right) \arrow[r] & \left(E^{'{t+1}}, F^{'t}, G^{'t}\right) \arrow[r] & \left(E^{'{t+1}}, F^{'{t+1}}, G^{'t}\right) \arrow[r] & \left(E^{'{t+1}}, F^{'{t+1}}, G^{'{t+1}}\right)
    \end{tikzcd}
\end{center}

\nS Let us suppose that $F'$ and $G'$ are given, and let us update $E'$. This reads\footnote{In HOSVD, the first line is: $\Abold \in E \otimes F \otimes G$}:\\
\begin{center}
\begin{tabular}{|clcl}
 &Given & & $\Abold \in E \otimes F' \otimes G'$ \\
 && & $F' \subset F$, $G' \subset G$\\
 && & a rank $m' < m$ \\
 && & \\
 &find && a subspace $E' \subset E$ \\
 & &&with $\dim E' = m'$ \\
 && & \\
 &such that & & $\|\Abold - \widehat{\Abold}\|$ minimum \\
 & && with $\widehat{\Abold}$ being the projection of $\Abold$ on $E' \otimes F' \otimes G'$ 
\end{tabular}
\end{center}

\nB Let us denote 
\[
 \begin{CD}
   \Abold @>\PC_{F' \otimes G'}>> \Abold_{\f'\g'}
 \end{CD}
\]
and, to keep notations shorter
\[
 \PC_X^{\perp}(\Abold) = \PC_{X^\perp}(\Abold)
\]
Then (the elementary changes are given in detail)
\begin{equation}
 \begin{aligned}
   \|\Abold - \widehat{\Abold}\|^2 &= \|\PC_{F' \otimes G'}(\Abold - \widehat{\Abold})  + \PC_{F' \otimes G'}^{\perp}(\Abold - \widehat{\Abold})\|^2   \\
   (a) &= \|\PC_{F' \otimes G'}(\Abold - \widehat{\Abold})\|^2 + \|\PC_{F' \otimes G'}^{\perp}(\Abold - \widehat{\Abold})\|^2   \\
   (b) &= \|\PC_{F' \otimes G'}(\Abold - \widehat{\Abold})\|^2 + \|\PC_{F' \otimes G'}^{\perp}(\Abold) - \PC_{F' \otimes G'}^{\perp}(\widehat{\Abold})\|^2   \\
   (c) &= \|\PC_{F' \otimes G'}(\Abold - \widehat{\Abold})\|^2 + \|\PC_{F' \otimes G'}^{\perp}(\Abold)\|^2   \\
   (d) &= \|\PC_{F' \otimes G'}(\Abold) - \PC_{F' \otimes G'}(\widehat{\Abold})\|^2 + \|\PC_{F' \otimes G'}^{\perp}(\Abold)\|^2   \\
   (e) &= \|\PC_{F' \otimes G'}(\Abold) - \widehat{\Abold}\|^2 + \|\PC_{F' \otimes G'}^{\perp}(\Abold)\|^2    
 \end{aligned}
\end{equation}
because (for each step)
\begin{description}
 \item[$(a)$] $\PC_{F' \otimes G'}(\Abold - \widehat{\Abold})$ and $\PC_{F' \otimes G'}^{\perp}(\Abold - \widehat{\Abold})$ are orthogonal
 \item[$(b)$] linearity of $\PC_{F' \otimes G'}^{\perp}$
 \item[$(c)$] $\PC_{F' \otimes G'}^{\perp}(\widehat{\Abold})=0$
 \item[$(d)$] linearity of $\PC_{F' \otimes G'}$
 \item[$(e)$] $\PC_{F' \otimes G'}(\widehat{\Abold})= \widehat{\Abold})$
\end{description}

\noindent The last term $\|\PC_{F' \otimes G'}^{\perp}(\Abold)\|^2 $ is not affected by the choice of $E'$. Hence, we have, because $\widehat{A} \in E' \otimes F' \otimes G'$
\begin{equation}
  \|\Abold - \widehat{\Abold}\|^2 \mbox{ minimum } \quad \Longleftrightarrow \quad \|\PC_{F' \otimes G'}(\Abold) - \widehat{\Abold}\|^2 \mbox{ minimum } 
\end{equation}
As $F',G'$ are given, let us denote by $A'_{\f'\g'}$ the matricization of $\PC_{F' \otimes G'}(\Abold)$
\[
 A'_{\f'\g'} \in E \otimes (F' \otimes G') \simeq \R^{n'p' \times m} 
\]
Then, the problem can be reformulated as
\begin{center}
\begin{tabular}{|clcl}
 &Given & & $A'_{\f\g} \in \R^{n'p' \times m} $ \\
 && & a rank $m' < m$ \\
 && & \\
 &find && a subspace $E' \subset E$ \\
 & &&with $\dim E' = m'$ \\
 && & \\
 &such that & & $ \|A'_{\f'\g'} - \widehat{\Abold}\|^2$  is minimal \\
 && &with $\widehat{\Abold}$ being the projection of $A'_{\f'\g'}$ on $E' \otimes (F' \otimes G')$ 
\end{tabular}
\end{center}

\nT{Link with PCA with latent variables} Written like this, we recognize the PCA of the projection of $A'_{\f'\g'} \in \R^{n'p' \times m}$, or the PCA with latent variables of matrix $A'$, matricization of $\Abold$ in $\K^{np \times m}$, with constraints set by $F'$ and $G'$.

\nT{Projection operators} Let us now detail the construction of $\Bbold = \PC_{F' \otimes G'}(\Abold)$. Let $\Abold = \sum_{i,j,k}\, a_{ijk}\;\ebold_i \otimes \fbold_j \otimes \gbold_k$. Let $F' \subset F$ and $G' \subset G$ defined as
\[
 F' = \span (\vbold_1, \ldots, \vbold_{n'}), \qquad G' = \span (\wbold_1,\ldots, \wbold_{p'})
\]
Let us denote by $\PC_{F'}$ the orthogonal projector in $F$ on $F'$, and $\PC_{G'}$ the orthogonal projector in $G$ on $G'$. Then, if the $(\vbold_j)_j$ are orthonormal in $F$ and the $(\wbold_k)_k$ orthonormals in $G$, we have
\begin{equation}
 \PC_{F'} = \sum_\alpha \vbold_\alpha \otimes \vbold_\alpha = VV^\t, \qquad \PC_{G'} = \sum_\beta \wbold_\beta \otimes \wbold_\beta = WW^\t
\end{equation}
Then, the projection $\Abold_{\textsc{f'g'}}$ of $\Abold$ on $E \otimes F' \otimes G'$ is given by (see section \ref{sec:alstruct:kron} for notation $\boxtimes$)
\begin{equation}
 \begin{aligned}
   \Abold_{\textsc{f'g'}} &= (\I \boxtimes \PC_{F'} \boxtimes \PC_{G'})\Abold \\
   &= \sum_{i,j,k}\, a_{ijk}\; \ebold_j \otimes \PC_{F'}\fbold_j \otimes \PC_{G'}\gbold_k \\
   &= \sum_{i,j,k}\, a_{ijk}\; \ebold_i \otimes VV^\t \fbold_j \otimes WW^\t \gbold_k \\
 \end{aligned}
\end{equation}
This is similar to the decomposition of $\Abold_{\f'\g'}$ on basis given by matrices $(\I, VV^\t,WW^\t)$. 

\notes Tucker decomposition and truncation has been proposed in \cite{Tucker1966} as a consistent mathematically oriented presentation, after several years of work and presentations made by Tucker and Levine (see \cite{Tucker1966} for a short report of this early history). Although an ALS procedure to go beyond truncated decomposition for a better approximation has been evoked in Tucker's paper, such an approach has been founded by Kroonenberg and de Leeuw seminal paper \cite{Kroonenberg1980} under the name of Tuckals3 for $3-$modes tensors (see as well \cite{Kroonenberg2008}, section 5.7). The extension to four modes tensors or multiway-array has been made by Lastovicka in \cite{Lastovicka1981}, including ALS procedure, and more generally to $n-$modes multiway-arrays including ALS for best approximation by a group in Groningen in 1986 \cite{Kapteyn1986}. Both approaches (decomposition and best approximation) have been later presented by two papers introducing HOSVD and HOOI, with a new and elegant approach relying on matricization and SVD as main tools \cite{Lathauwer2000a,Lathauwer2000b}. Since then, many surveys of Tucker decomposition and approximations have been issued, among which we can cite \cite{Kolda2009,Grasedick2010}. The presentation of Tucker best approximation as a a problem where the unknowns are vector subspaces of given dimension, like the PCA, can be found in \cite[chap. D]{Franc1992}.  Indeed, HOSVD, and HOOI, can be considered as natural extensions of PCA, even if, as the existence of at least two methods suggests, there are more than one natural extension. PCA is an iconic dimension reduction method equivalent to a truncation of the SVD. It has been enriched by several extensions, like PCA with metrics on row and column spaces (among which Correspondence Analysis \cite{Benzecri1973a,Cailliez1979}), and taking into account linear constraints on axis and components (PCA with instrumental variables). All these approaches have been derived for Tucker model, be it HOSVD or HOOI. Taking into account an inner product for each mode has been proposed through the so called \emph{schéma de dualité} \cite{Cailliez1979} or with tensor product in a series of works in the 80's. \cite{Jaffrenou1978} has established a link between Tucker's model with three modes and the so called \emph{schéma de dualité}. \cite{Polit1986} has generalized this approach to $n$ modes, and \cite{Aubigny1989} have developed the study with particular inner products on each mode.  See as well \cite{Mizere1981}. Recently, \cite{Iannacito2021} have shown that the correspondence between row and column point clouds in Correspondence Analysis can be extended to correspondence between point clouds attached to each mode. The analysis of a $n-$ways table with constraints on axis or components has been proposed e.g. in \cite{Kiers1988}. The geometric approach with duality scheme relies on an abstract algebra layer where modes and their dual spaces are involved. Such an approach has been simplified, as well as the approach with Kronecker product in the dutch school by a systematic use of the tensor product in \cite{Franc1992}, which proposes a unifying framework for $d-$ modes Tucker and CP model with metrics and constraints on each modes.

%
\section{Tucker Best rank $\rbold$ approximation for $d$ modes}\label{sec:blratuck:dmodes}
%

Extension of HOSVD to $d-$modes tensors is now straightforward.

\nS Let
\begin{equation}
 \Abold \in \bigotimes_{\mu=1}^d E_\mu
\end{equation}
and $\EC_\mu = (\ebold_{\mu,1}, \ldots,\ebold_{\mu,n_\mu})$ be an orthonormal basis of $E_\mu$. $\Abold$ can be decomposed on these basis as
\begin{equation}\label{eq:tucker:basis}
 \Abold = \sum_{i_1=1}^{n_1}\ldots \sum_{i_d=1}^{n_d}\alpha_{i_1\ldots i_d}\; \ebold_{1,i_i} \otimes  \ldots \otimes \ebold_{d,i_d}
\end{equation}
Let
\[
 \rbold = (r_1, \ldots,r_d) \in \N^d \qquad \mbox{with} \quad r_\mu \leq n_\mu
\]
The truncation of (\ref{eq:tucker:basis}) at $\rbold$ is
\begin{equation}\label{eq:tucker:rank_r}
 \widehat{\Abold}_{\rbold} = \sum_{i_1 \leq r_1}\ldots \sum_{i_d \leq r_d}\alpha_{i_1\ldots i_d}\, \ebold_{1,i_1} \otimes  \ldots \otimes \ebold_{d,i_d}
\end{equation}
$\Abold$ has \kw{Tucker rank}\index{rank!Tucker} or \kw{multilinear rank}\index{rank!multilinear} $  \rbold$ if the decomposition (\ref{eq:tucker:rank_r}) is exact, i.e. there exists a set of basis $(\EC_\mu)_\mu$ such that 
\begin{equation*}
 \Abold= \sum_{i_1=1}^{r_1}\ldots \sum_{i_d=1}^{r_d}\alpha_{i_1\ldots i_d}\, \ebold_{1,i_1} \otimes  \ldots \otimes \ebold_{d,i_d}
\end{equation*}
Tucker best rank $\rbold$ approximation of $\Abold$ is finding a tensor $\Abold_\rbold$ of rank $\rbold$ such that $\|\Abold - \Abold_\rbold\|$ is minimal. Having in mind that $(\ebold_{\mu,i_\mu})_{i_\mu \leq r_\mu}$ defines a vector subspace $E'_\mu \subset E_\mu$ of dimension $r_\mu$, finding best rank $\rbold$ approximation of a tensor $\Abold$ can be reformulated as \\
\begin{center}
 \ovalbox{
 \begin{tabular}{lll}
  given && $E_1, \ldots, E_\mu, \ldots, E_d$ \\
  && with $\dim E_\mu = n_\mu$ \\
  && $\Abold \in E_1 \otimes \ldots \otimes E_d$ \\
  && ranks $r_\mu \leq n_\mu$ for $\mu \in \llbracket 1,d \rrbracket$ \\
  &&\\
  find && $d$ subspaces $E'_\mu \subset E_\mu$ \\
  && with $\dim E'_\mu=r_\mu$ \\
  &&\\
  such that && $\|\Abold - \widehat{\Abold}_\rbold\|$ is minimal \\
  && with $\widehat{\Abold}_\rbold$ being the projection of $\Abold$ on $E'_1 \otimes \ldots \otimes E'_d$
 \end{tabular}
 }
\end{center}

\nP{$\rightarrow$} There is no known algorithm which converges to best Tucker rank $\rbold$ approximation. However, both HOSVD and ALS of constraints SVD can be extended in a straightforward way to $d-$modes tensors.

%
\section{$d-$modes HOSVD}\label{sec:blratuck:hosvddmodes}
%

$d-$modes HOSVD is running $d$ SVD on matricizations of the tensor $\Abold$ on each mode, as a straightforward extension of HOSVD for three modes (this is a continuation of section \ref{sec:blratuck:hosvd3}). 

\nS Let
\[
 \Abold \in E_1 \otimes \ldots \otimes E_d \simeq \R^{n_1 \times \ldots \times n_d}
\]
Let us select a mode $\mu$ and define
\begin{equation}
 \begin{CD}
   E_\mu @>A_\mu>> E_1 \otimes E_{\mu-1} \otimes E_{\mu+1}\ldots E_d
 \end{CD}
\end{equation}
(on the r.h.s. we have $\bigotimes_\nu E_\nu$ but $E_\mu$) i.e. $A_\mu$ is a matrix with $n_1\ldots n_{\mu-1}n_{\mu+1}\ldots n_d$ rows and $n_\mu$ columns:
\[
 A_\mu \in \R^{n_1\ldots n_{\mu-1}n_{\mu+1}\ldots n_d \times n_\mu}
\]
The SVD of $A_\mu$ reads
\begin{equation}
 A_\mu = U_\mu \Sigma_\mu V_\mu^\t \qquad \mbox{with} \quad 
  \begin{cases}
   U_\mu & \in \R^{n_1\ldots n_{\mu-1}n_{\mu+1}\ldots n_d \times n_\mu} \\
   \Sigma_\mu & \in \R^{n_\mu \times n_\mu} \\
   V_\mu & \in \R^{n_\mu \times n_\mu} \\
  \end{cases}
\end{equation}
Let us denote
\[
 Y_\mu = U_\mu \Sigma_\mu
\]
Then
\begin{equation}
 A_\mu = Y_\mu V_\mu^\t
\end{equation}

\nS Let $(\vbold_{\mu 1}, \ldots, \vbold_{\mu d})$ be the columns of $V_\mu$. They form an orthonormal system. $d-$modes HOSVD at rank $\rbold$ consists in running over modes $\mu$ and for each mode
\begin{enumerate}
 \item matricize $\Abold$ along mode $\mu$ as matrix $A_\mu$
 \item run the SVD of $A_\mu$ as $A_\mu = U_\mu\Sigma_\mu V_\mu^\t$
 \item select the columns of $V_\mu$ as a new basis for $E_\mu$: $\EC_\mu = (\vbold_{\mu,1}, \ldots,\vbold_{\mu,n_\mu})$ 
 \end{enumerate}
 Then, once it has been run over all modes,
 \begin{enumerate}
 \item decomposing $\Abold$ on these basis
 \item truncating the decomposition at $\rbold$
\end{enumerate}

%
\section{$d-$modes HOOI}\label{sec:blratuck:hooidmodes}
%

HOSVD simply is a series of SVD on matrices obtained by matricization along a mode. This provides subspaces $E'_\mu$ for each mode, independently. As for 3 modes, it is possible to go one step further towards the solution of best Tucker approximation at prescribed rank, with extension of HOOI to $d$ modes, which is an ALS process with constrained SVD. This is a continuation of section \ref{sec:blratuck:hooi}.

\nB Let us assume the state at step $t$ is 
\begin{equation}
 \Abold^t \in E_1^t \otimes \ldots \otimes E_d^t
\end{equation}
with $E^t_\mu \subset E_\mu$ and $\dim E^t_\mu=r_\mu$. We first select a mode, say $\mu$. Next step consists in 
\begin{enumerate}
 \item keeping all spaces $E^t_\nu$ fixed but $E_\mu^t=E_\mu$ which is the new unknown (it will be given by its basis)
 \item matricize $\Abold^t$ along mode $\mu$ as $A^t_\mu$
 \item run the SVD of $A^t_\mu$ as $A^t_\mu=U_t\Sigma_tV_t^\t$
 \item find the new basis as the columns of $V_t$,
 \item truncate the basis at the $r_\mu$ first vectors, which will span $E^{t+1}_\mu$
\end{enumerate}
i.e. optimize on subset $E'_\mu$ knowing all other subspaces $E'_\nu$ for $\nu \neq \mu$. Such a process is run over all modes. It is a step. HOOI consists in iterating this step until a stopping condition is met.

%
\chapter{Best low rank approximation of a tensor: Tensor Train}\label{chap:tt}
%

Let us have a $d-$ modes tensor
\[
 \Abold \in \bigotimes_\mu E_\mu, \qquad \mbox{with} \quad 
 \left\{
 \begin{array}{rcl}
   \dim E_\mu &=& n_\mu \\
   \mu & \in & \llbracket 1,d \rrbracket
 \end{array}
 \right.
\]
It can be written component wise (see notations in section \ref{sec:tens_prod:notations})
\[
 \begin{array}{lcl}
  \Abold &=& \displaystyle \sum_{\ibold \in \IC} \, a_\ibold \, \Ebold_\ibold \\
  &=& \displaystyle \sum_{i_1}\ldots \sum_{i_d}\, a_{i_1\ldots i_d}\, \ebold_{i_1}^{(1)} \otimes \ldots \otimes \ebold_{i_d}^{(d)}
 \end{array}
\]
In statistical physics or computational statistics, a tensor may represent a joint law, like
\[
a_{i_1\ldots i_d} = \P(X_1=i_1, \ldots,X_d = i_d)
\]
where the $X_\mu$ are random variable taking their values in a discrete set $\llbracket 1, n_\mu \rrbracket$. A standard issue is to compute some marginals of the joint law, among which the partition function or normalisation constant $Z$ defined by
\[
 \begin{array}{lcl}
  Z(\Abold) &=& \displaystyle \sum_{\ibold \in \IC} \, a_\ibold\\
  &=& \displaystyle \sum_{i_1}\ldots \sum_{i_d}\, a_{i_1\ldots i_d}
 \end{array}
\]
which requires $n^d$ sums with brute force (if each $n_\mu$ is equal to $n$). If $\Abold$ has rank one, i.e.
\[
 \Abold = \xbold_1 \otimes \ldots \otimes \xbold_d, \qquad \xbold_\mu \in E_\mu
\]
with
\[
 \xbold_\mu = \left(x_{1}^{(\mu)}, \ldots, x_{n_\mu}^{(\mu)}\right)
\]
then
\[
 \begin{array}{lcl}
  Z(\Abold) &=& \displaystyle \sum_{i_1}\ldots \sum_{i_d}\, x_{i_1}^{(1)} \ldots x_{i_d}^{(d)} \\
  &=& \displaystyle \left(\sum_{i_1} x_{i_1}^{(1)}\right) \ldots \left(\sum_{i_d} x_{i_d}^{(d)}\right)
 \end{array}
\]
which requires $dn$ sums and $d$ products (approximately), as for example for $d=3$
\[
\sum_{i,j,k} \, x_iy_jz_k = \left(\sum_i x_i\right)\left(\sum_j y_j\right)\left(\sum_k z_k\right)
\]
Marginalization is easy because variables are separate. Separation of variables is a \emph{grail} in many domains of mathematics, notably marginalisation, but as well for building analytical and numerical solutions of ODEs and PDEs, or disentangling interactions in models in statistical physics or many-body quantum systems, among many others. 

\nB Such a calculation works because variables are separate in a rank one tensor, and distributivity of multiplication over addition can be used. Tensor-Train format (TT) is a generalization of rank one tensors suitable for marginalisation, where each term of $\Abold$ is no longer a product of scalars, one per mode as in $a_{ijk}=x_iy_jz_k$, but a product of matrices, one per mode too. Then, variables (modes) are separate too, and distributivity of multiplication over addition works as well in the ring of matrices and can be used. For example, if we have, with $d=3$
\[
\forall \: i,j,k, \quad a_{ijk} = \ubold_i.G_j.\vbold_k \qquad \mbox{with} \quad 
\begin{cases}
  \ubold_i &  \in \R^{1 \times r} \\
  G_j &  \in \R^{r \times r} \\
  \vbold_k &  \in \R^{r \times 1} \\
\end{cases}
\]
then $a_{ijk} \in \R$ and
\[
\sum_{ijk} a_{ijk} = \left(\sum_i \ubold_i\right)\left(\sum_j G_j\right)\left(\sum_k \vbold_k\right)
\]

\nB This chapter is rather long (it is the longest of the book with 10 sections), and organized as follows. First, there is a series of sections on the Tensor train decomposition of a tensor, which could have been placed in the rest of section \ref{chap:rank3} but are presented here to propose a self-contained chapter on Tensor Trains. We can notice that tensor train designates both a format and a tensor: the vocabulary is not fully stabilized.
\begin{description}
\item[section \ref{sec:tt:presentation}] The TT-format of a tensor is presented first on a simple (and classical) example where, for $d=3$, coefficients $a_{ijk}$ can be decomposed as $a_{ijk}=x_i+y_j+z_k$
\item[section \ref{sec:tt:def}] gives the definition of Tensor Train format for a tensor. It is a generalization of separation of variables in a rank 1 tensor as a matrix form: each coefficient $a_{ijk}$ can be written as $a_{ijk}=u_i.G_j.v_k$ where $u_i,G_j,v_k$ are matrices.
\item[section \ref{sec:tt:bounds}] Some simple bounds of TT-rank by dimension analysis are established. We show that for a tensor in $(\R^n)^{\otimes d}$, the generic rank cannot be lower than $\sim \sqrt{n^{d-1}/d}$. Then, low TT rank tensors are scarce!
\item[section \ref{sec:tt:basis_change}] We show that the TT-rank is a tensor property, as it is conserved by a basis change.
\item[section \ref{sec:tt:elem}] The TT-decomposition of the sum $\Abold + \Bbold$ or Hadamard product $\Abold \odot \Bbold$ of two tensors are given from the TT-decomposition of $\Abold$ and $\Bbold$ only. This has far reaching consequences, as it implies that calculations involving combinations of those two operations on tensors can be done with their TT-decomposition only, and do not require to store the tensor in memory with all its coefficients. This is particularly useful when a good quality approximation of a tensor by a TT of low rank is available: when the order $d$ is large, calculations with full matrices are impossible., and can be done in TT format only.
\item[section \ref{sec:tt:cpformat}] Exact inequalities are given between the CP rank and the TT rank of a given format. It will be useful in chapter \ref{chap:quant} for showing that tensors which are a discretization on a Cartesian grid of a multivariate function have low TT rank. This is an occasion to present non unicity of TT decomposition beyond obvious multiplications by scalars.
\end{description}
Second, there is a series of sections dedicated to best TT-rank$-r$ approximation of a given tensor:
\begin{description}
\item[section \ref{sec:tt:svd}] This is a key section for usefulness of TT-format: an algorithm, proposed in \cite{Oseledets2009,Oseledets2011}, is presented for computing the best low TT-rank approximation of a tensor for prescribed TT-rank. This algorithm is fast and robust, and called TT-SVD. It is presented first on tensors in $E \otimes F \otimes G$ for sake of clarity and to avoid as much as possible technicalities.
\item[section \ref{sec:tt:best}] TT-SVD is presented in  its full generality, for $d-$modes tensors. It is the same algorithm that in section \ref{sec:tt:svd} to find the best low rank TT-approximation of a given tensor by alternating SVD on matricized forms of the tensor, in a way not so different from HOSVD, leading to separation of variables. Its main advantage is that storing a tensor $\Abold \in (\R^n)^{\otimes d}$ requires $n^d$ terms, whereas its best TT rank $r$ approximation requires $ndr^2$ terms only.
\item[section \ref{sec:tt:rebuild}] In this very short section, it is shown how to reconstruct a tensor $\Abold$ knowing its TT-decomposition, a useful operation often omitted in Tensor Train presentations.
\item[section \ref{sec:tt:qual}] It presents an estimation of the quality of the TT approximation knowing the quality of each SVD involved in TT-SVD algorithm. The take home message is that there is no free lunch: the dramatic decrease in required storage capacity (from $n^d$ to $dr^2n$) is paid by a dramatic decrease in approximation quality too (in $\O(\theta^d)$ if $\theta$ is the quality of each SVD, to keep it simple). 
\item[section \ref{sec:tt:diagram}] Finally, to finish on a visual understanding of TT, the sketch of Tensor Trains with diagrams is given.
\end{description}

\nT{Notations} Here are, beyon,d notations common to the whole book,  some notations for this chapter:\\
\\
\begin{center}
    \ovalbox{ 
       \begin{tabular}{ll}
       $f,g,h$ & maps $\llbracket 1, n\rrbracket \longrightarrow \R$ \\
       $(G(j))_j$ & family of matrices in $\R^{r \times r}$\\
       $r$ & TT rank of a tensor \\
       $\rbold = (r_1,\ldots,r_d)$ & multilinear TT rank of a tensor \\
       $(\ubold(i))_i$ & family of vectors in $\R^{1 \times r}$ \\
       $(\ubold,G,\vbold)_{\t\t}$ & TT decomposition of a tensor in $E \otimes F \otimes G$ \\
       $(\ubold, G_2, \ldots, G_{d-2},\vbold)_{\t\t}$ & TT decomposition of a tensor in $E_1 \otimes \ldots \otimes E_d$ \\
        $(\vbold(k))_k$ & family of vectors in $\R^{r \times 1}$ \\
        $Z(\Abold)$ & partition function of a tensor $\left(\sum_\ibold a_\ibold\right)$
       \end{tabular}
    }
\end{center}

%
\section{Presentation of TT format on an example}\label{sec:tt:presentation}
%

Let us start by introducing TT format on a simple example. The rigorous definition will come later, in section \ref{sec:tt:def}.

\nS Let $\Abold \in E \otimes F \otimes G$ with general term $a_{ijk}$ such that
\begin{equation}
 a_{ijk} = f(i) + g(j) + h(k)
\end{equation}
with
\[
 \begin{CD}
  f,g,h \: : \: \{1,\ldots,n\} @>>> \R
 \end{CD}
\]
where $\dim E = \dim F = \dim G = n$ for sake of simplicity. Then, one can write 
\begin{equation}
 \forall \: i,j,k, \quad a_{ijk} = 
 \begin{pmatrix}
   f(i) & 1
 \end{pmatrix}
 \begin{pmatrix}
  1 & 0 \\
  g(j) & 1& 
 \end{pmatrix}
 \begin{pmatrix}
   1 \\
   h(k)
 \end{pmatrix}
\end{equation}
where each coefficient $a_{ijk}$ is a product of three matrices
\begin{equation}
 a_{ijk} = 
 \underbrace{\begin{pmatrix}
   f(i) & 1
 \end{pmatrix}}_{1 \times r}
 \,
 \underbrace{\begin{pmatrix}
  1 & 0 \\
  g(j) & 1& 
 \end{pmatrix}}_{r \times r}
 \,
 \underbrace{\begin{pmatrix}
   1 \\
   h(k)
 \end{pmatrix}}_{r \times 1}
\end{equation}
with $r=2$. Let us denote
\begin{equation}
\left\{
 \begin{array}{lclcc}
   u(i) &=& \displaystyle \begin{pmatrix}
   f(i) & 1
 \end{pmatrix} & \in & \R^{1 \times r} \\
 G(j) &=& \begin{pmatrix}
  1 & 0 \\
  g(j) & 1& 
 \end{pmatrix} & \in & \R^{r \times r} \\
 v(k) &=& \begin{pmatrix}
   1 \\
   h(k)
 \end{pmatrix} & \in & \R^{r \times 1} 
 \end{array}
 \right.
\end{equation}
Then
\begin{equation}\label{eq:blra:tt:1}
 \forall \: i,j,k, \quad a_{ijk} = u(i).G(j).v(k)
\end{equation}
where $"."$ denotes a matrix product. This is an example of a tensor train format (\index{tensor train!see TT}\kwnind{TT-format}\index{TT!format}) or a  tensor-train decomposition (\kwnind{TT-decomposition}\index{TT!decomposition}) of $\Abold$ where each coefficient is a product of matrices, each matrix depending on one index only. This is an example of separation of variables $(i,j,k)$ where $a$ is read as a function on $\{1,n\}^3$
\begin{equation}
 \begin{CD}
  I \times J \times K @>a>> \R \\
  (i,j,k) @>>> u(i).G(j).v(k)
 \end{CD}
\end{equation}

\nT{Remark} Here we have made a slight but useful abuse of notations. A vector is an element of a vector space, and is represented by the list of its coordinates in a given basis. A vector can be vertical and written as a $r \times 1$ matrix, or horizontal and written as a $1 \times r$ matrix, and both should be distinguished: they are related by a transposition. However, they are not distinguished here to avoid the use of transposition, and a vector can be written horizontally of vertically. The matrix product can be sketched as (in our example, $r=2$)\\
\begin{center}
    \begin{tikzpicture}
      \node() at (-1,0){$a_{ijk} = $};
      \node (ud) at (0,0){};
      \node (uf) at (1.9,0){};
      \node () at (1,0.3){$u(i)$};
      \node () at (1,-0.3){$1 \times r$};
      \node (vd) at (4.3,2.1){};
      \node (vf) at (4.3,-0.1){};
      \node () at (4.8,1.5){$v(k)$};
      \node () at (4.9,0.5){$r \times 1$};
      \node () at (3,1.5){$G(j)$};
      \node () at (3,0.5){$r \times r$};
      \draw[thick] (ud)--(uf);
      \draw (2.1,0) rectangle (3.9,2);
      \draw[thick] (vd)--(vf);
    \end{tikzpicture}
\end{center}

\nS We then have, for the calculation of of the normalisation constant or partition function $Z(\Abold)$:
\begin{equation}
 \begin{array}{lcl}
   Z(\Abold) &=& \displaystyle \sum_{i,j,k} \, a_{ijk} \\
   &=& \displaystyle \sum_{i,j,k} \, u(i).G(j).v(k) \\
   &=& \displaystyle \underbrace{\left(\sum_iu(i)\right)}_{1 \times r} \, \underbrace{\left(\sum_jG(j)\right)}_{r \times r} \, \underbrace{\left(\sum_kv(k)\right)}_{r \times 1}
 \end{array}
\end{equation}
as there is the property of distributivity of multiplication over addition in the ring of matrices. This requires $3n$ matrix additions and two matrix multiplications, instead of $n^3$ additions. Then, tensor train format is suitable for marginalisation. The same economy is done on storage: storing $\Abold$ exactly requires $n$ terms for all $f(i)$, $g(j)$ or $h(k)$; hence $3n$ terms in total instead of $n^3$, plus the program to compute each $a_{ijk}$ knowing $f,g,h$.

\nT{Kolmogorov complexity} Let us recall that Kolmogorov complexity of a string is the ratio between the shortest program which has this string as output and stops over the length of the string itself. Random sequence have a complexity of one (this has been set as the definition of a random sequence by Kolmogrov and simultaneously by Chaitin). A tensor can be given by the set of its coefficients, which is written as a string in the memory of the computer. A TT format like
\[
\forall \: i,j,k, \quad a_{ijk} = u(i).G(j).v(k)
\]
can be considered as a short program to build all coefficients of the tensor. The length of the program is mostly in the storage of the inputs $\{(u(i))_i \cup (G(j))_j \cup (v(k))_k\}$ hence $nr + nr^2+nr=nr(r+2) \sim nr^2$. Then, when $n$ is large, the Kolmogorov complexity of $\Abold$ is about $\frac{nr^2}{n^3}=r^2/n^2$ which is very close to zero.   

\nT{Notation} If $a_{ijk} = u(i).G(j).v(k)$, I denote
\begin{equation}\label{eq:tt:presentation:notation}
    \Abold = (U,\Gbold,V)_{\t\t}
\end{equation}
where
\[
\left|
\begin{array}{clclclc}
     U & \mbox{is the } & n \times r & \mbox{matrix with } & \ubold_i & \mbox{in column} & i \\ 
     \Gbold & \mbox{is the } & n \times r \times r & \mbox{tensor with } & G_j & \mbox{in mode} & j\\ 
     V & \mbox{is the } & n \times r & \mbox{matrix with } & \vbold_k & \mbox{in column} & k \\ 
\end{array}
\right.
\]
which will be useful later in this section.

\nS To finish the example, let us note that
\begin{equation}
 \left\{
   \begin{array}{lcl}
    \displaystyle \sum_iu(i) &=& \displaystyle
    \begin{pmatrix}
     \sum_if(i),n
    \end{pmatrix} \\
     \displaystyle \sum_jG(j) &=& \displaystyle
    \begin{pmatrix}
     n & 0 \\
     \sum_jg(j) & n
    \end{pmatrix} \\
     \displaystyle \sum_kv(k) &=& \displaystyle
    \begin{pmatrix}
     n \\ \sum_kh(k)
    \end{pmatrix} \\
   \end{array}
 \right.
\end{equation}
which shows that
\begin{equation}
 \sum_{i,j,k} \left(f(i)+g(j)+h(k)\right) = n^2\left(\sum_if(i)+ \sum_jg(j) + \sum_kh(k)\right)
\end{equation}
(and not $\sum_if(i)+\sum_jg(j)+ \sum_kh(k)$).

\nS This can be extended to any tensor of the form 
\begin{equation}
 \Abold \in \bigotimes_\mu E_\mu
\end{equation}
with
\begin{equation}
 a_\ibold = \sum_\mu f_{\mu}(i_\mu)
\end{equation}
We then have
\begin{equation}
 a_{i_1\ldots i_d} = 
 \underbrace{\begin{pmatrix}
   f_1(i_1) & 1
 \end{pmatrix}}_{1 \times r}
 \,
 \underbrace{\begin{pmatrix}
  1 & 0 \\
  f_2(i_2) & 1& 
 \end{pmatrix}}_{r \times r}
 \,
 \underbrace{\begin{pmatrix}
  1 & 0 \\
  f_3(i_3) & 1& 
 \end{pmatrix}}_{r \times r}
 \, \ldots \,
 \,
 \underbrace{\begin{pmatrix}
  1 & 0 \\
  f_{d-1}(i_{d-1}) & 1& 
 \end{pmatrix}}_{r \times r}
 \underbrace{\begin{pmatrix}
   1 \\
   f_d(i_d)
 \end{pmatrix}}_{r \times 1}
\end{equation}
which can be written
\begin{equation}
 a_{i_1\ldots i_d} = u(i_1).G_2(i_2)\times \ldots \times G_{d-1}(i_{d-1}).v(i_d)
\end{equation}
and is an example of a tensor-train decomposition of a $d-$modes tensor. Let us note that marginalisation runs smoothly thanks to distributivity as
\begin{equation}
 Z(\Abold) = \displaystyle \underbrace{\left(\sum_{i_1}u(i_1)\right)}_{1 \times r} \, \underbrace{\left(\sum_{i_2}G_2(i_2)\right)}_{r \times r} \, \ldots \, \underbrace{\left(\sum_{i_{d-1}}G_{d-1}(i_{d-1})\right)}_{r \times r} \,\underbrace{\left(\sum_{i_d}v(i_d)\right)}_{r \times 1}
\end{equation}
which requires $nd$ additions instead of $n^d$. Storage as well is very parcimonious: it requires $nd$ terms, $n$ for each function $f_\mu$, and Kolmogorov complexity is $\sim d/n^{d-1}$.

\nS Marginalization is made easy for such a format with variable separation, and not limited to the normalisation constant. For example, coming back to $d=3$ and $a_{ijk}= u(i).G(j).v(k)$, the marginalisation over modes $2$ and $3$ is slicing along these modes and computing the sums of coefficients in each slice
\begin{equation}
    A_1(i) = \sum_{j,k=1}^n a_{ijk}
\end{equation}
which can be computed as
\begin{equation}
    \begin{array}{lcl}
        A_1(i) &=& \displaystyle \sum_{j,k=1}^n a_{ijk} \\
        &=& \displaystyle \sum_{j,k=1}^n u(i).G(j).v(k) \\
        &=& \displaystyle u(i).\sum_{j,k=1}^n G(j).v(k) \\
        &=& \displaystyle u(i)\left(\sum_{j=1}^n G(j)\right)\left(\sum_{k=1}^nv(k)\right) \\
    \end{array}
\end{equation}

%
\section{Definition}\label{sec:tt:def}
%

Let us recall first some vocabulary which is relevant for any compact format, like CP, Tucker for tensors. Let us phrase it for rank $r$ format of a matrix. Let $A \in \R^{m \times n}$. It is said of rank $r$ if it can be written as $A=UV^\t$ with $U \in \R^{m \times r}$ and $V \in \R^{n \times r}$. $UV^\t$ is the rank $r$ format of $A$ or rank $r$ decomposition of $A$. The best rank $r$ approximation of $A$ if $A$ is not of rank $r$ is finding a rank $r$ matrix $\widetilde{A}$ such that $\|A-\widetilde{A}\|$ is minimal. SVD is an algorithm for finding a best approximation of a matrix $A$ at rank $r$ or a decomposition if $A$ has rank $r$. This can be transposed to TT format of tensors, and is presented in the following table:\\
\marginpar{\emph{TT-format}}\index{TT!format} \marginpar{\emph{TT-decomposition}} \index{TT!decomposition} \marginpar{\emph{TT-approximation}}\index{TT!approximation} \marginpar{\emph{TT-SVD}}\index{TT!SVD}    
\begin{center}
 \begin{tabular}{l|l}
   \hline
   TT-format & a Tensor Train format of a tensor \\
   TT-decomposition & writes a given tensor $\Abold$ in a TT-format \\
   TT-approximation & finds a tensor $\Bbold$ in a given TT-format \\
   & $\qquad $which approximates a given tensor $\Abold$ at best \\
   TT-SVD & An algorithm for building a TT-decomposition or a TT-approximation \\
   \hline
 \end{tabular}
\end{center}
\nB This is developed in this section.

\nT{TT-format} Let $\Abold$ be a $d-$ modes tensor
\[
 \Abold \in \bigotimes_\mu E_\mu
\]
where
\[
\left\{
    \begin{array}{lcl}
        \mu &\in& \llbracket 1,d \rrbracket \\
        \dim E_\mu &=& n_\mu \\
        \Abold &=& (a_{i_1\ldots i_d})_{i_1\ldots i_d} \quad \mbox{with} \quad i_\mu \in \llbracket 1,n_\mu \rrbracket
    \end{array}
\right.
\]
It is said to be in \kwnind{TT-format}\index{TT!format} of rank $r$ if there exists a family $u(i_1)$ of $1 \times r$ matrix, a family $G_\mu(i_\mu)$ of $r \times r$ matrices for each $\mu$ with $2 \leq \mu \leq d-1$ and a family $v(i_d)$ of  $r \times 1$ matrices such that each coefficient $a_{i_1\ldots i_d}$ can be written as ($\times$ denotes a matrix product on r.h.s.)
\begin{equation}
 a_{i_1\ldots i_d} = u(i_1) \times G_2(i_2)\times \ldots \times G_{d-1}(i_{d-1}) \times v(i_d)
\end{equation}
$r$ is the \kwnind{TT-rank}\index{TT-rank} of $\Abold$. Sometimes, for sake of consistency, we will denote
\begin{equation}
 u(i_1) = G_1(i_1), \qquad v(i_d)=G_d(i_d)
\end{equation}
yielding
\begin{equation}\label{eq:tt:def:1}
 \begin{array}{lcl}
   a_{i_1\ldots i_d} &=& G_1(i_1) \times G_2(i_2)\times \ldots \times G_{d-1}(i_{d-1}) \times G_d(i_d) \\
   &=& \displaystyle \prod_{\mu=1}^d\, G_\mu(i_\mu)
 \end{array}
\end{equation}

\nT{Notation} I will denote
\begin{equation}
    \begin{array}{lcl}
         \Abold &=& (U, \Gbold_2, \ldots,\Gbold_{d-1},V)_{\t\t} \\
         && \mbox{or} \\
         &=& (\Gbold_1, \ldots,\Gbold_d)_{\t\t} 
    \end{array}i_1
\end{equation}
where the columns of $U$ (resp. $V$) are the $\ubold_{i_1}$ (resp. $\vbold_{i_d}$) and the slices of the $\Gbold_{\mu}$ the matrices  $(G(i_\mu))_\mu$.

\nS Let us note that such a decomposition is not unique. Let us show it with $d=3$ (the extension to any $d \geq 2$ is straightforward). Let
\[
 a_{i,j,k} = u(i).G(j).v(k)
\]
Then, if $\lambda,\mu,\nu \in \R$ with $\lambda\mu\nu=1$, we have
\begin{equation}
    a_{i,j,k} = u'(i).G'(j).v'(k) \qquad \mbox{with} \quad 
    \begin{cases}
     u'(i) &= \lambda\, u(i) \\
     G'(j) &= \mu\, G(j) \\
     v'(k) &= \nu \, v(k)
    \end{cases}
\end{equation}
Non unicity of TT-decomposition will be developed in section \ref{sec:tt:cpformat}.

\nT{TT format of matrices} Let us remark that, if $d=2$, $\Abold$ is a matrix and tensor-train format is a decomposition of a matrix as a product of two matrices of rank $r$. To see this, let 
\begin{equation*}
 A \in E \otimes F \simeq \L(F,E), \qquad \mbox{with} \quad 
 \begin{cases}
  \dim E &= m \\
  \dim F &= n
 \end{cases}
\end{equation*}
Tensor train decomposition of $A$ at rank $r$ is
\begin{equation}
 \forall \: i,j, \quad a_{ij} = u(i).v(j)
\end{equation}
with $u \in \R^{1 \times r}$ and $v \in \R^{k \times 1}$. Let us denote $U \in \R^{m \times r}$ and $W \in \R^{n \times r}$ with row $i$ of $U$ being $u(i) \in \R^r$ and row $j$ of $W$ being $v(j) \in \R^r$. Then
\begin{equation}
 A = UW^\t
\end{equation}
and TT decomposition is a rank $r$ decomposition of $A$.

\nT{Multilinear TT-rank} Let us give the component-wise expression of (\ref{eq:blra:tt:1}). It reads
\begin{equation}
 a(i,j,k) = \sum_{\alpha,\beta=1}^r\, u(i,\alpha)G(\alpha,j,\beta)v(\beta,k)
\end{equation}
which is easily extended to TT-decomposition for $d-$modes tensors, from (\ref{eq:tt:def:1})
\begin{multline}
 \forall \: i_1,\ldots,i_d, \qquad a_{i_1\ldots i_d} = \sum_{\alpha_1=1}^r \ldots \sum_{\alpha_{d-1}=1}^r \; G_1(i_1,\alpha_1)G_2(\alpha_1,i_2,\alpha_2)\ldots \\
 \ldots G_\mu(\alpha_{\mu-1},i_\mu,\alpha_\mu) \ldots \ldots G_{d-1}(\alpha_{d-2},i_{d-1},\alpha_{d-1})G_d(\alpha_{d-1},i_d)
\end{multline}
This shows that the rank $r$ may vary between modes and needs not to be unique over all matrices, as we can define
\begin{multline}
 \forall \: i_1,\ldots,i_d, \qquad a_{i_1\ldots i_d} = \sum_{\alpha_1=1}^{r_1} \ldots \sum_{\alpha_{d-1}=1}^{r_{d-1}} \; G_1(i_1,\alpha_1)G_2(\alpha_1,i_2,\alpha_2)\ldots \\
 \ldots G_\mu(\alpha_{\mu-1},i_\mu,\alpha_\mu) \ldots \ldots G_{d-1}(\alpha_{d-2},i_{d-1},\alpha_{d-1})G_d(\alpha_{d-1},i_d)
\end{multline}
The tuple
\[
\rbold = (r_1,\ldots,r_{d-1}) \in \N^{d-1}
\]
is called the \kw{TT multilinear rank}\index{rank!TT multilinear}\index{TT!multilinear rank} of $\Abold$ and is denoted $\rbold_{\t\t}(\Abold)$. However, for sake of simplicity of notations, we will keep here the choice $r_1=\ldots=r_{d-1}=r$. Extension of results to a more general family $\rbold=(r_1,\ldots,r_{d-1})$ will be straightforward.

%
\section{Bounds for TT-ranks}\label{sec:tt:bounds}
%

\nS Let us fix for sake of simplicity $\Abold \in \R^{n \times n \times n}$. This limitation will be dropped later. Let $r_{\t\t}^{\mathrm{max}}$ be the maximum TT-rank of a generic tensor $\Abold$, i.e. the smallest integer satisfying to
\begin{equation}
 \forall \:\: \Abold \in \R^{n \times n \times n}, \quad r_{\t\t}(\Abold) \leq r_{\t\t}^{\mathrm{max}}
\end{equation}
If
\begin{equation}
    \Abold = (u,G,v)_{\t\t}
\end{equation}
we have
\begin{equation}\label{eq:tt:cprank:1}
    \forall \: i,j,k,\quad a_{ijk} = u_i.G_j.v_k \qquad \mathrm{with} \quad 
    \begin{cases}
      i,j,k & \in \{1,\ldots,n\} \\
      u_i & \in \R^r \\
      G_j & \in \R^{r \times r} \\
      v_k & \in \R^r
    \end{cases}
\end{equation}
So, TT-decomposition requires $nr$ terms for $(u_i)_i$ and $(v_k)_k$, and $nr^2$ for $(G_j)_j$, or $R=2nr+nr^2=nr(r+2)$ in total. We can improve this by setting $\|u\|=\|G\|=\|v\|=1$ and multiplying the expression by  scalar. But here, we are interested by the order of magnitude per dimension $r,n$ and this will be omitted. Let us denote by $N=n^3$. The set $(u,G,v)_{\t\t} := ((u_i)_i, (G_j)_j,(v_k)_k)$ is a set of $R$ coordinates. Then, equation (\ref{eq:tt:cprank:1}) defines a continuous surjective map
\begin{equation}
    \begin{CD}
      \R^R @>\phi_{\t\t}>> \R^N \\
      (u,G,v)_{\t\t} @>>> \Abold
    \end{CD}
\end{equation}
which implies $R \geq N$, or $nr(r+2) \geq n^3$. If we simplify it as $nr^2 \geq n^3$, it yields
\begin{equation}
    r_{\t\t}^{\mathrm{max}} = n
\end{equation}
(the value without this simplification is $\sqrt{n^2+1}-1$). So, the expected generic TT-rank of a $3-$modes tensor in $\R^{n \times n \times n}$ scales with $n$. When $n$ is large, the set of tensors of low rank is negligible (having low TT-rank is a rare event, even when $d=3$). 

\nT{Extension to $d-$modes tensors} This can be extended to $d$ modes tensors, keeping the simplification that if $\Abold \in \bigotimes_\mu E_\mu$, we have $\dim E_\mu = n$ for any $\mu$, or $\Abold \in (\R^n)^{\otimes d} \equiv \R^{n \times \ldots n}$.  A tensor $\Abold \in (\R^n)^{\otimes d}$ has $n^d$ coefficients. A rank $r$ TT-decomposition requires about $ndr^2$ terms ($r^2$ terms for any matrix $G(i_\mu)$, $d$ modes $\mu$, $n$ such matrices for each mode $\mu$). We then have, by a simple dimension analysis
\begin{equation}
 ndr_{\t\t}^2 \sim n^d
\end{equation}
or
\begin{equation}
 r_{\t\t} \sim \sqrt{\frac{n^{d-1}}{d}}
\end{equation}
Hence, the expected generic TT-rank of a tensor $\Abold \in (\R^n)^{\otimes d}$ scales with $n^{(d-1)/2}$.

\nT{Observation} This is huge as soon as $d$ is significant. For example, if $n=100$ and $d=10$, $r_{\t\t} \simeq 3.2 \times 10^8$.

\nT{Non homogeneous dimensions and multilinear rank} Extension to $\dim E_\mu = n_\mu$ and multilinear TT-rank $\rbold = (r_1,\ldots,r_{d-1})$ is straightforward. Let $N=\prod_\mu n_\mu$ and 
\begin{equation}
    \Abold = (G_1,\ldots,G_d)_{\t\t}
\end{equation}
Let us illustrate a situation with $d=5$, $\nbold = (n_1,\ldots,n_5)$ and $\rbold = (r_1,\ldots,r_4)$. The equality
\begin{equation}
    a_{i_1i_2i_3i_4i_5} = u(i_1).G_2(i_2).G_3(i_3).G_4(i_4).v(i_5)
\end{equation}
can be sketched as the following matrix product where dimensions are shown\\
\begin{center}
    \begin{tikzpicture}
      \draw (0,0) -- (1.5,0);
      \draw (2,0) rectangle (3.5,1.5) ; 
      \draw (4,0) rectangle (5.5,1.5) ; 
      \draw (6,0) rectangle (7.5,1.5) ; 
      \draw (8,0) -- (8,1.5);
      \node () at (0.75,-0.5){$u$};
      \node () at (2.75,-0.5){$G_2$};
      \node () at (4.75,-0.5){$G_3$};
      \node () at (6.75,-0.5){$G_4$};
      \node () at (8,-0.5){$v$};
      \node () at (0.75,0.3){$r_1$};
      \node () at (2.3,0.75){$r_1$};
      \node () at (2.75,1.8){$r_2$};
      \node () at (4.3,0.75){$r_2$};
      \node () at (4.75,1.8){$r_3$};
      \node () at (6.3,0.75){$r_3$};
      \node () at (6.75,1.8){$r_4$};  
      \node () at (8.3,0.75){$r_4$};  
    \end{tikzpicture}
\end{center}
As there are $n_1$ vectors $u(i_1)$, $n_2$ matrices $G_2(i_2)$, $\ldots$, this yields that 
\begin{equation}
    R = n_1r_1 + n_2r_1r_2 + n_3r_2r_3 + n_4r_3r_4 + n_5r_4
\end{equation}
where $R$ is the number of elements in $(u,G_2,G_3,G_4,v)_{\t\t}$. More generally
\begin{equation}
    R = n_1r_1 + \sum_{\mu=2}^{d-1} \, n_\mu r_{\mu-1}r_\mu + n_dr_{d-1}
\end{equation}
and the generic rank $\rbold$ must be such that $R=N$. 

%
\section{Change of basis}\label{sec:tt:basis_change}
%

The existence of a rank $r$ TT decomposition is a property of the tensor, not of its array when decomposed on a basis: its rank is invariant by a change of basis. We show it for a $3-$modes tensor. 

\nS Let us first define a basis change. Let us have
$$
 \Abold \in E \otimes F \otimes G \quad = [a_{ijk}] \quad \mbox{in basis } \{(e_i)_i,(f_j)_j,(g_k)_k\}_{i,j,k}
 $$
 and let us define a new basis by (see section \ref{sec:dmodestens:basis})
 $$
 \begin{cases}
  e'_a &= \sum_im_{ia}e_i \\
  f'_b &= \sum_jn_{jb}f_j \\
  g'_c &= \sum_kp_{kc}g_k
 \end{cases}
 $$
 Then, the components of $\Abold$ in the new basis 
 \[
  \Abold = \sum_{abc} \, a'_{abc} \, e'_a \otimes f'_b \otimes g'_c
 \]
are expressed as
$$
a'_{abc} = \sum_{ijk}a_{ijk}m_{ia}n_{jb}p_{kc} 
$$
(see section \ref{sec:dmodestens:basis}). Let us now assume that there is a rank $r$ TT-decomposition of $\Abold$ as
\begin{equation}
 a_{ijk} = \sum_{\alpha,\beta=1}^ru_{i\alpha}g_{\alpha j\beta}v_{\beta k}
\end{equation}
Then
$$
  \begin{array}{lcl}
    a'_{abc} &=& \displaystyle \sum_{ijk}a_{ijk}m_{ia}n_{jb}p_{kc} \\
    &=& \displaystyle \sum_{i,j,k}\sum_{\alpha,\beta}u_{i\alpha}g_{\alpha j \beta}v_{\beta k}m_{ia}n_{jb}p_{kc} \\
    &=& \displaystyle \sum_{\alpha,\beta}\underbrace{\left(\sum_iu_{i\alpha}m_{ia}\right)}_{u'_{a\alpha}}\; \underbrace{\left(\sum_jg_{\alpha j \beta}y_{jb}\right)}_{g'_{\alpha b \beta}}\; \underbrace{\left(\sum_k v_{\beta k}p_{kc} \right)}_{v'_{\beta c}}\\
    &=& \displaystyle \sum_{\alpha,\beta=1}^r\;u'_{a\alpha}\, g'_{\alpha b \beta}\, v'_{\beta c} \\
    &&\\
    &=& u'_a.G'_b.v'_c, \qquad \mbox{with} \quad
    \begin{cases}
     u'_a & \in \R^{1 \times r} \\
     G'_b & \in \R^{r \times r} \\
     v'_k & \in \R^{r \times 1}
    \end{cases}
  \end{array}
$$
which shows that there is a rank $r$ TT decomposition of $\Abold$ expressed in new basis as well.

%
\section{Computations in TT-format}\label{sec:tt:elem}
%

Some elementary properties of TT decomposition are presented here. For sake of simplicity, they are presented on $3-$modes tensors, but generalization to $d-$modes tensors is straightforward, unless otherwise stated. These properties look ordinary, and cannot be qualified as shiny. They have, however, far reaching consequences. In short, they permit to work for addition and Hadamard product of tensors in the TT format only. Let us explain this. Let us have
\[
\Abold = (U_\a,G_\a,V_\a)_{\t\t}, \qquad \Bbold = (U_\b,G_\b,V_\b)_{\t\t}
\]
Then, it is possible to compute the TT format $(U_{\a+\b},G_{\a+\b},V_{\a+\b})_{\t\t}$ of $\Abold + \Bbold$ from the TT format of $\Abold$ and $\Bbold$ only, without requiring any of the coefficients of $\Abold$ or $\Bbold$. If $d$ is significant, it makes it possible to compute the sum of two tensors with following procedure:
\begin{center}
    \begin{tikzpicture}
      \node (AB) at (0,2) {$\Abold,\Bbold$};
      \node (TTB) at (0,0){$(U_\b,G_\b,V_\b)$};
      \node (TTA) at (0,0.5){$(U_\a,G_\a,V_\a)$};
      \node (TTAB) at (6,0.25){$(U_{\a+\b},G_{\a+\b},V_{\a+\b})$};
      \node (*) at (1,0.25){};
      \node (A+B) at (6,2){$\Abold + \Bbold$};
      \draw[->] (AB) -- (TTA) ;
      \draw[->] (*) -- (TTAB) ;
      \draw[->] (TTAB) -- (A+B) ;
      \draw[->,dashed] (AB) -- (A+B) ; 
    \end{tikzpicture}
\end{center}
The diagram is commutative. The same can be done for Hadamard product.

\nT{TT rank of the sum of two tensors} Here we show that
\begin{equation}
 \left.
   \begin{array}{lcl}
      \ttrank \Abold &=& r \\
      \ttrank \Abold' &=& r'
   \end{array}
 \right\}
 \quad \Longrightarrow \quad 
 \ttrank (\Abold + \Abold') \leq r+r'
\end{equation}

\begin{proof} Let $\Abold,\Abold'$ be in TT format given by
 $$
  \left\{
    \begin{array}{lclcl}
     a_{ijk} &=& u(i)G(j)v(k) &\qquad & G \in \R^{r \times r}\\
     a'_{ijk} &=& u'(i)G'(j)v'(k) & \qquad & G' \in \R^{r' \times r'}
    \end{array}
  \right.
  $$
  If $\Bbold = \Abold + \Abold'$, we have
  $$
  b_{ijk} = a_{ijk} + a'_{ijk}
  $$
  It is easy to check that
  $$
  b_{ijk} = u_s(i)G_s(j)v_s(k)
  $$
  with
  $$
  \left\{
  \begin{array}{lclcl}
    u_s(i) &=& 
    \begin{pmatrix}
    u(i) & u'(i) \\     
  \end{pmatrix} & \quad & \in \R^{r+r'} \\
  G_s(j) &=& 
  \begin{pmatrix}
   G(i) & 0 \\
   0 & G'(j)
  \end{pmatrix} &
  \quad & \in \R^{(r+r') \times (r+r')} \\
  v_s(k) &=& 
  \begin{pmatrix}
   v(k) \\
   v'(k)
  \end{pmatrix} & \quad & \in \R^{r+r'} 
  \end{array}
  \right.
  $$
  hence the result
\end{proof}

\nS Let us give here a hint for extension to $d-$modes tensors from the example of $4-$modes tensors. Let
\begin{equation}
\left\{
    \begin{array}{lcl}
     a_{ijk\ell} &=& u(i)G(j)G(k)v(\ell) \\
     a'_{ijk\ell} &=& u'(i)G'(j)G'(k)v'(\ell) 
     \end{array}
  \right. 
\end{equation}
Then
\begin{equation}
  b_{ijk\ell} = u_s(i)G_s(j)G_s(k)v_s(\ell)
\end{equation}
with
\begin{equation}
 \left\{
  \begin{array}{lclcl}
    u_s(i) &=& 
    \begin{pmatrix}
    u(i) & u'(i) \\     
  \end{pmatrix} & \quad & \in \R^{r+r'} \\
  G_s(j) &=& 
  \begin{pmatrix}
   G(i) & 0 \\
   0 & G'(j)
  \end{pmatrix} &
  \quad & \in \R^{(r+r') \times (r+r')} \\
  G_s(k) &=& 
  \begin{pmatrix}
   G(k) & 0 \\
   0 & G'(k)
  \end{pmatrix} &
  \quad & \in \R^{(r+r') \times (r+r')} \\
  v_s(\ell) &=& 
  \begin{pmatrix}
   v(\ell) \\
   v'(\ell)
  \end{pmatrix} & \quad & \in \R^{r+r'} 
  \end{array}
  \right.
\end{equation}
We then observe therefore that
\begin{equation}
\begin{pmatrix}
   G(i) & 0 \\
   0 & G'(j)
  \end{pmatrix}
  \begin{pmatrix}
   G(k) & 0 \\
   0 & G'(k)
  \end{pmatrix}
  =
  \begin{pmatrix}
   G(i)G(k) & 0 \\
   0 & G'(i)G'(k)
  \end{pmatrix}   
\end{equation}
Moreover, this proves that
\begin{equation}
 \left. 
    \begin{array}{lcl}
      \Abold &=& (U,\Gbold,V)_{\t\t} \\
      \Abold' &=& (U',\Gbold',V')_{\t\t} \\
    \end{array}
 \right\}
 \quad \Longrightarrow \quad 
 \Abold+\Abold' = (U\mid U' \; , \; \Gbold \| \Gbold' \; , \; V \mid V')_{\t\t}
\end{equation}
where $\mid$ is column-wise concatenation visualized here:\\
\\
\begin{center}
 \begin{tikzpicture}
   \draw (0,0) rectangle (3,3) ;
   \node () at (3.5,1.5){$=$};
   \draw (4,0) rectangle (7,3) ;
   \draw (5.5,3)--(5.5,0);
   \node () at (1.5,1.5){$U(\Abold+\Abold')$};
   \node () at (4.75,1.5) {$U(\Abold)$};
   \node () at (6.25,1.5) {$U(\Abold')$};
 \end{tikzpicture}
\end{center}
and $\Gbold\|\Gbold'$ is defined in the proof just above. This has an important consequence: addition of two tensors can be done without access to their components, but within their TT-format. For $3-$modes tensors, the advantage is not obvious. But $d-$modes tensors with $d$ large have $n^d$ elements, which beyond a certain value of $d$ for a given $n$ cannot even be stored. Storing tensors in their TT-format when they have low rank $r$ requires much less memory space. The addition of two tensors can be achieved within their TT-format only, without knowing all components, as
\begin{equation}
 (U,\Gbold,V)_{\t\t} + (U',\Gbold',V')_{\t\t}= (U\mid U' \; , \; \Gbold \| \Gbold' \; , \; V \mid V')_{\t\t}
\end{equation}


\nT{TT-rank of the Hadamard product of two tensors} Let $\Abold, \Abold'$ be two $3-$modes tensors with
\begin{equation*}
 \left\{
    \begin{array}{lclcl}
      a_{ijk} &=& u_i.G_j.v_k & \quad & \rank r\\
      a'_{ijk} &=& u'_i.G'_j.v'_k & \quad & \rank r'\\
    \end{array}
 \right.
\end{equation*}
Let us recall that $\Abold \odot \Abold'= (a_{ijk}a'_{ijk})_{i,j,k}$, i.e. Hadamard product is component-wise product. We then have
\begin{equation}
 \ttrank \Abold \odot \Abold' \leq (\ttrank \Abold)(\ttrank \Abold')
\end{equation}
which is similar to a result for the rank of Hadamard product between matrices
\begin{equation}
 \rank (A \odot B) \leq (\rank A)(\rank B)
\end{equation}

\begin{proof}
A direct proof can be given by working with indices, but we give here a component free demonstration. Let us recall first that, if $x, y \in \R$, then 
\begin{equation}\label{eq:tt:elem:1}
xy = x \otimes_\k y 
\end{equation}
i.e. a product between two scalars can be read as their Kronecker product (a scalar is a vector in a space of dimension 1). A second preliminary result is a classical formula involving matrix product and Kronecker product:
\begin{equation}
  (AB) \otimes_\k (CD) = (A \otimes_\k C)(B \otimes_\k D)
\end{equation}
(this is the \emph{mixed product property}, see section \ref{sec:alstruct:prodmat}) which can be extended to product of more than two matrices as in
\begin{equation}\label{eq:tt:elem:2}
  (ABC) \otimes_\k (DEF) = (A \otimes_\k D)(B \otimes_\k E)(C \otimes_\k F)
\end{equation}
provided dimensions are consistent for matrix product, as
\begin{equation}
 \begin{array}{lcl}
  (ABC) \otimes_\k (DEF) &=& A(BC) \otimes_\k D(EF) \\
  &=& (A \otimes_\k D)(BC \otimes_\k EF) \\
   &=& (A \otimes_\k D)(B \otimes_\k E)(C \otimes_\k F)
 \end{array}
\end{equation}
Then
\begin{equation}
\begin{array}{lclcl}
a_{ijk}a'_{ijk} &=& [u(i).G(j).v(k)].[u'(i).G'(j).v'(k)] &&\\
&=& [u(i).G(j).v(k)] \otimes_\k [u'(i).G'(j).v'(k)] && \mbox{by } \ref{eq:tt:elem:1}\\
&=& [u(i) \otimes_\k u'(i)].[G(j) \otimes_\k G'(j)].[v(k) \otimes_\k v'(k)] & \qquad & \mbox{by }\ref{eq:tt:elem:2} \\
\end{array}
\end{equation}
which can be written
\begin{equation}\label{tt:svd:equation7}
 a_{ijk}a'_{ijk} = u_\odot(i).G_\odot(j).v_\odot(k) \qquad \mbox{with} \quad 
 \begin{cases}
   u_\odot(i) &= u(i) \otimes_\k u'(i) \\
   G_\odot(j) &= G(j) \otimes_\k G'(j) \\
   v_\odot(k) &= v(k) \otimes_\k v'(k)
 \end{cases}
\end{equation}
We then observe that
\begin{equation}
\left|
\begin{array}{ccl}
u(i) \otimes_\k u'(i) & \mbox{is} & 1 \times rr' \\
G(j) \otimes_\k G'(j) & \mbox{is} & rr' \times rr' \\
v(k) \otimes_\k v'(k) & \mbox{is} & rr' \times 1 \\
\end{array}
\right.
\end{equation}
which yields the result. Such a small calculation which avoids the manipulation of indices can easily be extended to Hadamard product of tensors of any order $d$.

\end{proof}

\noindent As for addition, this shows that Hadamard product can be computed within the framework of TT-formats, without computing the components. Indeed, we have
\begin{equation}
 (U,\Gbold,V)_{\t\t} \odot (U',\Gbold',V')_{\t\t} = (U \star U',\Gbold \star \Gbold',V \star V')_{\t\t}
\end{equation}
where $\star$ is defined for $U,V$ or $\Gbold$ in equation (\ref{tt:svd:equation7}). Let us note that $\star$ is the Khatri-Rao product (see section \ref{sec:alstruct:prodmat}).

\nT{Rounding} However, the calculation of Hadamard product of $m$ tensors of rank $r$ each leads to manipulation of tensors of rank $r^m$: in decomposition $(U,\Gbold,V)$ of their Hadamard product, $U$ and $V$ will have $r^m$ columns, and $\Gbold \in \R^{r^m \times n \times r^m}$. Hence, TT-format does not scale linearly with the number of operands. A possible solution to address this issue is to use the so called \kw{rounding}. It is a way to build a low TT rank approximation of a tensor when this tensor itself is given in TT format. As an example, if we have, say, two tensors $\Abold,\Bbold$ in TT format of rank $r$, we can compute $\Cbold= \Abold \odot \Bbold$ in TT format. It will be a tensor of rank $r^2$. Applying the rounding on $\Cbold$ which is in TT format yields a tensor $\Cbold'$ in TT format, of rank $r$ if wished. Rounding is an essential tool to control the growth of the rank in operations in TT format.

%
\section{Links with CP format}\label{sec:tt:cpformat}
%

Here, when useful, $r_{\c\p}$ denotes the CP-rank of a tensor and $r_{\t\t}$ its TT-rank. Otherwise, $r$ may denote one or the other according to the context. $\rank \Abold$ denotes the CP-rank of $\Abold$ and $\mathrm{TT-}\rank \Abold$ its TT-rank.

\nB We first establish a link between the CP-rank and the TT rank of a tensor as
\[
r_{\c\p}(\Abold) \leq r_{\t\t}^2(\Abold)
\]
and show next that if $(U,G,V)_{\t\t}$ is the TT decomposition of $\Abold$ and the slices $G_j = G(:,j,:)$ commute, then $r_{\c\p}(\Abold) \leq r_{\t\t}(\Abold)$. For this, we exploit the fact that the decomposition $\Abold = (U,G,V)_{\t\t}$ is not unique. Finally, we show that
\[
r_{\t\t}(\Abold) \leq r_{\c\p}(\Abold)
\]
by exhibiting a rank $r$ TT-decomposition of $\Abold$ if its CP rank is $r$. This shows the main result of this section
\[
r_{\t\t}(\Abold) \leq r_{\c\p}(\Abold) \leq r_{\t\t}^2(\Abold)
\]

\nT{Links between CP-rank and TT-rank} Let us recall that the CP-rank of a tensor $\Abold$ is the smallest integer $r$ such that $\Abold$ can be decomposed as a sum of $r$ rank one tensors, or
\begin{equation}
 \Abold = \sum_{a=1}^r\, \xbold_a \otimes \ybold_a \otimes \zbold_a
\end{equation}
for a $3-$modes tensor (see section \ref{sec:rank3:def}). Let $\Abold$ be a $3-$modes tensor with TT-rank $r$. Then, its CP-rank is not larger than $r^2$.
\begin{proof}
Indeed, let 
\begin{equation}
 \Abold = \sum_{ijk}\, a_{ijk}\, \ebold_i \otimes \fbold_j \otimes \gbold_k
\end{equation}
with
\begin{equation}
 a_{ijk} = \sum_{\alpha,\beta=1}^r\, u_{i\alpha}g_{\alpha j\beta}v_{\beta k}
\end{equation}
Then,
\begin{equation}\label{eq:tt:def:2}
 \begin{array}{lcl}
   \Abold &=& \displaystyle \sum_{i,j,k}\, a_{ijk}\, \textbf{e}_i \otimes \textbf{f}_j \otimes \textbf{g}_k \\
   &=& \displaystyle \sum_{i,j,k}\sum_{\alpha,\beta=1}^r\, u_{i\alpha}g_{\alpha j\beta}v_{\beta k}\, \textbf{e}_i \otimes \textbf{f}_j \otimes \textbf{g}_k \\
   &=& \displaystyle \sum_{i,j,k}\sum_{\alpha,\beta=1}^r\, u_{i\alpha} \textbf{e}_i \otimes g_{\alpha j\beta}\textbf{f}_j \otimes v_{\beta k}\textbf{g}_k \\
   &=& \displaystyle \sum_{\alpha,\beta=1}^r \underbrace{\left(\sum_iu_{i\alpha}\textbf{e}_i\right)}_{\textbf{u'}_\alpha} \otimes \underbrace{\left(\sum_jg_{\alpha j\beta} \textbf{f}_j\right)}_{\textbf{g'}_{\alpha\beta}} \otimes \underbrace{\left(\sum_k v_{\beta k}\textbf{g}_k\right)}_{\textbf{v'}_\beta} \\
   &=& \displaystyle \sum_{\alpha,\beta=1}^r \, \textbf{u'}_\alpha \otimes \textbf{g'}_{\alpha\beta} \otimes \textbf{v'}_\beta
 \end{array}
\end{equation}
which shows that
\[
 \rank \Abold \leq r^2
\]
\end{proof}

\nT{Extension to $d-$modes tensors} Calculation (\ref{eq:tt:def:2}) can be extended to $d-$modes tensors with $d>3$. Let us assume that 
\begin{equation}
 \Abold \in E_1 \otimes \ldots \otimes E_d \qquad \mbox{with} \quad \Abold = \sum_{i_1} \ldots \sum_{i_d} \, a_{i_1\ldots i_d}\, \textbf{e}_{1i_1} \otimes \ldots \otimes \textbf{e}_{di_d} 
\end{equation}
and
\begin{equation}
  r_{\t\t}(\Abold) = r
\end{equation}
Then,
\begin{equation}
 \begin{array}{lcl}
   \Abold &=& \displaystyle \sum_{i_1 \ldots i_d} a_{i_1\ldots i_d}\, \textbf{e}_{1i_1} \otimes \ldots \otimes \textbf{e}_{di_d} \\
   &=& \displaystyle \sum_{i_1 \ldots i_d} \:\: \sum_{\alpha_1 \ldots \alpha_{d-1}=1}^r \, u_{i_1\alpha_1}g_{\alpha_1 i_2 \alpha_2}^{(2)} \ldots g_{\alpha_{d-2} i_{d-1} \alpha_{d-1}}^{(d-1)}v_{\alpha_{d-1}i_d}\, \textbf{e}_{1i_1} \otimes \ldots \otimes \textbf{e}_{di_d} \\
   &=& \displaystyle \sum_{\alpha_1 \ldots \alpha_{d-1}=1}^r \left(\sum_{i_1}u_{i_1\alpha_1}\textbf{e}_{1i_1}\right) \otimes \ldots \otimes \left(\sum_{i_\mu}g_{\alpha_{\mu-1}i_\mu \alpha_\mu}^{(\mu)}\textbf{e}_{\mu i_\mu}\right) \otimes \ldots \otimes \left(\sum_{i_d}v_{\alpha_{d-1}i_d} \textbf{e}_{di_d}\right)\\
   &=& \displaystyle \sum_{\alpha_1 \ldots \alpha_{d-1}=1}^r \, \textbf{u'}_{\alpha_1} \otimes \ldots \otimes \textbf{g'}_{\alpha_{\mu-1}\alpha_{\mu}}^{(\mu)} \otimes \ldots \otimes \textbf{v'}_{\alpha_d}
 \end{array}
\end{equation}
from which it can be seen that
\begin{equation}
 \rank \Abold \leq r^{d-1}
\end{equation}

\nS Let us note that if $r=1$, the tensors with TT-rank and CP-rank equal to one are the same. In other words
\begin{equation}
 \mathbb{TT}_1\left((\R^n)^{\otimes 3}\right) = \mathbb{CP}_1\left((\R^n)^{\otimes 3}\right)
\end{equation}

\nT{Non unicity of the decomposition} Let us start with an analogy with matrices. Let $A \in \R^{m \times n}$ be a matrix with a rank $r$ decomposition
\[
 A = UW, \qquad \mbox{with} \quad
 \begin{cases}
  U & \in \R^{m \times r} \\
  W & \in \R^{r \times n}
 \end{cases}
\]
Then, letting $P \in \R^{r \times r}$ be an invertible matrix, we can write 
\begin{equation*}
 \begin{array}{lcl}
   A &=& UW \\
   &=& U\I_rW \\
   &=& U(P^{-1}P)W \\
   &=& (UP^{-1})(PW) \qquad \mbox{with} \quad
   \begin{cases}
     UP^{-1} & \in \R^{m \times r} \\
     PW & \in \R^{r \times n}
   \end{cases}
 \end{array}
\end{equation*}
Hence, $A= (UP^{-1})(PW)$ is another rank $r$ decomposition of $A$, decomposition which is not unique. Among all decompositions, the one provided by the SVD is such that if $r'<r$, the truncation of the SVD at rank $r'$ yields the best approximation of $A$ at rank $r'$. A similar situation exists for TT-decomposition and approximation.

\nS Let 
\begin{equation}
 \Abold \in E \otimes F \otimes G \qquad \mbox{with} \quad \Abold=(U,\Gbold,V)_{\t\t}
\end{equation}
We then have
\begin{equation}
 \forall \: i,j,k, \quad a_{ijk}= u(i).G(j).v(k) \qquad \mbox{with} \quad 
 \begin{cases}
   u(i) & \in \R^{1 \times r} \\
   G(j) & \in \R^{r \times r} \\
   v(k) & \in \R^{r \times 1}
 \end{cases}
\end{equation}
Let $P,Q \in \R^{r \times r}$ be two invertible matrices. This can be written
\begin{equation}
 \begin{array}{lcl}
     a_{ijk} &=& u(i).G(j).v(k) \\
     &=& u(i).\I_r.G(j).\I_r.v(k) \\
     &=& u(i).(P^{-1}P).G(j).(Q^{-1}Q).v(k) \\
     &=& (u(i).P^{-1})(P.G(j).Q^{-1})(Q.v(k)) \\
 \end{array}
\end{equation}
So, we have
\begin{equation}
 a_{ijk} = u'(i).G'(j).v'(k) \qquad \mbox{with} \quad 
 \begin{cases}
  u'(i) &= u(i).P^{-1} \\
  G'(j) &= P.G(j).Q^{-1} \\
  v'(k) &= Qv(k)
 \end{cases}
\end{equation}
which is another TT-format for $\Abold$, with same rank.

\nT{Simultaneous diagonalisable matrices} Let $(M_i)_i$ be a set of matrices in $\R^{n\times n}$. They are said simultaneously diagonalisable if $(i)$ each one is diagonalisable and $(ii)$ there exists an invertible matrix $P \in \R^{n \times n}$ such that for each matrix $M_i$, $P^{-1}M_iP$ is diagonal (see as well section \ref{sec:rank3:bound} for a similat approach). Then, if the $(M_i)_i$ are simultneously diagonalisable, they commute pairwise, i.e.
\[
 \forall \:i,j, \quad M_iM_j=M_jM_i
\]
The converse is almost true: if the $(M_i)_i$ commute pairwise, then, if each of them is diagonalisable, they are simultaneously diagonalisable (if the $(M_i)_i$ commute pairwise, it may happen that one or several of them is not diagonalisable).

\nT{Simultaneous diagonalization of core matrices} Let us consider the set of matrices $(G_j)_j$ with $G_j := \Gbold(j)$. We assume that they are simultaneously diagonalisable. Then, there exists a matrix $P \in \R^{r \times r}$ such that
\begin{equation}
 \forall \: j, \qquad G'_j=P^{-1}G_jP = \Lambda_j, \qquad \Lambda_j = \Diag (\lambda_{j1}, \ldots,\lambda_{jr})
\end{equation}
Then, we have
\begin{equation}
g'_{\alpha j \beta} = 
 \begin{cases}
   \lambda_{i\alpha} & \mbox{if } \alpha=\beta \\
   0 & \mbox{if } \alpha \neq \beta 
 \end{cases}
\end{equation}
Thus,
\begin{equation}
    \begin{array}{lcl}
      a_{ijk} &=& u'_i.\Lambda_j.v'_k \\
      &=& \displaystyle \sum_{\alpha,\beta=1}^r\, u'_{i\alpha}g_{\alpha j \beta}v'_{\alpha k} \\
      &=& \displaystyle \sum_{\alpha=1}^r\, u'_{i\alpha}\lambda_{j\alpha}v'_{\alpha k}
    \end{array}
\end{equation}
and
\begin{equation}
 \begin{array}{lcl}
   \Abold &=& \displaystyle \sum_{i,j,k}\, a_{ijk}\, \ebold_i \otimes \fbold_j \otimes \gbold_k \\
   &=& \displaystyle \sum_{i,j,k} \sum_{\alpha=1}^r \, u'_{i\alpha}\lambda_{j\alpha}v'_{\alpha k} \; \ebold_i \otimes \fbold_j \otimes \gbold_k \\
   &=& \displaystyle \sum_{\alpha=1}^r \left(\sum_{i,j,k}  u'_{i\alpha}\lambda_{j\alpha}v'_{\alpha k} \; \ebold_i \otimes \fbold_j \otimes \gbold_k \right) \\
   &=& \displaystyle \sum_{\alpha=1}^r \underbrace{\left(\sum_{i}   u'_{i\alpha}\ebold_i \right)}_{\xbold_\alpha} \otimes  \underbrace{\left(\sum_j\lambda_{j\alpha}\fbold_j\right)}_{\ybold_\alpha} \otimes \underbrace{\left(\sum_k v'_{\alpha k}\gbold_k \right)}_{\zbold_\alpha} \\
   &=& \displaystyle \sum_{\alpha=1}^r \, \xbold_\alpha \otimes \ybold_\alpha \otimes \zbold_\alpha
 \end{array}
\end{equation}
which shows that, if the cores $G_j$ are simultaneously diagonalisable,
\begin{equation}
 r_{\c\p}(\Abold) \leq r
\end{equation}

\nS The converse shows that if $r_{\c\p}(\Abold)=r$, then $r_{\t\t}(\Abold) \leq r$.
\begin{proof}
Let us have 
\begin{equation}
    \Abold = \sum_{\alpha=1}^r \, \xbold_\alpha \otimes \ybold_\alpha \otimes \zbold_\alpha
\end{equation}
Let us have
\begin{equation}
    \xbold_\alpha = \sum_i x_{i\alpha} \, \ebold_i, \qquad \ybold_\alpha = \sum_j y_{j\alpha} \, \fbold_j, \qquad \zbold_\alpha = \sum_k z_{k\alpha} \, \gbold_k, 
\end{equation}
Then,
\begin{equation}
    \forall \: i,j,k, \quad a_{ijk} = \sum_{\alpha=1}^r x_{i\alpha}y_{j\alpha}z_{k\alpha}
\end{equation}
Let us consider the matrix 
\[
G_j = \Diag (y_{j\alpha})_\alpha \in \R^{n \times n}
\]
and the vectors
\begin{equation}
    \ubold_i = (x_{i\alpha})_\alpha, \qquad \vbold_k = (z_{k\alpha})_\alpha
\end{equation}
Then
\begin{equation}
    \forall \: i,j,k, \quad a_{ijk} = \ubold_i.G_j.\vbold_k
\end{equation}
For example, if $r=2$
\begin{equation}
    \begin{array}{lcl}
         \ubold_i.G_j.\vbold_k &=& 
         \begin{pmatrix}
             x_{i1} & x_{i2}
         \end{pmatrix}
         \begin{pmatrix}
             y_{j1} & 0 \\
             0 & y_{j2}
         \end{pmatrix}
         \begin{pmatrix}
             z_{k1} \\
             z_{k2}
         \end{pmatrix} \\
         &=& \begin{pmatrix}
             x_{i1} & x_{i2}
         \end{pmatrix}
          \begin{pmatrix}
             y_{j1}z_{k1} \\
             y_{j2}z_{k2}
         \end{pmatrix} \\
         & = & x_{i1}y_{j1}z_{k1} + x_{i2}y_{j2}z_{k2}
    \end{array}
\end{equation}
which shows that
\[
r_{\c\p}(\Abold) \leq r_{\t\t}
\]
\end{proof}

%
\section{Tensor-Train format and SVD: TT-SVD}\label{sec:tt:svd}
%

SVD is a decomposition of a matrix $A$ which provides its best approximation by a matrix of rank $k$ (see section \ref{sec:inter:svd}). In a similar way, TT-SVD is an algorithm which provides the best approximation of a tensor $\Abold$ by a tensor of same dimensions with prescribed rank $r$, possibly multilinear. Let us first develop it for $\Abold \in E \otimes F \otimes G$. It is written as
\[
\begin{CD}
  \left\{ 
     \begin{array}{c}
          \Abold \in E \otimes F \otimes G \\
          r \in \N
     \end{array}
  \right\}
  @>>>
  \widetilde{\Abold} = (U, \Gbold,V)_{\t\t} \quad \st \quad 
    \left\{ 
     \begin{array}{c}
          \mathrm{TT-rank} \: (\widetilde{\Abold})=r\\
          \|\Abold - \widetilde{\Abold}\| \: \mbox{minimal}
     \end{array}
  \right\}
\end{CD}
\]

\nS There is an immediate correspondence between TT format and SVD, presented here first for $d=3$. This is a key construction, shared with HOSVD,  because the algorithm TT-SVD benefits from two qualities when compared to CP or Tucker rank (HOSVD is selected for best Tucker approximation):
\begin{itemize}[label=$\rightarrow$]
    \item the target, TT-format, is parsimonious
    \item TT-SVD is a quick and robust algorithm
\end{itemize}
Here is a rough comparison between CP, Tucker and TT best approximation at rank $r$:
\begin{center}
\begin{tabular}{c|cc}
    \hline
    & dimension & numerical stability \\
    \hline
    CP & $rnd$ & unstable \\
    HOSVD & $\sim r^d$ & stable \\
    TT & $dnr^2$ & stable \\
    \hline
\end{tabular}
\end{center}
In CP and TT, the complexity scales linearly with $n$, but approximation of CP rank can be numerically unstable (see section \ref{sec:blracp:ill}). HOSVD and TT-SVD rely on a well known, well behaved algorithm: SVD, but best approximation as a Tucker model scales exponentially with $d$. As there is no free lunch, it will be shown in section \ref{sec:tt:qual} that the price to pay is that, apart from some well understood situations (see chapter \ref{chap:quant}), the low rank TT best approximation may be very poor as far as quality is concerned. 

\nS SVD provides a decomposition and, if truncated, a best approximation. It has a counter-part for TT format, which provides a decomposition and, if truncated, a best approximation. Decomposition through TT-SVD is presented first in the case of a $3-$modes tensor, then the best approximation through truncation of TT-SVD, and finally decomposition and truncation are presented for a $d-$modes tensor in section \ref{sec:tt:best}, with a full development on the case $d=4$, which is more challenging that $d=3$. 

\nS Let us have $\Abold \in E \otimes F \otimes G$ written in TT-format
\begin{equation}
 a_{ijk} = \sum_{\alpha,\beta}u_{i\alpha}g_{\alpha j\beta}v_{\beta k}
\end{equation}
For sake of simplicity in notations, we select $\dim E=\dim F = \dim G=n$. Let 
\begin{equation}
 \left\{
    \begin{array}{lcccl}
      U \in \R^{n \times r} & \mbox{with} & U(i,\alpha) &=& u_{i\alpha} \\
      \Gbold \in \R^{r \times n \times r} & \mbox{with} & \Gbold(\alpha,j,\beta) &=& g_{\alpha i\beta} \\
      V \in \R^{n \times r} & \mbox{with} & V(k,\beta) &=& v_{k\beta}
    \end{array}
 \right.
\end{equation}
or
\begin{equation}
 \Abold = (U,\Gbold,V)_{\t\t}, \qquad \Abold(i,j,k) = \sum_{\alpha,\beta=1}^r\: U(i,\alpha)\Gbold(\alpha,j,\beta)V^\t(\beta,k)
\end{equation}

\nS It is possible to write
\begin{equation}\label{tt:svd:eq1}
 \Abold(i,j,k) = \sum_{\alpha=1}^rU(i,\alpha)\Bbold(\alpha,j,k) \qquad \mbox{with} \quad 
 \begin{cases}
  \Bbold & \in \R^{r \times n \times n} \\
  \Bbold(\alpha,j,k) & = \sum_{\beta=1}^r\Gbold(\alpha,j,\beta)V^\t(\beta,k) 
 \end{cases}
\end{equation}
Let us now consider the matricization $A \in \R^{n \times n^2}$ of $A \in E \otimes FG$ where $FG$ stands for $F \otimes G$. Then, equation (\ref{tt:svd:eq1}) means that $\rank A=r$. If $B \in \R^{r \times n^2}$ is the matricization of $\Bbold$ with $B(\alpha,jk)=\Bbold(\alpha,j,k)$, we have 
\begin{equation}\label{tt:svd:eq2}
 \underbrace{A}_{n \times n^2} = \underbrace{U}_{n \times r}.\underbrace{B}_{r \times n^2}
\end{equation}
which is a rank $r$ decomposition of $A$. An immediate possibility to obtain such a decomposition is to run the SVD of $A$, with $A=U\Sigma W^\t$ and $B=\Sigma W^\t$.

\nS Once (\ref{tt:svd:eq2}) has been obtained, let us focus on $B \in \R^{r \times n^2}$ which is the matricization of tensor $\Bbold$ defined by (see equation \ref{tt:svd:eq1})
\begin{equation}\label{tt:svd:eq3}
  \Bbold(\alpha,j,k) = \sum_{\beta=1}^r\Gbold(\alpha,j,\beta)V^\t(\beta,k) 
\end{equation}
Let us consider double index $(\alpha,j)$ as a single index $aj$, and matrix $B'\in \R^{rn \times n}$ defined by a reshaping of $B$
\begin{equation}
 B'(\alpha j,k)= B(\alpha,jk)
\end{equation}
Then, (\ref{tt:svd:eq3}) can be rewritten as
\begin{equation}\label{tt:svd:eq4}
  B'(\alpha j,k) = \sum_{\beta=1}^r G(\alpha j,\beta)V^\t(\beta,k) \qquad \mbox{with} \quad B' \in \R^{rn \times n}
\end{equation}
and $G \in \R^{rn \times r}$ being the matricization of $\Gbold$ over third mode on columns. In general, we expect from its dimensions that $\rank B'=n$. However, equation (\ref{tt:svd:eq4}) means that $\rank B'=r$. It is immediate that SVD of $B'$ yields the expected decomposition.

\nT{Result} Let $\Abold \in \R^{n \times n \times n}$, and let us assume it has a TT-decomposition of rank $r$
\begin{equation}
 \Abold = (U,\Gbold,V)_{\t\t}
\end{equation}
Then, the matricization $A$ of $\Abold$ over the first mode as rows has rank $r$, and can be written
\begin{equation}\label{tt:svd:eq5}
 A = UB
\end{equation}
$B$ itself can be reshaped as $B'$, with
\begin{equation}
 B'(\alpha j,k)=B(\alpha,jk)
\end{equation}
and decomposed as
\begin{equation}\label{tt:svd:eq6}
 B'=GV^\t
\end{equation}
Both decompositions (\ref{tt:svd:eq5}) and (\ref{tt:svd:eq6}) have been obtained by SVD of $A$ and $B'$ respectively. Here is the algorithm with tracking of the dimensions (which is crucial):\\
\\
\begin{algorithm}[H]
\begin{algorithmic}[1]
\STATE \textbf{given} a tensor $\Abold \in E \otimes F \otimes G$ and a rank $r \in \N$
\STATE \textbf{matricize} $\Abold$ as $A \in \R^{n \times n^2}$
\STATE \textbf{compute} the SVD of $A$: $A=U\Sigma W^\t$; $U \in \R^{n \times n}$
\STATE $B \leftarrow \Sigma W^\t$; $B \in \R^{n \times n^2}$
\STATE \textbf{reshape} $B' \leftarrow B$; $B' \in \R^{n^2 \times n}$
\STATE \textbf{run} the SVD of $B'$: $B'=U'\Sigma' V^\t$;
\STATE $G \leftarrow  U'\Sigma'$; $G \in \R^{n^2 \times n}$
\STATE \textbf{reshape} $\Gbold \leftarrow G$; $\Gbold \in \R^{n \times n \times n}$
\RETURN $(U,\Gbold,V)$
\end{algorithmic}
\caption{TT-SVD algorithm for the decomposition of a $3-$modes tensors}
\label{alg:mds}
\end{algorithm}

\nS This yields a full TT-decomposition of $\Abold$. One sees that it is not really compact, because $\Gbold$, like $\Abold$, is in $\R^{n \times n \times n}$. The advantage is visible when a rank $r$ best approximation is provided. It is obtained by truncating the SVDs at rank $r$, as in\\
\\
\begin{algorithm}[H]
\begin{algorithmic}[1]
\STATE \textbf{given} a tensor $\Abold \in E \otimes F \otimes G$ and a rank $r \in \N$
\STATE \textbf{matricize} $\Abold$ as $A \in \R^{n \times n^2}$
\STATE \textbf{compute} the SVD of $A$: $A=U\Sigma W^\t$; $U \in \R^{n \times n}$
\STATE \textbf{truncate} the SVD at rank $r$; $U \in \R^{n \times r}$
\STATE $B \leftarrow \Sigma W^\t$; $B \in \R^{r \times n^2}$
\STATE \textbf{reshape} $B' \leftarrow B$; $B' \in \R^{rn \times n}$
\STATE \textbf{run} the SVD of $B'$: $B'=U'\Sigma' V^\t$;
\STATE \textbf{truncate} the SVD at rank $r$; $V \in \R^{n \times r}$
\STATE $G \leftarrow  U'\Sigma'$; $G \in \R^{rn \times r}$
\STATE \textbf{reshape} $\Gbold \leftarrow G$; $\Gbold \in \R^{r \times n \times r}$
\RETURN $(U,\Gbold,V)$
\end{algorithmic}
\caption{TT-SVD algorithm for approximation of a $3-$modes tensors}
\label{alg:mds-approx}
\end{algorithm}

\nB When extended to $d-$modes tensors with $d \geq 3$, this leads to an algorithm called \texttt{TT-SVD} defined by Oseldets in \cite{Oseledets2011}.

%
\section{Best TT approximation of a $d-$modes tensor}\label{sec:tt:best}
%

As for matrix rank or tensor rank, one can be led to approximate at best a given tensor $\Abold$ by a tensor $\Bbold$ which can be written as a TT-format tensor at rank $r$, whereas $\Abold$ itself is not in TT-format. This is analogous to best low rank approximation of a matrix, and can be called best low rank TT-format approximation. Let us set this problem for $d-$modes tensors first, and show how a best approximation can be built with truncated SVD iterated on modes to separate variables one after the other.

\nT{Setting the problem}The problem of best TT-approximation can be set as follows:\\
\begin{center}
 \ovalbox{
    \begin{tabular}{ll}
      Given & a tensor $\displaystyle \Abold \in \bigotimes_\mu E_\mu$ \\
      & a multilinear rank  $\rbold = (r_1,\ldots,r_{d-1})$ \\
      &\\
      find & a family $u(i_1)$ of $1 \times r_1$ matrices \\
      & a family $G_\mu(i_\mu)$ of $r_{\mu-1} \times r_\mu$ matrices for each $\mu$ with $2 \leq \mu \leq d-1$ \\
      & a family $v(i_d)$ of  $r_{d-1} \times 1$ matrices \\
      &\\
      such that & $\|\Abold - \Bbold\|$ minimal \\
      &\\
      with & $\displaystyle \forall \: \ibold \in \IC, \quad \beta_{i_1\ldots i_d} = u(i_1).G_2(i_2)\times \ldots \times G_{d-1}(i_{d-1}).v(i_d)$
    \end{tabular}
 }
\end{center}

\nB If $\Bbold=\Abold$, it is called TT-decomposition (the decomposition is exact).

\nT{How to build a TT best approximation on an example} Let
\[
 \Abold \in \R^{m \times n \times p \times q} \equiv E \otimes F \otimes G \otimes H
\]
with general component
\[
 \Abold(i,j,k,\ell)
\]
$\Abold$ has $mnpq$ components. The guideline is to reshape $\Abold$ as a matrix, and develop SVD of such a matrix to isolate one variable. If a best approximation is looked for, the SVD is truncated at selected rank.

\nP{$\rightarrow$} As a first step let us matricize $\Abold$ as 
\begin{equation}\label{tt:best:eq:1}
 A_1 \in E \otimes (F \otimes_\k G \otimes_\k H) \simeq \L(F \otimes_\k G \otimes_\k H \, , \, E) \simeq \R^{m \times npq}
\end{equation}
(see section \ref{sec:tenselem:reshape} for matricization of tensors using Kronecker product) with $m$ rows and $npq$ columns. In general case, $m \leq npq$. This step consist in deriving the SVD of $A_1$ truncated at rank $r_1$, as
\begin{equation}
 A_1 = U_1\Sigma_1V_1^\t = U_1W_1, \qquad \mbox{with} \quad
 \begin{cases}
 W_1 & = \Sigma_1V_1^\t \\
   U_1 & \in \R^{m \times r_1} \\
   W_1 & \in \R^{r_1 \times npq}
 \end{cases}
\end{equation}
It reads component-wise
\begin{equation}\label{eq:ttxmpl1}
 \Abold(i,j,k,\ell) = \sum_{a=1}^{r_1}U_1(i,a)W_1(a, jk\ell)
\end{equation}
(see section \ref{sec:inter:svd}). One can observe that variable $i$ has been isolated from set of variables $(j,k,\ell)$. Let us iterate this process on $W_1$ to isolate variable $j$.


\nP{$\rightarrow$} Let us now reshape $W_1$ as
\begin{equation}\label{tt:best:eq:2}
 \begin{CD}
   W_1 \in \R^{r_1 \times npq} @>>> W'_1 \in \R^{r_1n \times pq}
 \end{CD}
\end{equation}
The general component of $W'_1$ is
\begin{equation}
 W'_1(aj,k\ell)
\end{equation}
The second step consists in deriving the SVD of $W'_1$ truncated at rank $r_2$ as
\begin{equation}
 \begin{array}{lcl}
  W'_1 &=& U_2\Sigma_2V_2^\t \\
  &=& U_2W_2
 \end{array}
\end{equation}
with
\begin{equation}
 \begin{cases}
   W'_1 & \in \R^{r_1n \times pq}\\
   U_2 & \in \R^{r_1n \times r_2} \\
   W_2 & \in \R^{r_2 \times pq} 
 \end{cases}
\end{equation}
It can be written component-wise
\begin{equation}\label{eq:ttxmpl2}
 W_1(a,  jkl) = W'_1(aj, kl) = \sum_{b=1}^{r_2}U_2(aj, b)W_2(b, k\ell)
\end{equation}
One can observe that we have isolated $j$ from $(k,\ell)$ in $W_1$. Then, associating (\ref{eq:ttxmpl1}) and (\ref{eq:ttxmpl2})
\begin{equation}\label{eq:ttxmpl3}
 \begin{array}{lcl}
   \Abold(i,j,k,\ell) &=& \displaystyle \sum_{a}U_1(i,a)W_1(a, jk\ell) \\
   &=& \displaystyle \sum_{a}U_1(i,a)\left(\sum_{b}U_2(aj, b)W_2(b, k\ell)\right)\\
   &=& \displaystyle \sum_{a=1}^{r_1}\sum_{b=1}^{r_2}U_1(i,a)U_2(a, j, b)W_2(b, k\ell)
 \end{array}
\end{equation}
where $i$ and $j$ have been each isolated from $(k,\ell)$.


\nP{$\rightarrow$} As a third and last step, quite similar to second step, let us reshape $W_2$ as
\begin{equation}
 \begin{CD}
   W_2 \in \R^{r_2 \times pq} @>>> W'_2 \in \R^{r_2p \times q}
 \end{CD}
\end{equation}
with general component
\begin{equation}
 W'_2(bk,\ell)
\end{equation}
Then, third step is deriving the SVD of $W'_2$ truncated at rank $r_3$ as
\begin{equation}
 \begin{array}{lcl}
  W'_2 &=& U_3\Sigma_3V_3^\t \\
  &=& U_3W_3
 \end{array}
\end{equation}
with
\begin{equation}
 \begin{cases}
   W'_2 & \in \R^{r_2p \times q}\\
   U_3 & \in \R^{r_2p \times r_3} \\
   W_3 & \in \R^{r_3 \times q} 
 \end{cases}
\end{equation}
It can be written componentwise
\begin{equation}\label{eq:ttxmpl4}
 W_2(b,  k\ell) = W'_2(bk, \ell) = \sum_{c=1}^{r_3}U_3(bk, c)W_3(c, \ell)
\end{equation}

\nP{$\rightarrow$} Hence the result: associating (\ref{eq:ttxmpl3}) and (\ref{eq:ttxmpl4})
\begin{equation}
 \begin{array}{lcl}
   \Abold(i,j,k,\ell) &=& \displaystyle \sum_{a=1}^{r_1}\sum_{b=1}^{r_2}U_1(i,a)U_2(a, j, b)W_2(b, k\ell)\\
   &=& \displaystyle \sum_{a=1}^{r_1}\sum_{b=1}^{r_2}U_1(i,a)U_2(a, j, b)\left(\sum_{c=1}^{r_3}U_3(bk, c)W_3(c, \ell)\right) \\
   &=& \displaystyle \sum_{a=1}^{r_1}\sum_{b=1}^{r_2}\sum_{c=1}^{r_3} U_1(i,a)U_2(a, j, b)U_3(b, k, c)W_3(c, \ell)
 \end{array}
\end{equation}
If SVD are not truncated, this is a\footnote{\emph{a}, not \emph{the}!} TT decomposition of $\Abold$ where all variables $(i,j,k,\ell)$ have been separated. Best low rank TT approximation at rank $\rbold$ simply is cutting the SVD's at rank $r_\mu$.

\nT{Remark on notations} Here, we have kept the classical notations for SVD ($A=U\Sigma V^\t$). The matrices $U_k$ involved in the successive SVD are the components of the TT decomposition (TTD) or best approximation (TTBA). We have kept the notation $U$. In the literature of TTD, namely in Oseledets seminal papers, they are translated as matrices denoted $G$, in the following way. $U_1 \in \R^{m \times r_1}$ is a matrix of general term $U_1(i,a)$ with\footnote{Here the simplified notation $[m]= \llbracket 1,m \rrbracket$ is adopted.} $i \in [m]$ and $a \in [r_1]$. The $i-$th row of $U_1$ is denoted $G_1(i)$. As well, $U_2 \in \R^{r_1 \times n \times r_2}$. The matrix $G_2(j)$ is defined as the matrix of components $(U_2(a,j,b))_{a,b}$. And so forth for $G_3(k)= (U_3(b,k,c))_{b,c}$. Finally, $G_4(\ell) = (U_4(c,\ell))_{c}$ . We then have
\begin{equation}
 \forall \: i,j,k,\ell, \quad \Abold(i,j,k,\ell) = \underset{\in \R^{1 \times r_1}}{G_1(i)}.\underset{\in \R^{r_1 \times r_2}}{G_2(j)}.\underset{\in \R^{r_2 \times r_3}}{G_3(k)}.\underset{\in \R^{r_3 \times 1}}{G_4(\ell)}
\end{equation}

%
\section{Reconstruction}\label{sec:tt:rebuild}
%

\nS The reconstruction can be done by first wrapping up all matrices built within the procedure, tracking their dimensions and the way they have been built, as in the following table to be read top down\\
\\
\begin{center}
 \begin{tabular}{c|ccl}
  Step & Matrix & Dimensions & How it has been built \\
  \hline 
  Data set & $\Abold$ & $m \times n \times p \times q$ & given as data set \\
  \hline
  &$W$ & $m \times npq$ & matricization of $\Abold$ \\
  step & $U_1$ & $m \times r_1$ & SVD of $W$: $W=U_1W_1$ \\
  1 & $W_1$ & $r_1 \times npq$ & idem\\
  & $W'_1$ & $r_1n \times pq$ & reshaping $W_1$ \\
  \hline
  & $U_2$ & $r_1n \times r_2$ & SVD of $W'_1$: $W'_1=U_2W_2$ \\
  step & $W_2$ & $r_2 \times pq$ & idem \\
  2 & $W'_2$ & $r_2p \times q$ & reshaping $W_2$ \\  
  \hline
  step &$U_3$ & $r_2p \times r_3$ & SVD of $W'_2$: $W'_2=U_3W_3$ \\
  3 & $W_3$ & $r_3 \times q$ & idem\\
 \end{tabular}
\end{center}

\nS It is then possible to rebuild both tensor $\Abold$ and its dimensions, knowing the ranks, as in the following diagram (reading previous table bottom up). We are given $(U_1,U_2,U_3,W_3)$.\\
\\
\begin{center}
\begin{tikzpicture}
\node (3) at (0,4) {$\left. \begin{array}{l}W_3\\U_3\end{array}\right\}$} ;
\node (2p) at (3,4) {$W'_2=U_3W_3$} ;
\node (2)  at (3,2) {$\left. \begin{array}{l}W_2\\U_2\end{array}\right\}$} ;
\node (1p) at (6,2) {$W'_1=U_2W_2$} ;
\node (1)  at (6,0) {$\left. \begin{array}{l}W_1\\U_1\end{array}\right\}$} ;
\node (0) at (9,0) {$W=U_1W_1$} ;
\node (A) at (11,0) {$\Abold$} ;
\draw[->] (3)--(2p) ;
\draw[->, dashed] (2p)--(2) ;
\draw[->] (2)--(1p) ;
\draw[->, dashed] (1p)--(1) ;
\draw[->] (1)--(0) ;
\draw[->, dashed] (0)--(A) ;
\end{tikzpicture}
\end{center}
where a dashed arrow means a reshaping.

%
\section{Quality of the best approximation}\label{sec:tt:qual}
%

Let us present the estimate of the quality of the approximation on the example developed in section \ref{sec:tt:best}. Here, for sake of simplicity of notations, all local ranks $r_\mu$ will be assumed to be equal and denoted $r$.

\nS TT best approximation of $\Abold$ is obtained by a succession of matricizations, reshapings and SVD. Quality of approximation of each step is presented here.

\nP{$\rightarrow$} \emph{First matricization:} We recall that $A_1 \in \R^{m \times npq}$ is the matricization of $\Abold$ on its first mode (see equation (\ref{tt:best:eq:1})). Then
\begin{equation}
 \|A_1\|^2=\|\Abold\|^2
\end{equation}

\nP{$\rightarrow$} \emph{SVD of $A_1$:} Let us assume that truncation of the SVD is done at rank $r$. We have
\begin{equation}
   A_1 = \underbrace{\sum_{a \leq r} u_a \otimes w_a}_{U_1W_I} \: + \: \underbrace{\sum_{a > r} u_a \otimes w_a}_{E_1}
\end{equation}
Hence
\begin{equation}
 \|A_1\|^2 = \|U_1W_1\|^2+ \|E_1\|^2
\end{equation}
because the families $(u_a)_a$ and $(w_a)_a$ are orthonormal. A key observation which allows to track the quality of the approximation is that 
\begin{equation}
 \|U_1W_1\|^2=\|W_1\|^2
\end{equation}
\begin{proof}
Indeed, the columns of $U_1$ form an orthonormal system with $U_1^\t U = \I_r$. Then 
\begin{equation}
 \begin{array}{lcl}
   \|U_1W_1\|^2 &=& \Tr (U_1W_1)^\t(U_1W_1) \\
   &=& \Tr W_1^\t U_1^t U_1W_1 \\
   &=& \Tr W_1^t W_1 \\
   &=& \|W_1\|^2
 \end{array}
\end{equation}
\end{proof}
\noindent Then
\begin{equation}\label{tt:qual:eq:1}
 \begin{array}{lcl}
\|\Abold\|^2 &=& \|A_1\|^2 \\
&=& \|W_1\|^2+\|E_1\|^2
 \end{array}
\end{equation}
with
\begin{equation}
 \|E_1\|^2 = \sum_{a > r}\; \sigma_a^2
\end{equation}
where the $\sigma_a$ are the singular values of the SVD of $A_1$.


\nP{$\rightarrow$} \emph{Second SVD of $W_1$} The same procedure can be run with $\|W_1\|^2$. First $W'_1$ is build by a reshaping of $W_1$ (see equation (\ref{tt:best:eq:2})). Then $\|W'_1\|=\|W_1\|$. SVD of $W'_1$ reads
\begin{equation}
 W'_1 = U_2W_2+E_2
\end{equation}
with (omitting the details similar to first SVD of $A_1$)
\begin{equation}\label{tt:qual:eq:2}
 \begin{array}{lcl}
   \|W_1\|^2 &=& \|W'_1\|^2 \\
   &=& \|U_2W_2\|^2 + \|E_2\|^2 \\
   &=& \|W_2\|^2 + \|E_2\|^2 \\
 \end{array}
\end{equation}
with
\begin{equation}
 \|E_2\|^2 = \sum_{b > r} \; \sigma_{2,b}^2
\end{equation}
where the $\sigma_{2,b}$ are the singular values of the SVD of $W_1$. Finally
\begin{equation}
 \|\Abold\|^2 = \|W_2\|^2 + \|E_2\|^2 + \|E_1\|^2
\end{equation}


\nP{$\rightarrow$}\emph{Third SVD of $W_2$} Finally, the same procedure applied to the SVD of $W_2$ yields
\begin{equation}\label{tt:qual:eq:a}
 \|\Abold\|^2 = \|W_3\|^2 + \|E_3\|^3 + \|E_2\|^2 + \|E_1\|^2
\end{equation}

\nT{Global quality} One might think looking at equation (\ref{tt:qual:eq:a}) that errors are additive between SVD. But they are multiplicative, i.e. if $\theta$ is the average quality of a SVD, then the global qulity is $\theta^{d-1}$ and not $(d-1)\theta$ (here, $d=4$). To see this, 

\nP{$\rightarrow$}Let $\theta_1 \in [0,1] \subset \R$ be defined from equation (\ref{tt:qual:eq:1}) as the quality of the first SVD
\begin{equation}
 \|\Abold\|^2 = \|W_1\|^2 + \|E_1\|^2\qquad \mbox{with} \quad 
 \begin{cases}
   \|W_1\|^2  &= \theta_1\|\Abold\|^2\\
   \|E_1\|^2 &= (1-\theta_1)\|\Abold\|^2
 \end{cases}
\end{equation}

\nP{$\rightarrow$} Let us now define similarly $\theta_2$ as the quality of the second SVD from equation (\ref{tt:qual:eq:2})
\begin{equation}
 \|W_1\|^2 = \|W_2\|^2 + \|E_2\|^2 \qquad \mbox{with} \quad 
 \begin{cases}
   \|W_2\|^2 &= \theta_2\|W_1\|^2\\
   \|E_2\|^2 &= (1-\theta_2)\|W_1\|^2
 \end{cases}
\end{equation}
So
\begin{equation}
 \left\{
 \begin{array}{lcrcr}
   \|W_2\|^2 &=& \theta_2\|W_1\|^2 &=&  \theta_2 \theta_1\|\Abold\|^2\\
   \|E_2\|^2 &=& (1-\theta_2)\|W_1\|^2 &=& (1-\theta_2)\theta_1\|\Abold\|^2 \\
   \|E_1\|^2 &=& (1-\theta_1)\|\Abold\|^2 & &
 \end{array}
 \right.
\end{equation}

\nP{$\rightarrow$} Let us define $\theta_3$ as the quality of the third SVD, by
\begin{equation}
 \|W_2\|^2 = \|W_3\|^2 + \|E_3\|^2 \qquad \mbox{with} \quad 
 \begin{cases}
   \|W_3\|^2 &= \theta_3\|W_2\|^2\\
   \|E_3\|^2 &= (1-\theta_3)\|W_2\|^2
 \end{cases}
\end{equation}

\nP{$\rightarrow$} Finally
\begin{equation}
 \left\{
 \begin{array}{lcrcr}
   \|W_3\|^2 &=& \theta_3\|W_2\|^2 &=&  \theta_3 \theta_2 \theta_1\|\Abold\|^2\\
   \|E_3\|^2 &=& (1-\theta_3)\|W_2\|^2 &=& (1-\theta_3)\theta_2 \theta_1\|\Abold\|^2 \\
   \|E_2\|^2 &=& (1-\theta_2)\|W_1\|^2 &=& (1-\theta_2)\theta_1\|\Abold\|^2 \\
   \|E_1\|^2 &=& & & (1-\theta_1)\|\Abold\|^2 
 \end{array}
 \right.
\end{equation}

\nT{Result} Thus, the quality of approximation of TT aproximation is $\theta_1\theta_2\theta_3$ as
\begin{equation}
 \|W_3\|^2 = \theta_1\theta_2\theta_3\|\Abold\|^2
\end{equation}
The error term can be decomposed as
\begin{itemize}[label=$\rightarrow$]
 \item $1-\theta_1$ for the first SVD
 \item $(1-\theta_2)\theta_1$ for the second SVD 
 \item $(1-\theta_3)\theta_2 \theta_1$ for the third SVD
\end{itemize}
In the general case of $d$ modes, one has
\begin{equation}
 \|W_{d-1}\|^2 = \left(\prod_{\mu=1}^{d-1}\theta_\mu\right)\|\Abold\|^2
\end{equation}
where $\theta_\mu$ is the quality of the $\mu-$th SVD, and global quality $\theta$ is
\begin{equation}
 \theta = \prod_{\mu=1}^{d-1}\theta_\mu
\end{equation}

\nT{Interpretation} Deriving best TT approximation of a tensor $\Abold$ at a prescribed rank $\rbold$ is trading $(i)$ a linear scaling of the storage requirement with the order $d$ (about $ndr^2$ terms instead of $n^d$ for a $d-$modes tensor and dimension $n$ on each mode for rank $r$) and a possibility to make some elementary operations like addition and Hadamard product in TT format (see section \ref{sec:tt:elem}) with $(ii)$ an exponential decrease in quality with the order $d$ as $\theta=\prod_{\mu=1}^{d-1}\theta_\mu$. In brief, \emph{there is no free lunch}.

%
\section{Sketch with tensor diagrams}\label{sec:tt:diagram}
%

Let us recall that tensor diagrams \marginpar{\emph{tensor diagram}}\index{tensor!diagram} have been introduced in section \ref{sec:matelem:diag}. Here, I present TT decomposition with a tensor diagram on the example used in section \ref{sec:tt:best}
\begin{equation}
 \Abold \in E \otimes F \otimes G \otimes H
\end{equation}
with the following dimensions and indices:\\
\\
\begin{center}
 \begin{tabular}{l|cccc}
  \hline
  space & $E$ & $F$ & $G$ & $H$ \\
  dimension & $m$ & $n$ & $p$ & $q$ \\
  indices & $i$ & $j$ & $k$ & $\ell$ \\
  \hline
 \end{tabular}

\end{center}

\nP{$\rightarrow$} The tensor $\Abold$ can be sketched as 
\begin{center}
 \begin{tikzpicture}
   \node[tensor_box] (A) at (1,0) {$\Abold$} ; 
   \draw (0,0)--(A) node[pos=0.5, above] {$i$};
   \draw (A)--(2,0) node[pos=0.5, below] {$k$};
   \draw (A)--(1,1) node[pos=0.5, right] {$j$};
   \draw (A)--(1,-1) node[pos=0.5, left] {$\ell$};
 \end{tikzpicture}
\end{center}

\nP{$\rightarrow$} The first SVD $A = U_1W_1$ written component-wise
\begin{equation}
 A(i,j,k,l) = \sum_a\, U_1(i,a)W_1(a, j,k,\ell)
\end{equation}
can be sketeched by 
\begin{center}
 \begin{tikzpicture}
   \node[tensor_box] (U) at (1,1) {$U_1$} ; 
   \node[tensor_box] (W) at (3,1) {$W_1$} ; 
   \draw (U)--(W) node [pos=0.5, above] {$a$} ;
   \draw (U)--(1,0) node [pos=0.5, left] {$i$} ;
   \draw (W)--(3,2) node [pos=0.5, right] {$j$} ;
   \draw (W)--(4,1) node [pos=0.5, below] {$k$} ;
   \draw (W)--(3,0) node [pos=0.5, left] {$\ell$} ;   
 \end{tikzpicture}
\end{center}

\nP{$\rightarrow$} The second SVD $W'_1= U_2W_2$ written componentwise
\begin{equation}
 W_1(a,j,k,\ell) = \sum_b \, U_2(a,j,b)W_2(b, k, \ell)
\end{equation}
can be sketched by
\begin{center}
 \begin{tikzpicture}
   \node[tensor_box] (U) at (1,1) {$U_1$} ; 
   \node[tensor_box] (U2) at (3,1) {$U_2$} ;
   \node[tensor_box] (W2) at (5,1) {$W_2$} ;
   \draw (U)--(U2) node [pos=0.5, above] {$a$} ;
   \draw (U2)--(W2) node [pos=0.5, above] {$b$} ;
   \draw (U)--(1,0) node [pos=0.5, left] {$i$} ;
   \draw (U2)--(3,0) node [pos=0.5, left] {$j$} ;
   \draw (W2)--(6,1) node [pos=0.5, above] {$k$} ;
   \draw (W2)--(5,0) node [pos=0.5, left] {$\ell$} ;   
 \end{tikzpicture}
\end{center}

\nP{$\rightarrow$} And finally the third SVD $W'_2= U_3W_3$ written componentwise
\begin{equation}
 W_2(b,k,\ell) = \sum_c \, U_3(b,k,c)W_3(c, \ell)
\end{equation}
can be sketched by
\begin{center}
 \begin{tikzpicture}
   \node[tensor_box] (U) at (1,1) {$U_1$} ; 
   \node[tensor_box] (U2) at (3,1) {$U_2$} ;
   \node[tensor_box] (U3) at (5,1) {$U_3$} ;
   \node[tensor_box] (W3) at (7,1) {$W_3$} ;
   \draw (U)--(U2) node [pos=0.5, above] {$a$} ;
   \draw (U2)--(U3) node [pos=0.5, above] {$b$} ;
   \draw (U3)--(W3) node [pos=0.5, above] {$c$} ;
   \draw (U)--(1,0) node [pos=0.5, left] {$i$} ;
   \draw (U2)--(3,0) node [pos=0.5, left] {$j$} ;
   \draw (U3)--(5,0) node [pos=0.5, left] {$k$} ;
   \draw (W3)--(7,0) node [pos=0.5, left] {$\ell$} ;   
 \end{tikzpicture}
\end{center}

\notes The Tensor-Train format as presented here has been proposed in \cite{Oseledets2009,Oseledets2011} and quickly adopted by the community of numerical tensor analysis. Its main advantages for numerical calculations is the numerical stability of the SVD on matricized forms of involved tensors. The two main tools to "enter into the game" are the TT-SVD and the rounding. A tensor $\Abold$ being given, TT-SVD provides as output a TT format decomposition or best approximation at prescribed multilinear rank $\rbold$. It is algorithm 1 in \cite{Oseledets2011}. The rounding is deriving the TT-SVD of a tensor when the tensor is given in TT format only, ignoring its actual coefficients. It is fully developed in section 3 of \cite{Oseledets2011}, from the link between QR decomposition and SVD. The procedure is as follows (with Oseledets notations): in $A=UV^\t$, write QR decomposition of $U$ and $V$ as $U=Q_uR_u$ and $V=Q_vR_v$. Build $P=R_uR_v^\t$ and compute the reduced SVD $P=XDY^\t$. Then, truncated SVD of $A$ reads $\widehat{U}=Q_uX$ and $\widehat{V}=Q_vY$. The following step is to observe that, if $A$ is "horizontal", $U$ is small, and its QR decomposition can be computed directly. Third step is to derive QR decomposition of $V$ in TT format only. The full algorithm is presented as algorithm 2 in \cite{Oseledets2011}.

\nB Tensor Train format is a particular case of a more general and rich structure: tensor networks, which have been developed (with many guises) for modeling and simulating many body quantum systems. Accessible surveys on tensor networks from viewpoint of theoretical physics are \cite{Orus2014,Montangero18,Ran2020}. Tensor Trains are equivalent to Matrix Product States, which have been a landmark for modeling and simulating 1d entangled quantum many body systems. It has been naturally extended to 2d systems under the name of Projected Entangled Pair States (PEPS). A historical survey is presented in the introduction of \cite{Ran2020}. 

\nB A tensor network is a family of elementary tensors coupled by contraction on common indices. This is better seen with tensor diagrams. Here is a family with repeated indices (a repeated index is a common mode; for example, $U \in E \otimes F$, $G \in F \otimes H \otimes K$ and repeated index is $\alpha$ for common mode $F$)
\begin{center}
 \begin{tikzpicture}
   \node[tensor_box] (U) at (1,1) {$U$} ; 
   \node[tensor_box] (G) at (4,1) {$G$} ;
   \node[tensor_box] (V) at (7,1) {$V$} ;
   \draw (U)--(2,1) node [pos=0.5, above] {$\alpha$} ;
   \draw (3,1)-- (G)  node [pos=0.5, above] {$\alpha$} ;
   \draw (G)-- (5,1)  node [pos=0.5, above] {$\beta$} ;
   \draw (6,1)-- (V)  node [pos=0.5, above] {$\beta$} ;
   \draw (U)--(1,0) node [pos=0.5, left] {$i$} ;
   \draw (G)--(4,0) node [pos=0.5, left] {$j$} ;
   \draw (V)--(7,0) node [pos=0.5, right] {$k$} ;
 \end{tikzpicture}
\end{center}
and here is the tensor network
\begin{center}
 \begin{tikzpicture}
   \node[tensor_box] (U) at (1,1) {$U$} ; 
   \node[tensor_box] (G) at (3,1) {$G$} ;
   \node[tensor_box] (V) at (5,1) {$V$} ;
   \draw (U)--(G) node [pos=0.5, above] {$\alpha$} ;
   \draw (G)-- (V)  node [pos=0.5, above] {$\beta$} ;
   \draw (U)--(1,0) node [pos=0.5, left] {$i$} ;
   \draw (G)--(3,0) node [pos=0.5, left] {$j$} ;
   \draw (V)--(5,0) node [pos=0.5, right] {$k$} ;
 \end{tikzpicture}
\end{center}
which is a tensor Train. Extension to PEPS is as follows:\\
\\
\begin{center}
    \begin{tikzpicture}
    %
    %
      \node[tensor_box] (00) at (0,0){};
      \node[tensor_box] (10) at (1,0){};
      \node[tensor_box] (20) at (2,0){};
      \node[tensor_box] (30) at (3,0){};
      \node[tensor_box] (40) at (4,0){};
      \node[tensor_box] (01) at (0.5,1){};
      \node[tensor_box] (11) at (1.5,1){};
      \node[tensor_box] (21) at (2.5,1){};
      \node[tensor_box] (31) at (3.5,1){};
      \node[tensor_box] (41) at (4.5,1){};
      \node[tensor_box] (02) at (1,2){};
      \node[tensor_box] (12) at (2,2){};
      \node[tensor_box] (22) at (3,2){};
      \node[tensor_box] (32) at (4,2){};
      \node[tensor_box] (42) at (5,2){};
      %
      %
      \draw(00)--(10);
      \draw(10)--(20);
      \draw(20)--(30);
      \draw(30)--(40);
      \draw(01)--(11);
      \draw(11)--(21);
      \draw(21)--(31);
      \draw(31)--(41);
      \draw(02)--(12);
      \draw(12)--(22);
      \draw(22)--(32);
      \draw(32)--(42);
      %
      %
      \draw (00)--(01);
      \draw (10)--(11);
      \draw (20)--(21);
      \draw (30)--(31);
      \draw (40)--(41);
      \draw (01)--(02);
      \draw (11)--(12);
      \draw (21)--(22);
      \draw (31)--(32);
      \draw (41)--(42);
      %
      %
      \draw (00)--(0,-.7);
      \draw (01)--(0.5, 0.3);
      \draw (02)--(1, 1.3);
      \draw (10)--(1,-.7);
      \draw (11)--(1.5, 0.3);
      \draw (12)--(2, 1.3);
      \draw (20)--(2,-.7);
      \draw (21)--(2.5, 0.3);
      \draw (22)--(3, 1.3);      
      \draw (30)--(3,-.7);
      \draw (31)--(3.5, 0.3);
      \draw (32)--(4, 1.3);  
      \draw (40)--(4,-.7);
      \draw (41)--(4.5, 0.3);
      \draw (42)--(5, 1.3);  
    \end{tikzpicture}
\end{center}
Here it is as tensors with indices on a smaller grid:\\
\\
\begin{center}
    \begin{tikzpicture}
    %
    %
      \node[tensor_box] (00) at (0,0){$A$};
      \node[tensor_box] (10) at (2,0){$B$};
      \node[tensor_box] (01) at (0.5,1.5){$D$};
      \node[tensor_box] (11) at (2.5,1.5){$C$};
      %
      %
      \draw(00)--(10) node [pos=0.5, above] {$\alpha$};
      \draw(01)--(11) node [pos=0.5, above] {$\gamma$};
      %
      %
      \draw (00)--(01) node [pos=0.5, left] {$\delta$};
      \draw (10)--(11) node [pos=0.5, left] {$\beta$};
      %
      %
      \draw (00)--(0,-1) node [pos=0.5, right] {$i$};
      \draw (01)--(0.5, 0.5) node [pos=0.5, right] {$\ell$};
      \draw (10)--(2,-1) node [pos=0.5, right] {$j$};
      \draw (11)--(2.5, 0.5) node [pos=0.5, right] {$k$};
    \end{tikzpicture}
\end{center}
The tensors are, with indices
\[
\Abold_{i\alpha\delta}, \quad \Bbold_{j\alpha\beta}, \quad \Cbold_{k\beta\gamma}, \quad \Dbold_{\ell\gamma\delta}
\]
and
\[
\Tbold_{ijk\ell}=\sum_{\alpha\beta\gamma\delta}\, a_{i\alpha\delta}\, b_{j\alpha\beta} \, c_{k\beta\gamma} \, d_{\ell\gamma\delta}
\]
which is a PEPS model for $\Tbold$.

%% file: quantization.tex
%
\chapter{Discretization of functions}\label{chap:quant}
%

One of the most nagging questions in the study of tensor rank is: which characteristics of a tensor confer to it the property of being of low rank? 

\nB Indeed, typical CP rank of tensors in $(\R^n)^{\otimes d}$ scales with $n^{d-1}/d$, whereas real world tensors have a far smaller rank. This is not yet fully understood. Here, an answer is given for specific tensors: a tensor the coefficients of which are the evaluation of a multivariate continuous function on a grid is expected to be low rank. The explanation is rather simple and given here with all details: first, it is shown that tensors as estimates of polynomials on a grid are low rank, and second the Stone-Weierstrass theorem states that any continuous function can be approximated as close as prescribed by a polynomial for the uniform norm: polynomials are dense in the space of continuous functions.

\nB The idea behind is very simple: if $f$ and $g$ are real univariate functions, $f \otimes g$ is the  bivariate function defined by $(f \otimes g)(x,y)=f(x)g(y)$. It is a rank one function. A polynomial like $P(x,y)=x^2 + 2xy + y^2$ can be written $x^2 \otimes \ones + 2 x \otimes y + \ones \otimes y^2$ where $\ones$ is the constant function $\ones(x)=1$. Such a polynomial is rank 3 at most. Here, it is shown that the rank can be transferred to the tensor (here a matrix) which is the evaluation of $P$ on a grid, whatever its size.

\nB In this chapter, the polynomials represent finite rank functions (each polynomial is a finite sum of monomials). This observation is at the root of controlling the rank of tensors built by discretization of a multivariate function on a grid. This chapter recalls and associates different classical approaches for best approximation of multivariate functions (continuous, with $\ell^2$ or $\ell^\infty$ norm) like Stone-Weierstrass theorem for approximating continuous functions by polynomials, development and truncation of Fourier series with orthogonal polynomials which is equivalent to Tucker model in tensor space of multivariate functions, etc... It is organized as follows:

\begin{description}
\item[section \ref{sec:quant:prelim}] Two standard products of functions with which we will work are recalled.
\item[section \ref{sec:quant:compact}] CP, Tucker and TT format for a tensor are recalled.
\item[section \ref{sec:quant:mesh}] The definition of a mesh is given in 1D, and extended to several dimension, with definition of a product between Cartesian grids which will appear to be dual to the tensor product between functions.
\item[section \ref{sec:quant:quant}] The discretization of a multivariate function on a multivariate grid is presented.
\item[section \ref{sec:quant:prod}] This is a key section where it is shown that the rank of the tensor built by discretization of a multivariate function on a grid is equal to the CP-rank of the function with tensor product of functions being defined by $f \otimes g (x,y)=f(x)g(y)$.
\item[section \ref{sec:quant:poly}] The finite CP-rank of polynomials in the tensor space of multivariate functions (defined with with $\otimes$ as in section \ref{sec:quant:prod}) is presented.
\item[section \ref{sec:quant:SW}] The Stone-Weierstrass  theorem for approximation of continuous functions by polynomials for $\ell^\infty$ norm is presented.
\item[section \ref{sec:quant:ortho}] The decomposition of a function on a basis of orthonormal polynomials is presented. A link is established with Tucker model in tensor space of multivariate functions. 
\item[section \ref{sec:quant:elem}] Some elementary operations on tensors have been presented in chapter \ref{chap:tenselem}. Their natural extension (or not) to tensor spaces of multivariate functions (or tensor space of Fourier coefficients of their development on basis of orthonormal polynomials) is presented.  
\item[section \ref{sec:quant:nutshell}] Finally, navigation between different spaces and approximations in this chapter is presented in a nutshell as a diagrammatic sketch.
\end{description}

%
\section{Preliminaries}\label{sec:quant:prelim}
%

Let us observe that a vector $\abold = (a_1, \ldots, a_n) \in \R^n$ can be considered as a function
\[
 \begin{CD}
  \llbracket 1,n\rrbracket @>\abold>> \R
 \end{CD}
\]
which can be made explicit in $\abold_i := \abold(i)$. Let $\abold, \bbold$ be two functions on $\llbracket 1,m\rrbracket$ and $\llbracket 1,n\rrbracket$ respectively. Then, the function associated to $\abold \otimes \bbold$ 
\[
 \begin{CD}
  \llbracket 1,m\rrbracket \times \llbracket 1,n\rrbracket @>\abold \otimes \bbold>> \R
 \end{CD}
\]
is a such that
\[
 \abold \otimes \bbold \; (i,j) = \abold(i)\bbold(j)
\]
This can be extended to functions
\[
 \begin{CD}
  \R @>f,g>> \R
 \end{CD}
\]
with, if $x,y \in \R$
\begin{equation}\label{eq:quant:1}
 f \otimes g \; (x,y) = f(x)g(y)
\end{equation}
However, for $\otimes$ to be meaningful, $f$ and $g$ must belong to some explicit vector spaces of functions. For example, if they are polynomials of degree not greater than $n$, $f$ and $g$ belong to a finite dimensional vector space, and all the theory presented up to now can be deployed.  If the vector space they belong to is infinite dimensional, there are more technicalities, but the theory of tensor product between infinite dimensional vector spaces is well developed, especially between Hilbert spaces (like spaces of square integrable functions).

\nS This is consistent with the definition of the tensor product as a universal property. Let $\H$ be a functional space of real functions of a real variable and $f,g \in \H$. Let $\H \otimes \H$ denote the space of bivariate functions $\R \times \R \longrightarrow \R$. Then, we have 
\begin{equation}
 \begin{CD}
  \H \times \H @>>> \H \otimes \H \\
  (f,g) @>>> f \otimes g
 \end{CD}
\end{equation}
Even if it seems obvious, it is important to have in mind that\marginpar{\dbend} $fg \neq f \otimes g$ as $fg \in \H$:
\begin{equation}
 \begin{CD}
  \H \times \H @>>> \H \\
  (f,g) @>>> fg
 \end{CD}
\end{equation}
with
\begin{equation}
 fg(x) = f(x)g(x)
\end{equation}
To be concrete, if 
\begin{equation}
 f(x)=x^2, \qquad g(x) = x^4
\end{equation}
then
\begin{equation}
 fg(x) = x^6, \qquad f \otimes g (x,y) = x^2y^4
\end{equation}

\nS Along this line, a tensor $\Abold =\sum_{ijk}\, a_{ijk} \, \ebold_i \otimes \fbold_j \otimes \gbold_k$ can be read as a function
\begin{equation}
\begin{CD}
  I \times J \times K @>\Abold>> \R, 
  \qquad \mbox{with} \quad
  \begin{cases}
   I &= \{1,\ldots, m\} \\
   J &= \{1,\ldots, n\} \\
   K &= \{1,\ldots, p\} \\
  \end{cases}
\end{CD}
\end{equation}
with
\begin{equation}
 \Abold(i,j,k) = a_{ijk}
\end{equation}
(this depends obviously on the selected basis in $E,F,G$) and (\ref{eq:quant:1}) can be extended to multivariate functions as
\begin{equation}
 (f \otimes g \otimes h)(x,y,z) = f(x)g(y)h(z)
\end{equation}
or more generally
\begin{equation}
 \left(\bigotimes_\mu f_\mu\right)(x_1,\ldots,x_n) = \prod_\mu f_\mu(x_\mu)
\end{equation}
This is a rank one function where separation of the variables is highlighted by the operator $\otimes$. It is important to specify in which functional space the functions $f$ and $g$ live. It is not required at this preliminary stage. It is sufficient to say that they live in a functional space $\H$ which is closed by ordinary product, i.e.
\begin{equation}
 \left. 
    \begin{array}{ll}
      f &\in  \H \\
      g & \in \H
    \end{array}
 \right\}
 \quad \Longrightarrow \quad fg \in \H
\end{equation}
Functional spaces $\H$ such that any function $f \in \H$ is a finite linear combination of rank one functions will play a central role, like spaces of polynomials. 

\nB In this chapter, we study the extension of compact representation of such functions with CP, Tucker and TT formats.

%
\section{Compact formats in a nutshell}\label{sec:quant:compact}
%

I recall here the definitions of some compact formats which have been presented in these notes, and their equivalent in infinite dimensional function spaces. Those formats have been developed for tensors which are elements in a tensor product of finite dimensional spaces. Extension to infinite dimensional vector spaces is delicate, and will not be developed here. Here are some starting points. The r.h.s. of the arrow is the functional equivalent of the l.h.s. classically developed for finite dimensional vector spaces.
\begin{equation}
    \begin{CD}
      \displaystyle \Abold \in \bigotimes_{\mu=1}^d E_\mu @>>> \bm{\phi} \: : \: \R^d \longrightarrow \R \\
      \mbox{tensor} @>>> \mbox{multivariate function}
    \end{CD}
\end{equation}
As illustrative examples, we have often selected
\begin{equation}
 \Abold \in E \otimes F \otimes G
\end{equation}
with
\begin{equation}
 \Abold = \sum_{ijk} \, a_{ijk} \, \ebold_i \otimes \fbold_j \otimes \gbold_k
\end{equation}
We have presented some compact formats, illustrated here for $d=3$:

\nP{$\rightarrow$} CP-format
\begin{equation}
 \begin{CD}
   \displaystyle \Abold = \sum_{a=1}^r\, \ubold_a \otimes \vbold_a \otimes \wbold_a @>>> \displaystyle \bm{\phi}(x,y,z) = \sum_a \, f_a(x)g_a(y)h_a(z)
 \end{CD}
\end{equation}

\nP{$\rightarrow$} Tucker format
\begin{equation}
 \begin{CD}
   \displaystyle \Abold = \sum_{i \leq r}\sum_{j\leq r}\sum_{k \leq r}\, a_{ijk} \, \ubold_i \otimes \vbold_j \otimes \wbold_k @>>> \displaystyle \bm{\phi}(x,y,z) = \sum_{\nbold \prec \rbold} \, \widehat{\phi}_\nbold \, T_{n_1}(x)T_{n_2}(y)T_{n_3}(z)
 \end{CD}
\end{equation}
(the notations are developed in section \ref{sec:quant:ortho})

\nP{$\rightarrow$} TT format 
\begin{equation}
 a_{ijk} = \ubold(i).G(j).\vbold(k)
\end{equation}
TT format is rather used for low rank approximation of finite dimensional tensors built from multivariate functions. 

\notes One of the key question we have addressed is: a tensor $\Abold$ being given, find a rank $r$ as small as possible such that there exists a tensor of rank $r$ as close as possible to $\Abold$. This question can be transported as such for functional tensors: a multivariate function $\bm{\phi}$ being given, find a rank $\rbold$ tensor ($\rbold$ as small as possible) such that there exists a multivariate function of rank $\rbold$ as close as possible to $\bm{\phi}$. This requires that $\H$ to which $\bm{\phi}$ belongs is a normed space. This has been thoroughly studied for Banach spaces (complete normed spaces) and Hilbert spaces (Banach spaces endowed with an inner product inducing the norm). This functional analysis side of the approach will not be developed in this chapter. See \cite{Ryan2002,Hackbusch2012, Hackbusch2019} for tensor product in Banach spaces.

%
\section{Mesh}\label{sec:quant:mesh}
%

A function
\[
 \begin{CD}
  \R @>f>> \R
 \end{CD}
\]
can be given
\begin{itemize}[label=$\rightarrow$]
 \item $f \in \H$: as an element in a functional space $\H$, like $\H = \mathrm{C}^{\infty}([0,1])$, the space of smooth real functions from $[0,1]$ to $\R$
 \item by its values for any point in its domain, like $f(x)= x^2+2x+1$
\end{itemize}
Here, we associate to each function in a functional space $\H$ its values on a grid, which can be multidimensional if the function is multivariate. A mesh is a grid on which a real function is evaluated. Such an evaluation will be called a \kw{discretization}.

\nS A \kw{mesh} of size $n$ in $\R$ is a sequence
\begin{equation}
 \xbold = (x_1, \ldots,x_n) \qquad \mbox{with}\quad x_1 < \ldots < x_n
\end{equation}
Note that the ordering $x_1 < \ldots $ is not compulsory, but convenient for continuous functions. The evaluation of a function $f\: : \:\R \longrightarrow \R$ on mesh $\xbold$ is $(f(x_1), \ldots,f(x_n))$ which will be developed later.

\nS Let us consider 2 meshes $\xbold, \ybold$ in $\R$ of respective sizes $m,n$. They define a mesh denoted $\xbold \otimes_\m \ybold$ in $\R^2$ by
\begin{equation}
 \xbold \otimes_\m \ybold = 
 \begin{pmatrix}
  (x_1,y_1) & (x_1,y_2) & \ldots & (x_1,y_n) \\
  (x_2,y_1) & (x_2,y_2) & \ldots & (x_2,y_n) \\
  \vdots & \vdots &  & \vdots \\
  (x_m,y_1) & (x_m,y_2) & \ldots & (x_m,y_n) \\
 \end{pmatrix}
\end{equation}

\nS This can be extended to building $d-$dimensional meshes, with less visualization. For $d=3$, the mesh $\xbold \otimes_\m \ybold \otimes_\m \zbold$ is defined by\\
\begin{center}
 \begin{tikzpicture}
  \node (mesh) at (-3,1){
  $\xbold \otimes_\m \ybold \otimes_\m \zbold =$
  };
  \node (111) at (0,2){$(x_1,y_1,z_1)$};
  \node (m11) at (0,0){$(x_m,y_1,z_1)$};
  \node (mn1) at (4,0){$(x_m,y_n,z_1)$};
  \node (1n1) at (4,2){$(x_1,y_n,z_1)$};
  \node (11p) at (2,3){$(x_1,y_1,z_p)$};
  \node (1np) at (6,3){$(x_1,y_n,z_p)$};
  \node (m1p) at (2,1){$(x_m,y_1,z_p)$};
  \node (mnp) at (6,1){$(x_m,y_n,z_p)$}; 
  \draw[dashed] (111) -- (m11) ; 
  \draw[dashed] (111) -- (1n1) ; 
  \draw[dashed] (1n1) -- (mn1) ; 
  \draw[dashed] (mn1) -- (m11) ; 
  \draw[dashed] (11p) -- (m1p) ; 
  \draw[dashed] (11p) -- (1np) ; 
  \draw[dashed] (1np) -- (mnp) ; 
  \draw[dashed] (mnp) -- (m1p) ; 
  \draw[dashed] (111) -- (11p) ; 
  \draw[dashed] (1n1) -- (1np) ; 
  \draw[dashed] (m11) -- (m1p) ; 
  \draw[dashed] (mn1) -- (mnp) ;   
 \end{tikzpicture}
\end{center}
and so on for $d>3$.

\nS Meshes defined like this are \kw{Cartesian grids}.

\notes Meshes have been defined and studied for decades, for example in signal processing (see \cite{Pages2015}). They are called \kw{discretization grid} or simpler \kw{discretizers} or \kw{grids}.  One issue which has been thoroughly studied in many domains is to find an optimal mesh for discretizing a continuous function for a given purpose (like approximation, integration, etc ...). If the size of a one dimensional mesh is $n$, the size of equivalent $d-$dimensional mesh is $n^d$ which hits the wall of the curse of dimensionality. Therefore, techniques have been elaborated for sampling a (well designed) set of points only and not the full mesh, known as \kw{sparse grids}. This topic is out of scope of these notes although it has contributed to the use of tensor based methods for addressing numerical issues in the domain of multivariate functions. What we show next is that the rank of the tensor of values of a multivariate polynomial on a full mesh  does not depend on the size of the mesh but on the number of monomials only. This means that this tensor needs not to be stored in memory even for exact calculations.

%
\section{Discretization}\label{sec:quant:quant}
%

Let $f$ be a real function of a real variable
\begin{equation}
 \begin{CD}
  \R @>f>>\R
 \end{CD}
\end{equation}
and $\xbold$ a mesh in $\R$. The \kw{discretization} of function $f$ on mesh $\xbold$ is the vector
\begin{equation}
 \QS(f,\xbold) = (f(x_1), \ldots,f(x_n)) \quad \in \R^n
\end{equation}

\nS Let $\phi$ be a real function with two variables
\begin{equation}
 \begin{CD}
  \R^2 @>\bm{\phi}>>\R
 \end{CD}
\end{equation}
The discretization of $\bm{\phi}$ on mesh $\xbold \otimes_\m \ybold$ is the matrix
\begin{equation}
 \QS(\bm{\phi}\, , \, \xbold \otimes_\m \ybold) = 
 \begin{pmatrix}
  \phi(x_1,y_1) & \phi(x_1,y_2) & \ldots & \phi(x_1,y_n) \\
  \phi(x_2,y_1) & \phi(x_2,y_2) & \ldots & \phi(x_2,y_n) \\
  \vdots & \vdots &  & \vdots \\
  \phi(x_m,y_1) & \phi(x_m,y_2) & \ldots & \phi(x_m,y_n) \\
 \end{pmatrix}
 \quad \in \R^{m \times n}
\end{equation}

%
\section{Product and tensor product of functions}\label{sec:quant:prod}
%

Two elementary but key results on the rank of discretization of product of functions are given here.

\nT{Functions and products of functions} Let us consider two real functions $f,g$
\begin{equation}
 \begin{CD}
  f,g \: : \: \R @>>> \R
 \end{CD}
\end{equation}
We recall that $fg$ is defined by
\begin{equation}
 \begin{CD}
  x @>fg>> f(x)g(x)
 \end{CD}
\end{equation}
whereas $f \otimes g$ is defined by
\begin{equation}
 \begin{CD}
  (x,y) @>f\otimes g>> f(x)g(y)
 \end{CD}
\end{equation}

\nT{Discretization of products of functions} Let $f,g \in \H$ and $\xbold$ be a mesh in $\R$. It is immediate to see that
\begin{equation}
 \QS(fg,\xbold) = \QS(f,\xbold) \odot \QS(g,\xbold)
\end{equation}
where $\odot$ is the Hadamard product (see section \ref{sec:alstruct:prodmat}), as $(fg)(x_i)=f(x_i)g(x_i)$. Similarily, if $\xbold,\ybold$ are two meshes, we have
\begin{equation}
 \QS(f \otimes g \, , \, \xbold \otimes_\m \ybold) = \QS(f,\xbold) \otimes \QS(g,\ybold)
\end{equation}
Those identities can be extended to products of $d$ functions as
\begin{equation}
 \QS\left(\prod_\mu f_\mu \; , \; \xbold\right) = \bigodot_\mu \QS(f_\mu,\xbold)
\end{equation}
and
\begin{equation}
 \QS\left(\bigotimes_\mu f_\mu \; , \; \otimes_\m \xbold_\mu\right) = \bigotimes_\mu \QS(f_\mu,\xbold_\mu)
\end{equation}

\nT{Rank of a function of several variables} A function
\begin{equation}
 \begin{CD}
   \R^2 @>\phi>> \R
 \end{CD}
\end{equation}
is said to have rank one if it can be written
\begin{equation}
 \bm{\phi}(x,y)=f(x)g(y)
\end{equation}
or
\begin{equation}
 \bm{\phi} = f \otimes g
\end{equation}
i.e. there is separation of the variables. Let $\bm{\phi} \: : \: \R^2 \longrightarrow \R$. Then, if $\phi$ has rank one, any discretization $\QS(\phi,\xbold \otimes_\m \ybold)$ is a matrix of rank one, whatever the meshes and their sizes. 
\begin{proof}
This is immediate, because
\begin{equation}
 \QS(f \otimes g, \xbold \otimes_\m \ybold) = \QS(f,\xbold) \otimes \QS(g,\ybold)
\end{equation}
Let us show it with more details as an exercise. Let $A$ be the discretization of $\phi$ with
\begin{equation}
 A(i,j) = \phi(x_i,y_j)
\end{equation}
and $\ubold, \vbold$ the discretization of $f,g$ as
\begin{equation}
 \left\{ 
    \begin{array}{lcl}
     \ubold_i &=& f(x_i) \\
     \vbold_j &=& g(y_j)
    \end{array}
 \right.
\end{equation}
Then
\begin{equation}
 \begin{array}{lcl}
   A(i,j) &=& \phi(x_i,y_j) \\
   &=& f(x_i)g(y_j) \\
   &=& \ubold_i\vbold_j
 \end{array}
\end{equation}
and
\begin{equation}
 A = \ubold \otimes \vbold
\end{equation}
\end{proof}

\nS Here is our main result, given first for bivariate functions, and next for general case. Let
\begin{equation}
 \begin{CD}
  \R^2 @>\bm{\phi}>>\R
 \end{CD}
\end{equation}
with CP rank $r$
\begin{equation}
 \bm{\phi} = \sum_{a=1}^r\, f_a \otimes g_a
\end{equation}
or
\begin{equation}
 \bm{\phi}(x,y) = \sum_{a=1}^r\, f_a(x)g_a(y)
\end{equation}
Let $\xbold \otimes_\m \ybold$ be a mesh in $\R^2$. Then
\begin{equation}
 \rank \QS(\bm{\phi} \, , \, \xbold \otimes_\m \ybold) = \rank \bm{\phi}
\end{equation}
This is independent of the size of the mesh. 
\begin{proof}
The full proof can be given in three lines 
\begin{equation}
 \begin{array}{lcl}
   \QS(\bm{\phi} \, , \, \xbold \otimes_\m \ybold) &=& \displaystyle \QS\left(\sum_{a=1}^r\, f_a \otimes g_a \; , \; \xbold \otimes_\m \ybold\right) \\
   &=& \displaystyle \sum_{a = 1}^r \, \QS(f_a \otimes g_a \, , \, \xbold \otimes_\m \ybold) \\
   &=& \displaystyle \sum_{a = 1}^r \, \QS(f_a, \xbold) \otimes \QS(g_a, \ybold)
   \end{array}
\end{equation}
\end{proof}


\nS This can be extended to $d-$variate functions.  Let
\begin{equation}
 \begin{CD}
  \R^d @>\bm{\phi}>>\R
 \end{CD}
\end{equation}
with CP rank $r$
\begin{equation}
 \bm{\phi} = \sum_{a=1}^r\, f_{a1} \otimes \ldots \otimes f_{ad}
\end{equation}
or
\begin{equation}
 \bm{\phi}(x_1,\ldots,x_d) = \sum_{a=1}^r\, f_{a1}(x_1) \ldots f_{ad}(x_d)
\end{equation}
Let $\xbold_1 \otimes_\m \ldots \otimes_\m \xbold_d$ be a mesh in $\R^d$. Then
\begin{equation}
 \rank \QS(\bm{\phi} \, , \, \xbold_1 \otimes_\m \ldots \otimes_\m \xbold_d) = \rank \bm{\phi}
\end{equation}
with same demonstration.

\notes A function of CP rank 1 is a function with variable separation, which opens the way to solve ODEs or PDEs, either analytically (sometimes by algebra alone) or numerically. It is an immense domain in both abstract and applied mathematics. A function of low CP rank is sometimes called \emph{weakly separated}, although this may induce some confusion with separability in week topology. The observation that the discretization of a rank $r$ function is a tensor of same CP-rank is a key observation: it explains why some compact formats like TT format are so efficient for solving numerically PDEs by working with tensors of estimates of their values on a grid. Indeed, many differential operators are linear, and often such studies boil down to systems like $Ax=b$. Data sets made by discretization are by nature low rank. They are far from being randomly generated tensors, and have ranks much lower than the typical rank induced by their sizes.

%
\section{Rank of discretized  Polynomials}\label{sec:quant:poly}
%

Here, by a slight abuse of notation, $x$ will denote either a scalar $x \in \R$ or the monomial function associated to polynomial $P=X$ in $\R[X]$. 

\nB The set of univariate polynomials is denoted $\PS$\marginpar{$\PS$}  and the space of $d-$variate polynomials is denoted $\PS^{\otimes d}$\marginpar{$\PS^{\otimes d}$} for reasons which will become clear later.

\nS Let us start with an example in $\R^2$. Let 
\begin{equation}
 P(x,y) = x^2 +2xy + y^2
\end{equation}
Then, the rank of $P$ is 3. To see this, let us denote by $\ones_x$ the function 
\begin{equation}
 \begin{CD}
  \R @>\ones_x>>\R \\
  x @>>> 1
 \end{CD}
\end{equation}
i.e.
\begin{equation}
 \forall \: x \in \R, \quad \ones_x(x)=1
\end{equation}
Then $x$ can be identified with the bivariate function $x \otimes \ones_y$, which has rank one. Thus, one can write
\begin{equation}
 P = x^2 \otimes \ones_y + 2 x\otimes y + \ones_x \otimes y
\end{equation}
This is a CP-decomposition of $P$ as a sum of three rank one functions, hence $r_{\c\p}(P) \leq 3$. As a consequence
\begin{equation}
 \rank \QS(P, \xbold \otimes_\m \ybold) \leq 3
\end{equation}
because 
\begin{equation}
 \rank \QS(P, \xbold \otimes_\m \ybold) = \rank P
\end{equation}

\nS This can easily be extended to any polynomial. Let $\Pbold$ be a general $d-$variate polynomial written
\begin{equation}
 \Pbold(x_1,\ldots,x_d) = \sum_{\nbold} a_{n_1\ldots n_d}\, x_1^{n_1} \ldots x_d^{n_d}
\end{equation}
where the sum $\sum_{\nbold}$ is over a finite number of terms. Then, one can write
\begin{equation}
 \Pbold = \sum_{\nbold} a_{n_1\ldots n_d}\, x_1^{n_1} \otimes  \ldots \otimes x_d^{n_d}
\end{equation}
which is a CP decomposition of $\Pbold$. Hence
\begin{equation}
 \rank \Pbold \leq N
\end{equation}
where $N$ is the number of monomials in $\Pbold$. As a consequence
\begin{equation}
 \rank \QS(P, \otimes_\m \xbold_\mu) \leq N
\end{equation}

%
\section{Approximation of functions by Polynomials}\label{sec:quant:SW}
%

Here I present the Stone-Weierstrass theorem which states that any real continuous function on a nice subset of $\R^d$ can be approximated as close as it is wished by a multivariate polynomial for norm $\ell^\infty$. As polynomials lead to low rank quantized forms, any discretization of a real continuous function can be approximated as close as wished by a low rank tensor. The approximation is obtained for uniform convergence, which is related to $\ell^\infty$ norm. There is much more room in a ball of radius $\epsilon$ for $\ell^\infty$ norm than for $\ell^2$ norm for which best low rank approximations have been defined (see section \ref{sec:inter:some}). Hence a second effort must be made for having a desired approximation for $\ell^2$ norm. 

\nT{Stone-Weierstrass theorem} Let us first recall Stone-Weierstrass theorem in a general form before an application for approximating a continuous function by polynomials. Let $S$ be  a compact Hausdorff space\footnote{In topology, a Hausdorff space is a set $E$ such that if $x,y \in E$ with $x \neq y$, there exists a neighborhood $\NC(x)$ of $x$ and a neighborhood $\NC(y)$ of $y$ such that $\NC(x) \cap \NC(y)= \emptyset$.}. Denote by $C(S)$ the set of all real-valued continuous functions on $S$
\[C(S) = \left\{
 f \: : \: S \longrightarrow \R, \quad f \: \mbox{continuous}
 \right\}
\]
Let $A$ be a subset of $C(S)$ closed by multiplication, i.e.
\[
 f,g \in A \quad \Longrightarrow fg \in A
\]
such that any function $f \in A$ separates $S$, i.e. if $x,y \in S$
\[
 x \neq y \quad \Longrightarrow \quad \exists \: f \in A \: : \: f(x) \neq f(y)
\]
and that any constant function is in $A$. Then, $A$ is dense in $C(S)$ in the sense of maximum norm, i.e.
\begin{equation}
 \forall \: \epsilon >0, \quad \forall f \in C(S), \quad  \exists \: \varphi \in A \: : \: \forall \: x \in S, \quad |f(x)-\varphi(x)| < \epsilon
\end{equation}
As a consequence, for any function $f \in C(S)$, there is a sequence $(\varphi_n)_n$ with $\varphi_n \in A$ which converges uniformly to $f$.


\nT{Application to approximation of continuous functions by polynomials}  Here, Stone Weierstrass  theorem is used with $A$ being the set of polynomials on $S$. Let us present it with the example for real continuous functions on a product of closed interval
\begin{equation}
 \begin{CD}
  [a,b]^d @>f>> \R
 \end{CD}
\end{equation}
and generalization to more general situations (like convex bodies) is straightforward. Let $S = [a,b]^d$, and consider the set $C(S)$ of real continuous functions on $[a,b]^d$. The product of two polynomials is a polynomials, the polynomials separate the points in $S$ (if two points are different, there exists a polynomial which has different values on both points), and any constant function is a polynomial (of degree 0). Then, Stone-Weierstrass theorem applies, and 
\begin{equation}
 \forall \: f \in C(S), \quad  \forall \: \epsilon > 0, \quad \exists \: \Pbold \in \PS^{\otimes d} \: : \: \forall \: x \in S, \quad |f(x)-\Pbold(x)| < \epsilon 
\end{equation}

\notes Approximation of a continuous function by polynomials has a long history and has been thoroughly studied for decades for finding best numerical approximations of functions (see e.g. \cite{Rivlin1981,Canuto2006,Olver2010}). Weierstrass first established his theorem in 1885 for functions continuous on $[0,1]$ best approximated by univariate polynomials. His theorem has been extended to more general settings by Stone in 1937, and such a generalized setting has been presented here following \cite[p. 126]{Lax2002}. A classical proof of Weierstrass theorem has been published by Bernstein in 1912, relying on what is called nowadays Bernstein polynomials (see \cite[section 10.3]{Davidson2010} or \cite[section 1.1.1.]{Rivlin1981} for a demonstration). Borel elaborating on works by Chebychev has proved that, given a continuous function $f$ on compact interval $[a,b] \subset \R$, and an integer $n \in \N$, there exists a polynomial $P_n$ of degree $n$ closest to $f$ for $\ell^\infty$ norm and that this polynomial is unique (this is called the \emph{Chebychev approximation theorem}, see \cite[theorem 1.8]{Rivlin1981}). Chebychev polynomials $T_n(x)$ of degree $n$ are the ones on $[-1,1]$ with leading monomial $x^n$ that deviate least from 0 for norm $\ell^\infty$ (or, equivalently, with minimal $\ell^\infty$ norm). They have been derived by Chebychev in 1853, and if $x = \cos \theta$, $T_n(x)=\cos n\theta$ (see \cite{Rivlin1981,Rivlin1974} for further information). This leads to an efficient algorithm for finding $P_n$, polynomial of degree $n$ which approximates best function $f$ on $[-1,+1]$ with $\ell^{\infty}$ norm based on the following idea: find the $n+1$ roots $(x_k)_k$ with $x_k \in [-1,+1]$ of Chebychev polynomial of degree $n+1$ (called \emph{Chebyshev node}) and find the polynomial $P$ of degree $n$ such that $P_n(x_k)=f(x_k)$ (see \cite[section 5.8]{Press2007} for further details).

%
\section{Decomposition of a function on a basis of orthonormal polynomials}\label{sec:quant:ortho}
%

Any square integrable function on, say, $S = [a,b]$ with possibly $a=-\infty$ and $b=+\infty$ can be approximated as accurately as prescribed by an univariate polynomial for $\ell^2$ norm by a truncation of its Fourier development on a basis or orthonormal polynomials. 

\nS Here, $\H$ denotes the Hilbert space of $L^2$ integrable functions on $S$. The scalar product is defined by
\begin{equation}
 \langle f,g\rangle = \int_a^b\, w(x)f(x)g(x)\, dx
\end{equation}
where $w$ is called a \kw{weight function} and depends on $a$ and $b$. There exists a basis for $\H$ of orthonormal polynomial. If $[a,b] = [-1,+1]$
\begin{equation}
 w(x) = \frac{1}{\sqrt{1-x^2}}
\end{equation}
and the polynomial of degree $n$ is called the \kw{Chebyshev polynomial}\index{polynomial!Chebyshev} and denoted by $T_n(x)$. One has
\begin{equation}
 \forall \:\: m,n \in \N, \qquad 
 \int_{-1}^{+1}\, \frac{T_m(x)T_n(x)}{\sqrt{1-x^2}}\,  \, dx = 
 \begin{cases}
  1 & \mbox{if} \:\: m=n \\
  0 & \mbox{if} \:\: m \neq n
 \end{cases}
\end{equation}

\nS Any (measurable) function $f$ such that $\int_{-1}^{+1} (1-x^2)^{-1/2}|f(x)|^2\, dx < +\infty $ can be decomposed on this basis as
\begin{equation}
 f = \sum_{n \in \N}\, \widehat{f}_n\, T_n \qquad (\mbox{in }\H)
\end{equation}
with
\begin{equation}
f, T_n \in \H, \qquad \widehat{f}_n \in \R 
\end{equation}
or
\begin{equation}
 f(x) = \sum_{n \in \N}\, \widehat{f}\, T_n(x) \qquad (\mbox{in }\R)
\end{equation}
The coefficients $\widehat{f}_n$ are called \kw{Fourier coefficients}, by reference to Fourier series built on $\R/\Z$ and where the orthonormal basis is made in complex plane $\C$ ($\R/\Z$ is isomorphic with the unit circle in complex plane) by functions $\exp 2i\pi k \omega x$. 


\nT{Multivariate functions and Tucker approximation} Let us extend this to $d-$ variate real functions. Therefore, let $d \geq 1$, $S=[-1,+1]^d$, $\xbold \in S$ and $\H$ be the set of real functions on $S$ such that $\int_S (1-x_1^2)^{-1/2}\ldots (1-x_d^2)^{-1/2}|f(\xbold)|^2\, d\xbold < +\infty $. Let
\[
 \mbold = (m_1, \ldots,m_d), \quad \nbold = (n_1,\ldots,n_d) \in \N^d, \qquad \xbold = (x_1,\ldots,x_d) \in \R^d
\]
and define 
\begin{equation}
 \Tbold_\nbold(\xbold) = T_{n_1}(x_1) \ldots T_{n_d}(x_d)
\end{equation}
which can be written as
\begin{equation}
 \Tbold_\nbold = \bigotimes_{\mu=1}^d\, T_{n_\mu}
\end{equation}
It is easy to see that the system $(\Tbold_\nbold)_{\nbold \in \N^d}$ forms an orthonormal basis for $\H$, i.e.
\begin{equation}
 \forall \:\: \mbold,\nbold \in \N^d, \qquad 
 \int_S\, \frac{\Tbold_\mbold(\xbold)\Tbold_\nbold(\xbold)}{\prod_\mu\sqrt{1-x_\mu^2}}\,  \, d\xbold = 
 \begin{cases}
  1 & \mbox{if} \:\: \mbold = \nbold \\
  0 & \mbox{if} \:\: \mbold \neq \nbold
 \end{cases}
\end{equation}
and that each function $\bm{\phi} \in \H$ can be decomposed as
\begin{equation}
 \bm{\phi} = \sum_{\nbold \in \N^d}\, \widehat{\phi}_\nbold \, \Tbold_\nbold
\end{equation}
This can, be truncated at a given rank $\rbold = (r_1,\ldots,r_d)$  to have a finite rank approximation of $\bm{\phi}$ as 
\begin{equation}\label{eq:quant:ortho:1}
 \bm{\phi} = \sum_{\nbold \prec \rbold}\, \widehat{\phi}_\nbold \, \Tbold_\nbold + \Ebold
\end{equation}
where $\Ebold$ is an error term and $\rbold \prec \nbold$ means that $\forall \; \mu, \: r_\mu \leq n_\mu$. This is a finite rank \kwnind{Tucker approximation}\index{Tucker approximation} of $\bm{\phi}$. The Fourier coefficients $\widehat{\phi}_\nbold$ can be arranged in a $d-$modes tensor $\widehat{\bm{\phi}} \in \R^{r_1 \times \ldots \times r_d}$ with components $\widehat{\phi}_\nbold$.

\notes These notes on finite rank approximation of multivariate functions are elementary and basic. Approximation of multivariate functions is a flourishing domain (see e.g. several chapters in \cite{Olver2010}), and the aim of this section is to establish a natural link with low rank approximation of tensors. 

%
\section{Elementary tensor operations on multivariate functions}\label{sec:quant:elem} 
%

What is Tucker approximation useful for? A good reason for developing this approach is the observation that elementary operations on tensor generalize to multivariate functions sometimes naturally, and sometimes not. These operations are the backbone for implementing and understanding best low rank approximation of a tensor (a function) by a Tucker model (its Fourier development on a basis or orthonormal polynomials), like HOSVD. Let us recall here some elementary operations useful for that:\\
\begin{center}
 \begin{tabular}{l|c}
   \hline
   Elementary operation & see section \\
   \hline
   Permutation & \ref{sec:tenselem:permut} \\
   Symetrization & \ref{sec:tenselem:sym} \\
   slicing & \ref{sec:tenselem:slice} \\
   Kronecker product & \ref{sec:tenselem:reshape} \\
   reshaping & \ref{sec:tenselem:reshape} \\
   matricization & \ref{sec:tenselem:reshape} \\
   \hline
 \end{tabular}
\end{center}


\nT{Permutation} Let us transpose the example in section \ref{sec:tenselem:permut} on a $3-$modes tensor, hence a $3-$variate function. Let $\sigma \in \SF_3$ be defined as
\begin{equation}
 \sigma = 
 \begin{pmatrix}
  0 & 1 & 2 \\
  2 & 0 & 1
 \end{pmatrix}
 , \qquad
  \sigma^{-1} = 
 \begin{pmatrix}
  0 & 1 & 2 \\
  1 & 2 & 0
 \end{pmatrix}
\end{equation}
or
\begin{equation}
 \sigma(0)=2, \qquad \sigma(1)=0, \qquad \sigma(2)=1
\end{equation}
Then
\begin{equation}
 \sigma\left(\sum_{i,j,k} \, a_{ijk} \, \ebold_i \otimes \fbold_j \otimes \gbold_k\right) = \sum_{i,j,k} \, a_{ijk} \, \gbold_k \otimes \ebold_i \otimes \fbold_j 
 \end{equation}
 which is transported as
 \begin{equation}
  (\sigma(\bm{\phi}))(x,y,z) = \bm{\phi}(y,z,x)
 \end{equation}
 For example,
 \begin{equation}
  \bm{\phi}(x,y,z) = xy^2z^3 \quad \Longrightarrow \quad (\sigma(\bm{\phi}))(x,y,z) = x^3yz^2
 \end{equation}
 Note\marginpar{\dbend} that the variables $x,y,z$ play the role of the indices in a finite dimensional tensor, hence the permutation $\sigma^{-1}$ is applied to them (see equation \ref{eq:tenselem:permut:1}) 
 
 \nT{Slicing} Let us transport the example of section \ref{sec:tenselem:slice}. We fix $i$ and let indices $j,k$ running to yield $\Abold_{1,i}= \sum_{j,k}\, a_{ijk}\, \fbold_j \otimes \gbold_k$. This can be transported as: $(i)$ select first variable $(ii)$ fix $x$ and define the slice
 \begin{equation}
  \bm{\phi}_{1,x}(y,z) = \bm{\phi}(x,y,z)
 \end{equation}
For example, if $\bm{\phi}(x,y,z)=xy^2z^3$, then
\begin{equation}
 \bm{\phi}_{1,3/2}(y,z)= \frac{3y^2z^3}{2}
\end{equation}
(here, $x=3/2$).

\nT{Kronecker product, vectorization, reshaping} None of these elementary operations can be transported as easily as permutation or slicing. Let us show it on the example of the vectorization of a matrix. Let us recall that if $A \in \R^{m \times n}$ with vector $\abold_i \in \R^m$ as column $i$, $\vec A$ is a vector in $\R^{mn}$ built by concatenation of columns of $A$ from left to right. This relies on a function
\begin{equation}\label{eq:quant:elem:1}
 \begin{CD}
   (i,j) @>>> k
 \end{CD}
\end{equation}
where $(i,j)$ are the indices in $A$ and $k$ is the index in $\vec A$, i.e
\begin{equation}
 (\vec A)_k = A_{ij}, \qquad k = k(i,j)
\end{equation}
with
\begin{equation}
 k(i,j) = (j-1)m + i
\end{equation}
This is possible because there exists a bijection between $\{1,m\} \times \{1,n\}$ and $\{1,mn\}$ or even between $\N$ and $\N \times \N$. But an equivalent bijection cannot be explicitly calculated between $\R \times \R$ and $\R$ (Cantor has shown that there is a one to one function between $\R^2$ and $\R$; Peano curve is an example of a surjective and continuous curve from $[0,1]$ to $[0,1]^2$ but it is not injective, and there cannot be an homeomorphism between $\R \times \R$ and $\R$ ($\R^p$ and $\R^q$ are not homeomorphic if $p \neq q$). So, vectorization cannot be built easily this way. Such a limit exists for Kronecker product, matricization and reshaping, which all rely on the same function defined in equation (\ref{eq:quant:elem:1}).

\nS This is where Tucker best approximation enters into the game. Let $\bm{\phi}$ be a $d-$variate function in a functional space $\H^{\otimes d}$, and $\rbold = (r_1, \ldots,r_d)$. Then, a major quality of \kw{Tucker approximation} of $\bm{\phi}$ at rank $\rbold$ (see equation (\ref{eq:quant:ortho:1} and comment reproduced here)
\begin{equation}
 \bm{\phi} = \sum_{\nbold \prec \rbold}\, \widehat{\phi}_\nbold \, \Tbold_\nbold + \Ebold
\end{equation}
is that the Fourier coefficients $\widehat{\phi}_\nbold$ for $\nbold \prec \rbold$ can be arranged in a $d-$modes tensor $\widehat{\bm{\phi}} \in \R^{r_1 \times \ldots \times r_d}$ with components $\widehat{\phi}_\nbold$. Then, matricization and reshaping can be operated on tensor $\widehat{\bm{\phi}}$.  This leads to a two steps procedure: given a $d-$variate function $\bm{\phi} \in \H^{\otimes d}$
\begin{enumerate}
 \item select a (not necessarily low) rank $\rbold$ which defines a finite dimensional vector space of polynomials spanned by the basis of orthonormal polynomials up to degree $\rbold$, and project $\bm{\phi}$ on this space yielding $\widehat{\bm{\phi}} \in \R^{r_1 \times \ldots \times r_d}$,
 \item apply tools for best low rank approximation (CP, Tucker, TT) on finite dimensional tensor $\widehat{\bm{\phi}}$
\end{enumerate}
Note \marginpar{\dbend} that for step 1, finding best approximation by a polynomial with prescribed accuracy is not an easy task.

\nT{HOSVD} An example for this strategy is HOSVD of a multivariate function. It requires matricization along all modes, says $M_x$ for matrix along the first variable $x$, and SVD of each matrix, say $M_x$. None of these operation is simple for $d-$variate functions. Hence, the strategy is to have a fairly good finite dimensional approximation $\widehat{\phibold}$ of $\phibold$ via a polynomial approximation and implement HOSVD on $\widehat{\phibold}$. As new basis are given from SVD of matricization of $\widehat{\phibold}$ along each mode, Eckart-Young theorem applies, and solutions are nested. It is quite possible to have the polynomial approximation of $\phibold$ with $\ell^\infty$ norm, and best low rank approximations of $\widehat{\phi}$  with Frobenius norm $\ell^2$.

%
\section{In a nutshell}\label{sec:quant:nutshell} 
%

The translation of a $d-$variate function either into a tensor $\Xbold \in \R^{n \times \ldots \times n}$ (its evaluation on a mesh $\xbold$) or its decomposition on a basis of orthonormal polynomials permits to take advantages of convergence theorems in functional spaces. A major application is to build $\Xbold$ which yields an as accurate as possible numerical solution of a linear PDE. We will not enter in these notes into this domain, which is immense (see \cite{Canuto2006,Shen2011,Einsiedler2017}). In this chapter, we have traveled on the following diagram:\\
\begin{center}
 \begin{tikzpicture}
   \node (fun) at (0,2) {$\phibold \in H^{\otimes d}$} ;
   \node (dec) at (0,0) {$\displaystyle \phibold = \sum_{\nbold \in \N^d} \, \widehat{\phi}_\nbold\, \Tbold_\nbold $} ;
   \node (num) at (4,2) {$\QS(\phibold,\xbold) \in (\R^n)^{\otimes d}$} ; 
   \node (low) at (4,0) {$\displaystyle \widetilde{\phibold} = \sum_{\nbold \prec \rbold} \widehat{\phi}_\nbold \, \Tbold_\nbold $} ; 
   \draw[->] (fun) -- (dec) node [pos=0.5, left] {$(b)$} ;
   \draw[->] (fun) -- (num) node [pos=0.5, above, sloped] {$(a)$} ; 
   \draw[->] (dec) -- (low) node [pos=0.5, below, sloped] {$(c)$} ;
   \draw[->] (low) -- (num) node [pos=0.5, right] {$(e)$} ;;
   \draw[->] (fun) -- (low) node [pos=0.5, above, sloped] {$(d)$};
 \end{tikzpicture}
 \begin{equation}\label{sk:func}
  \mbox{}
 \end{equation}

\end{center}

\noindent The meaning of the arrows are as follows:
\begin{enumerate}[label=$(\alph*)$]
 \item We are given function $\phibold$ which belongs to a Hilbert space $\H^{\otimes d}$ (e.g. $L^2(S^d)$ with $S=[-1,+1]$) and build a $d-$modes tensor $\QS(\phibold, \xbold)$ which is the discretization of $\phibold$ a given $d-$modes mesh (see sections\ref{sec:quant:mesh} \& \ref{sec:quant:quant}) 
 \item the function $\phibold$ can be exactltly decomposed as a countable sum on a basis $(\Tbold_\nbold)_\nbold$ of $\H^{\otimes d}$, like a Fourier basis. The coefficients $\widehat{\phi}_\nbold$ are in the frequency space (see section \ref{sec:quant:ortho}).
 \item Approximation theorems in functional spaces lead to a finite rank approximation $\widetilde{\phibold}$ of $\phibold$ by a truncation of the decomposition in the frequency space at a given order $\rbold$; the Fourier coefficients can be arranged in a tensor $\widehat{\phibold} \in \R^{r_1 \times \ldots \times r_d}$ (see sections \ref{sec:quant:ortho} \& \ref{sec:quant:elem}).
 \item The Stone-Weierstrass theorem says that polynomials in $d$ variables are dense in the set of continuous functions on a compact space for uniform convergence norm (see section \ref{sec:quant:SW}).
 \item The evaluation map of $\widehat{\phibold}$ on the mesh $\xbold$ leads to a finite Tucker rank approximation of $\QS(\phibold,\xbold)$, which is independent of the size of the mesh.
\end{enumerate}

\nB Let us note that arrows $(c)$ and $(d)$ are two ways to approximate a multivariate by a polynomial of given degree; it does not mean they lead to the same polynomial.